\documentclass[10pt]{article}
\usepackage{amsmath,amssymb,amsfonts,amsthm,amscd,graphicx,psfrag,epsfig}
\usepackage{color}
\definecolor{blue}{rgb}{0,0,0.7}
\definecolor{red}{rgb}{0.75, 0, 0}
\usepackage{stmaryrd}
\usepackage[titletoc,toc]{appendix}
\usepackage{hyperref}
\hypersetup{colorlinks, linkcolor=blue, citecolor=cyan}
\usepackage[all]{xy}           
\usepackage{marginnote}  
\usepackage{float}
\usepackage{soul}

\newtheorem{theorem}{Theorem}[section]
\newtheorem{lemma}[theorem]{Lemma}
\newtheorem{proposition}[theorem]{Proposition}
\newtheorem{corollary}[theorem]{Corollary}
\newtheorem{conjecture}[theorem]{Conjecture}
\newtheorem{definition}[theorem]{Definition}

\newcommand{\bpf}{\begin{proof}}
\newcommand{\epf}{\end{proof}}
\newcommand{\bs}{\begin{split}}
\newcommand{\es}{\end{split}}
\newcommand{\be}{\begin{equation}}
\newcommand{\ee}{\end{equation}}
\newcommand{\bt}{\begin{theorem}}
\newcommand{\et}{\end{theorem}}
\newcommand{\bd}{\begin{definition}}
\newcommand{\ed}{\end{definition}}
\newcommand{\bp}{\begin{proposition}}
\newcommand{\ep}{\end{proposition}}
\newcommand{\bl}{\begin{lemma}}
\newcommand{\el}{\end{lemma}}
\newcommand{\bc}{\begin{corollary}}
\newcommand{\ec}{\end{corollary}}
\newcommand{\bcon}{\begin{conjecture}}
\newcommand{\econ}{\end{conjecture}}
\newcommand{\la}{\label}
\newcommand{\B}{{\rm B}}
\newcommand{\A}{{\rm A}}
\newcommand{\bE}{{\bf E}}
\newcommand{\Z}{{\mathbb Z}}
\newcommand{\R}{{\mathbb R}}
\newcommand{\Q}{{\mathbb Q}}
\newcommand{\C}{{\mathbb C}}
\newcommand{\G}{{\rm G}}

\newcommand{\hra}{\hookrightarrow}
\newcommand{\lra}{\longrightarrow}
\newcommand{\lms}{\longmapsto}
\newcommand{\bS}{{\Bbb S}}

\begin{document}

 \title{Spectral description of non-commutative local systems on surfaces and non-commutative cluster varieties}

 \author{Alexander Goncharov, Maxim Kontsevich}
\date{\it To the memory of Yuri Ivanovich Manin}

\maketitle

 \tableofcontents

 \begin{abstract} Let $R$ be a non-commutative field. We   prove that generic triples of flags in an $m-$dimensional $R-$vector space  are described by 
 flat $R-$line bundles on the honeycomb graph with $\frac{(m-1)(m-2)}{2}$ holes. \vskip 2mm
 
 Generalising this, we prove that the non-commutative stack ${\cal X}_{m, \bS}$ 
 of     framed flat $R-$vector bundles of rank $m$ 
on  a decorated surface $\bS$ contains  open dense substacks, identified with   stacks of      
  flat line bundles on   certain bipartite graphs $\Gamma$ on $\bS$.  \vskip 2mm
  
  We  introduce non-commutative cluster Poisson varieties related to bipartite ribbon graphs. They 
   carry a canonical non-commutative Poisson structure. The result above just means  that the   space ${\cal X}_{m, \bS}$  
   has a structure of a non-commutative cluster Poisson variety, equivariant under  the action of the mapping class group of $\bS$.   \vskip 2mm   
   
   For bipartite graphs on a torus, we get the non-commutative    dimer cluster integrable system.  \vskip 2mm
     
  We develop a parallel dual story of  non-commutative cluster ${\cal A}-$varieties related to bipartite ribbon graphs. They carry a canonical non-commutative 2-form. 
  The  dual non-commutative moduli space ${\cal A}_{m, \bS}$ of twisted decorated local systems on $\bS$ carries a cluster ${\cal A}-$variety structure, equivariant under the action of the mapping class group of $\bS$. The non-commutative cluster ${\cal A}-$coordinates on the space ${\cal A}_{m, \bS}$ are expressed as ratios of Gelfand-Retakh 
  quasideterminants. In the case $m=2$ this recovers the Berenstein-Retakh non-commutative cluster algebras related to surfaces. \vskip 1mm

For any split reductive group $\G$ with connected center, we prove 
    that all stacks of framed $\G-$Stokes data  carry a  cluster Poisson structure, equivariant under the wild mapping class group. 
  Therefore these stacks can be equivariantly quantized.  
  
  The similar stacks of decorated $\G-$Stokes data carry an equivariant cluster ${\cal A}-$variety structure.

    We introduce    admissible dg-sheaves,  and define non-commutative stacks of  Stokes data as    stacks  of admissible dg-sheaves of certain type.

   \end{abstract}

\section{Introduction}  

{\bf Summary}. Section \ref{SECT1.1a} is a selfcontained  exposition of several key results of the paper. All of them but 
Proposition \ref{PR1.7} are equipped with complete proofs. 
Computations proving Proposition \ref{PR1.7} are presented in Section \ref{sec4.2}. 
They show how  Gelfand-Retakh's quasideterminants \cite{GR} appear from the geometric considerations, making many  of their properties obvious.

In particular, we prove in Section \ref{SECT1.1a} the following. Take  a surface ${\cal S}$, possibly singular and non-orientable, obtained by gluing triangles and  deleting their vertices. 
Then  generic local system of $m$-dimensional 
vector spaces over a skew field $R$ on  ${\cal S}$, equipped with a flat section of the associated  local system of flags near each deleted vertex, can be described  
by a flat line bundle on a certain graph $\Gamma$. \vskip 2mm

One can  glue a surface from triangles in various  ways. Different gluing patterns lead to different graphs $\Gamma$. 
The main question not resolved in Section \ref{SECT1.1a} is how to relate   flat line bundles on  different graphs $\Gamma$   describing the same generic local system on the surface. 
To address this question   we introduce  in  Sections \ref{secziya2} $\&$ \ref{SSEECC5AA} 
non-commutative cluster varieties of two flavors, related to  bipartite ribbon graphs. They carry, respectively, a non-commutative Poisson structure and a closed non-commutative 2-form. In Sections  \ref{secc3x} $\&$ \ref{sec4} 
we define these  cluster  structures  on  moduli spaces of $R-$local systems on surfaces of two flavors, depending on  the type of flags near the vertices. 
The crucial role play  bipartite graphs on the surface called {\it ${\rm GL}_m-$graphs}  \cite{G}. A preview can be found in Section \ref{SECT1.1}. 
  \vskip 2mm

${\rm GL}_m-$graphs  appear in the description of local systems on surfaces out of the blue. 
In Section \ref{SECCT7}  we present another approach, which explains the role of ${\rm GL}_m-$graphs, and  proves the main results  in a  bit stronger form.  
We start from any collection ${\cal H}$ of smooth cooriented curves  with disjoint Legendrians on an oriented  surface.  
We assign to ${\cal H}$ a  finite type  stack 
${\cal M}({\cal H}, {\bf 1})$. We show explicitly that  deformations of the collection ${\cal H}$ leaving its components  smooth and their Legendrians disjoint  
 lead to equivalent stacks. 

Here is an example.   Let ${\cal S}$ be a surface $S$ with a finite set  of punctures. Draw $m$  small  simple non-intersecting loops around each puncture $s$. 
Denote by  ${\cal H}_m$  their union.   
Then   ${\cal M}({\cal H}_m, {\bf 1})$ is the the stack ${\cal X}_m(S)$ of  
 $m-$dimensional local systems on $S$  with a complete flag near each puncture, invariant under the  monodromy around $s$. 
On the other hand,  zig-zag paths of 
any bipartite graph $\Gamma$ on $S$ provide   a 
collection ${\cal Z}_\Gamma$ of cooriented curves on $S$. We show that  ${\cal M}({{\cal Z}_\Gamma}, {\bf 1})$ contains as an open subspace the space   
${\rm Loc}_1(\Gamma)$ of flat line bundles on 
the graph $\Gamma$. 
By the very definition of a ${\rm GL}_m-$graph,  the collection of its zig-zag paths  is isotopic to the  collection of loops ${\cal H}_m$. 
Such an isotopy provides an open embedding  $$
{\rm Loc}_1(\Gamma)\subset {\cal M}({{\cal Z}_\Gamma}, {\bf 1}) \stackrel{\sim}{=} {\cal M}({{\cal H}_m}, {\bf 1}) = {\cal X}_m(S).
$$
So we embed the non-commutative torus ${\rm Loc}_1(\Gamma)$ as an open part  in ${\cal X}_m(S)$. The union of these tori for all ${\rm GL}_m-$graphs on $S$ 
is a cluster Poisson atlas in the sense of Section \ref{secziya2}, equivariant under the   mapping class group action.  See the preview  in Section \ref{SEC1.2}. 
\vskip 2mm

In Section \ref{sec8.4} we construct a wide supply of stacks with a similar cluster description. We start from    a split reductive group $\G$ over $\Q$, and pick 
 any collection $\{\beta_s\}$ of elements of the cyclic envelope of the Artin braid semigroup ${\rm Br}^+_\G$, assigned to the punctures $s$ of $S$. We introduce a 
 stack\footnote{Here stack means traditional stack used in algebraic geometry.}  ${\cal X}_\G(S, \{\beta_s\})$ of $\G-$local systems on $S$ with a pattern  of flags near each puncture $s$  determined 
by  the $\beta_s$. When $\beta_s=e$ is the unit for all $s$, 
and $\G={\rm GL}_m$, we get the stack ${\cal X}_m(S)$. When $S= S^2 -\{0, \infty\}$ and $\beta_0=e$, 
we recover the cluster Poisson varieties introduced in \cite{FG3}. Using this and  \cite{GS19} we show that, 
skipping  $S^2-\{0\}$, any stack ${\cal X}_\G({S}, \{\beta_s\})$ carries a cluster Poisson structure equivariant under the action of the wild mapping class group for the pair $({S}, \{\beta_s\})$. 
We  get similar results  for any  surface ${\cal S}$ glued from triangles,  decorated in addition to the $\{\beta_s\}$ at the punctures $s$ 
by the elements of the braid semigroup ${\rm Br}^+_\G$ at the  corners. 

These stacks are of  independent interest. 
For example, any stack of Stokes data for holomorphic $\G-$bundles with meromorphic connections on a Riemann surface  is holomorphically 
equivalent to one of them. 
However unless $\G$ is of type $\A_1$,  many stacks ${\cal X}_\G(S, \{\beta_s\})$ do not appear this way. 

We describe  non-commutative analogs of stacks  ${\cal X}_{{\rm GL}_m}({\cal S}, \{\beta_s\})$ and prove that they 
carry a  non-commutative cluster Poisson structure,  equivariant under the   wild mapping class group action.  \vskip 2mm

Section \ref{SEC6} contains an application of the constructions of Section \ref{secziya2} to non-commutative cluster integrable systems. It is entirely independent of the rest of the paper. 

\subsection{Spectral description of generic triples of flags} \la{SECT1.1a}

Let $R$ be a {\it skew field}, i.e. a non-commutative division algebra. 
 An $m$-dimensional $R-$vector space is a  free rank $m$  left 
$R$-module. If $m=1$, we call it an 
 $R-$line. 

\vskip 2mm
A {\it flag} ${\cal F} $ in an $m$-dimensional $R$-vector space ${V}$ is  a filtration  
by  $R-$vector spaces
\be \la{FLAG}
V = {\cal F}^{0} \supset {\cal F}^{1} \supset {\cal F}^{2} \supset \ldots \supset 
{\cal F}^{m-1} \supset {\cal F}^{m} = 0, 
\quad {\rm codim}~{\cal F}^{i} = i.
\ee
We associate to a flag ${\cal F}$ a collection of $m$ lines 
${\rm gr}^i{\cal F}:= {\cal F}^{i}/{\cal F}^{i+1}$, $i=0, ..., m-1$.  

The group of automorphisms of the $R$-vector space ${V}$ acts transitively on the set of all flags in $V$. \\

\begin{figure}[t]
\centerline{\epsfbox{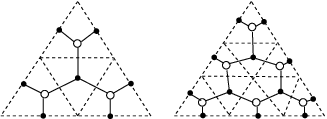}}
\caption{The bipartite graphs $\Gamma_2$ and $\Gamma_3$, shown by solid lines.}
\label{ncls10}
\end{figure}

\begin{figure}[ht]
\centerline{\epsfbox{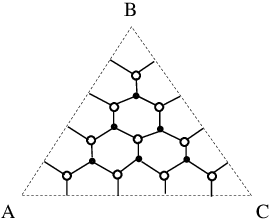}}
\caption{A bipartite graph $\Gamma_4$ for a  generic triple of 
flags $({\cal A}, {\cal B}, {\cal C})$ in 
a $4$-dimensional space. }
\label{ncls3aa}
\end{figure}

In Section \ref{SECT1.1a} we describe    generic  triples of flags in an $m-$dimensional $R-$vector space via flat line bundles on the bipartite graph $\Gamma_m$, 
shown on Figures \ref{ncls10}, \ref{ncls3aa}. 
It is a self-contained part   of the  paper, 
  describing  the simplest, and at the same time crucial result with all details.  
  Let us start with generic pairs of flags. 
 
\paragraph{1. Generic pairs of flags.}  A pair of flags $({\cal A}, {\cal B})$   in an $m$-dimensional $R-$vector space $V$    is  
{\it generic} if for 
any  integers $a,b\geq 0$ such that $a+b=m$ the   projections  
$V \lra V/{\cal A}^a$ and $V \lra V/{\cal B}^b$ induce an isomorphism 
$$
V \stackrel{\sim}{\lra}  V/{\cal A}^a \oplus V/{\cal B}^b. 
$$ 
Given a generic pair of flags (${\cal A}, {\cal B}$), 
let us set 
$$
{\rm L}^{b+1}:=  {\cal A}^{a} \cap {\cal B}^{b}, ~~~~a+b=m-1, ~~a,b \geq 0.
$$
The  embeddings ${\rm L}^a \hra V$, $a=1, ..., m$,  give rise to
 a canonical isomorphism 
$$
{\rm L}^1 \oplus \ldots \oplus {\rm L}^{m} \stackrel{\sim}{\lra} V. 
$$

\bl Assigning to a generic pair of  flags $({\cal A}, {\cal B})$ in an $m$-dimensional 
$R$-vector space $V$ the ordered collection of $m$ lines 
$({\rm L}^1, \ldots , {\rm L}^{m})$ we get an equivalence of  groupoids: 
\be 
\begin{split}
&\mbox{Groupoid of generic pairs of flags in $m$-dimensional spaces $(V, {\cal A}, {\cal B})$} 
\stackrel{\sim}{\lra} \\
&~~~~~~~~~~\mbox{Groupoid of ordered collections of $m$ lines $({\rm L}^1, \ldots , {\rm L}^{m})$}.\\
\end{split}
\ee
\el

\begin{proof}
The inverse functor $({\rm L}^1, \ldots , {\rm L}^m) \lra (V, {\cal A}, {\cal B})$ is given by  
$$
V:= {\rm L}^1 \oplus \ldots \oplus {\rm L}^{m}, ~~~ {\cal A}^a:= {\rm L}^1\oplus \ldots \oplus  {\rm L}^{m-a}, ~~~ 
{\cal B}^b:= {\rm L}^{b+1}\oplus \ldots \oplus  {\rm L}^{m}.
$$
\end{proof}

\paragraph{2. Spectral description of the groupoid of generic triples of flags.}

A triple of  flags $({\cal A}, {\cal B}, {\cal C})$ in an $m$-dimensional $R$-vector space $V$  is {\it generic} if for any   non-negative integers  $(a,b,c)$ with $a+b+c=m$,   
the following map is an isomorphism:\footnote{Note that a  flag     
in   $V$ gives rise to the dual flag  
 in   $V^*$, but the dual to a generic triple of flags may not be generic.}  
 \be\la{THFMI}
V \lra V/{\cal A}^a \oplus V/{\cal B}^b \oplus V/{\cal C}^c, \ \ \ \ \ \ a+b+c=m.
\ee
Let us  describe generic triples of  flags  in an $m$-dimensional $R$-vector space.  

Let us assign  a triple of flags $({\cal A}, {\cal B}, {\cal C})$ to the vertices of a triangle $t$,  labeled      by ${\rm A}, {\rm B}, {\rm C}$. 

We introduce a 
bipartite graph $\Gamma_m$,  shown in Figures \ref{ncls10}-\ref{ncls3aa} for   $m\leq 4$. 
Recall that a {\it bipartite graph} is a graph with vertices of two types, $\bullet$-vertices and $\circ$-vertices, 
such that each edge has a vertex of each type.

To define the graph $\Gamma_m$, take an $m$-triangulation of the triangle $t$, shown by punctured lines on  Figure \ref{ncls10}. It  
subdivides each side  of the triangle  into 
$m$  equal segments, and tessellates the triangle  into   
  little triangles of two types: the  
$\bullet$-triangles, and the $\circ$-triangles.  
 The bipartite graph $\Gamma_m$ is the dual graph for this 
tessellation. \vskip 1mm

Recall  that a flat line bundle $L$ on a graph $\Gamma$ is given by the following data:

 \begin{enumerate}
\item  A line $L_v$ for every vertex $v$ of $\Gamma$. 

\item  An isomorphism $t_{w,  v}: L_v\to L_w$ 
for each edge of $\Gamma$ with the ends $v,w$, with $t_{v, w}= t_{w, v}^{-1}$. 
 \end{enumerate}

\begin{figure}[ht]
\centerline{\epsfbox{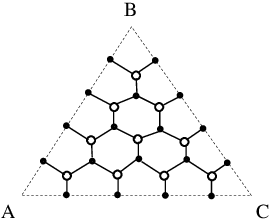}}
\caption{The bipartite graph $\overline \Gamma_4$ includes the $\bullet-$vertices on the sides of the triangle.}
\label{ncls3a}
\end{figure}

\bt \la{MT1}
There is a canonical equivalence between the following two groupoids:
\begin{enumerate}
\item Groupoid of generic triples of flags in an $m$-dimensional $R$-vector space.

\item Groupoid of $R$-line bundles with  connections on the bipartite graph $\Gamma_m$, see Figure \ref{ncls3a}.
\end{enumerate}
\et

 \begin{proof}  
1.  {\it The functor ${\cal L}: \{ \mbox{Generic triples of flags}\} \lra \{\mbox{Flat line bundles on the graph $\Gamma_m$}\}$.}

\noindent
Given a generic triple of  flags $({\cal A}, {\cal B}, {\cal C})$,   
 let us define a flat $R$-line bundle     ${\cal L}({\cal A}, {\cal B}, {\cal C})$ on the    graph $\Gamma_m$.  
 
\begin{enumerate}

\item  

 The $\circ-$vertices of the graph $\Gamma_m$ are parametrized by the triples of non-negative integers $(a,b,c)$ such that $a+b+c=m-1$. 
We assign to each $\circ-$vertex a one dimensional subspace:
\be
{\rm L}^\circ_{a,b,c}:= {\cal A}^{a} \cap {\cal B}^{b} \cap {\cal C}^{c}, \ \ \ \ \ \ a+b+c=m-1. 
\ee

\item   The $\bullet-$vertices of the  graph $\Gamma_m$ are parametrized by the triples of non-negative integers $(a,b,c)$ such that $a+b+c=m-2$. We assign to each $\bullet-$vertex 
 a two dimensional subspace: 
 \be
{\rm P}_{a,b,c}:= {\cal A}^{a} \cap {\cal B}^{b} \cap {\cal C}^{c}, \ \ \ \ \ \ a+b+c=m-2. 
\ee  
\end{enumerate}
 The plane ${\rm P}_{a,b,c}$  assigned to a $\bullet-$vertex $b$  contains three lines  assigned to the   $\circ-$vertices  incident to $b$:
  \be \la{CANISOa}
\begin{gathered}
    \xymatrix{
    &{\rm L}^\circ_{a,b+1,c} \ar[d]  &\\
 &  {\rm P}_{a,b,c}    &\\
 {\rm L}^\circ_{a+1,b,c}  \ar[ru]   &&   {\rm L}^\circ_{a,b,c+1} \ar[lu]   }
\end{gathered}
 \ee
 Indeed, the embedding 
${\cal A}^{a+1}\hra {\cal A}^{a}$ induces the embedding $ 
 {\cal A}^{a+1}\cap {\cal B}^{b} \cap {\cal C}^{c} \hra {\cal A}^{a} \cap {\cal B}^{b} \cap {\cal C}^{c}.
 $   
Condition (\ref{THFMI}) just means that for each plane ${\rm P}_{a,b,c}$, the three lines (\ref{CANISOa})   are disjoint. We define a line ${\rm L}^\bullet_{a,b,c}$ as the kernel of the natural map from the sum of the three lines in (\ref{CANISOa}) to the plane:
 \be
 {\rm L}^\bullet_{a,b,c}:= {\rm Ker}\Bigl({\rm L}_{a+1,b,c} \oplus   {\rm L}_{a,b+1,c} \oplus   {\rm L}_{a+1,b,c} \lra  {\rm P}_{a,b,c} \Bigr).\ee
 Then there are three maps, induced by  projections of the subspace $ {\rm L}^\bullet_{a,b,c}$ onto the   summands:
   \be \la{CANISO}
\begin{gathered}
    \xymatrix{
    &{\rm L}^\circ_{a,b+1,c} &\\
 & \ar[ld]_\sim  \ar[u]^\sim {\rm L}^\bullet_{a,b,c}    \ar[rd]^\sim  &\\
 {\rm L}^\circ_{a+1,b,c}     &&   {\rm L}^\circ_{a,b,c+1}  }
\end{gathered}
 \ee 
 These three maps are isomorphisms if and only if   the three lines in the plane    ${\rm P}_{a,b,c}$ are disjoint.

So we get a flat line bundle   on the graph $\Gamma_m$: its fibers at the vertices are given by the lines ${\rm L}^\circ$ and ${\rm L}^\bullet$, and the 
parallel transport along an edge  $\bullet \to \circ$ is given by a map ${\rm L}^\bullet \lra {\rm L}^\circ$ in (\ref{CANISO}). \\ 

{\bf Remark.} If  $R=K$ is commutative, the moduli space of flat line bundles on a graph $\Gamma$ is identified with  $H^1(\Gamma, K^*)$. 
So  flat line bundles on the graph $\Gamma_m$ are described by  their monodromies around the 
${m-1 \choose 2}$  holes in the graph. These are literally the canonical coordinates from \cite{FG1}, called the {\it triple ratios}.  

If  $R$ is non-commutative, we no longer have  canonical  coordinates  describing  generic triples. 
For example,  the isomorphism classes of   flat $R$-line bundles on  the graph $\Gamma_3$ are described   by the elements of $R^*$ considered modulo the conjugation. 
For $m>3$ we   need to pick   a base point and paths to get a collection of elements of $R^*$ defined up to a common conjugation.  \\

\begin{figure}[t]
\centerline{\epsfbox{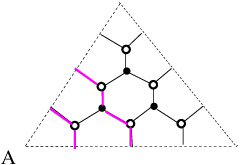}}
\caption{Vertices of the graph $\Gamma_3$ on the distance $\leq 2$ from the vertex ${\rm A}$ are on the two (purple) zig-zags.}
\label{ncls3z1}
\end{figure}

2. {\it  The reconstruction functor    ${\cal R}:  \{\mbox{Flat line bundles on $\Gamma_m$}\} \lra \{\mbox{generic triples of flags}\}$.} Given a flat line bundle ${\cal L}$ on the graph $\Gamma_m$, we get a map 
  from the direct sum over all internal $\bullet-$vertices $\{\bf b\}$ of the graph $\Gamma_m$ to the direct sum over all $\circ-$vertices $\{\bf w\}$ of the graph $\Gamma_m$: 
\be \la{Kastq}
{\rm K}_{\cal L}: \bigoplus_{\mbox{internal $\bullet-$vertices ${\bf b}$}} {\cal  L}_{\bf b}  \lra \bigoplus_{\mbox{ $\circ-$vertices ${\bf w}$}} {\cal L}_{\bf w}. 
\ee
We define a vector space $V_{\cal L}$ as the cokernel of this map:
\be \la{kocq}
V_{\cal L}:= {\rm Coker}({\rm K}_{\cal L}). 
\ee 
Let us define three flags in this vector space, assigned to the vertices of the triangle   $t$. Given a vertex $\A$,   
consider a function $d_{\A} : \{\mbox{vertices   of the graph $\Gamma_m$}\} \lra \{1, ..., m\}$, given by the distance from a vertex    to the   $\A$, see Figure \ref{ncls3z1}.\footnote{That is the distance to $\A$ from the zig-zag strand  containing the vertex on the graph $\Gamma_m$ and parallel to the side $\B{\rm C}$.} The   function $m-d_{\A}$ induces a decreasing filtration on the  complex (\ref{Kastq}).\footnote{That is a filtration on each of the spaces, preserved by  
the operator $K_{\cal L}$.} So it induces a  filtration on  
${\rm Coker}({\rm K}_{\cal L})$. We define 
  the flag ${\cal A}^\bullet_{\cal L}$ as the   space (\ref{kocq}) with this filtration. Explicitly: 
\be \la{FA}
{\cal A}^a_{\cal L}:= {\rm Coker}\Bigl( \bigoplus_{d_{\A}({\bf b})\leq  m-a}  {\cal L}_{\bf b}  \lra \bigoplus_{d_{\A}({\bf w})\leq m-a}  {\cal L}_{\bf w} \Bigr), \ \ \ \ \ \  a=0, ..., m-1.
\ee
 This way we get a functor
 \be
 {\cal R}: {\cal L} \lra ({\cal A}^\bullet_{\cal L}, {\cal B}^\bullet_{\cal L}, {\cal C}^\bullet_{\cal L}).
 \ee
By the  definition, the  functors ${\cal L}$ and ${\cal R}$ provide the equivalence of  categories.    Theorem \ref{MT1} is proved. \end{proof}

\paragraph{3. Reconstruction functor revisited.} 
Denote by $\Gamma_m'$  the graph  obtained from the graph $\Gamma_m$ by adding   $\bullet-$vertices at the endpoints of  external edges, see Figure \ref{ncls3a}, 
pictured  on the sides of the triangle.

A flat  line bundle ${\cal L}$ on   $\overline \Gamma_m$ provides three $m-$dimensional $R-$vector spaces assigned to the sides of the triangle as follows. 
For example, take the side $\A\B$. Let  $V^\bullet_{{\rm A}{\rm B}}$   be   the direct sum of the fibers of ${\cal L}$   at the 
  $\bullet-$vertices on the   side ${{\rm A}{\rm B}}$, and  $V^\circ_{{\rm A}{\rm B}}$    the direct sum of the fibers of ${\cal L}$   at the 
  $\circ-$vertices closest to the   side ${{\rm A}{\rm B}}$. The parallel transport along the edges from $\bullet-$vertices on the side $\A\B$ to the closest  $\circ-$vertices provides a canonical  isomorphism $V^\bullet_{{\rm A}{\rm B}}  = V^\circ_{{\rm A}{\rm B}}$. So we set 
  $$
  V_{\A\B}:= V^\bullet_{{\rm A}{\rm B}}  = V^\circ_{{\rm A}{\rm B}}.
  $$

 Denote by $V^\bullet_{\cal L}$ (respectively $V^\circ_{\cal L}$) the direct sum of the lines ${\cal L}_{\bf b}$ over all $\bullet-$vertices ${\bf b}$ (respectively the lines ${\cal L}_{\bf w}$ over the $\circ-$vertices ${\bf w}$) 
of the graph $\Gamma_m$ on Figure \ref{ncls3aa}.  
 Then there is a canonical map, induced by the embedding $V^\circ_{{\rm A}{\rm B}} \hra V^\circ_{{\cal L}}$, followed by the projection to  the quotient $V_{\cal L}$, defined in (\ref{kocq}):
\be
\psi_{\A\B}: V_{\A\B} \lra V_{\cal L}.
\ee

\bl \la{Prr} One has $\psi_{\A\B}=\psi_{\B\A}$. The map $\psi_{\A\B}$ is an isomorphism. \el

\begin{proof}  The first claim  is the tautology. 
 The second  is equivalent to the statement that 
$$
 V^\circ_{\A\B} \cap  {\rm K}_{\cal L}(V_{\cal L}^\bullet) =0.
 $$
Let  $\widetilde V^\circ_{\cal L}$ be the direct sum of ${\cal L}_{\bf w}$ over all but the ones closest to the side $\A\B$ $\circ-$vertices. So
$
 V^\circ_{\cal L} = V_{\A\B} \oplus \widetilde V^\circ_{\cal L}.
 $ 
 Note that the map  $K_{{\cal L}, {\rm C}}: V^\bullet_{\cal L}\stackrel{\sim}{\lra} \widetilde V^\circ_{\cal L}$ given by the parallel transport from each $\bullet-$vertex  along the edge  towards the vertex ${\rm C}$ is an isomorphism.  There are increasing filtrations of the spaces  $\widetilde V^\circ_{\cal L}$ and $V^\bullet_{\cal L}$, given by the codistance from the $\circ$/$\bullet-$vertices to the vertex ${\rm C}$ of the triangle, such that the  map ${\rm gr}K_{\cal L}$ of their associated graded induced by the map $K_{\cal L}$ is the map $K_{{\cal L}, {\rm C}}$. 
 This implies the claim. \end{proof}

Therefore we have three canonical isomorphisms $\psi_{\A\B}$, $\psi_{\B{\rm C}}$, $\psi_{{\rm C}\A}$:
   \be \la{CANISO**}
\begin{gathered}
    \xymatrix{
 {\rm V}_{\A\B}   \ar[dr]_{\psi_{\A\B}}   &&  \ar[dl]^{\psi_{\B{\rm C}}}  {V}_{\B{\rm C}}  \\
 & {\rm V}_{\cal L}    &\\
 &{\rm V}_{{\rm C}\A}  \ar[u]_{\psi_{{\rm C}\A}} &\\}
\end{gathered}
 \ee 
 
\begin{figure}[t]
\centerline{\epsfbox{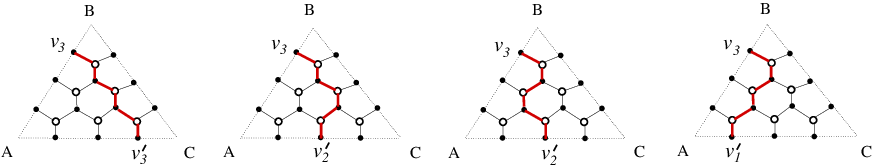}}
\caption{There are four ${\rm A}-$paths from the vertex $v_3$ on the side $AB$ to the side $AC$.}
\label{ncls3abc}
\end{figure}

Our next goal is to calculate  the composition
\be \la{psi1}
\psi_{\A{\rm C}}^{-1}\circ \psi_{\A\B}:   V_{\A\B} \stackrel{}{\lra} V_{\cal L} \stackrel{}{\lra} V_{\A{\rm C}}. 
\ee
To present the answer,  we define  explicitly an isomorphism  
 \be \la{AB1*}
\begin{split}
&\varphi_{{\rm A}, {\rm B} \to {\rm C}}:  V_{{\rm A}{\rm B}} \lra V_{{\rm A}{\rm C}}. \\
\end{split}
\ee 
Namely, we define an {\it  ${\rm A}-$path} as a path on the graph from a vertex on the  side $\A\B$ to a vertex on the   side $\A{\rm C}$, 
which  always  go towards $\A{\rm C}$. So given a vertex $v_i$ on 
the side ${\rm A}{\rm B}$ on the distance $i$ from   ${\rm A}$, 
and a vertex $v_j'$ on the side 
${\rm A}{\rm C}$ on the distance $j$  from  ${\rm A}$, there are ${i-1}\choose{j-1}$  ${\rm A}-${paths}  $\gamma: v_i\lra v_j'$, see   Figure \ref{ncls3abc}. 
Then we set:
\be \la{AB1}
\begin{split}
&\varphi_{{\rm A}, {\rm B} \to {\rm C}} := \sum_{\gamma} (-1)^{i-1}t_\gamma,\\
\end{split}
\ee
where the sum is over all    
${\rm A}-$paths $\gamma: v_i \lra v_j'$ from  ${\rm AB}$ to   ${\rm AC}$, and 
$t_\gamma: {\cal L}_{v_i}\lra {\cal L}_{v'_j}$ is the parallel transport in the local system ${\cal L}$ along the path  $\gamma$.  

\bp \la{CR} The map (\ref{psi1}) coincides with the map (\ref{AB1*}):
$$
\psi_{\A{\rm C}}^{-1}\circ \psi_{\A\B} = \varphi_{{\rm A}, {\rm B} \to {\rm C}}.
$$
\ep

 \begin{proof} Let $f$ be a $\bullet-$vertex of the graph $\Gamma_m$. Denote by $a,b,c$ the  $\circ-$vertices next to $f$, so that they are in the natural  bijection with the vertices of the triangle, e.g. the $a$ is north of $f$ and matches   $\A$. Pick a non-zero vector $f\in {\cal L}_f$. Denote by $e_a, e_b, e_c$ its parallel transports to the fibers at the vertices $a,b,c$.  Denote by $e'_a, e'_b, e'_c$ their projections to ${V}_{\cal L}$. Then  $e'_a=-e'_b-e'_c$. Note that $e_b$ and $e_c$ are parallel transports of the 
 vector  $e_a$ to the fibers of ${\cal L}$ at the  vertices $b$ and $c$. 
 Then apply the same procedure to $e_b$ and $e_c$ and so on. In the end we express the vector $e'_a$ as a linear combination of the vectors at the vertices on 
 the side $\A{\rm C}$.  Applying this to the $\circ-$vertices closest to the side $\A\B$ we prove the formula. 
 \end{proof}

The space $V_{\A\B}$ 
carries  flags ${\cal F}^\bullet_{{\rm A}{\rm B}}$ and ${\cal F}^\bullet_{{\rm B}{\rm A}}$,   assigned to the orientations   of the side. 
Namely, let ${\bf b}_1, ..., {\bf b}_m$ be the  $\bullet-$vertices on the side ${{\rm A}{\rm B}}$, oriented $\A\to \B$. 
Then the subspaces of the flag ${\cal F}^\bullet_{{\rm A}{\rm B}}$ are  
\be
{\cal F}^a_{{\rm A}{\rm B}}:= {\cal L}_{{\bf b}_1} \oplus \ldots \oplus  {\cal L}_{{\bf b}_{m-a}}. 
\ee

It is clear from the definitions that the following three flags in the space $V_{\cal L}$ coincide 
\be \la{CDF}
 \psi_{\A\B}({\cal F}^\bullet_{\A\B}) = \psi_{\A{\rm C}}({\cal F}^\bullet_{{\rm A}{\rm C}}) = {\cal A}_{\cal L}^\bullet.
\ee
 \vskip 2mm
 
So we can define a variant ${\cal R}'$ of the reconstruction functor ${\cal R}$ as follows.
 
 Take the three spaces $V_{\A\B}, V_{\B{\rm C}}, V_{{\rm C}\A}$ assigned to the sides of the triangle. There are 
 six explicitly defined isomorphisms between them: the  isomorphism $\psi_{\A, \B\to {\rm C}}$  given   by  (\ref{AB1}),  and the  ones obtained by permuting $\{\A, \B, {\rm C}\}$. 
 Proposition \ref{CR} implies that 
\be
\begin{split}
 & \varphi_{{\rm A}, {\rm C} \to {\rm B}} \circ \varphi_{{\rm A}, {\rm B} \to {\rm C}}   = {\rm Id}.\\
 & \varphi_{{\rm B}, {\rm C} \to {\rm A}} \circ 
  \varphi_{{\rm C}, {\rm A} \to {\rm B}} \circ \varphi_{{\rm A}, {\rm B} \to {\rm C}} = {\rm Id}.\\
\end{split}
\ee
Similar identities are obtained by  permuting $\{\A, \B, {\rm C}\}$.
Therefore we can identify     the three vector spaces, and  transform the flags to one of them, getting a configuration of three flags. 
These flags   are  in generic position.  Finally, it is clear from the construction that we have the following Proposition.

\bp \la{Pr} The two reconstruction functors ${\cal R}$  and ${\cal R}'$  coincide. 
\ep

  The  functor ${\cal R}'$ is a variant of the snake algorithm \cite[Section 9.7]{FG1} in the commutative case. We use it  to  prove in Theorem \ref{T1.9} an analog of the Laurent phenomenon in the non-commutative case.

 \paragraph{4. Spectral description of   decorated triples of flags.}  
A {\it decorated flag}  is 
a flag ${\cal F}^\bullet$ as in (\ref{FLAG})  with 
a choice of a non-zero element $f_i\in  {\rm gr}^i{\cal F}:= {\cal F}^{i-1}/{\cal F}^{i}$ in each successive quotient. 

\begin{figure}[ht]
\centerline{\epsfbox{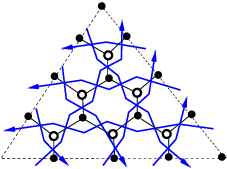}}
\caption{Nine zig-zag paths on the bipartite graph $\Gamma_3$. We push zig-zags a bit out of the vertices.}
\label{ncls3z}
\end{figure}

A {\it zig-zag path}, or simply a {\it zig-zag},  on  a ribbon graph $\Gamma$ is a path turning at every vertex either left or right, so that the left and right turns alternate. 
Zig-zags on a bipartite ribbon graph are oriented so that they turn right at     $\circ-$vertices. 
The  graph $\Gamma_m$ carries $3m$ oriented zig-zags,   see Figure \ref{ncls3z}. 

\bt \la{MT1A}
There is a canonical equivalence between the following two groupoids:
\begin{enumerate}
\item Groupoid of generic triples of decorated flags in an $m$-dimensional $R$-vector space.

\item Groupoid of $R$-line bundles with  connections on the bipartite graph $\Gamma_m$, trivialised on every zig-zag.
\end{enumerate}
\et

Theorem \ref{MT1A} is immediately deduced from Theorem \ref{MT1} and the description of the reconstruction functor, which  turns  
 trivializations of the local system on zig-zags   into decoration vectors of the flags. \\

\begin{figure}[ht]
\centerline{\epsfbox{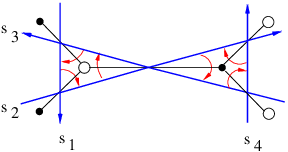}}
\caption{Writing $\Delta_E = s_1/s_2$ etc. where $s_i$ are flat sections on zig-zags, the counterclockwise products   are: $s_1/s_2 \cdot s_2/s_3\cdot s_3/s_1=1$ near the $\circ-$vertex, and 
 $s_2/s_3\cdot s_3/s_4\cdot s_4/s_2=1$ near the $\bullet-$vertex.}
\label{nclsm}
\end{figure}

The notable difference is that in the decorated setting  there is a coordinate description. Indeed, for each   edge $E$ of the graph $\Gamma_m$ there are   two zig-zags $\gamma_1, \gamma_2$ 
containing $E$. 
They provide two  trivialisations $s_{\gamma_1}$ and $s_{\gamma_2}$ of the restriction of  a flat line bundle to $E$. We  assume that   the oriented 
zig-zag $\gamma_1$ goes along     the edge $E$ in the direction $\circ\to \bullet$.  We define an edge coordinate $\Delta_E$ as their ratio:
\be
s_{\gamma_1} = \Delta_E\cdot s_{\gamma_2}, \ \ \ \ \Delta_E \in R^*.
\ee
Elements $\{\Delta_E\}$   satisfy the following monomial relations. For each internal vertex  of $\Gamma_m$, the  
  counterclockwise product of the  elements $\Delta_E$ over the edges   incident to the vertex 
is equal to $1$, see Figure \ref{nclsm}.  There are no other relations between  elements $\{\Delta_E\}$. 

  Recall that, given an $n \times n$   matrix $M$  over a non-commutative skew field 
   $R$,  Gelfand and Retakh \cite{GR} defined quasideterminants $M_{ij}$,   $1 \leq i,j\leq n$.  If  $R$ is commutative, then  the quasideterminant  $M_{ij}$ is, up to a sign,  the  ratio of the determinant of $M$ by a minor obtained by 
  deleting the row $i$ and the column $j$ of $M$.  See 
 (\ref{GRREF}) for examples of quasideterminants of $2\times 2$ matrices.

  \bp \la{PR1.7}
 Given a generic triple of decorated flags, the invariants $\{\Delta_E\}$ at the edges of the graph $\Gamma_m$ are ratios of two Gelfand-Retakh quasideterminants. 
\ep
See a detailed discussion of the simplest $m=2$ case in Section \ref{sec4.2}. \\

In the commutative case, assuming   that the space $V$ carries a volume form,  the  invariants $\Delta_E$ are   ratios of  the ${\cal A}-$coordinates  $A_f$ \cite{FG1}, assigned to
 the faces  of the graph $\Gamma_m$.
Precisely, each edge $E$ is shared by two faces $f^+$ and $f^-$, where $f^+$ is to the left of the $\circ \to \bullet$ oriented edge $E$, and 
 $
\Delta_E=  {A_{f^+}}/{A_{f^-}}.
 $

\paragraph{5. A generalization for   surfaces glued from triangles.} Theorem \ref{MT1} and \ref{MT1A} have the following immediate generalizations. 
Take a finite collection of triangles,  glue them edge-to-edge  in an arbitrary way, and delete the vertices of the glued triangles. 
  The obtained  "surface" ${\cal S}$  is \underline{singular }along the edges   where more then two triangles meet. It might be \underline {non-orientable}. 
  Its boundary is   the union of the edges which are not glued to other edges. It has {\it corners}, also referred to as  {\it special points},  and {\it punctures} provided by the deleted vertices. 
  The punctures correspond to those deleted vertices whose neighborhoods are topologically non-trivial. We call the union of the corners and  punctures {\it marked points}.  These are precisely the points which has been deleted. 
  It comes with a  triangulation ${\cal T}$.
  
     Let us consider the following groupoids.

\begin{enumerate}

\item Groupoid ${\cal X}_{m, {\cal S}}$ of local systems of $m$-dimensional $R$-vector spaces on  ${\cal S}$, equipped with a     flag 
in a fiber near each special point $s$, invariant under the   monodromy around  $s$, such that  for each triangle of ${\cal T}$  
 the  triple of flags obtained by moving the flags near the vertices to one point is generic. 

\item Groupoid ${\cal A}_{m, {\cal S}}$ of local systems of $m$-dimensional $R$-vector spaces on  ${\cal S}$, equipped with a  \underline{decorated}    flag 
in a fiber near each special point $s$, invariant under the   monodromy $s$, such that  for each triangle of ${\cal T}$  
 the  triple of flags obtained by moving the flags near the vertices to one point is generic. 
 
 \end{enumerate}
 
 We call the local systems with the data as in 1) {\it framed} local systems, and as in 2) the {\it decorated} ones. \vskip 2mm

If ${\cal S}=t$ is just a  triangle,   ${\cal X}_{m, {t}}$ is the groupoid of triples of   flags in generic position  an $m$-dimensional $R$-vector space. 
Indeed, any local system on a triangle is trivial, so the triple of flags at the vertices  is the only data left. 
Similarly, ${\cal A}_{m, {t}}$ is the groupoid of triples of  decorated flags in generic position in  an $m$-dimensional $R$-vector space. 

\vskip 3mm
The graphs $\Gamma'_m$ are glued  into a bipartite graph $\Gamma_{m}({\cal T})$: gluing the triangles along the sides we glue the matching black vertices on the sides.

\bt \la{TH2.2} Denote by ${\cal S}$ such a possibly  singular and possibly unorientable triangulated surface.

1) The groupoid ${\cal X}_{m, {\cal S}}$    is equivalent to the groupoid of flat line  bundles 
on the  graph $\Gamma_m({\cal T})$. 

2) The groupoid ${\cal A}_{m, {\cal S}}$   is equivalent to the groupoid of flat line bundles 
on the  graph $\Gamma_m({\cal T})$, trivialized on each zig-zag path.  \et

\begin{proof} 1) 
Take triangles with the vertices $(a,b,c)$ and $(b',c',d)$ and glue them along the
 $bc-$sides, matching the $b$ and $c$ vertices, getting a rectangle $abcd$. The ${\rm GL}_m-$graphs $\Gamma_{abc}$ and $\Gamma_{b'c'd}$ 
assigned to these triangles are glued into a ${\rm GL}_m-$graph $\Gamma_{abcd}$ 
for the rectangle $abcd$. Given flat  line bundles ${\cal L}_{abc}$ and ${\cal L}_{b'c'd}$ on the graphs $\Gamma_{abc}$ and $\Gamma_{b'c'd}$, Theorem \ref{MT1} provides generic triples of flags $(A,B,C)$ and $(B',C',D)$ in  
$m-$dimensional $R-$vector spaces $V_{abc}$  and respectively $V_{b'c'd}$. Each of the two  spaces is decomposed  into a direct sum of $m$ lines, provided by the pair of flags $BC$ and $B'C'$. The sets of these lines are  identified.
So there  is a collection of  isomorphisms of $R-$vector spaces $i_{bc}: V_{abc} \lra V_{b'c'd}$ which identifies 
the lines and  the flags $\{B, C\}$ and $\{B', C'\}$.  These  isomorphisms form a right $(R^\times)^m-$torsor. 
Indeed, an isomorphism $i_{bc}$  is defined uniquely up to the composition with an $R-$linear map $\varphi'$ 
preserving  the lines assigned to the flags $B', C'$. 
So we get an $(R^\times)^m-$torsor of 
configurations of four flags $A,B,C,D$. 
Given such an isomorophism $i_{bc}$, there is a canonical way to 
glue the flat line bundles ${\cal L}_{abc}$ and ${\cal L}_{b'c'd}$ into a flat line bundle ${\cal L}_{abcd}$ on the graph $\Gamma_{abcd}$.  The obtained  line bundles  form a  $(R^\times)^m-$torsor. Conversely, we can  start from gluing the flat line bundles  ${\cal L}_{abc}$ and ${\cal L}_{b'c'd}$. This way we get a canonical equivalence between the two torsors.  Performing this equivalence for each pair of glued triangles, we prove the claim. 

2) The same arguments plus Theorem \ref{MT1A} provide flat line bundles on the graphs $\Gamma_m({\cal T})$ trivialized on zig-zag paths. \end{proof}

\paragraph{6. The non-commutative Laurent phenomenon.} 

Take a generic framed local system ${\cal V}$ of $m$-dimensional $R-$vector spaces on ${\cal S}$. 
Denote by ${\cal L}$ the associated flat line bundle on the graph $\Gamma_m({\cal T})$. Take an edge $e$ of the triangulation ${\cal T}$. 
Recall the  reconstruction functor ${\cal R}'$. The local system ${\rm V}:= {\cal R}'({\cal L})$ is isomorphic to ${\cal E}$ by Theorem \ref{MT1} and Proposition \ref{Pr}. Its fiber  ${\rm V}'_{e}$  over  $e$ is  the direct sum 
of the fibers ${\cal L}_{\bf b_i}$ of the flat line bundle ${\cal L}$ over the $\bullet-$vertices ${{\bf b}_i}$ on the edge $e$:
$$
{\rm V}'_{e}  =  \bigoplus_{{\bf b}_i}{\cal L}_{{\bf b}_i}.
$$
So the monodromy of  ${\rm V}$ along any closed path $\alpha$ on ${\cal S}$ interseting  the edge $e$ provides  a linear map 
$$
{\rm M}_\alpha: \bigoplus_{{\bf b}_i}{\cal L}_{{\bf b}_i} \stackrel{}{\lra} \bigoplus_{{\bf b}_i}{\cal L}_{{\bf b}_i}.
$$
It is described by the matrix elements $\mu_{ij}: {\cal L}_{{\bf b}_i} \lra {\cal L}_{{\bf b}_j}$.

Given a decorated local system  ${\cal F}$, the lines ${\cal L}_{{\bf b}_i}$ carry non-zero vectors. Therefore the monodromy is described  by a matrix $\mu_{ij}$ with the entries in $R^\times$. Recall that in this case there are elements $\Delta_E\in R^\times$ assigned to the   edges $E$ of the bipartite graph $\Gamma_m({\cal T})$, called the ${\cal A}-$coordinates of ${\cal F}$.

\bt \la{T1.9} The  matrix elements $\mu_{ij}$ of the monodromy of a framed local system ${\cal V}$ on ${\cal S}$  along a loop $\alpha$  are finite signed sums 
of the monodromies of the associated flat line bundle ${\cal L}$ on the graph $\Gamma_m({\cal T})$ along 
certain paths from ${\bf b}_i$ to ${\bf b}_j$ on the graph, homotopic to $\alpha$.

The matrix elements of the monodromy of a decorated local system ${\cal F}$ on ${\cal S}$ are finite signed sums 
of the products of the associated ${\cal A}-$coordinates or their inverses:
\be \la{LFE}
\mu_{ij} = \sum_{\gamma} \pm \Delta_{E_{\gamma, 1}}^{\varepsilon_1} \cdots \Delta_{E_{\gamma, k}}^{\varepsilon_k}, \qquad \varepsilon_j \in \pm1.
\ee
The sum is over certain paths $\gamma$ on the graph from ${\bf b}_i$ to  ${\bf b}_j$, and $E_{\gamma, j}$ are the consecutive edges of the   $\gamma$. 
The sign $\varepsilon$ in $\Delta_E^\varepsilon$ is determined by the direction of the edge $E$ on the  path $\gamma$: from  $\circ \to \bullet$, or $\bullet \to \circ$. 
\et

 \begin{proof}  The loop $\alpha$   intersects 
  edges $e_1, e_2, ..., e_n$ of  ${\cal T}$, ordered by the intersection points  $e_i \cap \alpha$ on the  $\alpha$. The   reconstruction functor ${\cal R}'$ for the triangle provides  explicit isomorphisms defined by formula (\ref{AB1}):
$$
{\rm V}_e \stackrel{{\cal R'}}{\lra} {\rm V}_{e_1}  \stackrel{{\cal R'}}{\lra} \ldots \stackrel{{\cal R'}}{\lra} {\rm V}_{e_n}  \stackrel{{\cal R'}}{\lra}  {\rm V}_e.
$$
The monodromy ${\rm M}_\alpha$ is their composition. This proves the first claim. The second follows from this since for the decorated local system each map ${\cal R}'$ is described explicitly by a matrix, whose entries are given by the signed sums over paths, and each path  provides a monomial in ${\cal A}-$coordinates or their inverses. 
\end{proof}

Formula (\ref{LFE}) is the non-commutative analog of the Laurent phenomenon for the monodromies. \\

When ${\cal S}$ is smooth and orientable, it is  the surface glued from triangles into a decorated surface. Recall that 
a {\it decorated surface} $\bS$ is a smooth oriented surface  with 
a finite set of {\it punctures} inside and  {\it special points} on the boundary, considered  modulo isotopy. 
We assume that each boundary component  has special points.  
A {\it marked point} is either a special   point or a puncture. \vskip 2mm

We described generic framed/decorated local systems on a decorated surface $\bS$ via flat line bundles on the graph $\Gamma_m({\cal T})$ for an ideal triangulation ${\cal T}$ of $\bS$. The mapping class group $\Gamma_\bS$ of $\bS$ acts by automorphisms of the space of local systems on $\bS$, preserving the data at the marked points. 
We want to see how its action  affects  the above description. For this we need, for each $\gamma \in \Gamma_\bS$,  to relate the flat line bundles on the graphs 
$\Gamma_m({\cal T})$ and $\Gamma_m(\gamma({\cal T}))$ corresponding to  a generic local system on $\bS$. Recall that a flip of a triangulation ${\cal T}$ is a new triangulation ${\cal T}'$ obtained by flipping the diagonal in one of the rectangles of ${\cal T}$. Since any two ideal triangulations are related by a sequence of flips, we have to solve the problem for a single flip. 
However if $m>2$, a flip  ${\cal T} \to {\cal T}'$ results in a complicated transformation of the bipartite graphs 
$\Gamma_m({\cal T}) \to \Gamma_m({\cal T}')$.  In \cite{G} this transformation was decomposed into  a composition of $(m-1)^2$ elementary 
moves, known as  {\it two by two moves},  of bipartite graphs of special kind, called ${\rm GL}_m-$graphs.  
So to relate the flat line bundles on the graphs $\Gamma_m({\cal T})$ and $\Gamma_m({\cal T}')$ assigned to a generic framed/decorated local system on $\bS$ we generalize 
the main constructions of Section \ref{SECT1.1a} to ${\rm GL}_m-$graphs. 
We discuss this in  Sections \ref{SECT1.1} - \ref{SEC1.2}.

\subsection{Describing noncommutative local systems on  surfaces} \la{SECT1.1} 
\paragraph{1. Spectral description of   local systems on   decorated surfaces via  {\it ${\rm GL_m}-$graphs}.}  For a decorated surface $\bS$  
  there are two more flexible variants  of Theorem \ref{TH2.2} 
  described below. 

Let us recall few more  facts about bipartite ribbon graphs, see also  Section \ref{SEC2B}.  
A bipartite ribbon graph $\Gamma$ gives rise to a decorated surface $\bS_\Gamma$ homotopy equivalent to $\Gamma$, obtained by gluing the ribbons. 
The graph $\Gamma$ is embedded to the surface  $\bS_\Gamma$. 
A bipartite ribbon graph $\Gamma$ gives rise to the conjugate bipartite ribbon graph $\Gamma^\ast$, 
obtained by altering the cyclic order of the edges at each $\bullet$-vertex. 
The surface $\Sigma_\Gamma$ associated with the ribbon graph $\Gamma^\ast$ is called the {\it spectral surface}.  
The collection ${\cal Z}_\Gamma$ of oriented zig-zags of $\Gamma$, being pushed to  $\bS_\Gamma$, is well defined up to an isotopy. 
It allows to reconstruct the  graph 
$\Gamma$. 
There are elementary transformations of bipartite ribbon graphs, called the {two by two move}. \vskip 2mm

\begin{figure}[t]
\centerline{\epsfbox{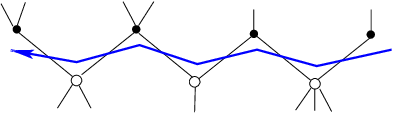}}
\caption{A zig-zag on a bipartite graph, naturally oriented,  and  pushed  out of the   vertices.}
\label{zig}
\end{figure}
 
 A special type of bipartite graphs on  $\bS$, called {\it ${\rm GL_m}-$graphs}, $m\geq 2$,  
was introduced in \cite{G}.   They are  characterized   the way their zig-zags relate to the marked points.  
Precisely, (on a finite cover of $\bS$,) each zig-zag path goes around a single marked point on $\bS$, and for each marked point $s$ on $\bS$ there are exactly $m$ zig-zag paths going around $s$. 
There is also a  technical condition: no parallel bigons, see Figure \ref{ncls12}.

Any {\it ideal triangulation} ${\cal T}$ of $\bS$, i.e., a triangulation with vertices 
at the marked points, gives rise to a  {${\rm GL_m}-$graph}  $\Gamma_m({\cal T})$  on $\bS$.
The {${\rm GL_m}-$graphs}  $\Gamma_m({\cal T})$ and $\Gamma_m({\cal T'})$ assigned to the triangulations ${\cal T}$ and ${\cal T}'$ related by a flip of an edge can be connected by a sequence of two by two moves \cite[Section 1.2]{G}.
\vskip 2mm

A flat connection  on an  $R$-vector bundle induces a flat connection  on the 
bundle of flags in its fibers.

\bd \la{9.9.11.1}
A  framed $R$-local system on a decorated surface $\bS$  
is a local system of finite dimensional $R$-vector spaces on $\bS$ 
with a  framing, given by  
 flat sections of the restriction of the associated local system of flags 
to the neighborhoods of the marked points.
\ed

\bt \la{MTI} 
 Let $\Gamma$ be a ${\rm GL_m}-$graph on a decorated surface $\bS$,    homotopy equivalent to $\bS$, 
and $\Sigma_\Gamma$ the spectral surface. 
 Consider   the following two  groupoids:  

\begin{enumerate}

\item 
Groupoid ${\cal X}_{m, \bS}$ of  $m$-dimensional   framed $R$-local systems on $\bS$;

\item  Groupoid ${\rm Loc}_1(\Sigma_\Gamma)$ of flat  $R$-line bundles with connection on the 
 spectral surface $\Sigma_\Gamma$. 
\end{enumerate}

Then  there is a canonical open embedding\footnote{In particular, they  are birationally equivalent.} 
\be \la{11.9.11.1}
\varphi_\Gamma:   {\rm Loc}_1(\Sigma_\Gamma)\hra {\cal X}_{m, \bS}. 
\ee  
\et

The spectral surface $\Sigma_\Gamma$ is homotopy equivalent to 
the graph $  \Gamma^*$. So local systems on $\Sigma_\Gamma$ are the same thing as   local systems on the graph $\Gamma^*$. Note also that $\Gamma =\Gamma^*$. 
So there are canonical identifications
\be
{\rm Loc}_1(\Sigma_\Gamma) = {\rm Loc}_1(\Gamma) = {\rm Loc}_1(\Gamma^*).
\ee

 Theorem \ref{MTI}  is proved  in Section \ref{SECCT7}. The main idea is 
  outlined in Section \ref{SEC1.2}. It is based on the interpretation of each of the moduli spaces as a  moduli space of 
{\it admissible dg-sheaves}.   We use only a combinatorial description of admissible dg-sheaves, and elementary linear algebra, valid over a skew field.\footnote{We do not use  constructions involving tensor products of vector spaces. The latter are not defined as $R-$vector spaces.}

An elementary proof of a slightly weaker version of Theorem \ref{MTI} is given in Theorem \ref{MTII}. Namely, we restrict to bipartite graphs which can be obtained by 
compositions of  two by two moves from the graphs $\Gamma_m({\cal T})$ related to ideal triangulations, 
and   give a simple  construction   of a  functor assigning to 
a  generic  local system on $\bS$ a flat  line bundle on the spectral surface $\Sigma_\Gamma$, 
providing  a birational equivalence   of the groupoids 1) and 2). Theorem \ref{MTII} is all we need to describe explicitly  the action of the mapping class group $\Gamma_\bS$ on the space of all framed local systems, discussed in the end of Section \ref{SECT1.1a}. \vskip 2mm


We show in Section \ref{SEC5.1} that 
{any} two by two move of   bipartite ribbon graphs $\mu: \Gamma \to \Gamma'$ gives rise to a 
non-commutative birational isomorphism $\mu$ intertwining $\varphi_{\Gamma}$ and $\varphi_{\Gamma'}$:
\be \la{11.9.11.2}
\mu: {\rm Loc}_1( \Gamma) \stackrel{\sim}{\lra} {\rm Loc}_1( \Gamma').
\ee
Therefore the natural action of the mapping class group of  $\bS$   on the non-commutative moduli space ${\cal X}_{m, \bS}$ can be presented
 by compositions of  birational isomorphisms (\ref{11.9.11.2}).\vskip 2mm 

Theorems \ref{MTI} $\&$ \ref{MTII} are  closely related to the Gaiotto-Moore-Neitzke  abelianization  of vector bundles with flat connection on a Riemann 
surface via flat line bundle on the associated spectral curve  \cite{GMN}.\footnote{Following the request of the referee, let us elaborate this. 
Gaiotto-Moore-Neitzke \cite{GMN}  start from the spectral cover $\pi: \Sigma_{\rm an}\to C$ related to a point of the Hitchin base of a Riemann surface $C$, and  assign to it a spectral network 
on $C$.
Given a flat line bundle ${\cal L}$ on $\Sigma_{\rm an}$, they define, using the spectral network,  an $n-$dimensional  flat bundle ${\cal E}$ on  $C$.
The construction  in general involves infinite formal sums. When $n>2$, it is not clear (and probably not always true) that  the topological spectral covers  corresponding to ideal webs 
can be realized as the analytic ones $\Sigma_{\rm an}$. However when the ideal web comes from an ideal triangulation of $S$, \cite{GMNs} defines the  spectral network bypassing the spectral curve, and the flat bundle ${\cal E}$ is given then by a finite sum construction. It is not clear how to do this for   more general ideal webs. In particular, it is not clear 
how to prove using this approach that flips of triangulations give rise to non-commutative cluster transformations, or to describe the action  of the mapping class group of  $\bS$   on the non-commutative moduli space ${\cal X}_{m, \bS}$. }\vskip 2mm 

When $R=K$ is    commutative  and $\Gamma$ is a bipartite ribbon graph,   
 ${\rm Loc}_1(\Gamma)(K)  = H^1(\Gamma, K^*)$. 
Any loop $\gamma$ on $\Gamma$ defines a   function   on ${\rm Loc}_1( \Gamma)$  given by the monodromy  around the loop. 
Given  a ${\rm GL_m}$-graph $\Gamma \subset \bS$, the  birational isomorphism (\ref{11.9.11.1}) is nothing else but  
the cluster Poisson  coordinate system on the  space ${\cal X}_{{\rm GL_m}, \bS}$ of  framed ${\rm GL_m}$-local systems on $\bS$  \cite{G}.  The ${\rm PGL}_m$-version of the 
  story was done in \cite{FG1}. \vskip 2mm

Going back to a skew field  $R$,  
   the monodromy of a flat $R$-line bundle over a circle 
is described by an element of   $R^*$ modulo conjugation. So  it does not give rise to a well defined 
  rational function on   non-commutative   space. The analogs of cluster Poisson coordinate systems are  non-commutative 
birational isomorphisms 
(\ref{11.9.11.1}). 
The analogs of cluster Poisson transformations  
are  non-commutative birational isomorphisms  
(\ref{11.9.11.2}) 
intertwining the ones  
(\ref{11.9.11.1}).   \vskip 2mm

Let us introduce  the non-commutative analog of the Poisson structure  
 in a   more general setup of non-commutative cluster Poisson varieties assigned to arbitrary bipartite ribbon graphs.  Here 
 Theorem \ref{MTII} is sufficient for the development of the theory, and all its applications discussed in the paper. 

\paragraph{2. Non-commutative cluster Poisson varieties.}  In the commutative case, a cluster Poisson 
variety is defined by starting from a quiver, and considering a collection of quivers related by sequences 
of  mutations. We assign to each quiver a cluster Poisson torus, and glue them by using  birational isomorphisms between the tori assigned to mutations. 

A special class of quivers and mutations  is provided by bipartite ribbon graphs and two by two moves between them \cite{GK}.  A cluster 
Poisson torus assigned to a bipartite graph $\Gamma$ is given by ${\rm H}^1(\Gamma, {\Bbb G}_m) = {\rm Loc}_1(\Gamma)$.    It has the following
 Poisson structure. Recall that the graphs $\Gamma$ and $\Gamma^*$ coincide - 
only their ribbon structures are different. On the other hand, the graph $\Gamma^*$ is homotopy equivalent to the spectral surface $\Sigma_\Gamma$. So the intersection pairing 
on ${\rm H}_1(\Sigma_\Gamma, \Z)$ induces a Poisson structure on 
  ${\rm H}^1(\Gamma, {\Bbb G}_m)$. 
Birational transformations (\ref{11.9.11.2}) in the commutative case are  cluster Poisson transformations. \vskip 2mm

In the non-commutative case, the analogs of   quivers are   bipartite ribbon graphs. The analogs of mutations are   two by two moves. The non-commutative  cluster Poisson torus 
assigned to a bipartite graph $\Gamma$  is  the moduli space ${\rm Loc}_1(\Gamma)$. 
Consider a  collection of   bipartite ribbon graphs $\{\Gamma'\}$ obtained by a sequence of two by two moves from the original graph $\Gamma$. 
Then groupoids ${\rm Loc}_1(\Gamma')$  are related by   
birational transformations   (\ref{11.9.11.2}). We show in the  end of Section \ref{SEC5.1} that they satisfy the pentagon relations. Therefore we get a non-commutative analog of a 
cluster Poisson variety ${\cal X}_{|\Gamma|}$ by gluing them according to these birational isomorphisms. \vskip 2mm

To introduce a non-commutative analog of the Poisson structure, consider    
 the free abelian group $L(\Gamma)$ generated by the loops on the graph $\Gamma$. It has a Lie algebra structure, provided by the bipartite ribbon structure on $\Gamma$. 
Since the graphs $\Gamma$ and $\Gamma^*$ coincide, there is another Lie bracket    on the same vector space, provided by the bipartite ribbon 
structure of  $\Gamma^*$. We denote by    $L(\Gamma^*)$ the latter Lie algebra.   

The commutative 
algebra given by the symmetric algebra   of $L(\Gamma^*)$ has a unique commutative Poisson 
algebra structure extending the Lie bracket on $L(\Gamma^*)$. It is the non-commutative analog of the Poisson structure on $H^1(\Gamma, {\Bbb G}_m)$. 
We prove that  it is preserved by non-commutative birational isomorphisms (\ref{11.9.11.2}). 
\vskip 2mm

  We conclude that  the    moduli spaces of $m$-dimensional  $R$-local systems 
on decorated surfaces carry a non-commutative cluster Poisson variety structure. 

\vskip 2mm 
Another class of interesting examples   are    non-commutative Grassmannians, 
parametrising the $m$-dimensional subspaces in a standard
 $N$-dimensional $R$-vector space, and their strata.  
 
 
 \vskip 2mm
Finally,   a few words about "non-commutative spaces". We say that $X$ is a non-commutative space if for any skew field $R$ there is a set $X(R)$ of $R$-points, depending functorially on $R$. We say that $X$ is a non-commutative stack if there is a non-commutative 
space ${\cal X}$ such that for each skew field $R$ the group ${\rm GL}_n(R)$ acts on ${\cal X}(R)$, and $X(R):={\cal X}(R)/{\rm GL}_n(R)$.  
For example,   for  the non-commutive space $ {\rm Loc}'_m(S)$ of $m$-dimensional local systems on a surface $S$, trivialized at a point $x\in S$, the 
  ${\rm Loc}'_m(S)(R)$ is the set of  $m$-dimensional $R$-local systems on  $S$, trivialized $x$. For the  non-commutative stack 
  ${\rm Loc}_m(S)$  we have ${\rm Loc}_m(S)(R):= {\rm Loc}'_m(S)(R)/{\rm GL}_m(R)$. 
  Non-commutative spaces/stacks  in the paper, e.g. the ones  ${\cal X}_m(\bS)$/${\cal A}_m(R)$ of framed/decorated 
   $m$-dimensional 
  local systems on a decorated surface $\bS$,  are understood this way. 
  
Quite often  
one can define points of a non-commutative space $X$ with values in the   $N \times N$ matrices ${\rm Mat}_N(K)$ over a field $K$, and 
 there are commutative spaces $X_N$ such that $X({\rm Mat}_N(K))=X_n(K)$. 
For example, $  {\rm Loc}'_m(S)({\rm Mat}_N(K)) :=   {\rm Loc}'_{mN}(S)(K)$. Then  a 
 non-commutative space gives rise to a collection of classical moduli spaces parametrized by positive integers $N$.  
The non-commutative differential forms, and symplectic or Poisson structures on such a non-commutative space $X$ are defined so that they 
provide  the usual differential forms, symplectic or Poisson structures on $X({\rm Mat}_N(K))$.

\paragraph{3. An application: Non-commutative cluster integrable systems.}   
In Section \ref{SEC6}, in the special case when $\Gamma$ is a bipartite graph on  a torus, we define a non-commutative generalisation of the 
  dimer models,  and  a non-commutative analog of the dimer cluster integrable system constructed 
in \cite{GK}.  After the definition of cluster Poisson varieties is given, the construction of the Hamiltonians and the proof that they Poisson commute follows literally the commutative case. 
In the simplest example we recover the non-commutative integrable system from \cite{K}. 

Since the dimer cluster integrable systems is equivalent to the relativistic Toda   systems 
studied by Fock and Marshakov  \cite{FM}, this also gives a non-commutative analog of the latter for ${\rm GL}_m$. 

\paragraph{4. Spectral description of   \underline{twisted} local systems on   decorated surfaces.}  
Let $T'\bS$ be the complement to the zero section of the 
tangent bundle on $\bS$.

\bd A {\it twisted local system on $\bS$} is  
a local system on $T'\bS$ with the monodromy $-1$ around a  
loop generating  $\pi_1(T'_s\bS)$. 
\ed
 
We consider twisted framed local systems of $R$-modules of rank $m$ on $\bS$.
 Notice that the $-1$ monodromy around a  
loop generating  $\pi_1(T'_s\bS)$ preserves the framing  flags.  

\vskip 2mm

\bt \la{MTIa} 
 Let $\Gamma$ be a ${\rm GL_m}-$graph on a decorated surface $\bS$,    homotopy equivalent to $\bS$, 
and $\Sigma_\Gamma$ the spectral surface. 
 Consider   the following two  groupoids:  

\begin{enumerate}

\item 
Groupoid of  $m$-dimensional \underline{twisted}  framed $R$-local systems on $\bS$;

\item  Groupoid of flat \underline{twisted} $R$-line bundles with connection on the 
 spectral surface $\Sigma_\Gamma$. 
\end{enumerate}

Then  there is a canonical open embedding of the second groupoid to the first. 
\et

  Theorem \ref{MTI} is the  "untwisted" version of Theorem \ref{MTIa}, with the   underlined words "twisted"  deleted. 
We can identify untwisted and twisted local systems on $\Sigma_\Gamma$ since the obstruction lies in $H^2(\Sigma_\Gamma, \Z/2\Z)=0$.

\paragraph{5. Non-commutative cluster ${\cal A}-$varieties.} In the commutative set-up, cluster varieties come in pairs, as cluster Poisson varieties, 
and cluster $K_2-$varieties, also referred to as cluster ${\cal A}-$varieties  \cite{FG2}. Moreover, there are two dual variants of  moduli spaces of local systems on decorated surfaces, which carry 
an equivariant cluster variety structure  of either  Poisson or  $K_2-$flavor. We show that in the non-commutative set-up the story is similar.  
Let us discuss now the non-commutative cluster ${\cal A}-$varieties. \vskip 2mm

A bipartite ribbon graph $\Gamma$ gives rise to a non-commutative torus ${\cal A}_\Gamma$ - the moduli space of 
 twisted flat $R-$line bundles on the associated spectral surface $\Sigma_\Gamma$, trivialized at the boundary.  It comes with a canonical collection of ${\cal A}-$coordinates $\{\Delta_E\}$,  
parametrised by the edges of the graph $\Gamma$, and satisfying monomial relations parametrised by the vertices of $\Gamma$.\footnote{Note that  boundary components of the spectral surface $\Sigma_\Gamma$ are in natural bijection with zig-zags on the  
  original surface $\bS_\Gamma$.  So this description match the one   in subsection 3, except that now we can talk   about the twisted flavor of the story.}  
For a two by two move 
$\mu: \Gamma_1 \lra \Gamma_2$, we define in Section \ref{SSEECC5} a non-commutative birational isomorphism  
\be \la{AMu}
\mu: {\cal A}_{\Gamma_1} \lra {\cal A}_{\Gamma_2}.
\ee
It serves the role of a non-commutative cluster ${\cal A}-$transformation. 
We show that  transformations (\ref{AMu}) satisfy  the pentagon relations. Therefore we can glue them into a non-commutative space  ${\cal A}_{|{\Gamma}|}$, 
which we call a non-commutative cluster ${\cal A}-$variety. 

We prove that a  non-commutative cluster ${\cal A}-$variety  carries a canonical non-commutative 
$2-$form $\Omega$. Precisely, given a bipartite ribbon graph $\Gamma$ we introduce a non-commutative 2-form $\Omega_\Gamma$  
on the    torus ${\cal A}_\Gamma$. It has a  simple expression in the cluster ${\cal A}-$coordinates. 
Birational isomorphisms  (\ref{AMu}) preserve them:
$$
\mu^*\Omega_{\Gamma_2} = \Omega_{\Gamma_1}.
$$
 This  feature of non-commutative cluster ${\cal A}-$varieties is   quite non-trivial.  
 Its proof  Section \ref{sec9} leads to a  version of the Steinberg relation for non-commutative 2-forms, see Theorem \ref{T11.4}. 
 In the commutative case, we recover  the   2-form associated with cluster algebras/cluster ${\cal A}-$varieties \cite{GSV}, \cite{FG2}. \vskip 2mm

   Just as in the commutative case, there is a canonical map of non-commutative cluster varieties 
$$
  p: {\cal A} \lra {\cal X}.
$$
  It has Hamiltonian properties similar to the ones in the commutative case, where it is the projection along the null-foliation of the canonical 2-form 
  on $\Omega$ on ${\cal A}$. \vskip 2mm
  
We stress that in sharp contrast with the commutative case,   non-commutative cluster Poisson varieties do not carry any   cluster coordinates. 
  Therefore cluster Poisson considerations require new language. 
  
  On the other hand, the non-commutative cluster ${\cal A}-$varieties  come with a canonical collection of functions. Although we refer to them  as    
  cluster ${\cal A}-$coordinates, unlike in the commutative case,  cluster coordinates in each cluster are not independent; they  satisfy certain monomial relations.

\paragraph{6. Spectral description of the non-commutative   \underline{decorated}  local systems on $\bS$.}  
In the commutative case there is the dual moduli space ${\cal A}_{{\rm SL_m}, \bS}$ parametrised 
twisted ${\rm SL_m}$-local systems on $\bS$ equipped with a {\it decoration}: a choice of an invariant decorated flag 
near every marked point of $\bS$. It has a cluster  
 ${\cal A}$-variety structure \cite{FG1}, which is a
 geometric incarnation  of the Fomin-Zelevinsky upper cluster algebras   \cite{FZI}. 
There are cluster ${\cal A}$-coordinate systems on ${\cal A}_{{\rm SL_m}, \bS}$ parametrised by 
ideal triangulations   \cite{FG1} or, more generally, ${\rm SL_m}$-graphs   \cite{G}. 

The noncommutative analog   is the moduli space ${\cal A}_{m, \bS}$ parametrised twisted 
decorated rank $m$ $R$-local systems on $\bS$.  
Note that there is no non-commutative analog of ${\rm SL_m}$. 
Its  generic points    parametrise  twisted 
 flat line bundles on the spectral surface $\Sigma$ which  are  
  trivialized at the boundary   of ${\Sigma}$. 
The trivilalisations allow   to construct rational functions $\Delta_E$ on the moduli space ${\cal A}_{m, \bS}$ 
parametrised by the edges $E$ of the graph $\Gamma$. 
These functions are analogs of the cluster ${\cal A}$-coordinates. 
As we stressed above, the functions $\Delta_E$   
do not form a coordinate system: they satisfy non-trivial monomial relations. There is 
 no natural way to select a subset of the generators.  A similar phenomenon occurs already for the commutative moduli space ${\cal A}_{{\rm GL_m}, \bS}$ \cite{G}.

\paragraph{7. Berenstein-Retakh non-commutative cluster algebras related to surfaces. } In Section \ref{sec4} we show that in  the   case of two dimensional twisted decorated non-commutative local systems 
the obtained coordinate description recovers the non-commutative cluster functions 
decribed as  ratios of $2\times 2$ 
quasideterminants by A. Berenstein and V. Retakh in \cite{BR}. 
In our approach   formulas for these functions were obtained by geometric considerations, 
which made their properties evident. 

\vskip 2mm
In general   the functions $\Delta_E$ on the moduli spaces ${\cal A}_{m,  \bS}$ are ratios of the Gelfand-Retakh quasideterminants \cite{GR}. 
This way we get a geometric approach to  quasideterminants.

\subsection{Non-commutative cluster varieties and stacks of admissible dg-sheaves} \la{SEC1.2}

Given  a bipartite ribbon graph $\Gamma$,  we take  stacks 
${\rm Loc}_1(\Gamma')$, where   graphs $\Gamma'$  are obtained from $\Gamma$ by any sequences of two by two moves, and glue them    
  via  birational isomorphisms (\ref{11.9.11.2}) into a space  ${\cal X}_{|\Gamma|}$.   

Then a  question arises: can we assign to 
each bipartite ribbon graph $\Gamma$ a finite stack   ${\cal M}^\circ_{|\Gamma|}$,  whose equivalence class   is invariant under   
 two by two moves, and which contains   ${\rm Loc}_1(\Gamma)$ as a Zariski open  substack? 
Equivalently, we want to have a Zariski open embedding  of   ${\cal X}_{|\Gamma|}$ into an an priori defined stack: 
\be \la{ADS}
{\cal X}_{|\Gamma|} \subset {\cal M}^\circ_{|\Gamma|}.
\ee
To approach this question, we introduce stacks of admissible dg-sheaves on manifolds.    This allows to answer the  question in many interesting situations, discussed below, although not   
in full generality yet.

\paragraph{1. Admissible dg-sheaves.} 

Given a cell complex ${\cal C}$, we define in Section \ref{SEC7.1}  a dg-category of {\it dg-sheaves}, a dg-version of the derived category of complexes of  
sheaves constructible for the stratification of  ${\cal C}$. 

 Recall that the microlocal support of a constructible complex of sheaves on a manifold $X$ is a  complex of sheaves 
on a conical Lagrangian in  the cotangent bundle $ {\rm T}^*X$ \cite{KS}.   
 A cooriented hypersurface in a manifold $X$ gives rise to a Lagrangian in ${\rm T}^*X$, and a Legendrian in its spherical bundle  
$
{\rm ST}^*X$. 
 They are given, respectively, by the    conormal vectors/rays in the direction specified by the coorientation. \vskip 2mm

 Let ${\cal H}$ be a  collection   of  \underline{smooth} cooriented hypersurfaces  in a  manifold  $X$ with   \underline{disjoint} Legendrians.

We consider a subcategory of 
  ${\cal H}$-{\it supported} dg-sheaves. These are dg-sheaves whose   microlocal support   is given by complexes supported  at   the zero section    and 
   Lagrangians assigned to the hypersurfaces of ${\cal H}$. \vskip 2mm

Then we consider  a dg-subcategory of   {\it ${\cal H}$-admissible} dg-sheaves.  
Its key features  are the following.

 1. Admissible dg-sheaves on a manifold $X$ are concentrated in the degrees $[0,1]$.  
 
 2. The microlocal support of an {${\cal H}$-admissible} dg-sheaf is supported on the union of the zero section of ${\rm T}^*X$  and  
 a   Lagrangian  ${\rm T}^*_{\cal H}X$,   given by the union of conormal bundles to the cooriented hypersurfaces of ${\cal H}$. 
  For each component ${\cal C}$ of ${\cal H}$, the  microlocalisation along  ${\rm T}^*_{{\cal C}}X$can be identified with a local system on ${\cal C}$.

  3. Stacks of {  ${\cal H}$-admissible} dg-sheaves are invariant under   deformations of   ${\cal H}$ on $X$, 
 called below {\it admissible isotopies},  which  keep  each 
 component  of  ${\cal H}$ smooth, 
 and their  Legendrians disjoint.  \vskip 2mm
 
 We demonstrate property 3) directly in Section \ref{Sec3} in the dimension 2.   We stess   that   we allow in the process of deformation non-transversal intersections, e.g. three lines intersecting in a point.  
 Note   that a smooth Hamiltonian isotopy of  a Legendrian curve can develop a cusp, which we do not allow. \vskip 2mm

Note that according to \cite{GKS},    the derived category  of all constructible sheaves with a given   microlocal support    is invariant under 
   Legendrian isotopies. However this  does not imply that admissible isotopies preserve ${\cal H}$-admissibility, and   {not all} Legendrian isotopies 
 preserve  ${\cal H}$-admissibility. \vskip 2mm

    Since  local systems on the conormal bundle to a cooriented hypersurface are pull backs of   local systems on the hypersurface, 
we  abuse terminology by referring to them  as local systems on the hypersurface. 
Therefore below we refer to the microlocal support on ${\rm T}^*_{\cal H}X$  as the {\it Legendrian microlocal support}.     \vskip 2mm

Since ${\cal H}$-admissible dg-sheaves are concentrated in   degrees $[0, 1]$,  not all complexes of constructible sheaves with   Legendrian 
microlocal support  in  ${\rm T}^*_{\cal H}X$ 
 are  ${\cal H}$-admissible.  Stacks of ${\cal H}-$admissible dg-sheaves are of finite type, while stacks of complexes of constructible sheaves with the Legendrian 
microlocal support  in  ${\rm T}^*_{\cal H}X$ are not:  the latter contain   local systems on $X$ placed  in arbitrary degrees. \vskip 2mm

\begin{figure}[ht]
\centerline{\epsfbox{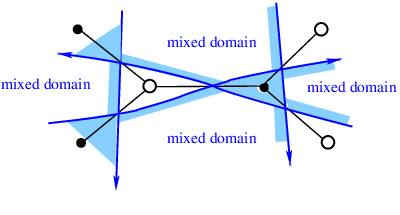}}
\caption{Zig-zag paths for a bipartite graph $\Gamma$ on a surface $\bS$ cut  the surface   into $\bullet$, $\circ$, and mixed domains. 
Each zig-zag path is coorientated   by a shaded area to the right of the path.}
\label{ncls104a}
\end{figure}

Let   $\Gamma$ be a bipartite ribbon graph. Denote by      ${\cal Z}_\Gamma$ the collection   
    of oriented zig-zags on the decorated surface  $\bS_\Gamma$ associated to $\Gamma$, well defined up to isotopy. Since the surface  $\bS_\Gamma$ is oriented,  zig-zags are  cooriented. 
The complement  $\bS_\Gamma - {\cal Z}_\Gamma$ is a union of domains of three types,  depending on the orientation types of their boundary segments, see Figure \ref{ncls104a} $\&$  Section \ref{SSec4}.

  \bl \la{LE2}
  The stack ${\rm Loc}_1(\Gamma)$ of   flat  line bundles   on  a bipartite ribbon graph  $\Gamma$ is canonically equivalent to     
 the stack   of ${\cal Z}_\Gamma$-admissible  dg-sheaves  on $\bS_\Gamma$ satisfying the following   conditions: 
 
 1. They are acyclic on mixed domains - note that  this is an open condition.\footnote{The condition that  a complex   does not have cohomology in two neighboring degrees is an open condition. In particular the conditions that a complex is acyclic,   or has  a single non-trivial cohomology group, are open conditions. The latter are precisely the conditions we impose in Lemma \ref{LE2}.}

 2. Their Legendrian 
   microlocal support is given by 
 flat $R-$line bundles at every curve of 
${\cal Z}_\Gamma$.
  \el
Our construction of the functor from the sheaves on the surface to flat line bundles on the graph $\Gamma$ uses dg-sheaves: the connection on the line bundle on   $\Gamma$ is  given by "homotopies between morphisms", 
absent in the classical description of the derived category. In the commutative setting, another construction of the functor follows from \cite{STWZ}.\footnote{Namely, as was pointed out  by the referee, in \cite{STWZ} a Lagrangian $L$ with the topology of the spectral surface  of $\Gamma$ is constructed; it fills the Legendrian of  the zig-zag curves. Let $\widetilde L$ be its Legendrian lift to $J^1(S) \subset T^*S \times R^1$. By Guillermou's  quantization theorem \cite{G12}, \cite{GKS}, ${\rm Loc}_1(L)$ can be identified with the sheaves which are zero at $S \times -\infty$, and pure of microlocal rank one along the conormal to $\widetilde L$. Restricting the sheaves to $S \times \infty$ we get the  functor.}\vskip 2mm

  Denote by ${\cal M}({\cal H}, {\bf 1})$ the stack of ${\cal H}-$admissible  dg-sheaves  on a surface satisfying  property 2) of Lemma \ref{LE2} for every curve of ${\cal H}$.   It is a finite type stack. Then By Lemma \ref{LE2} there is  an embedding
 \be
 {\rm Loc}_1(\Gamma)\subset {\cal M}({\cal Z}_\Gamma, {\bf 1}).
  \ee
  
 We define   ${\cal M}^\circ_{|\Gamma|}$ as the component 
of   ${\cal M}({\cal Z}_\Gamma, {\bf 1})$ containing ${\rm Loc}_1(\Gamma)$. We   show that in many interesting cases we   get this way   Zarisky open embeddings (\ref{ADS}).

  \paragraph{2. Examples of stacks of ${\cal H}-$admissible dg-sheaves.} 
  
  They   include  the   stacks of    Stokes data, which are the Betti side of  stacks of    flat connections with   irregular singularities on Riemann surfaces, see Sections    \ref{sec9.1} - \ref{SECT9.2}.
This way we  get   non-commutative     stacks of Stokes data stacks as well. \vskip 2mm

Let us explain   how to get the Betti stack of framed flat connections   in the   case of regular singularities.

\begin{figure}[ht]
\centerline{\epsfbox{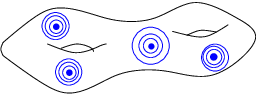}}
\caption{A collection of curves ${\cal I}_{\bS,  {3}}$ on a genus two   surface $\bS$ with four punctures.}
\label{ds}
\end{figure}

  For each marked point $s$ of $\bS$,   in a small punctured at    $s$ disc/half disc,      consider  $m$ 
    non intersecting and non selfintersecting   circles/arcs ending on the boundary,  cooriented out of $s$. Denote by ${\cal I}_{\bS,  {m}}$   their union, see Figure \ref{ds}.   
    Its isotopy class is determined uniquely by the surface $\bS$ and the integer $m$. 
  
The following observation is an immediate corollary of the definition of ${\cal I}_{\bS,  {m}}-$admissible dg-sheaves. 

 \bl \la{LE1} The  stack ${\cal X}_{\bS, m}$ of  $m$-dimensional  $R$-local systems on a decorated surface $\bS$   with 
  complete filtrations near     marked points   is canonically equivalent to the stack of   ${\cal I}_{\bS,  {m}}-$admissible  dg-sheaves on $\bS$ satisfying the following two conditions:  
  
  1.    They are acyclic near the marked points. 
   
   2.  Their microlocal support is 
 \mbox{ the union of the zero section and    flat $R-$line bundles at every curve of 
  ${\cal I}_{\bS,  {m}}$.}
  \el

\begin{figure}[t]
\centerline{\epsfbox{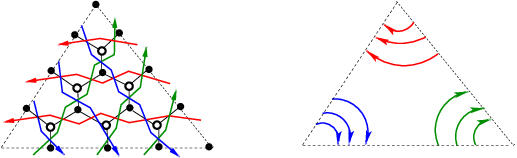}}
\caption{Zig-zag paths on the bipartite graph $\Gamma_3$ on the left can be admissibly deformed to   concentric arcs around the vertices, shown  on the right.}
\label{stokes20}
\end{figure}

  Lemmas \ref{LE2} - \ref{LE1} provide    a non-commutative cluster Poisson structure on the stack ${\cal X}_{\bS, m}$.  Indeed, 
 by Lemma \ref{LE1},    ${\cal X}_{\bS, m}$   is   realized as  the  substack ${\cal M}_0({\cal I}_{\bS, m}, {\bf 1})\subset {\cal M}({\cal I}_{\bS, m}, {\bf 1})$  
given by   condition 1):  
$$
{\cal X}_{\bS, m} = {\cal M}_0({\cal I}_{\bS, m}, {\bf 1}) \subset {\cal M}({\cal I}_{\bS, m}, {\bf 1}).
$$  
Consider a bipartite ${\rm GL_m}-$graph  $\Gamma$ on $\bS$.   By its very definition \cite{G}, there  is an admissible isotopy of 
    the  collection ${\cal Z}_\Gamma$ of its zig-zags 
  to the collection  of   curves ${\cal I}_{\bS, m}$, see Figure \ref{stokes20}: 
\be \la{ISOi}
i_\Gamma: {\cal Z}_\Gamma \lra {\cal I}_{\bS, m}.
\ee 
By  property 3),  isotopy (\ref{ISOi})   provides an equivalence of the related  stacks of ${\cal H}$-admissible dg-sheaves
$$
i_{\Gamma*}: {\cal M}({\cal Z}_\Gamma, {\bf 1}) \stackrel{\sim}{\lra} {\cal M}({\cal I}_{\bS, m}, {\bf 1}).
$$
By Lemma \ref{LE2},  the stack ${\rm Loc}_1(\Gamma)$ is realized as  the Zariski open substack ${\cal M}_{{\rm mix}=0}({\cal Z}_\Gamma, {\bf 1})$ of 
  ${\cal M}({\cal Z}_\Gamma, {\bf 1})$ defined by   condition 1):
$$
{\rm Loc}_1(\Gamma) = {\cal M}_{{\rm mix}=0}({\cal Z}_\Gamma, {\bf 1}) \subset {\cal M}({\cal Z}_\Gamma, {\bf 1}).
$$
 Indeed, the acyclicity on mixed domains   in Lemma    \ref{LE2} is an open condition. Acyclicity on mixed domains implies  acyclicity near  marked points since marked points are contained in the mixed domains, see Figure \ref{stokes20}  when $\bS$ is a triangle. The general case reduces to this since we can get a ${\rm GL}_m-$graph on $\bS$ by amalgamating the graphs $\Gamma_m$ on triangles of an ideal triangulation of $\bS$. 
 Therefore   $i_{\Gamma*}$ induces an open embedding 
 $j_\Gamma: {\cal M}_{{\rm mix}=0}({\cal Z}_\Gamma, {\bf 1}) \subset {\cal M}_0({\cal I}_{\bS, m}, {\bf 1})$. Summarising, we get  the commutative diagram
\be
\begin{gathered}
    \xymatrix{
  {\rm Loc}_1(\Gamma) = {\cal M}_{{\rm mix}=0}({\cal Z}_\Gamma, {\bf 1})   \ar[d]^{\cap}_{j_\Gamma}        \ar[r]^{~~~~~~~~~~~\subset}& {\cal M}({\cal Z}_\Gamma, {\bf 1}) \ar[d]_{i_{\Gamma*}}^{\sim}  \\
   {\cal X}_{\bS, m} = {\cal M}_0({\cal I}_{\bS, m}, {\bf 1})   \ar[r]^{~~~~~~~\subset}&  {\cal M}({\cal I}_{\bS, m}, {\bf 1}).}
\end{gathered}
 \ee
In particular,    for any  ${\rm GL}_m-$graph $\Gamma$ on $\bS$, we have,  as was promised in (\ref{ADS}), 
 \be
 {\cal M}^0_{|\Gamma|} = {\cal X}_{\bS, m}. 
 \ee
 For example,  the bipartite graph $\Gamma_m$ on Figure \ref{stokes20} allows to recover the whole stack of all triples of flags as the stack 
${\cal M}^0_{|\Gamma_m|}$, 
rather than describing its generic part, as 
Theorem \ref{MT1} does.  

The stack ${\rm Loc}_1(\Gamma)$ is a non-commutative Poisson torus. So  
  open embeddings $j_\Gamma: {\rm Loc}_1(\Gamma)\subset {\cal X}_{\bS, m}$ for   bipartite  ${\rm GL_m}-$graphs  $\Gamma$ on $\bS$ define a non-commutative cluster Poisson   structure on the stack  ${\cal X}_{\bS, m}$. 
The action of the mapping class group $\Gamma_\bS$ of $\bS$ preserves the class of  ${\rm GL_m}-$graphs on $\bS$. Any two of them  are 
 related by a sequence of two by two moves \cite{G}, 
providing admissible isotopies of  their  zig-zags. 
Therefore  we get a $\Gamma_\bS-$equivariant   non-commutative cluster Poisson structure on the stack ${\cal X}_{\bS, m}$.\\



 In   Section \ref{Sect9.5} we describe   stacks of framed Stokes data via  admissible dg-sheaves, and show that  
 \be \la{1111}
\mbox{ \it Stacks of (framed) Stokes data are shadows of non-commutative (framed) stacks of Stokes data:}
\ee  
\bt \la{SSDNCa}
 For any stack of framed Stokes data ${\cal S}$ on a decorated surface   $\bS$, there is a non-commutative stack ${\cal S}^{\rm nc}$  over ${\rm Spec}(\Z)$, whose restriction to commutative fields is the  stack ${\cal S}$.  
Assume that  $\bS \not = S^2-\{\infty\}$. Then the      stack  ${\cal S}^{\rm nc}$ carries a non-commutative cluster Poisson structure equivariant under the wild mapping class group. 
\et

The argument is based on the observation that, similarly to Lemma \ref{LE1},   stacks of Stokes data can be
 realized as stacks of ${\cal L}-$admissible dg-sheaves on $\bS$ vanishing near the punctures,  
where ${\cal L}$ is the union of certain local collection of curves ${\cal L}_p$ around the punctures $p$. 

 When $\bS=S^2-\{0, \infty\}$, with     one regular and one irregular point, 
these stacks are  identified birationally with   cluster Poisson varieties assigned in \cite{FG3}  to the cyclic closures of  elements of the braid semigroup. Thanks to    \cite[Section 3.10]{FG3}, 
their quivers are described 
by bipartite graphs. The general case 
is   reduced to this one plus the regular triangle case. 
The case when  $\bS  = S^2-\{\infty\}$ requires a different technology, and is not considered in the paper. Yet the result is true in  this case as well.

\paragraph{3. Stacks of framed $\G-$Stokes data,   their cluster Poisson structure $\&$ quantization.} 
Stacks of $\G-$Stokes data 
 provide the   
Betti realization of  stacks of  $\G-$bundles with meromorphic connections with arbitrary  singularities on a Riemann surface $\Sigma$. 
We introduce stacks of {\it framed} $\G-$Stokes data, and 
  prove that they carry  a cluster Poisson structure, equivariant under  the wild mapping class group action.

\bt \la{SSDNC}
Assume that $\bS \not = S^2-\{\infty\}$.  For any split reductive group $\G$ over ${\Q}$ with connected center, the stack of \underline{framed} $\G-$Stokes data on $\bS$ has a cluster Poisson structure, equivariant under the action of the wild mapping class group. 
Stacks of \underline{plain} $\G-$Stokes data carry   compatible Poisson structures. 
 \et

 The Poisson nature of the moduli spaces of Stokes data   was studied by P. Boalch \cite{Bo11}.

 Theorems \ref{SSDNC} - \ref{TTSD}  are valid for a much larger class of {\it ideal} Legendrian links ${\cal L}$, see Section \ref{SECT9.2}.   \vskip 2mm
 
 Theorem  \ref{SSDNC}  
is proved in   Section \ref{Sect9.6} using the following  key ingredients.   \vskip 1mm

1. If $\bS=S^2-\{0, \infty\}$, with   one regular and one irregular point,  
the  stacks of $\G-$Stokes data are  identified birationally with the  cluster Poisson varieties assigned in \cite{FG3}  to the cyclic closures of certain elements of the braid semigroup ${\rm Br}^+_\G$ related to $\G$. 

2. When $\bS=t$ is a triangle, the cluster Poisson structure of the related space ${\cal P}_{\G, t}$ of configurations of triples of flags with pinnings was described in \cite{GS19}.  The key point is its  cyclic invariance. 

3. If   
$\bS = S^2-\{\infty\}$,   the stacks of Stokes data are  given by certain  cyclically ordered configurations of flags.
   The cluster Poisson structure for   sufficiently generic cyclic configurations    follows from    \cite{GS19}.  \vskip 1mm
 
 The general case 
is   built from these ones,  except  for   single irregular point -  that is $\bS= S^2-\{\infty\}$ -  of rather degenerate type,  which requires different technique, and is omitted.\\

Combining Theorem \ref{SSDNC} with the cluster quantization machine \cite{FG4},  we get the following.
  
  \bt \la{TTSD} Under the same assumptions  as in Theorem \ref{SSDNC}, the   stacks of $\G-$Stokes data, plain or framed,  admit a  cluster Poisson quantization, 
  equivariant under the  wild mapping class group action. 
\et

Precisely, we apply   the cluster quantization to the cluster Poisson space  of  {framed} $\G-$Stokes data. Then we use the fact that forgetting the framing we get a map 
to the space of  {plane} $\G-$Stokes data, which is  finite Galois map over a generic point. Its Galois group, given by a product of the Weyl groups of $\G$,   acts by cluster Poisson transformations. This allows us to   pass to its invariants.

 Theorem \ref{TTSD} completes the program of cluster quantization of the moduli spaces of flat connections on Riemann surfaces.  
 In the regular case this was done in \cite{FG1} for $\G = {\rm PGL}_m$, and in \cite{GS19} in general.  \vskip 2mm

 \paragraph{4. Cluster ${\cal A}-$varieties via stacks of admissible dg-sheaves.} The  analogs of  cluster ${\cal A}-$varieties in this set-up     are given by   stacks of  (twisted) 
 ${\cal H}-$admissible  dg-sheaves with  the following  extra condition: 
 \be
 \mbox{\it Local systems on all   components of the microlocal support but the zero secion are trivialized.}
 \ee

In particular, given a bipartite ribbon graph $\Gamma$ we recover non-commutative cluster ${\cal A}-$varieties discussed in Section \ref{SECT1.1}. Indeed, 
  zig-zags on the   surface $\bS_\Gamma$ match   boundary curves on the spectral surface $\Sigma_\Gamma$. So the trivialization of the local systems on zig-zags on $\bS_\Gamma$, 
  describing   the microlocal support 
  of an admissible dg-sheaf,   amounts to the trivialization of the local system on the spectral surface $\Sigma_\Gamma$ on the boundary  of    $\Sigma_\Gamma$. 
  The cluster ${\cal A}-$coordinates are assigned to the edges of $\Gamma$.
  They are  
obtained  by comparing   trivialisations  of the local systems on the two boundary components by using    
   the parallel transport along the canonical up to isotopy path connecting two boundary components.

 \paragraph{Remark.} The stacks of Stokes data are the geometric analogs of (families of)     
 wildly ramified  representations of the absolute Galois group, see   \cite{DMR}. So Theorem \ref{TTSD} should be useful in various flavors of the geometric  quantum Langlands correspondence.

\paragraph{5. ${\rm CY}_3$-categories assigned to stacks of Stokes data.}   Description of   the   cluster Poisson structure on stacks  of Stokes data via 
webs of zig-zags for  certain bipartite graphs on the surface $S$ 
has another  important application. 
 Recall that any bipartite graph gives rise to a canonical quiver with potential, which determines a combinatorially defined 3d CY $A_\infty-$category. 
The categories  assigned to bipartite graphs related by   two by two moves are equivalant. 
Therefore we arrive to a combinatorially defined  3d Calabi-Yau $A_\infty-$category  assigned to  any   topological Stokes data $({\cal L}, d)$. This category  should describe  
the Fukaya category  
of the family of open complex CY threefolds  assigned to the same    data $({\cal L}, d)$. 
This story in the regular case is described  in \cite{G}.  The $m=2$ case is established in \cite{BS}, \cite{S}.

\paragraph{6. Comments.} Sections \ref{SECCT7} and \ref{Sect9.5} are very closely related to the  series of  works of Shende-Treumann-Williams-Zaslow \cite{STZ14}, \cite{STWZ},  \cite{STZ16}, where they studied  extensively, in the commutative setting,  moduli spaces of  rank one microlocal sheaves, including the ones related to  bipartite surface graphs. 
They usually work with   all constructible complexes of sheaves with given microlocal  support  rather than with  finite type  stacks of  admissible ones. 
  The restrictions on constructible sheaf categories to cut down to finite type stacks appear in 
\cite[Proposition 5.20]{STZ14}, \cite[Propoisition 5.20]{STWZ}. 
The admissibility condition appears in \cite[Proposition 5.3]{STWZ}.\footnote{We are grateful to the referee for the last two comments.}  
  We use the very existence of the category of sheaves with given microlocal support \cite{KS}. Other than that, our exposition is selfcontained and very  elementary. 

 Our work also relates to some aspects of the   Gao-Shen-Weng  program  \cite{GSW}. In particular, non-commutative Plucker relations appear in their approach as well.

\paragraph{7. Acknowledgments.} We are grateful to A. Berenstein and V. Retakh for useful discussions during the Summer 2011 in IHES, when   the project  was started, and all results of Section \ref{SECT1.1a} were obtained. The hospitality of IHES was crucial  for the developing of the project. 
We are very grateful to the referee, who read the paper very carefully, 
and made many comments, incorporated in the text, which improved and clarified the exposition.

A.G. was 
supported by the  NSF grants DMS-1059129, DMS-1564385, DMS-1900743,  DMS-2153059. He was a Fellow of the Simons Foundation in 2019 and 2023. 
He is grateful to the NSF, Simons Foundation, and  IHES  for   hospitality and support.

 \section{Non-commutative cluster Poisson varieties} \la{secziya2}

\subsection{Two by two moves for non-commutative flat line bundles} \la{SEC5.1}

 Definition \ref{AFLBG} for the action of   two by two moves on flat line bundles  is dictated by the example provided by configurations of four lines in a two 
 dimensional space given  in Section \ref{SECT1.1a}. \vskip 2mm
 
Recall that a  ribbon graph  is a graph with  a cyclic order 
of the edges at every vertex. We allow ribbon graphs with external, that is 1-valent,  vertices. A ribbon graph $\Gamma$ is the same thing as a graph on a decorated  
 surface $\bS$ which is homotopy equivalent to the graph.  
 
A {\it face path} on    $\Gamma$ is a path 
turning   right at every vertex. There is a bijection  
\be \la{ZZFPa}
\{\mbox{face paths on $\Gamma$}\}~ \longleftrightarrow ~\{\mbox{marked points on  $\bS$}\}. 
\ee
 Face intervals correspond to  special points on  $\bS$. Face loops correspond to punctures of $\bS$.

A {\it zig-zag path} on   $\Gamma$ is a path 
turning at every vertex either left or right, so that the left and right turns alternate.  
A zig-zag path   on a bipartite ribbon graph $\Gamma$ is oriented so that 
it turns right at the $\circ$-vertices, and turns left at the $\bullet$-vertices, see Figure \ref{ncls1a}.

 \paragraph{1. Two by two moves.} 

\begin{figure}[ht]
\centerline{\epsfbox{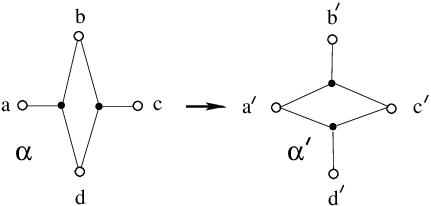}}
\caption{A two by two move $\alpha \to \alpha'$. 
}
\label{ncls6a}
\end{figure}

\bd \la{DEF6.1} Let $\Gamma$ be a bipartite ribbon graph containing a subgraph $\alpha$ 
  on the left of Figure \ref{ncls6a}. 
A {\em two by two move $\mu_\alpha:\Gamma \to \Gamma'$ at $\alpha$ }
produces a new bipartite graph $\Gamma'$, 
obtained by replacing  $\alpha$  by  
the subgraph $\alpha'$ on the right  of Figure \ref{ncls6a}. 
\ed

Given a two by two move $\mu_\alpha:\Gamma \to \Gamma'$, let us define 
a birational isomorphism of the moduli spaces of non-commutative flat line bundles on the graphs: 
\be \la{birat}
\mu_\alpha: {\rm Loc}_1(\Gamma) \lra{\rm Loc}_1(\Gamma').   
\ee

\begin{figure}[ht]
\centerline{\epsfbox{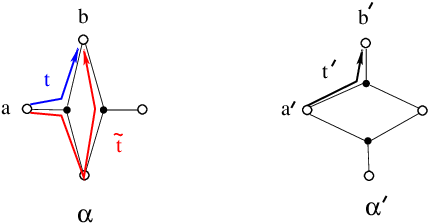}}
\caption{The two by two move action on the parallel transform: $t':= t + \widetilde t = (1+M_b)t$. }
\label{ncls7aa}
\end{figure}

Consider a two by two move at  $\alpha$ shown   on 
Figure \ref{ncls6a}:
$$
\mu_\alpha:\Gamma \to \Gamma'. 
$$
Let $\alpha_0$ be the graph $\alpha$ without the 
$\circ$-vertices, and similarly $\alpha'_0$. 
By the  
definition,   
 $\Gamma - \alpha_0 = \Gamma' - \alpha'_0$. Denote by $ b\bullet a$ the arc $t$ on the graph $\alpha$ 
going from the 
$\circ$-vertex $a$ to the  $\bullet$-vertex, and then to the 
$\circ$-vertex $b$ - see  Figure \ref{ncls7aa}.  We use a similar notation for the graph $\alpha'$.

\begin{figure}[ht]
\centerline{\epsfbox{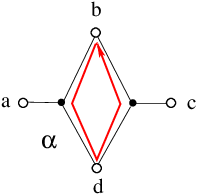}}
\caption{$M_b$ is the negative of the counterclockwise monodromy around the loop $\alpha$.}
\label{ncls7x}
\end{figure}

 Starting from a flat line bundle $L$  on $\Gamma$, we produce 
a new flat line bundle $L'$ on $\Gamma'$ as follows:

Let $M_v:L_v \to L_v$ be the \underline{minus} of the operator of  counterclockwise  monodromy  
around the simple loop on   $L_\alpha$, acting 
at  the fiber at 
 $v$, see  Figure \ref{ncls7x}. 

\begin{itemize}
\item The restriction of   $L'$ to $\Gamma' - \alpha'_0$ is defined as  
the restriction of  $L$ to $\Gamma - \alpha_0$.  

\item The restriction of $L'$ to $\alpha'$, denoted by $L_{\alpha'}$,
  is defined in Definition \ref{AFLBG}.  
\end{itemize}

\bd \la{AFLBG}  Given a flat line bundle $L_\alpha$ on the graph $\alpha$, we define a  flat line bundle    
$L_{\alpha'}$ on the graph $\alpha'$ by setting   
the parallel transports along the 
arcs  
between the $\circ$-vertices on $\alpha'$ to be:
\be \la{10.12.11.1b}
\begin{split}
&t_{b' \bullet a'} := (1+M_b)t_{b\bullet a}, \quad  t_{c' \bullet b'} := t_{c\bullet  b}(1+M_{b})^{-1}, \\
&t_{d'\bullet c'} := (1+M_d)t_{d\bullet c}, \quad  t_{a'\bullet d'} := t_{a\bullet d}(1+M_{d})^{-1}. \\
\end{split}
\ee

We glue the flat line bundle $L_{\alpha'}$ 
to the restriction of the flat line bundle $L$  to $\Gamma - \alpha_0$. 
\ed

Let $\widetilde t_{b, a}$ be ``the other way around $\alpha$''  parallel transport from $a$ 
to $b$, see Figure \ref{ncls7aa}. Then 
\be \la{IDENT}
t_{b'\bullet a'}\stackrel{(\ref{10.12.11.1b})}{=} (1+M_b)  t_{b\bullet a} = t_{b\bullet a}(1+M_a) = 
t_{b\bullet a} + \widetilde t_{b, a}. 
\ee
It is a sum of the two different simple paths 
from $a$ to $b$ on $\alpha$.

\bl \la{MONL}
Monodromies of the local systems $L_\alpha$ and $L_{\alpha'}$ at each of the four 
$\circ$-points  
coinside: 
$$
M_a=M'_a, ~~~~M_b =M'_b, ~~~~M_c=M'_c, ~~~~M_d =M'_d.
$$  
\el

\begin{proof} 

We have
\be \begin{split}
&M_a' = t_{a'\bullet d'}  t_{d'\bullet c'}  t_{c' \bullet b'}  t_{b' \bullet a'} = \\
&t_{a\bullet d}(1+M_{d})^{-1} (1+M_d)t_{d\bullet c}  t_{c\bullet  b}(1+M_{b})^{-1} (1+M_b)t_{b\bullet a} = \\
&t_{a\bullet d}t_{d\bullet c}  t_{c\bullet  b}t_{b\bullet a} = M_a.\\
\end{split}
\ee
Next, using $M_a=M_a'$ we have 
\be \la{SID}
M_b' =  t_{b' \bullet a'} M_a' 
t^{-1}_{b' \bullet a'} \stackrel{(\ref{10.12.11.1b})}{=}  t_{b \bullet a} (1+M_a) M_at^{-1}_{b' \bullet a'}= t_{b \bullet a} M_a(1+M_a) t^{-1}_{b' \bullet a'}.  
\ee
The first two equalities in (\ref{IDENT})  imply:
$$
t^{-1}_{b' \bullet a'}  \stackrel{(\ref{IDENT})}{=}  (t_{b \bullet a}(1+M_a))^{-1} = (1+M_a)^{-1}t^{-1}_{b \bullet a}.
$$
Substituting this to (\ref{SID}), we get:
$$
M_{b'}\stackrel{(\ref{SID})}{=}    t_{b \bullet a} M_a(1+M_a)(1+M_a)^{-1} t^{-1}_{b \bullet a} 
= t_{b \bullet a} M_a t^{-1}_{b \bullet a} = M_b. 
$$
The last two equalities   follow   
from the cyclic shift by two symmetry of the picture. 
\end{proof}

\bc
i) We can rewrite Definition \ref{AFLBG}  
as follows:
\be \la{10.12.11.1bx}
\begin{split}
&t_{b' \bullet a'} := (1 + M_b)t_{b\bullet a}, \quad  (1 +M_{c'}) t_{c' \bullet b'}:= t_{c\bullet  b}, \\
&t_{d'\bullet c'} := (1 + M_d)t_{d\bullet c}, \quad  (1 + M_{a'})t_{a'\bullet d'} := t_{a\bullet d}. \\
\end{split}
\ee

ii) One has $\mu_{\alpha'} \circ \mu_\alpha = {\rm Id}$. 
\ec

\begin{proof} i) Indeed, using Lemma \ref{MONL}, 
$$
t_{c'\bullet  b'} = t_{c\bullet  b}(1+M_{b})^{-1} = t_{c\bullet  b}(1+M_{b'})^{-1} 
~~\implies~~t_{c'\bullet  b'}(1+M_{b'}) = t_{c\bullet  b}~~
\implies~~(1+M_{c'})t_{c'\bullet  b'} = t_{c\bullet  b}.
$$  

The identity $(1+M_{a'})t_{a'\bullet d'} := t_{a\bullet d}$ follows by the cyclic shift by two symmetry.

ii) Since the cyclic shift by one of $\alpha'$ is isomorphic to $\alpha$, 
 this implies $\mu_{\alpha'} \circ \mu_\alpha = {\rm Id}$. 
\end{proof}

\paragraph{2. Configurations of lines and two by two moves.}  
To  prove  that   transformations (\ref{birat}) satisfy the pentagon relations we need the following comparison result. 
 Our  construction from Section \ref{SECT1.1a}, see Theorems \ref{MT1} and  \ref{TH2.2},  provides a birational  equivalence of groupoids, see Figure \ref{ncls6a} for the notation:
\be \la{10.12.11.4az}
\begin{split}
&\{\mbox{Quadruples of lines 
in a $2$-dimensional $R$-vector space $V$,   at the $\circ$-vertices of $\alpha$}\} \\
&~~~~~~~~~~~~~~~~~~~~~~~\longleftrightarrow \{\mbox{Flat $R$-line bundles on the graph $\alpha$}\}. \\
\end{split}
\ee

\bp \la{10.12.11.5X}
The transformation in Definition \ref{AFLBG} is characterised by the condition that it commutes with   
equivalence (\ref{10.12.11.4az}), providing a commutative diagram  \be \la{12.9.11.1}
\begin{gathered}
    \xymatrix{
    &\{\mbox{Line bundles ${\cal L}$ on $\alpha$}\}   \ar[rd]_{(\ref{10.12.11.4az})}  &      {\xrightarrow{\makebox[2.5cm]{Definition \ref{AFLBG}}} } &\ar[ld]^{(\ref{10.12.11.4az})} \{\mbox{Line bundle ${\cal L}'$ on $\alpha'$}\}&\\
    &&  (A,B,C,D) \subset V  &&\\   }
\end{gathered}
 \ee
\ep

 \begin{proof}  
We start with a configuration of four generic lines $(A, B, C, D)$ 
in $V$. We assign them  to the  $\circ$-vertices $(a,b,c,d)$ of the 
graph $\alpha$, and to the  $\circ$-vertices $(a',b',c',d')$ of the graph $\alpha'$. 

We abuse notation by denoting by the same letter $a$ a vector in the line $A$, etc.

The two $\bullet$-vertices provide us five vectors in the lines $(A,B,C,D)$, 
defined uniquely up to a common left $R^*$-factor 
by the condition that the sum of the three vectors at the $\circ$-vertices incident to any $\bullet$-vertex is zero. 
These are the vectors illustrated on the left of Figure \ref{ncls6}:
$$
(a, b, b_1, c, d), \quad
~~\mbox{So $b_1= Xb, ~~X\in R^*$}.
$$
The counterclockwise monodromy around   $\alpha$ acts 
on any of them by the left multiplication by $X$.  
Indeed, 
\be \la{MONCOOR}
t_{a\bullet b}(b) = a, \quad 
t_{d\bullet a}( a) =d, \quad t_{c\bullet  d}(d) = c, \quad t_{b_1\bullet  c}( c) =b_1= Xb.
\ee
There are similar five vectors in $V$ related to the graph $\alpha'$, illustrated on the right in Figure \ref{ncls6}:
$$
(a', b', b_1', c', d'),   ~~
\mbox{So $b'_1= Yb', ~~Y\in R^*$}.
$$  
The counterclockwise monodromy around   $\alpha'$ acts  
on them as the left multiplication by $Y$. 

\begin{figure}[ht]
\centerline{\epsfbox{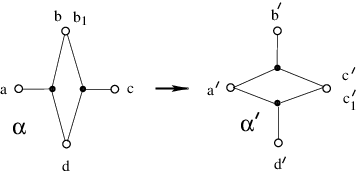}}
\caption{Calculating the two by two move for twisted configurations of lines in $V_2$.}
\label{ncls6}
\end{figure} 

One summarises definition of these vectors by the following relations, one per each $\bullet$-vertex:
\be \la{1.8.15.1}
a+b+d = 0, \quad d+c +b_1=0; \qquad a'+b'+c' = 0, \quad a'+d'+c_1'=0.
\ee
Notice that $b_1=Xb$ and $c'_1= Yc'$.   
So the equations read now as
$$
a+b+d = 0, \quad d+c+Xb=0; \qquad a'+b'+c' = 0, \quad a'+d'+Yc' = 0.
$$
These vectors are related as follows.

\bl \la{OMMO} One has 
\be
a'=a, ~~~~~~b'= (1-X)b, ~~~~~~ c'= -c,~~~~~~~~d'= (1-X^{-1})d.
\ee\el

\begin{proof}
Each pair of the vectors $(a', a), \ldots ,  (d',d)$ is collinear, so $a'=\lambda_aa$, $b'=\lambda_bb$, $c'=\lambda_cc$, $d'=\lambda_dd$.  
Let us solve these equations. First, let us  
express $c, d$ via $a,b$, and  respectively $c', d'$ via $a',b'$:
 $$
 d=-a-b, \quad c= (1-X)b + a, \quad c'= -a'-b', \quad d'= 
 - (1-Y) a'  + Yb'.
 $$
Next, the identities $\lambda_cc = c'$ and $\lambda_dd = d'$ give:
 $$
\lambda_c((1-X)b + a)= -\lambda_aa-\lambda_bb, \qquad 
\lambda_d(a+b)= (1-Y)a - Y\lambda_bb.
 $$
Set $\lambda_a=1$.  Then we get 
  the following equations:
$$
\lambda_b= 1- X,  ~~\lambda_c= -1,  ~~\lambda_d= 1-Y, ~~ -Y\lambda_b = \lambda_d.
$$
From the first and last equations we see that $Y=X^{-1}$, and therefore   we finally get 
$$
\lambda_a= 1, ~~\lambda_b= 1-X, ~~\lambda_c= -1,~~\lambda_d= 1-X^{-1}.
$$
\end{proof}

Lemma \ref{OMMO}  together with (\ref{MONCOOR}) implies
\be \la{10.12.11.1}
\begin{split}
&t_{b'\bullet a'} = (1-X)t_{b\bullet a}, \quad  t_{c'\bullet  b'} = (1-X^{-1})^{-1}t_{c\bullet  b_1}, \\
&t_{d'\bullet c'_1} = (1-X)t_{d\bullet c}, \quad  t_{a'\bullet  d'} = (1-X^{-1})^{-1}t_{a\bullet  d}. \\
\end{split}
\ee

Recall that $M_a$ is the {\it negative}  of the   counterclockwise  monodromy around $\alpha$. 
\bl \la{L5.6}
Formulas (\ref{10.12.11.1})    are equivalent to  formulas
 (\ref{10.12.11.1b}). 
\el

\begin{proof} 
The first line in   (\ref{10.12.11.1b}) is equivalent to the second one   
by  the cyclic shift by two argument. The    right formula in each line  in (\ref{10.12.11.1b}) 
is equivalent to the   left one, applied to the two by two move 
$\mu_{\alpha'}: \Gamma'\to \Gamma$. So we need to check only one of the four formulas in (\ref{10.12.11.1b}). 

 In the top left formula in (\ref{10.12.11.1}) the variable $X$, by its  definition,  is 
 the action of the counterclockwise monodromy along  $\alpha$ on the vector $b$. 
So it can be written  as $t_{b'\bullet a'} = (1+M_b)t_{b\bullet a}$, which is just the top  left formula in (\ref{10.12.11.1b}). 
  \end{proof}

Proposition \ref{10.12.11.5X} is proved \end{proof}

\paragraph{3. The pentagon relation.} 

\begin{figure}[ht]
\centerline{\epsfbox{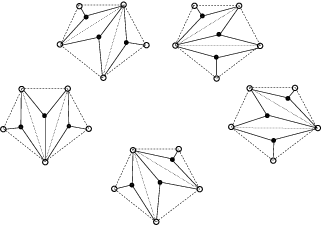}}
\caption{The pentagon relation. }
\label{ncls8}
\end{figure}

\bt The fifth power of 
 "two by two move followed by the non-trivial   symmetry of the triangulated pentagon" 
 for the 
two by two moves on Figure \ref{ncls8} acts as the identity on flat line bundles. 
\et 

\begin{proof} Let $\{\Gamma_i\}$, where $i \in \Z/5\Z$, be the sequence of bipartite graphs on Figure \ref{ncls8}. 
They are assigned to triangulations of the pentagon. The two by two move $\Gamma_i\to \Gamma_{i+1}$ corresponds to a
 flip of the triangulation. 
For each of the graphs $\Gamma_i$  there is an equivalence of groupoids
\be \la{10.12.11.40azaz}
\{\mbox{Generic flat line bundles on the graph $\Gamma_i$}\} \longleftrightarrow 
\ee
$$\{\mbox{Generic $5$-tuples of lines 
in a two dimensional space $V$,   assigned to $\circ$-vertices of a pentagon.}\}
$$ 
Thanks to Lemma \ref{L5.6}  
 the mutation of flat 
line bundles corresponding to a two by two move intertwines these equivalences. 
Indeed, a two by two move corresponding to a flip of a triangulation of the pentagon keeps
 intact one of the triangles of the triangulation, and does not affect
 the restriction of the flat line bundle to the graph inside of that triangle.  
\end{proof}

\subsection{Non-commutative cluster Poisson varieties} \la{SEC5S}


In this Section  we define a non-commutative analog of the Poisson algebra assigned to a bipartite ribbon graph   in \cite{GK}.  
The output  is still  a {\it commutative} Poisson algebra. We prove that two by two moves give rise to homomorphisms of Poisson algebras. 

\paragraph{1. The algebra ${\cal O}_\Gamma$.} Let $\Gamma$ be a graph. An oriented loop $\alpha$ on $\Gamma$ is a continuous 
map $\alpha: S^1 \to \Gamma$ of an oriented circle. 
A vertex of the loop $\alpha$  is a point $v \in S^1$ which is mapped to 
a vertex of $\Gamma$.

Denote by $l_\Gamma$  the set of all oriented loops on  $\Gamma$, considered modulo homotopy. 
Let 
$$
L_\Gamma:= \Z[l_\Gamma]
$$
 be the 
free abelian group generated by the set $l_\Gamma$. 
Let ${\cal O}_\Gamma$ be its symmetric algebra:
\be \la{OALG}
{\cal O}_\Gamma:= S^*(L_\Gamma).
\ee

Let ${\rm Loc}_{N, \C}(\Gamma)$ be the moduli space of all $N$-dimensional 
complex vector bundles with connections on $\Gamma$. 
Denote by ${\cal O}({\rm Loc}_{N, \C}(\Gamma))$ the algebra of regular functions 
on this space. 
There is a canonical map of commutative algebras 
\be \la{CAH}
r: {\cal O}_\Gamma \lra {\cal O}({\rm Loc}_{N, \C}(\Gamma)).
\ee
Namely, a loop $\alpha$ 
on $\Gamma$ give rise to a function $M_\alpha$ 
on ${\rm Loc}_{N, \C}(\Gamma)$ given by the $1/N$ times the trace of the monodromy 
 along the loop $\alpha$. Homotopic loops give   the same function. 
Since ${\cal O}_\Gamma$ is  the free commutative algebra generated by the set 
$l_\Gamma$, the map $\alpha \lms M_\alpha$ extends uniquely to a homomorphism
(\ref{CAH}). The trivial loop maps to the unit.

\paragraph{2. The  Goldman   bracket  for a ribbon graph $\Gamma$.}  Let us define a Lie bracket $\{\ast, \ast\}_{\rm G}$ on ${L}_\Gamma$. 

A valency $m$ vertex $v$ of a ribbon graph $\Gamma$ gives rise to a free abelian group ${\mathbb A}_v$ of rank $m-1$, given by  
$\Z$-linear combinations of the oriented out of $v$ edges $\stackrel{\to}{E}_i$ at $v$, of total degree zero: 
$$
{\mathbb A}_v = \{\sum n_i \stackrel{\to}{E}_i ~|~ \sum n_i=0, ~n_i\in \Z\}. 
$$
It is 
generated by  the ``oriented paths'' $\stackrel{\to}E_i - \stackrel{\to}E_j$, 
where $- \stackrel{\to}E_j$ is the edge $E_j$ oriented towards $v$. 

\begin{figure}[ht]
\centerline{\epsfbox{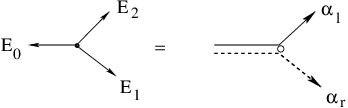}}
\caption{One has  $\delta_v(\alpha_r, \alpha_l)=\frac{1}{2}$.\label{di8}}
\end{figure}

\bl Given a vertex $v$ of $\Gamma$, there is a unique skew symmetric bilinear form 
$$
\delta_v: {\mathbb A}_v\wedge {\mathbb A}_v \lra \frac{1}{2}\Z
$$
such that for any triple  $E_0, E_1, E_2$ of edges as in Figure \ref{di8} one has 
\be \la{4.18.10.111}
\delta_v( \alpha_r\wedge \alpha_l) = \frac{1}{2}, \qquad \alpha_r:= \stackrel{\to}{E_1} - 
\stackrel{\to}{E_0}, \quad \alpha_l:= \stackrel{\to}{E_2} - \stackrel{\to}{E_0}. 
\ee
\el

Take two loops $\alpha$ and $\beta$ on a ribbon graph $\Gamma$. 
Let $\alpha_v\in {\Bbb A}_v$ be the element induced by the 
loop $\alpha$ at the vertex $v$. Set $\delta_v(\alpha, \beta):= \delta_v(\alpha_v, \beta_v)$. 
Given a pair $(v, w)$, where $v$ is a vertex of 
$\alpha$, $w$ is a vertex of $\beta$, and $\alpha(v) = \beta(w)$, consider 
a new loop $\alpha\circ_{v=w} \beta $ obtained by starting at the vertex $\beta(w)$, 
going around the loop $\beta$, and then around the loop $\alpha$. 
We define a bracket $\{\alpha,  \beta\}_{\rm G}$ 
by taking the sum over all such pairs of vertices $(v, w)$:
$$
\{\alpha,  \beta\}_{\rm G}:= \sum_{{v=w}} \delta_v(\alpha, \beta) \cdot \alpha\circ_{v=w} \beta. 
$$ 
This is a ribbon graph version of the Goldman bracket.  
It is easy to check the following. 

\bt
The   bracket $\{\ast, \ast\}_{\rm G}$ provides $L_\Gamma$ with a Lie algebra structure. 
So it extends via the Leibniz rule to a Poisson bracket $\{\ast, \ast\}_{\rm G}$ 
on the commutative algebra ${\cal O}_\Gamma$. The map (\ref{CAH}) is a map of Poisson algebras. 
\et

\paragraph{3. The Lie bracket $\{\ast, \ast\}$ on the space ${L}_\Gamma$  for a {bipartite} ribbon graph $\Gamma$.}  
Let $\Gamma$ be a bipartite ribbon graph.  
We apply the previous construction to the conjugate graph $ \Gamma^*$.

\bd \la{DGBRG}
i) The Lie bracket $\{\ast, \ast\}$ on  $L_\Gamma$ is induced, via the isomorphism 
$L_{\Gamma} = L_{ \Gamma^*}$, by the 
Lie bracket  on $L_{ \Gamma^*}$, for the conjugate graph $\Gamma^*$. 

ii) It 
 induces, via the Leibniz rule, a Poisson algebra structure $\{\ast, \ast\}$ 
on ${\cal O}_\Gamma$. 
\ed

Let $\alpha, \beta$ be two oriented loops on a bipartite surface graph $\Gamma$. 
Then 
$$
\{\alpha,\beta\} = \sum_v{\rm sgn}(v)\delta_v(\alpha, \beta) \cdot \alpha \circ_v \beta. 
$$
The sum is over all vertices $v \in\alpha\cap \beta$, 
and ${\rm sgn}(v)=1$ for the $\circ$-vertex, and $-1$ for the $\bullet$-vertex.

\begin{figure}[ht]
\centerline{\epsfbox{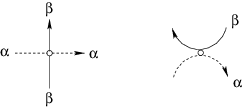}}
\caption{On the left: $\delta_v(\alpha, \beta)=1$. On the right: $\delta_v(\alpha, \beta)=0$.\label{di10}}
\end{figure}
To calculate the bracket $\{\alpha,\beta \}$ between  two loops $\alpha, \beta$ on $\Gamma$ 
we have the following two cases:

\begin{figure}[ht]
\centerline{\epsfbox{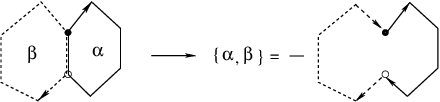}}
\caption{One has $\{\alpha, \beta\}=- \alpha \circ \beta$, where 
$\alpha \circ \beta$ is the loop on the right.\label{di9}}
\end{figure}

\begin{itemize}

\item A vertex $v \in \alpha \cap \beta$ is an isolated intersection point. Then its contribution is 
$\pm 1, 0$ depending on the geometry of the intersection, as depicted on Figure \ref{di10}.

\item A vertex $v \in \alpha \cap \beta$ lies on an edge $vw \subset \alpha \cap \beta$. Then its contribution is 
$\pm 1, 0$ depending on the geometry of the intersection, as depicted on Figures \ref{di9}, \ref{di9a}.

\begin{figure}[ht]
\centerline{\epsfbox{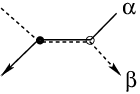}}
\caption{The total  contribution of the $\circ$ and $\bullet$-vertices to $\{\alpha, \beta\}$ is zero.\label{di9a}}
\end{figure}

\end {itemize}

Changing the orientation of the surface amounts to changing the 
sign of the Poisson bracket. Interchanging black and white we change the 
sign of the Poisson bracket. 
\vskip 2mm

\begin{figure}[ht]
\centerline{\epsfbox{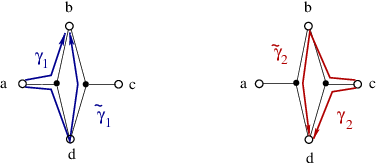}}
\caption{The sum of the contributions of the edges 
$\bullet b$ and $\bullet d$ to $\{\gamma_1, \widetilde \gamma_2\} + 
\{\widetilde \gamma_1, \gamma_2\}$ is zero. 
Thus  $\{\gamma_1, \widetilde \gamma_2\} + 
\{\widetilde \gamma_1, \gamma_2\} =0$. One has $\{\widetilde \gamma_1, \widetilde \gamma_2\} =0$ since the paths are enter and exit at the common edges 
$b\bullet $ and $d\bullet $ in the same direction. So
$\{\gamma_1+ \widetilde \gamma_1, \gamma_2+ \widetilde \gamma_2\} = 
\{\gamma_1, \gamma_2\} $. }
\label{ncls7b}
\end{figure}

Let us consider the  formal completion $\widehat L_{\Gamma} $ of $L_{\Gamma}$  by the length of the loop on the graph $\Gamma$, 
counted as the total number of edges in the loop. Then we make sense of  $(1+M)^{-1}\gamma$ as $(1-M+M^2-M^3 + ...) \gamma$. 

\bt \la{2.12}
Let $\Gamma \to \Gamma'$ be a two by two move of bipartite ribbon graphs. 
Then the map induced by Definition \ref{AFLBG} is a Lie algebra homomorphism 
of completed Lie algebras: 
$$
\mu_{\Gamma \to \Gamma'}: \widehat L_{\Gamma} \lra \widehat L_{\Gamma'}.
$$
\et

\begin{figure}[ht]
\centerline{\epsfbox{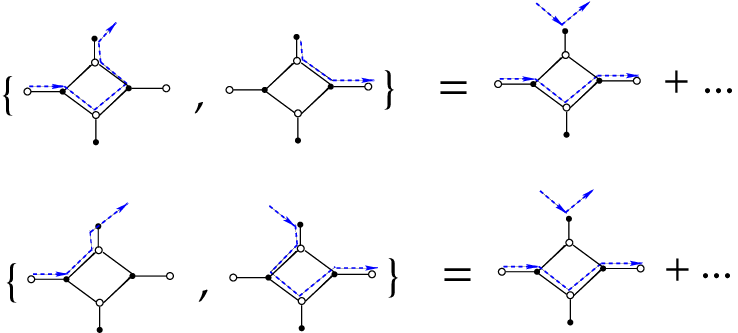}}
\caption{Calculating the Poisson brackets $\{M_b\gamma_{b \bullet a}, X_{c\bullet b}\}$ and $ \{\gamma_{b \bullet a},  M_c X_{c\bullet b}\}$ we get  the same result.}
\label{nclsid}
\end{figure}

\begin{proof} Denote by $\gamma_{b \bullet a}$ a loop on $\Gamma$ containing the segment $t_{b \bullet a}$, etc. Consider the following three identities:
\be \la{threeid}
\begin{split}
&\{\gamma_{b'\bullet a'}, \gamma_{d'\bullet c'}\} = \{\gamma_{b\bullet a}, \gamma_{d\bullet c}\},\\
&\{\gamma_{b' \bullet a'}, \gamma_{c' \bullet b'}\}= \{\gamma_{b \bullet a}, \gamma_{c \bullet b}\},\\
&\{\gamma_{b' \bullet a'}, \gamma_{c' \bullet b'}\}= \{\gamma_{b \bullet a}, \gamma_{c \bullet b}\}.\\
\end{split}
\ee
Theorem \ref{2.12} follows from them. Indeed, we have twelve identities like this to check. 
Using the skew-symmetry of the Poisson bracket, we reduce their number to six. These six can be reduced to the  three identities (\ref{threeid}): two  identities are reduced  by the rotation by $180^\circ$ of the graph, and  one by changing the plane orientation. Relation 3) in (\ref{threeid}) is reduced to 1) by changing the plane orientation. 
Let us check  identities 1) and 2). 

1) We have the following, where the first equality is Definition (\ref{10.12.11.1b}), and the second is explained on Figure \ref{ncls7b}:
$$
\{\gamma_{b'\bullet a'}, \gamma_{d'\bullet c'}\} = \{(1+M_b) \gamma_{b\bullet a}, (1+M_d) \gamma_{d\bullet c}\} 
\stackrel{\mbox{{\rm Fig} \ref{ncls7b}}}{=} \{\gamma_{b\bullet a}, \gamma_{d\bullet c}\}. 
$$

2) The second boils down to
\be
\begin{split}
&\{(1+M_b) \gamma_{b \bullet a}, (1+M_{c})^{-1}\gamma_{c \bullet b}\} \stackrel{?}{=} 
\{\gamma_{b \bullet a}, \gamma_{c \bullet b}\}.\\
\end{split}
\ee
So we need to prove that 
\be
\begin{split}
&\{M_b \gamma_{b \bullet a}, (1-M_{c}^{-1}+ M_c^{2}-...) \gamma_{c \bullet b}\} \stackrel{?}{=} 
\{\gamma_{b \bullet a}, (M_{c}^{-1}- M_c^{2}+ M_c^3 - ...) \gamma_{c \bullet b}\}.\\
\end{split}
\ee
This is checked by using the following general identity, where $X_{c\bullet b}$ is a loop containing the segment $c \bullet b$: 
$$
\{M_b\gamma_{b \bullet a}, X_{c\bullet b}\}= \{\gamma_{b \bullet a},  M_c X_{c\bullet b}\}.
$$
This identity proved by a calculation illustrated on Figure \ref{nclsid}.
\end{proof}

\bc A two by two move gives rise to   a birational isomorphism of Poisson algebras:
$$
\mu_{\Gamma \to \Gamma'}: \widehat {\cal O}_{\Gamma} \lra \widehat {\cal O}_{\Gamma'}.
$$
\ec

 \paragraph{Conclusion.} {\it Groupoids ${\rm Loc}_1(\Gamma)$ assigned to 
bipartite graphs $\Gamma$ and related by 
birational transformations (\ref{birat}) form a non-commutative analog of a 
cluster Poisson variety, with the Poisson bracket on the algebra (\ref{OALG}) given by Definition \ref{DGBRG}.}

\section{Dimers   and non-commutative cluster integrable systems} \la{SEC6}

In Section \ref{SEC6} we apply the general results about non-commutative cluster Poisson varieties from previous Sections to produce non-commutative cluster integrable systems, which in the commutative case were studied in \cite{GK}. Section \ref{SEC6} will not be used anywhere in the rest of the paper.

\paragraph{1. Dimer covers.} A {\it dimer cover} of a bipartite graph 
$\Gamma$ 
is a subset of edges of $\Gamma$ 
such that each vertex is the endpoint of a unique edge in the dimer cover. 
Given a dimer cover $M$, we assign to each edge $E$  of $\Gamma$ 
the weight $\omega_M(E)\in  \{0,1\}$ such that 
$\omega_M(E) = 1$ if and only if the edge $E$ belongs to $M$. 
Let $[E]$ be the edge $E$ oriented as  $\bullet \to \circ$. 
Then we get  a $1$-chain 
$$
[M]:= \sum_E \omega_M(E)[E]. 
$$
Here the sum is over all edges. We define ${\rm sgn}(v):
=-1$ for the $\bullet$-vertex $v$, and ${\rm sgn}(v): = 1$ for the $\circ$-vertex. 
The  condition that $M$ is a dimer cover  just means that 
$$
d [M] = \sum_{v}{\rm sgn}(v)[v].
$$
Here the sum is over all vertices $v$ of $\Gamma$. 
So if $M_1,M_2$ are dimer covers,  $[M_1]-[M_2]$ is a $1$-cycle.

\paragraph{2. Zig-zag paths revisited.}   
Recall that a zig-zag   $\gamma$ on a bipartite ribbon graph  $\Gamma$ is oriented so that 
it turns right at the $\circ$-vertices, and turns left at the $\bullet$-vertices.   
For every edge  there are exactly two 
zig-zag paths containing the edge. They traverse the edge  the opposite ways.

\paragraph{3. A reference dimer cover $\Phi$ for a bipartite graph on a torus \cite{GK}.}  We consider from now   bipartite graphs $\Gamma$ on a torus ${\rm T}$. 
A bipartite  graph $\Gamma$ on a torus is  {\it minimal} if it does not have parallel bigons, see the left picture on Figure \ref{ncls12}. 
Homology classes of  oriented zig-zag paths on a bipartite graph $\Gamma\subset {\rm T}$ 
form a collection of vectors $e_1, ..., e_n$ of the rank two   lattice 
$H_1({\rm T}, \Z)$. We order them cyclically by the angle they 
make with a given direction. Their sum is zero. So 
they are oriented sides of a convex polygon $N \subset H_1({\rm T}, \Z)$,  called the {\it Newton polygon} 
of   $\Gamma$. 
The   set $\{e_1,\dots,e_n\}$ is identified with the set of primitive edges of the boundary $\partial
N$ of $N$. Alternatively, $N$ is the convex hull of 
$0$,  $e_1$, $e_1+e_2$, ..., $e_1+ ... + e_{n-1}$.

Let $X$ be the set of circular-order-preserving maps 
\be \la{alpha1}
\alpha: \{e_1,\dots,e_n\} \to \R/\Z.
\ee
Let $E$ be an edge of $\Gamma$. 
Denote by $z_r$ (respectively $z_l$) 
the zig-zag path containing $E$ in which the $\circ$-vertex of $E$
precedes (respectively follows) the $\bullet$-vertex. 
Then a map $\alpha$ provides a function 
\be \la{alpha2}
\varphi_\alpha = \varphi:  \{\mbox{Edges of $\Gamma$}\} \lra \R, \qquad \varphi(E):= 
\alpha_r-\alpha_l.
\ee
Here we use the identification of  zig-zag loops with their homology classes $e_i$. Then 
$\alpha_r$ (respectively $\alpha_l$) is the evaluation  of the map
$\alpha$
on the zig-zag
path $z_r$ 
(respectively $z_l$) crossing the
edge $E$, and $\alpha_r-\alpha_l$ denotes the
length of the counterclockwise arc on $S^1=\R/\Z$ from $\alpha_l$ to $\alpha_r$.

\bt \la{4.17.10.1}  \cite[Theorem 3.3]{GK}   
For any map $\alpha \in X$, the function
$\varphi_\alpha$ satisfies
$$
d \varphi_\alpha = \sum_{v}{\rm sgn}(v)[v].
$$ 
\et

Take a map $\alpha$ which sends all $e_i$ to $0$ and assigns $1$ to the arc between $e_i$ and $e_{i+1}$,  
and $0$ to all other arcs. Then $\varphi_{\alpha}$ is a dimer cover of $\Gamma$. 
We denote it by $\Phi$. It corresponds   to a vertex of the 
Newton polygone, perhaps a degenerate one, i.e. located on a side.

\paragraph{4. The commuting Hamiltonians.} Since $M$ and $\Phi$ are perfect matchings, 
the $1$-cycle $[M]-[\Phi]$ is a union of 
non-intersecting and nonselfintersecting loops. 
Considered modulo homotopy on the graph $\Gamma$, 
it defines an element 
$$
(M):= [M]-[\Phi] \in L(\Gamma). 
$$

\bd Let $\Gamma$ be a bipartite graph on a torus. 
The non-commutative partition function ${\mathcal P}_{\Phi}$
 is an element of ${L}(\Gamma)$ given as the sum 
over all dimer covers
$M$ of $\Gamma$  
\be \la{partsum}
{\mathcal P}_{\Phi}:= \sum_{M}{\rm sgn}(M) (M).
\ee
\ed

The element $[M]-[\Phi]$ defines a class in $H_1({\rm T}, \Z)$,  
denoted by $h_{\rm T}(M)$.  The set of   homology classes   for different dimer covers $M$ 
of $\Gamma$ 
coincides with the set of internal points of a shift of the Newton polygon $N$. 
Take the sum over the dimer covers $M$ 
with the same homology class  $a$:
\be \la{decomph}
H_{a} := \sum_{h_{\rm T}(M)=a}(M)\in L(\Gamma).
\ee
Furthermore, for any nonzero  integer $n$, we define a new element 
$$
H_{na}\in L(\Gamma).
$$  
Namely, if $n>0$, it is obtained by 
travelling  $n$ times each loop of $H_{a}$. If $n<0$, we reverse 
the orientations of the loops of $H_{a}$, and then travel each obtained loop $|n|$ times. 

\bd 
The  elements $H_{na}$ are the {\it Hamiltonians of the dimer system}. 
\ed

\paragraph{Remark.} A map $\alpha$ related to another vertex of the Newton polygon leads to another collection  
of Hamiltonians which differ from $\{H_{a}\}$ 
by a sum of zig-zag loops, which lies in the center of the Poisson algebra ${\cal O}(\Gamma)$. 
Therefore the  Hamiltonian flows 
do not depend on the choice of $\alpha$,  and on each symplectic leaf any two collections of Hamiltonians 
differ by a common factor.

\bt \la{completeintegrabilitytheorem}
Let $\Gamma$ be a minimal bipartite graph on a torus ${\rm T}$. Then 
the Hamiltonians $H_{na}$ commute under the Lie bracket 
 on 
${L}(\Gamma)$. Therefore they Poisson commute in ${\cal O}(\Gamma)$. 

\et

\begin{proof} We prove that $\{H_a, H_b\}=0$. The proof almost literally follows 
the proof of \cite[Theorem 3.7]{GK}. For the convenience of the reader 
we indicate the main steps. The claim that 
$\{H_{na}, H_{mb}\}=0$ follows from this. 

Take a pair of matchings $(M_1, M_2)$ on $\Gamma$. 
Let us assign to them  another pair of matchings 
$(\widetilde M_1, \widetilde M_2)$ on $\Gamma$. 
Since $[M_1]-[M_2]$ is a $1$-cycle, and $M_1, M_2$ are dimer covers, 
it is a disjoint union of non-itersecting and 
non-selfintersecting loops. 
So $[M_1]-[M_2]$  is a disjoint union  of: 

\begin{enumerate}

\item homologically trivial on the torus loops,

\item  homologically  non-trivial  loops,  

\item edges shared by both matchings. 
\end{enumerate}

For every edge $E$ of each homologically trivial loop 
we switch the label of $E$: if $E$ belongs 
to a matching $M_1$ (respectively $M_2$), we declare that it will belong to a matching 
$\widetilde M_2$ (respectively $\widetilde M_1$). 
For all other edges we keep their labels intact. 
Switching the labels of a homologically trivial loop we do not change 
 the homology classes $[M_1]-[\Phi]$,  $[M_2]-[\Phi]$.\footnote{Contrary to this, switching the labels 
of homologically non-trivial loops in $[M_1]-[M_2]$, we do change the homology classes $[M_i]-[\Phi]$.} 
Clearly $\widetilde {\widetilde M}_i = M_i$.

\bl \la{4.16.10.1} \cite[Lemma 3.8]{GK}
One has 
\be \la{4.16.10.3}
\{(M_1), (M_2)\} + \{(\widetilde M_1), (\widetilde M_2)\} =0.
\ee
\el

\begin{proof} 
Let us write (\ref{4.16.10.3}) 
as a sum of the local contributions corresponding to the vertices $v$ of $\Gamma$, 
and break it in three pieces as follows: 
$$
\varepsilon(\mu_1, \mu_2) + \varepsilon(\widetilde \mu_1, \widetilde \mu_2) = 
(\sum_{v \in {\mathcal E}_1} + \sum_{v \in {\mathcal E}_2} + \sum_{v \in {\mathcal E}_3})
\Bigl(\delta_v(\mu_1, \mu_2) + \delta_v(\widetilde \mu_1, \widetilde \mu_2)\Bigr).
$$
Here 
${\mathcal E}_1$ (respectively ${\mathcal E}_2$ and ${\mathcal E}_3$) is the set of vertices which belongs to 
the homologically trivial loops (respectively homologically non-trivial loops, double edges).

For any $v \in {\mathcal E}_1$ we have 
$\delta_v(\mu_1, \mu_2) + \delta_v(\widetilde \mu_1, \widetilde \mu_2)=0$ since the elements 
$l_v(\mu_1)\in {\mathbb A}_v$ provided by $\mu_i$ coincides with the 
 $l_v(\widetilde \mu_2)$, and similarly $l_v(\mu_2) = l_v(\widetilde \mu_1)$. 
So the first sum is a sum of zeros.

The third sum is also a sum of zeros by the skew symmetry of the bracket. 

It remains to prove that 
$$
\sum_{v \in {\mathcal E}_2} 
\Bigl(\delta_v(\mu_1, \mu_2) + \delta_v(\widetilde \mu_1, \widetilde \mu_2)\Bigr)=0.
$$

Let $\gamma$ be an oriented loop on $\Gamma$. The {\it bending} 
$b_v(\gamma; \varphi)$ of the function $\varphi$ at a vertex 
$v$ of $\gamma$ is 
$$
b_v(\gamma; \varphi)=\sum_{E \in R_v}\varphi(E) - \sum_{E \in L_v}\varphi(E)\in \R. 
$$
Here $R_v$ (respectively $L_v$) 
is the set of all edges sharing the vertex $v$ which are on the right
(respectively left) of the oriented path $\gamma$. Lemma \ref{4.16.10.1} now follows
from Lemma \ref{Euclidean} below. 
\end{proof}

\bl\la{Euclidean} \cite[Lemma 3.9]{GK} For any simple topologically nontrivial loop $\gamma$, and for any  $\alpha\in X$, see (\ref{alpha1}), the corresponding function $\varphi_\alpha$ 
satisfies,  
$$
\sum_{v \in \gamma}b_v(\gamma; \varphi_\alpha)=0. 
$$
\el

\end{proof}

 \section{Non-commutative cluster ${\cal A}-$varieties} \la{SSEECC5AA}

 \subsection{Non-commutative cluster ${\cal A}-$varieties  from bipartite ribbon graphs} \la{SSEECC5}
 
  {\it Convention.} Dealing with \underline{twisted} local systems on a decorated surface $\bS$,  we alter $\bS$ by cutting out a little disc  near each marked point $m$ on $\bS$. Abusing notation, we  
 denote the obtained surface   by $\bS$.  The tangent vectors to each   boundary 
component of $\bS$,  positively oriented for the boundary orientation of $\bS$, 
form a homotopy (semi)circle $S^1_m$ in the bundle   of non-zero tangent vectors to $\bS$. 
We perform a similar procedure with the marked points on the spectral curve $\Sigma$. 
Talking about a flat section of a twisted local system near a marked point $m$  of $\bS$ or $\Sigma$, we mean a flat section of its restriction to the 
  (semi)circle $S^1_m$.

 The definition of a framed local system on a decorated surface in Definition \ref{9.9.11.1}    
 easily extends to the case of a twisted framed local system on $\bS$. So by the monodromy of a twisted local system near a puncture we mean 
the monodromy of its restriction to the oriented homotopy circle near the puncture.\\ 
  
 \begin{figure}[ht]
\centerline{\epsfbox{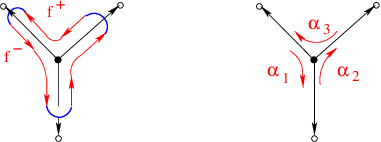}}
\caption{An   edge $E$ of a ribbon graph $\Gamma$ determines  two  faces
  of $\bS_\Gamma$, denoted $f_E^+$ and $f_E^-$.   
  There is a  canonical homotopy class of path    $p_E$, connecting boundary components $f_E^+$ and $f_E^-$, shown by a blue arc. }
\label{ncls00}
\end{figure} 

Pick a  vector  field tangent to the boundary of each face of a decorated surface, following the boundary orientation,  see Figure \ref{ncls00}. 
A {trivialization of   a twisted flat line bundle ${\cal L}$ on the surface  at the boundary } is given by a choice of a non-zero section   of the fiber of ${\cal L}$ over the tangent boundary   vector field.

\noindent 
Recall the surface $\bS_\Gamma$  for a ribbon graph $\Gamma$, and the spectral surface $\Sigma_\Gamma$ for a bipartite ribbon graph $\Gamma$. 

  Given a bipartite ribbon graph $\Gamma$, consider the groupoid, defined in either of the following two ways: 
\be
\begin{split}
{\cal A}_\Gamma:= &\{\mbox{Twisted   flat line bundles on  the spectral 
  surface ${\Sigma_\Gamma }$,   trivialized
at   the boundary}\}.\\
  =&\{\mbox{Twisted   flat line bundles on  the  
  surface ${\bS_\Gamma }$,   trivialized
at the zig-zag strands}\}.\\
\end{split}
\ee
These definitions are equivalent since  there are canonical identifications: 
$$
\Gamma =\Gamma^*, \ \ \  \Sigma_\Gamma = \bS_{\Gamma^*}, \ \ \ \mbox{  
$\{$zig-zags on $\bS_\Gamma$$\}$ =  $\{$boundary components of $\Sigma_\Gamma\}$}.
$$

Let us give  a coordinate description of the non-commutative tori 
${\cal A}_\Gamma$ in terms of the ${\cal A}-$coordinates.

 \paragraph{1. ${\cal A}-$coordinates on   ribbon graphs.} 
  
\bd \la{DDD} Let $\Gamma$ be a ribbon graph. The {\em ${\cal A}-$coordinates} on $\Gamma$ are   the  elements $\{\Delta_E\in R^*\}$ assigned to 
 the  \underline{oriented} edges $E$ of $\Gamma$,    satisfying  the following monomial relations:

  \begin{itemize}
 \item Let $E$ be an oriented edge, and   $\overline E$ the same edge  with  the opposite orientation. Then
    \be \la{MEQi}
  \Delta_{E}\Delta_{\overline E}=-1.
  \ee  
 \item Let   $E_1, ..., E_n$ be  the edges    incident to any   vertex $v$ of $\Gamma$, oriented out of the vertex $v$, whose order 
is compatible with  their cyclic order. Then     
   \be \la{MEQ}
  \Delta_{E_1}\Delta_{E_2} \ldots \Delta_{E_n}=-1.
  \ee 
\end{itemize}
 \ed

Calling the elements $\Delta_{E}$ "coordinates" is an abuse of terminology since they are not independent.   We   use the name {\it ${\cal A}-$decorated ribbon graph} for a ribbon graph    with a collection of ${\cal A}-$coordinates.

\bl \la{CRTH}  Given a   
  ribbon graph $\Gamma$, the  ${\cal A}-$coordinates on   $\Gamma$ parametrise 
 the isomorphism classes of twisted   flat line bundles on  the  
  surface ${\bS_\Gamma }$,   trivialized
at   the boundary. 
\el

\begin{proof}  
We start from a twisted line bundle on  the oriented decorated surface $\bS_\Gamma$.   
An   edge $E$ of $\Gamma$ determines  two  faces $f_E^+$ and $f_E^-$  of $\bS_\Gamma$.   
  There is a  canonical homotopy class of path    $p_E$, connecting boundary components $f_E^+$ and $f_E^-$, 
see Figure \ref{ncls00}.  
  The parallel transport along the   $p_E$   acts on the trivializations $s_E^\pm$ at the boundary components   $f_E^\pm$ 
as follows, defining   elements $\Delta_{E}\in R^*$:
 \be \la{para}
 {\rm par}_{p_E}: s^+_{E} \lra \Delta_{E} s^-_{E}, \ \ \ \ \Delta_{E}\in R^*. 
 \ee 
 Going around any edge, or any vertex along the cyclic order at this vertex, as   on Figure \ref{ncls00},   we rotate the tangent vector by $2\pi$,    getting both relations  (\ref{MEQi}) and (\ref{MEQ}). 
 For example, for a 3-valent vertex we get 
 \be
\begin{split}
& s^+_{E_1} \lra \Delta_{E_1} s^-_{E_1} = \Delta_{E_1} s^+_{E_2}   \lra \Delta_{E_1} \Delta_{E_2} s^-_{E_2} = \Delta_{E_1} \Delta_{E_2}s^+_{E_3}  \\
&  \lra \Delta_{E_1} \Delta_{E_2} \Delta_{E_3} s^-_{E_3} = -\Delta_{E_1} \Delta_{E_2} \Delta_{E_3} s^+_{E_1}. \\
\end{split}
\ee
Here the arrows are induced  by the parallel transports (\ref{para}) ${\rm par}_{p_{E_i}}$, $i \in \Z/3\Z$. The  $=$ signs mean canonical identifications provided by the parallel transport along the arcs $\alpha_i$   between $E_i$ to $E_{i+1}$, see Figure \ref{ncls00},  
identifying $s_{E_i}^-$ with $s_{E_{i+1}}^+$.  The $-$ sign in the end results from the fact that going around the circle amounts to the monodromy $-1$.

Conversely,  picking a section $s$ over a tangent vector to a boundary component, and using   elements $\Delta_{ E}$ satisfying the monomial relations 
  (\ref{MEQi}) and (\ref{MEQ}), we recover 
the twisted flat line bundle on $\bS_\Gamma$,   trivialized at the boundary. Changing $s$  we get an isomorphic object. \end{proof} 

\paragraph{2. Conventions on ${\cal A}-$coordinates on \underline{bipartite}  ribbon graphs.} {\it Given a bipartite ribbon graph, its edges  are orientated $\circ\lra \bullet$}. 

Using this convention,  the set of ${\cal A}-$coordinates on    any bipartite ribbon graph can be given by a collection of   elements $\{\Delta_E\}$ assigned to {\it nonoriented} edges $E$, 
assuming their default $\circ\lra \bullet$ orientation,  subject to the following relations, where {\it antycyclic} means the order opposite to the cyclic order:
 \be \la{CONV} 
\begin{split}
&\mbox{ \it The cyclic product of the elements on the edges sharing   any $\circ-$vertex is equal to $-1$.}\\
&\mbox{ \it The antycyclic product of the elements on the edges sharing   any $\bullet-$vertex $v$ is   $(-1)^{{\rm val}(v)-1}$.}\\
\end{split}
\ee 
Indeed, given a $\bullet-$vertex $v$, we have,  using $\Delta_{\overline E_i} = - \Delta^{-1}_{E_i}$:
\be
-1=\prod_{i=1}^{{\rm val}(v)}\Delta_{\overline E_i} = \prod_{i=1}^{{\rm val}(v)}(-\Delta^{-1}_{  E_i})= 
(-1)^{{\rm val}(v)}\Bigl(\prod^{i=1}_{{\rm val}(v)} \Delta_{ E_i} \Bigr)^{-1}.
\ee
We apply this convention to the bipartite ribbon graph $\Gamma^*$, since the spectral surface where the twisted flat line bundles live, is  $\bS_{\Gamma^*}$.  
  
Note that  the bipartite ribbon graphs  $\Gamma$ and $\Gamma^*$  are identical as  graphs. 
  Thinking about the graph $\Gamma$, e.g. drawing Figures,    we   depict the cyclic order on   $\Gamma$ as the  counterclockwise one, and 
   use below  and in Section \ref{sec4.2}  a different convention: 
 \be \la{CONVa} 
\begin{split}
&\mbox{ \it The counterclockwise product of the elements on the edges sharing   a  $\circ-$vertex is equal to $-1$.}\\
&\mbox{ \it The counterclockwise product of the elements on the edges sharing   a $\bullet-$vertex $v$ is   $(-1)^{{\rm val}(v)-1}$.}\\
\end{split}
\ee

\paragraph{3. The action of the two by two moves on ${\cal A}-$coordinates on bipartite ribbon graphs.}  
Recall the two by two moves of bipartite ribbon graphs, shown on Figure \ref{Amut}, see also  Definition \ref{DEF6.1} below. 
 
   \begin{figure}[ht]
\centerline{\epsfbox{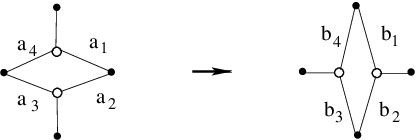}}
\caption{The  relations between the coordinates $a_i$ and $b_i$ are  $a_ia_{i+1}=(b_{i}b_{i+1})^{-1}$, $\forall i \in \Z/4\Z$.   The cases $i=1,3$  are checked using Figure \ref{AB}, the cases  $i=2,4$ are similar. }
\label{Amut}
\end{figure} 
We define  a  birational isomorphism  ${\cal A}_\Gamma \lra {\cal A}_{\Gamma'}$ corresponding to a two by two move  $\Gamma \lra \Gamma'$,  
  serving the role of an elementary non-commutative ${\cal A}-$cluster transformation.  
Let us assign the coordinates $\{a_1, a_2, a_3, a_4\}$ and $\{b_1, b_2, b_3, b_4\}$ to the internal edges on the graphs on Figure \ref{Amut}. 
Let us set 
\be
A_i = a_ia_{i+1}a_{i+2}a_{i+3}, ~~~~B_i = b_ib_{i+1}b_{i+2}b_{i+3}, \ \  \ \ i \in \Z/4\Z.
\ee
Consider the following transformation of    the coordinates  illustrated on Figure \ref{Amut}. 
\be \la{FRM1}
 \begin{split}
 &  b_1  =    (1 +A_3^{-1} ) a_3 \ = \ \ \ 
 a_3  (1 + A_4^{-1} );  \\
  &   b_2 =  (1+ A_4)^{-1}a_4 = \ \ \ 
 a_4(1+A_1)^{-1};\\
&  b_3  =   (1 +A_1^{-1} )a_1\ = \ \ \ 
a_1  (1 + A_2^{-1} ).\\
  & b_4 =  (1+A_2)^{-1}a_2 = \ \ \ 
 a_2(1+A_3)^{-1};\\
  \end{split}
 \ee
 
 To check the second equality in each line, note that  $a_ia_{i+1}a_{i+2}a_{i+3} a_{i+4} $ can be written in two   ways:
\be
\begin{split}
&a_{i}A_{i+1}  =   A_{i}a_{i}, \ \ \ \ \forall i \in \Z/4\Z;\\
\end{split}
\ee 

\bl \la{CORL} The square of the two by two transformation (\ref{FRM1}) is the identity transformation.
\el

\begin{proof} This   can be proved,  by using the calculation  of the coordinate transformations for the space ${\rm Conf}^{\rm df}_4(V_2)$ of four decorated flags in Section \ref{sec4.2}. 
Since the direct check is   easy, let us elaborate it. 

Going from $b$'s to $a$'s, the analog of the first formula in (\ref{FRM1}) is the following formula:
$
 a_2=b_4(1+B_1^{-1}). 
$ 
Since $B_1^{-1}=A_3$, it is equivalent to $b_4=a_2(1+A_3)^{-1}$, which is indeed the case by (\ref{FRM1}).

Similarly,   the analog of the second formula in (\ref{FRM1}) is the following: 
$
 a_3=b_1(1+B_2)^{-1}. 
$
Since $B_2^{-1}=A_4$, it is equivalent to $b_1=a_3(1+A_4^{-1})$, which is indeed the case by (\ref{FRM1}).
\end{proof}

 \bl We have for any $i \in \Z/4\Z$:
\be \la{RAB}
\begin{split}
&b_{i}b_{i+1}= (a_{i } a_{i+1})^{-1},   \\ 
&A_i = B_{i+2}^{-1}.\\
\end{split}
\ee
\el

\begin{proof} 
We have 
\be
\begin{split}
 &b_3b_4a_3a_4 = a_1(1+A_2^{-1} )  (1+A_2)^{-1}  a_2 a_3a_4=  a_1A_2^{-1}a_2  a_3a_4 = A_2^{-1}A_2=1.\\
 & b_4b_1a_4 a_1 = a_2(1+A_3)^{-1} (1+A_3 ^{-1})a_3a_4 a_1=  a_2A_3^{-1}a_3a_4 a_1 =A_3^{-1}A_3=1.\\
\end{split}
\ee
Since   transformation formulas (\ref{FRM1}) are   invariant under the cyclic shift $i \lms i+2$, the first formula in (\ref{RAB}) follows. The second formula follows from the first. \end{proof}

\begin{figure}[ht]
\centerline{\epsfbox{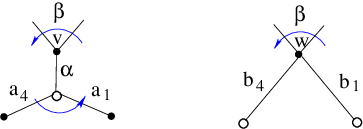}}
\caption{Using conventions (\ref{CONVa}), we have: $a_4a_1 \alpha=-1$,  $\alpha  \beta= (-1)^{{\rm val}(v)-1}$, $b_4b_1\beta= (-1)^{{\rm val}(v)}$.  
Therefore $  a_4a_1  = (b_4b_1)^{-1}$.}
\label{AB}
\end{figure}

Figure \ref{AB} shows that  formula  (\ref{RAB}) is consistent with the condition that the counterclockwise product of the coordinates around a   $\circ-$vertex is $ 1$, and around a 
$\bullet-$vertex is  $1$. 
 
In Section \ref{sec4.2} we   define  ${\cal A}-$coordinates for  ${\rm GL_2}-$graphs on a decorated surface $\bS$,  
parametrising the moduli space ${\cal A}_{2, \bS}$. The change of these ${\cal A}-$coordinates   under a flip of a triangulation, 
  studied in    
Section \ref{sec4.2}, motivated the  definition of the  ${\cal A}-$cluster coordinate transformations in (\ref{FRM1}).  

\bt \la{Th5.4} The coordinate transformation (\ref{FRM1}) for a two by two move satisfies the pentagon relation. 
\et

The proof can be obtained by a long   calculation. 
 Instead, we give a   simple proof based  on Theorem \ref{Th4.5}  in the  end of Section \ref{sec4.2}. \vskip 2mm
 
 Theorem \ref{Th5.4} allows to define the non-comutative cluster variety ${\cal A} = {\cal A}_{|\Gamma|}$ assigned to a bipartite ribbon graph $\Gamma$. Its clusters are parametrised by the bipartite ribbon graphs 
which can be obtained from $\Gamma$ by two by two moves and  the elementary transformations given by  shrinking  a two-valent vertex. 
We assign to each cluster the non-commutative torus ${\cal A}_\Gamma$, and glue them using the 
coordinate transformations (\ref{FRM1}) corresponding to the two by two moves. In Section \ref{NONCTF} we will show that there is a canonical non-commutative 2-form $\Omega_{\cal A}$ 
on the cluster variety ${\cal A}$, preserved by the two by two moves of the graphs.

\subsection{The  canonical 2-form on a  non-commutative cluster ${\cal A}-$variety.} \la{NONCTF}

  \paragraph{1. The non-commutative analog of $d\log (a_1) \wedge \ldots \wedge d\log (a_n)$.} Recall the cyclic envelope of the tensor algebra ${\rm T}(V)$ of a graded vector space $V$. It is  spanned by the elements 
$(v_1 \otimes v_2 \otimes \ldots \otimes v_n)_{\cal C}$, which are  cyclically invariant: 
 $$
(v_1 \otimes v_2 \otimes \ldots \otimes v_n)_{\cal C}  = (-1)^{|v_1| \cdot (|v_2| + \ldots + |v_n|)}(v_2 \otimes v_3 \otimes \ldots \otimes v_n\otimes v_1)_{\cal C}.
 $$
Below, if it can not lead to a misunderstanding,  we use a shorthand 
$$
 v_1   v_2   \ldots   v_n := (v_1 \otimes v_2 \otimes \ldots \otimes v_n)_{\cal C}.
 $$
 Given a graded algebra $A$,  denote by $\Omega_A$ the $\Q-$vector space\footnote{Characteristic $p$ would be fine as well, but we do not need this in the paper.}    generated by the symbols $da$ where $a \in A$, modulo the following relations: 
 $d(\alpha a+\beta b) = \alpha da+\beta db$, $\alpha, \beta \in \Q$,  and the Leibniz rule  $d(ab) = da\cdot b + (-1)^{|a|}a \cdot db$. Here $|da| = |a|+1$, where $|\cdot |$ denotes the degree in $A$. Denote by ${\cal C}{\rm T}(A \oplus \Omega_A)$ the cyclic envelope of the tensor algebra of the 
 graded $\Q-$vector space ${\rm T}(A \oplus \Omega_A)$.  
 
\bd Given any graded algebra $A$, and any invertible elements $a_i \in A^\times$, 
  the  cyclic product  
  \be
  \{a_1, \ldots,  a_n\}:= da_1   \ldots da_n a_n^{-1} \ldots a_1^{-1} \in {\cal C}{\rm T}(A \oplus \Omega_A)
    \ee  
 is the noncommutative analog of $d\log (a_1) \wedge \ldots \wedge d\log (a_n)$.
  \ed
  
For example,  the non-commutative 2-form  generalizing $d\log (a) \wedge d\log (b)$  is given by: 
\be \la{ab}
\{a,b\} = da ~db ~b^{-1}a^{-1}.
\ee
Note that $\{\alpha a, \beta b\}=\{a,b\}$ for any $\alpha, \beta\in \Q^\times$. In particular, we  often use below that $\{\pm a, \pm b\}=\{a,b\}$. 

  \bt For any elements $a_0, \ldots , a_n$ we have  the Hochshild cocycle property: 
 \be \la{F5}
 \{a_1, ..., a_n\} - \sum_{i=0}^{n-1} (-1)^i\{a_0, \ldots , a_{i-1}, a_ia_{i+1}, a_{i+2}, \ldots , a_n\}  +(-1)^n  \{a_0, ..., a_{n-1}\}=0.
 \ee
\et

For example: 
\be \la{F2a}
\begin{split}
&\{a\}- \{ab\} + \{b\}=0,\\
&\{b, c\} - \{ab, c\} + \{a, bc\} - \{a,b\} = 0.\\
\end{split}
\ee
 
 \begin{proof} Each of the following terms appears in the sum (\ref{F5}) twice, with the opposite signs
  $$
 d a_0 \ldots  da_{i-1}  d a_i ~a_{i+1}  da_{i+2} \ldots  da_n a^{-1}_n ... a_0^{-1}.
 $$
 \end{proof}

\bl i) One has 
\be \la{F2}
\begin{split}
&\{a,b\} =  -\{b^{-1}, a^{-1}\}.\\
\end{split}
\ee

ii) Suppose that $abc= \pm 1$. Then    (\ref{F1}) implies that   $\{a, b\}$ is cyclic invariant:  
\be  \la{F1}
\begin{split}
&\{a,b\} + \{b,c\}+ \{c,a\} = 3\cdot  \{a,b\}.\\
\end{split}
\ee
 iii) We have
$
 \{x, 1+x\} = \{1+x, x\}.
$
\el

\begin{proof} i) Recall the identity
$
c~dc^{-1} = - dc ~ c^{-1}.
$ 
Using this, we have 
$$
\{a,b\} = da~db ~b^{-1}a^{-1} =  a^{-1}da~db ~b^{-1} \stackrel{}{=} da^{-1}~a~b~db^{-1} =  - db^{-1}da^{-1} ~ab = -\{b^{-1}, a^{-1}\}.
$$

ii) Since $\{a,b\}=\{a,-b\}$, we  may assume  that $abc=1$. Then we have to show that
\be
\begin{split}
&\{a, b\} + \{b, (ab)^{-1}\} + \{(ab)^{-1}, a\} = 3 \{a, b\}.\\
\end{split}
\ee
Let us write the second identity in (\ref{F2a}) for $c=  (ab)^{-1}$: 
\be 
\begin{split}
&\{a, b\} - \{a, a^{-1}\} + \{ab, (ab)^{-1}\} - \{b, (ab)^{-1}\} = 0.\\
\end{split}
\ee
By   (\ref{F2}), we have  $\{a, a^{-1}\}= - \{a, a^{-1}\}$. So $\{a, a^{-1}\}=0$. So 
 \be \la{F3}
\begin{split}
&   \{a, b\} = \{b, (ab)^{-1} \}.\\
\end{split}
\ee
Using (\ref{F2}) we have: $\{a,b\}=-\{b^{-1}, a^{-1}\} \stackrel{(\ref{F3})}{=}  -\{a^{-1}, ab\} = \{(ab)^{-1}, a\}$.  
 This and (\ref{F3}) imply the claim. \\
 
 iii) One has  
  $
 \{x, 1+x\}= dx dx (1+x)^{-1}x^{-1} = dx dx x^{-1}(1+x)^{-1} = \{1+x, x\}.
 $  \end{proof}

 \paragraph{2. The   non-commutative 2-form  $\Omega_\Gamma$ for an ${\cal A}-$decorated bipartite  ribbon graph $\Gamma$.} Recall that an edge of a   graph $\Gamma$ is called an {\it external edge} if one of its vertices is univalent.  

For simplicity, we do not consider graphs which have a   component given by a single  external edge. \\

 We start from a  trivalent  ribbon graph $\Gamma_v$ with a single-vertex  $v$,     decorated by $a_1, a_2, a_3$ in the order compatible with the cyclic structure at $v$, such that $a_1a_2a_3=\pm 1$.
 Then we set, using notation (\ref{ab}):
 \be
 \Omega_{v}:= \{a_1, a_2\}.
  \ee
 It depends only on the  
   cyclic order of the edges. Indeed,    $\{a_1, a_2\} = \{a_2, a_3\} = \{a_3, a_1\}$ since $a_1a_2a_3=\pm 1$.

\begin{figure}[t]
\centerline{\epsfbox{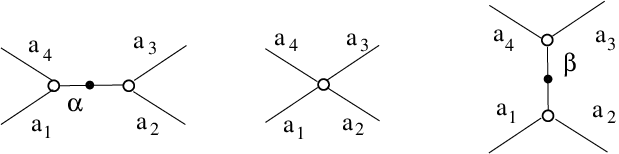}}
\caption{The cocycle condition guarantees that the 2-form $\Omega_\Gamma$ is well defined.}
\label{nclsI}
\end{figure}

Given a   ribbon graph $\Gamma_v$ with a single vertex $v$ of valency $>3$, we expand this vertex by adding two-valent vertices of the opposite color, as shown on Figure \ref{nclsI}, 
   producing a  bipartite ribbon graph $\Gamma'_v$   with trivalent vertices of the original color, and two-valent vertices of the opposite color. We set 
    \be
 \Omega_{v}:=   \sum_{x \in \Gamma'_v} \Omega_{x}.
 \ee   
By default, the contribution of the $2-$valent vertex is zero.\footnote{This is consistent with $\{a, a^{-1}\}=0$ since the decorations of a $2-$valent vertex are $a, \pm a^{-1}$.} The  sum is over trivalent  vertices  of the   graph $\Gamma'_v$. 
 To see that it does not depend on the choice of   $\Gamma'_v$, it is sufficient to check this for the    graphs on Figure \ref{nclsI}.  
This boils down  to the cocycle condition for  $\{a, b\}$:
 \be
 \begin{split}
& \{a_2, a_3\} + \{\alpha, a_4\} = \{a_1, a_2\} + \{\beta, a_3\}   \longleftrightarrow \\
& \{a_2, a_3\} + \{a_1, a_2a_3\} = \{a_1, a_2\} + \{a_1a_2, a_3\}.\\
 \end{split}
 \ee 
 Now we  introduce the non-commutative 2-form assigned to an any decorated bipartite ribbon graph $\Gamma$.

 Recall the ${\cal A}-$decorated   ribbon graph, see  Definition \ref{DDD}. 
 
 \bd    
  Given an ${\cal A}-$decorated  {bipartite} ribbon graph  $\Gamma$,  we set 
     \be
 \Omega_\Gamma:=  \sum_{w \in \Gamma} \Omega_{w} - \sum_{b \in \Gamma} \Omega_{b}.
 \ee   
 Here the first sum is over all $\circ-$vertices   of $\Gamma$, and the second over all $\bullet-$vertices.    
 \ed

 We assign to each   edge ${E}$ of an  ${\cal A}-$decorated ribbon graph $\Gamma$, decorated by an element $a_E$,  the following 3-form 
  \be
 \omega_{E}:= (a_{E}^{-1}da_{E})^3. 
  \ee

\bt \la{T111}
 For an  ${\cal A}-$decorated bipartite ribbon graph $\Gamma$, the  $d\Omega_\Gamma$ is a sum is over all external edges $E$ of   $\Gamma$:
 \be
d\Omega_\Gamma = \sum_{E} \omega_E.
\ee 
 The 2-form $\Omega_\Gamma$ assigned to an ${\cal A}-$decorated   ribbon graph $\Gamma$ without external edges is closed: $d\Omega_\Gamma=0$. \et

\begin{proof} Clearly the first claim implies the second. 

\bl
 Let $abc=1$. Then 
 \be \la{TRC}
 3 d \{a,b\} = -(a^{-1}da)^3 -  (b^{-1}db)^3 - (c^{-1}dc)^3 .
 \ee 
\el

\begin{proof}
Straightforward calculation, using $c^{-1}dc = - d(c^{-1}) c$ shows that the right hand side is 
$$
-(a^{-1}da)^3 -  (b^{-1}db)^3 +  \Bigl((ab)^{-1}d(ab)\Bigr)^3 = 3 da db d(b^{-1}) a^{-1} + 3 d(a^{-1}) da db b^{-1} =  3 d \{a,b\}. 
$$
 \end{proof}
 
 Therefore for each $\circ-$vertex $w$ of $\Gamma$ the differential of the 2-form $\Omega_w$ assigned to $w$ is  
  \be
  d\Omega_w = -\sum_{w \in {E}_i} \omega_{E_i}.
   \ee
      Here the sum is over all edges $E_i$ sharing the vertex $w$. 
      
      On the other hand, for each $\bullet-$vertex $b$ of $\Gamma$ the differential of the 2-form $\Omega_b$ assigned to $b$ is  
  \be
  d\Omega_b =  \sum_{b \in {E}_i} \omega_{E_i}.
   \ee  
   Here the sum is over all edges sharing the vertex $b$. Then evidently
  \be
 d\Omega_\Gamma =  \sum_{v \in \Gamma} d\Omega_v = -\sum_{E\subset \Gamma} \omega_{E} + \sum_{E\subset \Gamma} \omega_{{E}} = 0. 
  \ee 
\end{proof}

For the convenience of the reader, we reproduce  the picture illustrating  the change of the coordinates under the two by two move, which we refer to in Theorem \ref{THEOREM4.12}. 
 
    \begin{figure}[ht]
\centerline{\epsfbox{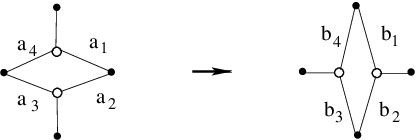}}
\caption{The  right picture is obtained from the left one by the clockwise rotation by $90^\circ$,  followed by the cyclic shift $a_i\lms b_{i+1}$,  $i \in \Z/4\Z$.}
\label{Amut*}
\end{figure} 
     
   \bt \la{THEOREM4.12}
The 2-form $\Omega_\Gamma$ assigned to a decorated bipartite ribbon graph $\Gamma$ is invariant under the two by two moves. 
This just means that  the following identity holds, using the notation (\ref{FRM1}):
 \be \la{LRHS}
 \{a_4, a_1\}-\{a_1, a_2\} +  \{a_2, a_3\}-\{a_3, a_4\}  \ \stackrel{}{= }\  \{b_1, b_2\}-\{b_2, b_3\} + \{b_3, b_4\}-\{b_4, b_1\}.   
 \ee
 \et
 
\begin{figure}[ht]
\centerline{\epsfbox{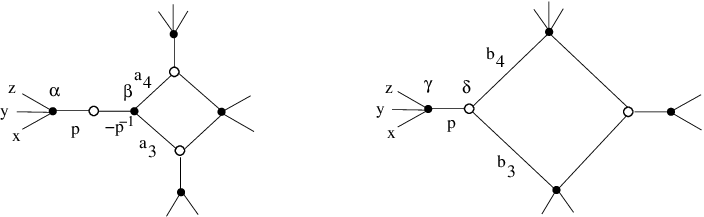}}
\caption{Contributions of  vertices $\alpha$ and $\gamma$ to, respectively,  $\Omega_\Gamma$ and $ \Omega_{\Gamma'}$ coincide: $\Omega_\alpha = \Omega_\gamma$. We use the identity $a_3a_4 = (b_3b_4)^{-1}$ to see that the decorations of  right edges of  vertices $\alpha$ and $\gamma$ are the same: $p$.  We use sign convention (\ref{CONVa}). The signs of ${\cal A}-$decorations do not affect the $2-$forms since $\{\pm a, \pm b\}=\{a,b\}$.}
\label{nclsom}
\end{figure}

 \begin{proof} Denote by $\Gamma$ and $\Gamma'$ the two ribbon graphs related by a two by two move. Then we claim that 
 \be \la{LHSRHS}
 \Omega_\Gamma - \Omega_{\Gamma'} =  \Bigl(\{a_4, a_1\}-\{a_1, a_2\} +  \{a_2, a_3\}-\{a_3, a_4\}\Bigr)  \ - \  \Bigl(\{b_1, b_2\}-\{b_2, b_3\} + \{b_3, b_4\}-\{b_4, b_1\}\Bigr).  
\ee
  Indeed, the contribution of the left $\bullet-$vertex on Figure \ref{Amut*} is equal to $-\Omega_\alpha-\Omega_\beta$, where $\alpha$ and $\beta$ are the two $\bullet-$vertices on 
  the left graph $\Gamma$ on Figure \ref{nclsom}. The contribution of the $\bullet-$ and $\circ-$vertices $\gamma$ and $\delta$ 
  on the right graph $\Gamma'$ on Figure \ref{nclsom} is equal to $-\Omega_\gamma+\Omega_\delta$. Since, as Figure \ref{nclsom} shows, $\Omega_\alpha = \Omega_\gamma$, 
  the total contribution of the left vertices on Figure \ref{Amut*} to $ \Omega_\Gamma - \Omega_{\Gamma'} $ is equal to $-\Omega_\beta - \Omega_\delta$.  
  Note that by defolt the contribution of thje two-valent vertices are zero. 
  Similar analysis for the top, right and bottom vertices on Figure \ref{Amut*} shows that their total contribution to the left hand side of (\ref{LHSRHS}) 
  is equal to the right hand side of (\ref{LHSRHS}). On the other hand, the contributions of all other matching vertices of the graphs $\Gamma$ and $\Gamma'$ evidently coincide.

Therefore it remains to prove the identity (\ref{LRHS}). This is done by a tedious   calculation in Section \ref{sec9}.  

Note  that, as we show there, the identity 
follows from a more elementary identity from Theorem \ref{T11.4}. The latter  is the non-commutative analog of the Steinberg relation $d\log (1+X) \wedge d\log(X)$. \end{proof} 
 
 \paragraph{Remarks.} 1. The right picture on  Figure \ref{Amut*} can be obtained from the left one  by the clockwise rotation by $90^\circ$, followed by the
  cyclic shift    $a_i \lms b_{i+1}$. The right hand side of relation (\ref{LRHS}) is obtained from the left hand side by the same cyclic shift.  
  
  2. Figure \ref{nclsom} shows that the identity $a_3a_4 = (b_3b_4)^{-1}$ is crucial to check that the invariance of the 2-form $\Omega_\Gamma$ under two by two moves 
  is equivalent to the identity (\ref{LRHS}).   
  
\paragraph{A variant.}    Take any    ribbon graph $\Gamma$. Assign to its oriented edges ${\bf E}$  elements $a_{\bf E}$ such that 
denoting by $\overline {\bf E}$ the edge ${\bf E}$ with the opposite orientation, we have 
$ 
a_{\bf E} a_{\overline {\bf E}} =1. 
$ 
For each vertex $v\in \Gamma$ we define  $\Omega_{v}$ using the edges are oriented out of the vertex, and 
  set  take the   sum  over all  vertices $v\in \Gamma$ of valency $>1$:
        \be \la{2F}
 \Omega_\Gamma:=  \sum_{v \in \Gamma} \Omega_{v}.
 \ee

\section{Spectral description of  framed   local systems on  surfaces} \la{secc3x}

In Section \ref{SEC2B} we start from recalling some material from  \cite{G}: a 
  class of bipartite ribbon graphs called ${\rm GL}_m-$graphs, and  their spectral surfaces. Then we  describe non-commutative   framed local systems on a decorated surface $\bS$ 
via  framed flat line bundles on the spectral surface $\Sigma$ of a bipartite graph on $\bS$. 
Namely, we show that the main construction of Section \ref{SECT1.1a} can be generalized 
to any  ${\rm GL}_m-$graph. \vskip 1mm
 
  In Sections \ref{SEC2} and \ref{SEC3.2} we  consider the \underline{twisted} variant: describe non-commutative  \underline{twisted} framed local systems on $\bS$ 
via \underline{twisted} framed flat line bundles on the spectral surface $\Sigma$.
  
  In Section \ref{SEC2} we consider the key example, when   $\bS$ is a triangle $t$. 

In Section \ref{SEC3.2} we consider the general case. 


\subsection{Spectral description of non-commutative framed local systems} \la{SEC2B}

\paragraph{1. The ${\rm GL_m}$-graphs and spectral surfaces \cite{G}.}
Given a  bipartite ribbon graph $\Gamma$, the 
{\it conjugate ribbon graph $\Gamma^*$}  introduced in \cite{GK} is obtained 
by changing the cyclic order of the edges at every $\bullet$-vertex. 
There is a   bijection:
\be \la{ZZFP}
\{\mbox{zig-zag paths on $\Gamma$}\}~ \longleftrightarrow ~\{\mbox{face paths on $ \Gamma^*$}\}. 
\ee

The conjugate ribbon graph $\Gamma^*$ provides a  topological surface $\Sigma$ with holes, 
whose boundaries are the face paths on $\Gamma^*$,  
described as follows. 
Let us add to the  ribbon graph $\Gamma^*$ 
its punctured faces $F_{\gamma_i}$ bounding the face paths $\gamma_i$. 
Precisely,  if $\gamma_i$ is a face loop, then $F_{\gamma_i}$ is a  disc bounded by this loop, punctured inside. 
If $\gamma_i$ is a  face interval, then $F_{\gamma_i}$ is a disc; the half of its boundary is 
identified with 
the path $\gamma_i$, and the rest becomes boundary of $\Sigma$. It is handy to puncture it at the marked point. 
We arrive at a  topological surface, called the 
{\it spectral surface} \cite{G}:
\be \la{2.13.12.1}
{\Sigma} = \Gamma^* \cup\cup_{\gamma_i}F_{\gamma_i}. 
\ee

Let $\bS$ be a decorated surface 
with  marked points $\{s_1, ..., s_n\}$. 
A {\it strand} on $\bS$ is either an oriented 
loop, or 
an oriented path on $\bS - \{s_1, ..., s_n\}$ connecting  boundary points.  
Let $\gamma$ be a strand such that $\bS-\gamma$ has two connected components: 
\be \la{6.20.11.1}
\bS - \gamma = S^\circ_{\gamma} \cup {S^\circ_{\gamma}}'.
\ee
They inherit orientatations from $\bS$.  
The domain $S^\circ_{\gamma}$ is the one whose orientation 
induces the original (clockwise) orientation of 
 $\gamma$. 
 We denote by $S_\gamma$ the closure of  $S_\gamma^\circ$. 

\begin{figure}[ht]
\centerline{\epsfbox{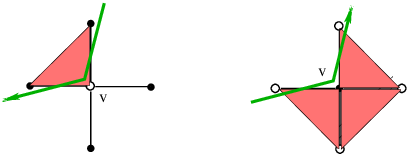}}
\caption{A (green) oriented  zig-zag strand $\gamma$  near a vertex $v$, deformed slightly off the vertex.}
\label{ncls1a}
\end{figure}

\vskip 3mm
Let $\bS$ be a decorated surface associated with a bipartite ribbon graph $\Gamma$, so that $\Gamma \subset \bS$. 
Given a zig-zag path $\gamma$ on $\Gamma$, the shape of the domain 
$S_\gamma$ near a vertex $v$ of $\Gamma$ depends  on   the vertex type. Namely, the intersection $S_\gamma\cap U_v$ with a neighborhood 
$U_v$ of  $v$ is a single sector if $v$ is a $\circ$-vertex, and a union of ${\rm val}(v)-1$ sectors if $v$ is a $\bullet$-vertex, see Figure \ref{ncls1a}. 

Pushing a bit zig-zag paths out of the vertices as  on Figure \ref{ncls1a}, we get {\it zig-zag strands} on $\bS$. 

Let  $\{\gamma\}$ be the collection  of zig-zag strands on $S$. 
We assume     that  for each zig-zag strand $\gamma$

\begin{itemize}

\item  $\bS-\gamma$ has two connected components, and the domains $S_\gamma$ are discs. \end{itemize}

Using   bijection (\ref{ZZFP}), we  construct the spectral surface ${\Sigma}$ as follows. Denote by $F_\gamma$ a copy of the disc $S_\gamma$. 
We attach  the  disc $F_\gamma$ to every zig-zag path $\gamma$ on  $\Gamma$. Since the $F_\gamma$ are discs, 
we  get a surface homeomorphic to  surface (\ref{2.13.12.1}). Canonical maps $F_\gamma \lra S_\gamma$ provide    a  ramified map:
\be \la{9.8.11.10}
\pi: {\Sigma} \lra {\bS}.
\ee

\bd 
The degree $m$ of the map $\pi$ is called the {\it rank} of the bipartite ribbon graph $\Gamma$. 
\ed

The map $\pi$ has the following properties: 
\begin{itemize}
\item The  
ramification points of $\pi$ are the $\bullet$-vertices of $\Gamma$ of valency $\geq 3$:  
the ramification index of a $\bullet$-vertex equals its valency minus $2$: this is clear from Figure \ref{ncls1a}. 

\item For every marked point $s \in {\bS}$ there is a bijection
\be \la{9.8.11.1}
\varphi_s: \{1, ..., m\} \lra \pi^{-1}(s).
\ee
 \end{itemize}

Indeed, the set $\pi^{-1}(s)$ is identified canonically with the set of zig-zag loops/arcs going around the point $s$. Each of the latter is the  boundary of a disc/angle containing 
$s$, which are ordered by the inclusion. 

\bd\la{Def0} Let $\Gamma \subset \bS$ be a bipartite graph associated with a decorated surface $\bS$.

(i) A zig-zag strand $\gamma$ on $\bS$ is   \underline{ideal } if  
    the domain  $S^\circ_{\gamma}$ in (\ref{6.20.11.1})  contains a single marked point $s$. 
We say that { $\gamma$ is associated to $s$}. 

\begin{figure}[t]
\centerline{\epsfbox{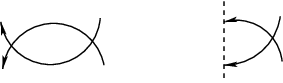}}
\caption{No parallel bigons (left) and half-bigons (right; the boundary of $\bS$ is punctured). }
\label{ncls12}
\end{figure}

(ii) 
   $\Gamma$ is a
\underline{strict ${\rm GL}_m$-graph} if its rank is $m>1$, all zig-zag strands are ideal, and: 

\begin{itemize}

\item There are $m$  zig-zag strands associated with each marked point.
\item There are no parallel bigons and parallel half-bigons, see Figure \ref{ncls12}.

\end{itemize}

(iii)   $\Gamma$   is a
\underline{${\rm GL}_m$-graph} if for some finite unramified cover $\rho: \widetilde \bS\to \bS$,  
$\rho^{-1}({\Gamma})$ is a strict ${\rm GL}_m$-graph.
\ed

For example, the bipartite ribbon graph $\Gamma_m({\cal T})$ assigned to an ideal triangulation ${\cal T}$ of $S$ just before Theorem \ref{TH2.2} is  a  ${\rm GL}_m$-graph. 
 \vskip 2mm

For a strict ${\rm GL}_m$-graph $\Gamma$ on $\bS$,  
the {\it codistance} $\langle x, s\rangle$ from a point $x\in \bS$ to a special point $s$ is the number 
of zig-zag strands $\gamma$ assocated with $s$ such that $x$ belongs to the \underline{open} disc $S^\circ_\gamma$. 
For example, the codistance  to a special point   $s$ from a point $x$ very close to  $s$ is $m$.

\bl \la{1.21.12.1d}  Let $\Gamma$ be a strict ${\rm GL}_m$-graph  on   
 $\bS$. 
Then 
  zig-zags associated with a marked point $s$ form $m $ concentric 
(half-) circles, going counterclockwise 
around $s$.

For any point $x \in \bS- \Gamma$, 
the sum of codistances 
from $x$ to the marked points  is equal to  ${\rm rk}(\Gamma)$:
\be \la{DF}
 \sum_s\langle x, s\rangle = {\rm rk}(\Gamma).
\ee
\el

\begin{proof} Formula (\ref{DF}) is equivalent to the definition of  rank as the degree of the covering  $\pi$. 

Alternatively, formula (\ref{DF}) is evident for a point $x$ close to a marked point: in this case all codistances but one are $0$, and the    non-zero one is  $m$. 
When   $x$ moves,  crossing an edge $E$ in an internal point, we alter codistances as follows: one of them  increases by $1$,   another decreases by $1$,   the rest stay intact.   \end{proof}


\paragraph{2. Spectral description of non-commutative framed local systems.} Given a decorated surface $\bS$, recall the notion of a {\it framing} on a  local system of finite dimensional vector spaces over a skew field $R$, see Definition \ref{9.9.11.1}. 
A framing is a choice of a flat section of the associated local system of flags near each marked point. A local system on $\bS$ is framed if it is equipped with a framing. \vskip 2mm

The monodromy of  an $m$-dimensional framed $R$-local system ${\cal V}$ on $\bS$ 
around a puncture $p$  is  an operator $N$ 
in    
the fiber ${\cal V}_x$  over a nearby point $x$  preserving the framing   
flag ${\cal F}$ in ${\cal V}_x$. The   triple $({\cal V}_x, {\cal F}, N)$ 
determines $m$ operators, called  the {\it monodromy eigenvalue operators}  at  $p$:
 $$
N_i: {\rm gr}^i{\cal F}\lra {\rm gr}^i{\cal F}, ~~ i=0, ..., m-1.
$$
If $R$ is commutative, these are the eigenvalues 
of $N$, which determine 
the semisimple part of $N$. 
  
The points of ${\Sigma_\Gamma}$ projecting to the marked points of ${\bS}$ are   
the {\it marked points of ${\Sigma_\Gamma}$}. 

\bt \la{MTII} 
 Let $\Gamma$ be a ${\rm GL_m}-$graph on a decorated surface $\bS$,    homotopy equivalent to $\bS$. We assume that it can be related by a sequence of two by two moves to a 
 ${\rm GL_m}-$graph $\Gamma_m({\cal T})$ assigned to an ideal triangulation of $\bS$.\footnote{We do not know whether  any  ${\rm GL_m}-$graph on a decorated surface $\bS$ has this property. This is the case when $\bS$  is a polygon by \cite[Section 3]{G}. Since the collection of  ${\rm GL_m}-$graphs with this property is equivariant under the action of the mapping class group of $\bS$, Theorem \ref{MTII} is sufficient for all our purposes.}   Then  the following two  groupoids are canonically birationally equivalent

\begin{enumerate}

\item 
Groupoid of  $m$-dimensional    framed $R$-local systems on $\bS$;

\item  Groupoid of flat   $R$-line bundles with connection on the 
 spectral surface  $\Sigma_\Gamma$. 
\end{enumerate}

The equivalence identifies the monodromy eigenvalue operators of  the connection on ${\cal V}$  
around a puncture $p$ with the monodromy operators  
of  
${\cal L}$ around 
the punctures on $\Sigma_\Gamma $ over $p$, so that 
the eigenvalue operator $N_i$ is the monodromy 
of  ${\cal L}$ around the point $\varphi_p(i)\in \pi^{-1}(p)$, see 
(\ref{9.8.11.1}).\et

 \begin{proof} Flat line bundles on the spectral surface $\Sigma_\Gamma$ are the same thing as flat line bundles on the graph $\Gamma$. 
 Follow the  proof of Theorem \ref{MT1}, adjusted to the   graph $\Gamma$,  
 let us define a functor 
\be \la{FUNCTL}
{\cal L}: \{ \mbox{Generic $m-$dimensional framed local systems ${\cal V}$ on $\bS$}\} \lra \{\mbox{Flat line bundles on   $\Gamma$}\}. 
\ee
 Given a vertex $v$ of   $\Gamma$ and a marked point $s$ on $\bS$  with the codistance $\langle v, s \rangle>0$, there is a unique 
 zig-zag $\gamma$ containing   $v$ such that $s \in S_\gamma$. Since $S_\gamma$ is a punctured disc, and the flag ${\cal F}^\bullet(s)$ near $s$  is monodromy invariant,  
 there is a canonical parallel transport of the  
 flag ${\cal F}^\bullet(s)$  to a flag ${\cal F}^\bullet_{v, s}$ in the fiber ${\cal V}_v$ of   ${\cal V}$ at $v$.    
  
\begin{enumerate}

\item  

Let $\bf w$ be a  $\circ-$vertex  of the graph $\Gamma$. Denote by $a_1, ..., a_n$  the codistances 
from the $\bf w$  to the marked points $s_1, ..., s_n$ on $\bS$.  Then
$$
a_1+...+a_n=m-1, \ \ \ \ \ a_i \geq 0.
$$
 We assign to the $\circ-$vertex  $\bf w$ a one dimensional subspace in the vector space ${\cal V}_{\bf w}$:
\be \la{OC1}
{\rm L}^\circ_{\bf w}:= {\cal F}_{{\bf w}, s_1}^{a_1} \cap \ldots \cap {\cal F}_{{\bf w}, s_n}^{a_n}. 
\ee

\item   Let $\bf b$ be a $k-$valent  $\bullet-$vertex of the  graph $\Gamma$. Denote by $b_1, ..., b_n$  the codistances 
from the $\bf b$  to the marked points $s_1, ..., s_n$ on $\bS$.  Then
$$
b_1+...+b_n=m-k+1, \ \ \ \ \ b_i \geq 0.
$$  We assign to the $k-$valent  $\bullet-$vertex  $\bf b$  
 a $(k-1)-$dimensional subspace in the vector space ${\cal V}_{\bf b}$: 
 \be \la{OC2}
{\rm P}_{\bf b}:= {\cal F}_{{\bf b}, s_1}^{b_1} \cap \ldots \cap {\cal F}_{{\bf b}, s_n}^{b_n}. 
\ee  
\end{enumerate}

The condition that for a configuration of flags the intersections (\ref{OC1}) and (\ref{OC2}) have given dimensions (1 or 2) is an open condition. 
To check that the collection of flags satisfying this condition is non-empty,  we note that this follows by Theorem \ref{MT1} for the bipartite graphs $\Gamma_m({\cal T})$.
Since two by two moves provide birational isomorphisms, the claim follows.  \vskip 2mm

 The plane ${\rm P}_{\bf b}$     contains the $k$ lines  ${\rm L}^\circ_{\bf w}$ assigned to the   $\circ-$vertices  $\bf w$ incident to $\bf b$. Indeed,   
  the  codistances $(a_1, ..., a_n)$ and  $(b_1, ..., b_n)$ from the vertices $\bf w$ and $\bf b$  are related as follows: 
  $0 \leq a_i- b_i\leq 1$, and 
 $a_j=b_j+1$ for exactly $k$ indices $j$. For any such an index $j$ 
  the embedding 
${\cal F}^{b_j+1}\hra {\cal F}^{b_j}$ induces the embedding $ 
 i_{{\bf w}\to {\bf b}}: {\rm L}^\circ_{\bf w} \lra {\rm P}_{\bf b}$. 
 
 We  define a line ${\rm L}^\bullet_{\bf b}$ as the kernel of the sum of  these maps  
  over all $\circ-$vertices ${\bf w}$ incident to the  ${\bf b}$:  \be
 {\rm L}^\bullet_{\bf b}:= {\rm Ker}\Bigl(\bigoplus_{{\bf w} - {\bf b}}{\rm L}^\circ_{\bf w}   \stackrel{i_{{\bf w}\to {\bf b}}}{\lra}  {\rm P}_{\bf b} \Bigr).
 \ee
Then there is a canonical map, induced by  projections of the subspace $ {\rm L}^\bullet_{\bf b}$ onto the   summands:
    \be \la{ABOV}
 {\rm L}^\bullet_{\bf b}  \lra  \bigoplus_{{\bf w}}{\rm L}^\circ_{\bf w}.
 \ee
 So we get a flat line bundle ${\cal L}_{\cal V}$  on the graph $\Gamma$. Its fibers at the vertices are given by the lines ${\rm L}^\circ$ and ${\rm L}^\bullet$, and the 
parallel transport along an edge  $\bullet \to \circ$ is given by a map ${\rm L}^\bullet \lra {\rm L}^\circ$ in (\ref{ABOV}). \vskip 2mm 

{\bf Remark.} If  $R=K$ is commutative, the moduli space of flat line bundles on a graph $\Gamma$ is identified with  $H^1(\Gamma, K^*)$. 
The latter include  the  monodromies around the 
  holes in the graph. These are   the coordinates from \cite{G}, generalizing the coordinates \cite{FG1}  to ideal bipartite graphs on $\bS$. 

If  $R$ is non-commutative, we no longer have  canonical  coordinates.  \vskip 2mm

   Let us show that the functor ${\cal L}$ is an equivalence. 
   Pick an ideal triangulation ${\cal T}$ of $\bS$. Recall the
 ${\rm GL}_m-$graph $\Gamma_m({\cal T})$ related to ${\cal T}$. By Theorem \ref{TH2.2}, applied to the surface ${\cal S}$ 
  glued from the triangles of ${\cal T}$,  which is just the decorated surface $\bS$, 
  the functor $(\ref{FUNCTL})$ is an equivalence of groupoids, where "generic" means that 
 the triple of flags for each triangle of ${\cal T}$ is generic in the sense of Section \ref{SECT1.1a}. 
Then the claim for a graph $\Gamma$ related to $\Gamma_m({\cal T})$ by a sequence of two by two moves follows from  results of Section \ref{SEC5.1}. 
\end{proof}


  \subsection{Spectral description of generic  twisted triples of flags} \la{SEC2}

 Let us replace a triangle $t$   by an oriented circle $S^1$ with three marked points $\A, \B, {\rm C}$. 
Twisted framed local systems   on $(S^1; \A, \B, {\rm C})$ are the same thing as   the  flat line bundle on $S^1$ with the monodromy $-1$, and  flags   $({\cal A}, {\cal B}, {\cal C})$ 
 in the fibers over the tangent vectors $v_\A, v_\B, v_{\rm C}$ to $S^1$ at the points   $\A,\B, {\rm C}$, following the orientation of $S^1$. 
We call such data {\it twisted triples of flags}. 
There is a canonical one-dimensional local system $\varepsilon$ on $S^1$ with the monodromy $-1$. Multiplying by $\varepsilon$, we identify 
the usual and twisted triples of flags. So there is almost no difference between the two notions. The distinction between the two will be important when we move to 
\underline{decorated} twisted 
configurations of flags.\vskip 2mm

Take a generic triple of  flags $({\cal A}, {\cal B}, {\cal C})$ in an $m$-dimensional $R$-vector space. We assign the flags  to the vertices of a triangle $t$. 
The bipartite ribbon graph $\Gamma_m$ on   Figure \ref{ncls3a} assigned to $t$ gives rise to the spectral surface $\Sigma_m$.  
Namely, 
each zig-zag $\gamma$ on  $\Gamma_m$ determines a smaller triangle $t_\gamma$,  given by the component of $t - \gamma$   
  containing a single   vertex of  $t$. Take the disjoint union of all triangles $t_\gamma$. Then, for each edge $E$ of the graph $\Gamma_m$,   glue 
the triangles $t_{\gamma'}$ and $t_{\gamma''}$ assigned to the  zig-zags $\gamma'$ and $\gamma''$ containing $E$ over the edge $E$. We get the surface $\Sigma_m$ and the projection 
$\Sigma_m \lra t$, ramified over the $\bullet-$vertices. To avoid confusion, we denote by $F_\gamma$ the triangle $t_\gamma$ when it sits on the spectral surface.

\bt \la{UHH}
There is a canonical equivalence between the following two groupoids:
\begin{enumerate}
 \item Groupoid of generic triples of flags in an $m$-dimensional $R$-vector space.
 
\item  Groupoid of flat \underline{twisted} $R$-line bundles  with connection  on the 
 spectral surface $\Sigma_m$ of  
  the   graph $\Gamma_m$. 
\end{enumerate}
\et

\begin{proof} $1\lra 2$. A triple of flags $({\cal A}, {\cal B}, 
{\cal C})$ in an $m$-dimensional $R$-vector space $V$ gives rise to $3m$ lines
\be \la{LINES}
{\rm gr}^a{\cal A}, ~~ {\rm gr}^b{\cal B}, ~~ {\rm gr}^c{\cal C},~~ a,b,c=0, ..., m-1. 
\ee
Given a  zig-zag $\gamma$ on the   graph $\Gamma_m$   going around of the vertex $\A$, we assign to the  face $F_\gamma$ of the spectral surface the line 
${\rm gr}^{a_\gamma-1}{\cal A}$, where $a_\gamma \in \{1, ..., m\}$ is the codistance from $\gamma$ to   $\A$. 
Abusing notation, we denote by ${\rm gr}^{a_\gamma-1}{\cal A}$ the
 twisted flat line bundle on  
  $F_\gamma$ whose restriction to   oriented  tangent vectors   to the oriented side $\B{\rm C}$ is given by the line ${\rm gr}^{a_\gamma-1}{\cal A}$. We do the same for the vertices $\B, {\rm C}$. 
Let us glue these twisted flat line bundles on the  $3m$ faces $F_\gamma$ to a twisted line bundle on the spectral surface. \\

Note that the edges of the graph $\Gamma_m$ are parametrised   by the following data: 

i) A triple  of non-negative  integers $(a,b,c)$ such that $a+b+c=m-1$,  plus

 ii) A choice of two of the three numbers $(a,b,c)$.  \vskip 2mm
 
The faces of the graph $\Gamma_m$ are labeled   by the triples of non-negative integers $(a,b,c)$ with  $a+b+c=m$. \\

Given  a data  $(\underline{a}, \underline{b}, c)$   assigned to an edge $E$,  where the chosen pair  in ii) is the pair $(a,b)$,  
  the two faces $F_{\gamma'}$ and $F_{\gamma''}$  of the spectral surface sharing the edge $E$  are labeled by the triples $(a,b+1,c)$ and $(a+1, b,c)$, accordantly. 
  Then the  faces $F_{\gamma'}$ and $F_{\gamma''}$ carry the line bundles 
$$
\L_{\gamma'}:= {\rm gr}^{a}{\cal A}, \ \ \ \ \ \ \L_{\gamma''}:= {\rm gr}^{b}{\cal B}.
$$
We assign to the edge $E$ the line 
$$
\L_E:= {\cal A}^{a}\cap  {\cal B}^b\cap   {\cal C}^c. 
$$
Then there are canonical  isomorphisms,  obtained by projecting along ${\cal A}^{a+1}$ and ${\cal B}^{b+1}$:
 \be \la{ISO2F}
 {\rm gr}^{a}{\cal A}  \stackrel{\sim}{\longleftarrow} {\cal A}^{a} \cap {\cal B}^{b} \cap  {\cal C}^c\stackrel{\sim}{\longrightarrow}
 {\rm gr}^{b}{\cal B}.  
 \ee
 They allow to  glue the twisted flat line bundles over the faces  $F_{\gamma'}$ and $F_{\gamma''}$ along their common edge $E$. 
 
Let us   extend  the obtained twisted flat line bundle on $\Sigma_m-\{\mbox{the vertices of $\Gamma_m^*$\}}$ to the spectral surface.   
In order to do this, we make the following   observations. \vskip 2mm

For any triple $(a,b,c)$ such that $a+b+c=m-1$ there are three canonical isomorphisms:
\be \la{CANISO}
\begin{gathered}
    \xymatrix{
    & {\rm gr}^b{\cal B}  &\\
 &  \ar[ld]_{\sim} {\cal A}^{a} \cap {\cal B}^{b} \cap {\cal C}^{c}  \ar[u]^{\sim}  \ar[rd]^{\sim}  &\\
  {\rm gr}^a{\cal A}    &&   {\rm gr}^c{\cal C} 
        }
\end{gathered}
 \ee
The first is  the composition  ${\cal A}^{a} \cap {\cal B}^{b} \cap {\cal C}^{c} \hra {\cal A}^{a} \to {\cal A}^{a}/{\cal A}^{a+1}$. The others are similar.

\bl \la{HEXTR} There are the following relations between the isomorphisms (\ref{CANISO}):

\begin{enumerate}
\item   For any integers $a,b,c \geq 0$  such that $a+b+c=m-1$ there is a commutative triangle, i.e. the  composition of the  three isomorphism is the identity:

\be \la{TRIANGLE}
\begin{gathered}
    \xymatrix{
    &&   {\rm gr}^b{\cal B}  \ar[rd]   &&\\
 &{\rm gr}^{a}{\cal A}   \ar[ru]    &   {\xleftarrow{\makebox[2cm]{  }} }   &      {\rm gr}^{c}{\cal C}  &\\
        }
\end{gathered}
 \ee
 
 \item For any  integers $a,b,c \geq 0$  such that $a+b+c=m-2$ there is a sign-commutative   hexagon,  i.e.   the   composition of the following  six isomorphisms/their inverces   is equal to $-1$:
\be \la{HEXAGON}
\begin{gathered}
    \xymatrix{
    && \ar[ld]_{\sim} {\rm gr}^{a+1}{\cal A}  \ar[rd]^{\sim}  &&\\
 &{\rm gr}^{b}{\cal B}    &     &   {\rm gr}^{c}{\cal C}  &\\
    &\ar[rd]_{\sim} {\rm gr}^{c+1}{\cal C} \ar[u]^{\sim} && \ar[u]_{\sim}  \ar[ld]^{\sim}{\rm gr}^{b+1}{\cal B}  &\\
    && {\rm gr}^{a}{\cal A}  &&    }
\end{gathered}
 \ee

\end{enumerate}
\el

\begin{proof} 1) Starting with a pair of isomorphisms
 ${\rm gr}^a{\cal A}\stackrel{\sim}{\longleftarrow}  {\cal A}^{a} \cap {\cal B}^{b} \cap {\cal C}^{c} \stackrel{\sim}{\longrightarrow} {\rm gr}^b{\cal B}$, and composing the inverse of the 
 first one with the  second, we get  an isomorphism ${\rm gr}^a{\cal A}\stackrel{\sim}{\lra} {\rm gr}^b{\cal B}$. The other two are similar. 
This makes the first claim is evident. 

2) Intersecting with the plane ${\cal A}^{a} \cap {\cal B}^{b} \cap {\cal C}^{c}$, we reduce the claim to the following.
 Let $A, B, C$ be three 
different lines  in a plane ${\rm P}$. 
Consider the quotient lines
$$
\overline A:= {\rm P}/A, ~~~ \overline B:= {\rm P}/B, ~~~ \overline C:= {\rm P}/C.
$$ 
Projecting the line $A$ 
to the quotient line $\overline B$ we get an isomorphism 
$
A \lra \overline B.
$
There are six such line isomorphisms: 
\be \la{MOREGON}
\begin{gathered}
    \xymatrix{
    &   \overline B   &  & \overline A    &&    \overline C    & \\
A   \ar[ru]     &&      \ar[lu]       C   \ar[ru]   &&   \ar[lu]  B    \ar[ru] && \ar[lu]    A \\
        }
\end{gathered}
 \ee
 We compose them to a map $A  \to A $, where going against 
an arrow means inverting the isomorphism. 
We claim that the composition is the minus identity map.  Indeed, 
choose a triple of non-zero vectors $x  \in A, y  \in B, z  \in C$ whose sum is zero: 
 \be \la{abc}
x +y +z  =0. 
 \ee
Such a triple exists, and is defined uniquely up to a multiplication from the left by a non-zero scalar. 
Condition (\ref{abc}) just means that $x +z $ projects to zero in $\overline B$.  
So the map  $A  \lra C$ provided by the diagram $A \stackrel{\sim}{\lra} \overline B \stackrel{\sim}{\longleftarrow} C $ is given by $x  \lms - z$. 
Therefore the composition 
$A  \lra C  \lra B  \lra A $ is given by 
$x  \lms -z  \lms y  \lms -x $. \end{proof}

The part i) of Lemma \ref{HEXTR} allows to extend the twisted flat line bundle to the $\circ-$vertices. Indeed, although the composition of the three maps on Figure \ref{ncls01} is the identity map, 
the tangent vectors are rotated by $2\pi$, thus resulting the $-{\rm Id}$ map. 
  \begin{figure}[ht]
\centerline{\epsfbox{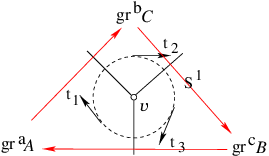}}
\caption{Let $a+b+c=m-1$, and the line ${\rm gr}^a{\cal A} $ sits over the tangent vector $t_1$, etc..  By (\ref{TRIANGLE}), the composition ${\rm gr}^a{\cal A} \to {\rm gr}^b{\cal B} \to {\rm gr}^c{\cal C}\to {\rm gr}^a{\cal A}$ is the identity map.  The   tangent vector $t_1$ rotates by 
 $t_1 \to t_2\to t_3\to t_1$ by $2\pi$.}
\label{ncls01}
\end{figure} 

Similarly, the part ii) of of Lemma \ref{HEXTR} allows to extend the twisted flat line bundle to the $\bullet-$vertices. Indeed, although the composition of the six maps (\ref{HEXAGON})  is $-{\rm Id}$, 
the tangent vectors are rotated by $4\pi$, thus still resulting the $-{\rm Id}$ map. \\

$2 \lra 1$.  {\it  The   functor    ${\cal R}:  \{\mbox{Twisted flat line bundles on $\Sigma_m$}\} \lra \{\mbox{generic twisted triples of flags}\}$.} 
Take a twisted flat line bundle ${\cal L}$ on the spectral surface $\Sigma_m$ of the bipartite ribbon graph ${\Gamma_m}$. It  contains the graph $\Gamma_m^*$.  
For each vertex $v\in \Gamma_m^*$,   restrict ${\cal L}$ to  the positively oriented tangent vectors to an oriented circle  $S^1_v\subset T_v\Sigma-\{0\}$, getting a flat line bundle   ${\cal L}_{S^1_v}$ on  
$S^1_v$ with the monodromy $-1$. 

Then the flat connection on ${\cal L}$ provides  a natural map 
\be \la{Kast}
{\rm K}_{\cal L}: \bigoplus_{\mbox{internal $\bullet-$vertices ${\bf b}$}} {\cal  L}_{S^1_{\bf b}}  \lra \bigoplus_{\mbox{ $\circ-$vertices ${\bf w}$}} {\cal L}_{S^1_{\bf w}}. 
\ee
We identify all      
circles $S^1_{\bf w}$  with a   circle $S^1$, and  
 define the vector space $V_{\cal L}$ as the cokernel of this map:
\be \la{koc}
V_{\cal L}:= {\rm Coker}({\rm K}_{\cal L}). 
\ee 
It is a flat   bundle on the circle $S^1$ with the monodromy $-{\rm Id}$. 
Let us define three flags in the fibers of this local system over the points  $\A, \B, {\rm C}$.   
The distance  functions $d_{\A}, d_\B, d_{\rm C}$, see the proof of Theorem \ref{MT1},   induce  three filtrations  in the fibers of    complex (\ref{Kast}) 
over the points $\A, \B, {\rm C}$. They induce   filtrations on  
${\rm Coker}({\rm K}_{\cal L})$, providing a twisted triple of flags.    The   functor  ${\cal R}$ is the inverse to the one defined before.    
   \end{proof}

  \paragraph{Remark.}  Let us compare   constructions  of the connection in the proofs of  Theorems \ref{MT1} and   \ref{UHH}. 
The key point is that given  three vectors $x,y,z$ in a two dimensional vector space   with $x+y+z=0$, we can define a connection either by $x \to y \to z \to x$, as we did proving Theorem \ref{MT1}, or by   
$x \to -y \to z \to -x$, getting a twisted connection.

\subsection{Spectral description twisted local systems on surfaces} \la{SEC3.2}

 Below we use the notation $\Sigma$ for the spectral surface $\Sigma_\Gamma$ since the bipartite graph $\Gamma$ is fixed.

\bt \la{9.9.11.3} 
 Let $\Gamma$ be a ${\rm GL}_m-$graph on a decorated surface $\bS$,  
and $\Sigma$ the corresponding spectral surface. 
 Then the following two  groupoids are canonically birationally equivalent:  

\begin{enumerate}

\item 
Groupoid of   $m$-dimensional \underline{twisted} framed $R$-vector bundles with flat connections   ${\cal V}$ on $\bS$.

\item  Groupoid of flat \underline{twisted} $R$-line bundles  with connection ${\cal L}$ on the 
 spectral surface $\Sigma$. 
\end{enumerate}
The equivalence identifies the monodromy eigenvalue operators of  the connection on ${\cal V}$  
around a puncture $p$ with the monodromy operators  
of  
${\cal L}$ around 
the punctures on $\Sigma $ over $p$ just as in Theorem \ref{MTII}.
\et

\begin{proof}   
Let ${\cal V}$ be an $m$-dimensional twisted framed $R$-local system on $\bS$.  
Let us construct a twisted 
flat line bundle   ${\cal L}$ on  $\Sigma$. 
Let 
${\cal F}(s)$ be a flat section of the local system of flags in ${\cal V}$ defining a framing near a marked point $s\in \bS$. 
For each zig-zag  $\gamma$ going around $s$,  
it provides a flat section over the punctured disc $S_\gamma$. Recall the punctured face $F_\gamma$ of the spectral surface $\Sigma$ assigned to $\gamma$. 
By the very construction of the spectral surface, there is a canonical identification:
$$
F_\gamma\stackrel{}{=} S_\gamma.
$$ 
We use it to  transport 
the flat section of the local system of flags  on $S_\gamma$ to   $F_\gamma$, and denote it by ${\cal F}(s)$. 
We assign to a zig-zag $\gamma$ a number $a_\gamma \in \{1, ..., m\}$ -  the codistance from a point of the open disc  $S_\gamma^\circ$ close to   $\gamma$ to the 
marked point associated with  $\gamma$. 
Then there is a twisted flat line bundle  on the face $F_\gamma$ of $\Sigma$:
\be \la{LLFFa}
{\cal L}_{F_{\gamma}}:= {\rm gr}^{a_\gamma-1}{\cal F}(s) 
\ee
 Let us glue  them into a   twisted flat line bundle  ${\cal L}$ on $\Sigma$.

\bl Let $s_1, \ldots , s_n$ be the  marked points    on   $\bS$. Let $y$ be a point of $\Sigma$, and $x:=\pi(y)\in \bS$. 

Let $0 \leq a_j \leq m$  be the codistance 
from $x$  to  $s_j$.  Then we have:
\begin{enumerate}

\item 
   Let $y \in \Sigma -   \Gamma^* $. Then  
\be \la{12.11.7.1a}
a_1+\ldots +a_{n} =m. 
\ee

\item   Let $y \in \Gamma^* - \{\bullet-$vertices$\}$. Then 
\be \la{12.11.7.1b}
a_1+\ldots +a_{n} =m-1. 
\ee

\item   Let $y$ be a $k-$valent $\bullet$-vertex  of  $ \Gamma^*$. Then 
\be \la{12.11.7.1c}
a_1+\ldots +a_{n} =m-k+1. 
\ee
\end{enumerate}

\el

\begin{proof} 1. When a point $x$ is   inside of a face, the claim follows from  (\ref{DF}). 

2. When a point $x'$ moves and hits an inside point $x$ of an edge $E$, just one of the codistances $a_i$ decreases by $1$. Namely, it is the one corresponding to the 
marked point $s_\gamma$ surrounded by the zig-zag $\gamma$ containing $E$, such that the open domain $S^\circ_\gamma$ containes $x'$. 
Same thing happens when $x$ is a $\circ-$vertex. 
Indeed, take a small disc $U_x$ containing $x$. Let $\gamma_1, \ldots , \gamma_k$ be the zig-zags containing  $x$. The   domains 
$S^\circ_\gamma \cap U_x$ are disjoint, and their closures cover  $U_x$, as shown on the left  in Figure \ref{ncls1a}.  

3. When a point $x'$ moves and hits a $\bullet-$vertex $x$, the $k-1$ of the codistances $\{a_i\}$ decrease by $1$. They are the ones corresponding to the zig-zags $\gamma$ 
containing the vertex $x$ such that $x' \in S_\gamma^\circ$.
\end{proof}

{\it Gluing over an edge $E$   of  $\Gamma$}. Let $\gamma_1$ and $\gamma_2$ be the   two zig-zags  passing through  $E$, and 
  $s_1$ and $s_2$ the related marked points of $\bS$. 
Faces $F_{\gamma_1}$ and $F_{\gamma_2}$ on $\Sigma$ intersect along the edge (identified with)  $E$. 
We denote by ${\cal F}^{a_i}$  the codimension $a_i$  subspace of the flag ${\cal F}(s_i)$. 
 Then, thanks to  (\ref{12.11.7.1b}),  $
a_1  + a_2 +  \ldots + a_n=m-1 
 $, and so there is a  twisted flat line bundle near the edge $E$ on $\bS$: 
 \be \la{LEa}
 {\cal L}_E:=  {\cal F}^{a_{ 1}}  \cap  {\cal F}^{a_{ 2}}  \cap  {\cal F}^{a_{3}} \cap \ldots \cap {\cal F}^{a_{n}}, \ \ \ \ \ \ a_1  + a_2 +  \ldots + a_n=m-1. 
\ee
Note that there are canonical projections  
$$
{\cal F}^{a_{1}} \lra {\rm gr}^{a_{1}}{\cal F}(s_1), \ \ \ \ \ \ 
{\cal F}^{a_{2}} \lra {\rm gr}^{a_{2}}{\cal F}(s_2).
 $$
Restricting them to   intersection (\ref{LEa})   we get   isomorphisms of   twisted flat line bundles  near the edge $E$:
\be \la{9.16.17.1}
{\cal L}_{F_{\gamma_1}} \stackrel{(\ref{LLFFa})}{=} {\rm gr}^{a_{1}}{\cal F}(s_1) \stackrel{\sim}{\longleftarrow}    {\cal L}_E  \stackrel{\sim}{\longrightarrow}   {\rm gr}^{a_{2}}{\cal F}(s_2)
 \stackrel{(\ref{LLFFa})}{=}  {\cal L}_{F_{\gamma_2}} .
\ee
Here the very  left and right isomorphisms come from the fact that the codistances from a generic point of $S^\circ_{\gamma_1}$ close to the zig-zag $\gamma_1$ are given by 
$(a_1+1, a_2, \ldots , a_n)$, and for $S^\circ_{\gamma_2}$ they are $(a_1, a_2+1, \ldots , a_n)$. 
 Using isomorphisms (\ref{9.16.17.1}), we glue  twisted flat line bundles (\ref{LLFFa})  on  faces $F_\gamma$ into a twisted flat line bundle     
 $$
 {\cal L}^\times  \mbox{on} \  \Sigma- \{\mbox{vertices of    $\Gamma^*$}\}.
 $$

 \bp\label{MONODROMY}
The twisted flat line bundle ${\cal L}^\times$ extends  to a twisted flat line bundle ${\cal L}$ on  $\Sigma$. 
\ep

\begin{proof}

Take a small  circle $S^1_v$ around a vertex $v\in \Gamma^*$ on the spectral surface $\Sigma$. Restrict the   ${\cal L}^\times$ to a section of the 
bundle of non-zero tangent vectors on $S^1_v$.  We get a flat line bundle ${\cal L}_v^\times$ over the circle $S^1_v$. 

\bl \la{MONODROMY1}

  The monodromy of the flat line bundle ${\cal L}_v^\times$ around every $\circ$-vertex $v$ is equal to $-1$. 
 
\el

    \begin{figure}[t]
\centerline{\epsfbox{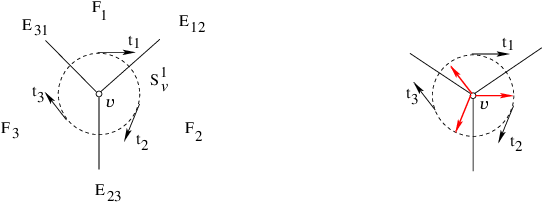}}
\caption{On the left: the edges $E_{i, i+1}$ and faces $F_i$ of $\Gamma^*$ near a $\circ$-vertex $v$.}
\label{ncls0}
\end{figure} \begin{proof}
  Denote by $E_{i, i+1}$ the edge of  $\Gamma^*$ shared by the faces $F_i$ and $F_{i+1}$, see Figure \ref{ncls0}. 
Take a non-zero tangent vector $t_i$ at a point of   $S_v^1 \cap F_i$ in the   clockwise direction.  
Denote by ${\cal L}_{t_i}$ the fiber of the twisted local system ${\cal L}_{F_i}$ over $t_i$. To define the parallel transform 
 ${\cal L}_{t_i} \to {\cal L}_{t_{i+1}} $, we  move tangent vectors $t_i$ and $t_{i+1}$ towards each other to the tangent vector $t_{i, i+1}$ at the 
 point $S^1_v \cap E_{i, i+1}$, and use 
   isomorphisms (\ref{9.16.17.1}). Denote by $\widetilde t_i$ the tangent vectors at the vertex $v$ parallel to the tangent vectors $t_i$, see the red vectors on the right picture on Figure \ref{ncls0}. 
   We have canonical isomorphisms  ${\cal L}_{t_i} = {\cal L}_{\widetilde t_i}$ provided by the parallel transform of tangent vectors $t_i\to \widetilde t_i$. After these  identifications, 
   the parallel transform of the vector $t_i$ around the circle $S_v^1$ amounts to the parallel transform ${\cal L}_{\widetilde t_i}$ for the   rotation of the  tangent vector $\widetilde t_i$ by  $360^\circ$,    
 resulting the  $-1$ monodromy.    
\end{proof}
 
Lemma \ref{MONODROMY1}  implies Proposition \ref{MONODROMY} a $\circ$-vertex. \\

Now let ${\bullet}$ be a  $\bullet$-vertex of $\Gamma$.  Let $v_{1}, \ldots , v_{k}$ be the $\circ$-vertices incident to ${\bullet}$, 
whose order is compatible with the cyclic structure at ${\bullet}$. 
A  zig-zag  $ \ldots v_i{\bullet}v_{i+1} \ldots $, where  $i \in \Z/k\Z$,   passing via the vertex ${\bullet}$ 
 determines  a   {\it $\bullet$-cosector}, see Figure \ref{ncls1}:
$$
C_{v_i{\bullet}v_{i+1}}:= U_{\bullet} \cap S_\gamma.
$$
The flat line bundle over this  $\bullet$-cosector     is denoted by 
$
{\cal L}_{ C_{v_i {\bullet}v_{i+1}}}.
$ 
Next, for each $\circ$-vertex $v_i$  consider  the  domain $C_{v_i}$ given by the intersection  of $U_{v_i}$ with an angle formed by the two edges at $v_i$ 
on the left and on the right of the edge ${\bullet}v_i$, see Figure \ref{ncls1}. Denote the twisted flat line bundle over this domain  by ${\cal L}_{{ C}_{v_i}}$. 
\begin{figure}[ht]
\centerline{\epsfbox{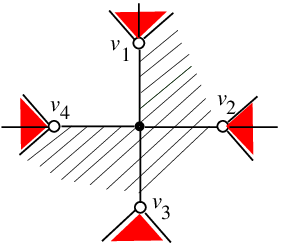}}
\caption{The shaded domain is the cosector $ C_{v_4 {\bullet} v_1}$. The  red ones are  domains $C_{v_i}$. }
\label{ncls1}
\end{figure} 

For each $i=1, ..., k$ the twisted flat structure of ${\cal L}$ provides   a pair of  
isomorphisms 
$$
{\cal L}_{C_{v_i}} \longrightarrow {\cal L}_{C_{v_i {\bullet} v_{i+1}}} \longleftarrow {\cal L}_{C_{v_{i+1}}}. 
$$
They combine into a diagram of $2k$ line isomorphisms. For example, for $k=3$ we get 
\be \la{MOREGON1}
\begin{gathered}
    \xymatrix{
    &     {\cal L}_{ C_{v_1 {\bullet} v_2}}   &  &    {\cal L}_{  C_{v_2 {\bullet} v_3}}   &&     {\cal L}_{ C_{v_3 {\bullet} v_1}}  \ar[rd]  & \\
{\cal L}_{C_{v_1}}  \ar[ur]   &&      \ar[ul]    {\cal L}_{C_{v_2}}  \ar[ur] &&  \ar[ul]  {\cal L}_{C_{v_{3}}}  \ar[ur]  &&  \ar[ul]   {\cal L}_{C_{v_1}}\\
        }
\end{gathered}
 \ee
The same argument as in the proof of Lemma \ref{HEXTR}(ii)  shows that the composition is  $(-1)^{k}$.

Moving around the circle $S^1_{\bullet}$ on the spectral surface $\Sigma$  we rotate the projection to $\bS$ of the tangent vector 
to the circle by $({\rm val}({\bullet})-1)2\pi$. The composition of $2{\rm val}(\bullet)$ line isomorphisms in (\ref{MOREGON1})  is $(-1)^{{\rm val}( {\bullet})}$. Since ${\rm val}({\bullet})-1 + {\rm val}({\bullet})$ is odd,  the result is multiplied by $-1$. So we can extend the twisted flat line bundle ${\cal L}^\times$ across the $\bullet$-vertex. Proposition \ref{MONODROMY} is proved. 
\end{proof}


To  prove that the constructed map is a birational equivalence of groupoids, we proceed similarly to the proof of  equivalence in Theorem \ref{UHH}.  
  Theorem \ref{9.9.11.3} is proved.
\end{proof}

\section{Spectral description of   twisted decorated local systems}\la{sec4}

\subsection{The spectral description} \la{sec4.1}

We use the same conventions 
about  twisted  local systems on $\bS$ as in Section \ref{SEC3.2}.  We recall decorated flags, defined in Section \ref{SECT1.1}.2.

\bd \la{9.9.11.2}
A \underline{\em twisted decorated} $R$-local system on a decorated surface $\bS$  
is a twisted local system  of finite dimensional $R-$vector spaces on $\bS$ 
with  
a flat section of the associated local system of decorated flags over each marked boundary component.
\ed

Let $\Gamma$ be a ${\rm GL_m}$-graph on  $\bS$, homotopy equivalent to $\bS$. Denote by 
$\Sigma$ its spectral surface. 

\bt \la{12.9.11.1z} Given a ${\rm GL_m}$-graph $\Gamma$ on $\bS$, homotopy equivalent to $\bS$, there is a birational equivalences of groupoids 
of non-commutative local systems:
\be \la{eq3}
\begin{split}
&\{\mbox{$m$-dimensional  \underline{twisted decorated}  local systems ${\cal V}$ on $\bS$}\} \longleftrightarrow \\
&\{\mbox{Flat twisted  line bundles ${\cal L}$ on   
  ${\Sigma}$, trivialized 
on the boundary}\}.\\
\end{split}
\ee 
\et

\begin{proof}
Equivalence (\ref{eq3}) is obtained by applying the 
equivalence from Theorem \ref{9.9.11.3} to    twisted  framed  local systems on $\bS$, equipped with an extra data: a decoration.\footnote{Note that a decorated local system is unipotent.} 
 The description of the eigenvalue monodromy operators  in Theorem \ref{9.9.11.3} implies that 
the twisted flat   line bundle ${\cal L}$ on ${\Sigma}$ assigned to  ${\cal V}$ 
  has the following property: 

\begin{itemize}
\item The restrictions of ${\cal L}$ to punctured discs around the   
 marked points in $\pi^{-1}(s)$ 
are identified with the twisted  flat line bundles ${\rm gr}^a{\cal A}$ 
associated to the decoration of ${\cal V}$ near $s$.  
\end{itemize}
So a choice of a decoration on   ${\cal V}$ near a marked point  $s$ is equivalent to a choice of trivializations of the restriction of ${\cal L}$ to the homotopy circles near  the marked points $\pi^{-1}(s)$.  
 \end{proof}

\paragraph{\it A coordinate description.} 
Recall that $\Gamma =   \Gamma^*$ as graphs, and  faces of 
$ \Gamma^*$ match zig-zags on $\Gamma$. 
 By Theorem \ref{CRTH},  the isomorphism classes of \underline{twisted} flat line bundles on the spectral surface ${\Sigma}$, trivialised on the marked boundary,  
are described 
by  a collection of invariants $\Delta_\bE\in R^\times$ at the oriented edges $\bE$ of $\Gamma$,
subject to the  monomial relations  (\ref{MEQ}).

Since  the edges have canonical orientation $\circ\to \bullet$, given a ${\rm GL}_m$-graph $\Gamma$ on   $\bS$, we get a collection of   
elements $\Delta_E\in R^\times$ ,  assigned to the   edges $E$ 
of $\Gamma$. They describe generic 
twisted decorated  local systems of $m-$dimensional $R-$vector spaces on $\bS$. They are  non-commutative counterparts of the ratios of the 
cluster ${\cal A}$-coordinates \cite{FG1}, introduced in the context of ${\rm GL}_m$-graphs $\Gamma$ in \cite{G}. 

Precisely, in the commutative case, for the twisted 
${\rm SL}_m-$local systems, the  ${\cal A}$-coordinates $\{A_{\rm F}\}$ are assigned to the faces $\rm F$ of the graph $\Gamma$. 
Then $\Delta_E = A_{{\rm F}_+}/A_{{\rm F}_-}$, where $\rm F_+$ and ${\rm F}_-$ are the two faces sharing the edge $E$, with the face ${\rm F}_+$ on the right of the oriented $\circ\to \bullet$ edge.   However the functions $\{\Delta_E\}$ 
do not form a coordinate system: they satisfy  monomial relations (\ref{MEQ}).  \vskip 2mm

In Section \ref{sec4.2a} we elaborate the $m=2$ example.

\subsection{Gluing   spectral surfaces for   ${\rm GL_2}-$graphs from hexagons} \la{sec4.2a}

 The spectral surface $\Sigma$ of a ${\rm GL}_m-$graph $\Gamma$ was discussed in Section \ref{SEC3.1}. 
 Below we revisit the  case of a ${\rm GL}_2-$graph from a more elementary  perspective, and use this to describe 
 non-commutative twisted decorated ${\rm GL}_2-$local systems on a decorated surface $\bS$ via flat line bundles on the spectral surface. 
 Then in Section \ref{sec4.2} we elaborate this description in detail in coordinates, relating to the Gelfand-Retakh quasideterminants and algebraic relations between them.

\paragraph{1. The spectral surface for the triangle.} Let $\Gamma$ be the ${\rm GL_2}$-graph associated to the triangle $t$, see Figure \ref{bgc11}. 
The spectral surface of $\Gamma$ 
is a hexagon $h_t$ providing a $2:1$ cover of the triangle $t$, ramified at  
the central $\bullet$-vertex. 
The map $h_t\to t$ can be realized by 
the map $z \to z^2$ sending a regular hexagon in the complex plane $z$ 
with the  vertices at the sixth roots of unity  onto a triangle with vertices at 
the cubic roots of unity. The vertices of the hexagon $h_t$ are labeled by elements of the set $\{1, 2\}$, 
so that going around the hexagon the labels alternate. 
For  each vertex $v$ of $t$ there is a bijection
$$
\{\mbox{The vertices of the hexagon $h_t$ 
  over $v$}\} \stackrel{\sim}{\lra} \{1,2\}.
$$   
The sides of the hexagon inherit a ``$1\to 2$'' orientation.

\paragraph{2. Gluing the spectral surface from hexagons.} 
Let ${\cal T}$ be an ideal triangulation of $\bS$. For each triangle $t$ of
 ${\cal T}$ choose a hexagon $h_t$ over  $t$ with the vertices labeled as above. 
Let us glue the hexagons  
$h_{t}$ into a complete smooth topological surface ${\bf \Sigma}$. 
Given an edge $E$ of  ${\cal T}$, there are 
two triangles $t$ and $t'$ sharing $E$. 
Let $h$ and $h'$ be the hexagons covering these triangles. 
We identify the vertices of the hexagons  with the same labels,  
projecting onto the same vertex of ${\cal T}$. Then we glue   two edges 
 whose vertices are identified. For example, gluing  two hexagons assigned to a triangulation of a rectangle we get an annulus.
Denote by ${\bf S}$ a complete surface obtained by filling the punctures on  $\bS$. 
Gluing hexagons $h_t$ we get a  surface ${\bf \Sigma}$ with the following properties: 
\begin{itemize}

\item There is a $2:1$ map 
$
\pi: {\bf \Sigma} \to {\bf S}
$ 
with  one ramification point in every triangle of ${\cal T}$. 

\item For every marked point $s \in {\bf S}$ there is a bijection $\pi^{-1}(s)\lra \{1,2\}$. 
Moving around a triangle of ${\cal T}$ we get a permutation $\{1,2\} \lra \{2,1\}$. 
\end{itemize}

The spectral surface $\Sigma$ is obtained by deleting the vertices of the triangles from ${\bf \Sigma}$. 

\paragraph{3. The  twisted  decorated flat line bundle  on ${\Sigma}$.} 
Let ${\cal V}$ be a twisted decorated local system of two dimensional $R$-vector spaces over ${\bS}$. 
So there is an invariant decorated flag ${\cal A}_s$ near each marked point $s$.  
For every ideal triangle $t$ 
there are decorated flags ${\cal A}$, ${\cal B}$, 
 ${\cal C}$ assigned to the vertices of $t$. We assign  
to the vertices of the hexagon $h_t$ the following six trivialised lines
 \be \la{12.10.11.1}
\begin{split}
&A:=  {\cal A}^1, \ \ \ \ \ \ \ B:=  {\cal B}^1, \ \ \ \ \ \ \ C:=  {\cal C}^1,\\
&\overline A:=  {\cal A}/{\cal A}^1, ~~~ \overline B:=  {\cal B}/{\cal B}^1, ~~~ \overline C:=  {\cal C}/{\cal C}^1, \\
\end{split}
\ee 
 so that 
going around the hexagon we get the lines $(A, \overline B, C, \overline A, B, \overline C)$. 
We denote by $a_1$ the vector trivializing the line $\A$, by $a_2$ the vector trivializing the line $\overline \A$ etc. 

\begin{figure}[t]
\centerline{\epsfbox{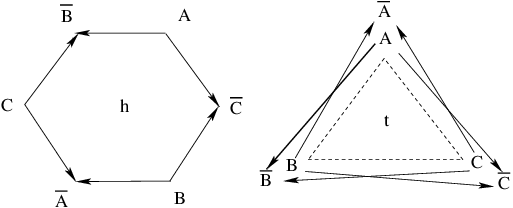}}
\caption{The twisted flat line bundle   on a hexagon assigned to a triple of flags.}
\label{bgc7}
\end{figure}

The map $h_t \to t$ sends the vertices labeled by $A, \overline A$ to the 
vertex $a$, and so on, see Figure \ref{bgc7}.  
We get a twisted flat line bundle ${\cal L}_{h_t}$ over the hexagon, trivialised at the 
 tangent vectors to the   boundary near   the vertices by the vectors provided by the decorations. For example, the line bundle near the vertex $\A$ of the hexagon 
 is trivialised by the vector $a_1\in \A$. The  parallel transports along the  oriented edges are defined by 
the following maps/their inverses:
 \be \la{12.9.11.100*}
\varphi_{a_1, b_2}: A  \lra \overline B.
\ee
Given   an ideal triangulation  of $\bS$, and  
gluing the hexagons $h_t$ into the surface ${\Sigma}$,  
the twisted flat line bundles ${\cal L}_{h_t}$ are glued into 
a twisted  flat line bundle ${\cal L}$  
on ${\Sigma}$, trivialized near the marked points. 

\vskip 2mm

 We   construct    ${\cal A}-$coordinates  
by assigning to each   edge $E$ of the hexagonal tessellation of ${\Sigma}$ an  $\Delta_E\in R^\times$ as follows. 
Let $A$ and $\overline B$ be the lines over the vertices of $E$,   
   trivialized by non-zero vectors $a_1 \in A$ and $b_2 \in \overline B$, provided by the decorations.  
The parallel transport $\varphi_E: A \to \overline B$ along the $E$  
acts by 
 $
\varphi_E(a_1) = \Delta_E b_2. 
 $
For every hexagon $h_t$,
 the alternated product of the six functions $\Delta_E$ assigned to its  consecutive edges is equal to $-1$. 
There are no more relations between the functions $\{\Delta_E\}$.

Let us  elaborate this, to  see what happens when we alter a ${\rm GL}_2-$graph $\Gamma$ by a two by two move.

\subsection{Coordinates for   twisted decorated non-commutative 2-dimensional local systems on $\bS$} \la{sec4.2}

\paragraph{1. Twisted  triples of decorated flags.} We define them as   collections of decorated flags 
$$
\{{\cal F}_i\}, \ i \in \Z/6\Z, \ \ \mbox{such that} \ \  {\cal F}_{i+3} = -{\cal F}_i \ \ \forall i \in \Z/6\Z. 
$$
 
\bl
A twisted triple  of decorated flags $\{{\cal F}_i\}$ is the same thing as a twisted decorated local system on the triangle $t$. 
\el

\begin{proof} We consider   $t$ as a disc with three marked points $x_1, x_2, x_3$. Pick a non-zero tangent vector $v_i$ at the point $x_i$, which follow  the boundary orientation. 
There is the unique twisted 
 local system  ${\cal L}$ on the disc. We identify the fibers ${\cal L}_{v_i}$ of ${\cal L}$ at the vectors $v_i$ by   the parallel transport from $v_1$ to $v_2$ and   to $v_3$ 
along the boundary circle, following the circle orientation, depicted as counterclockwise on the pictures. 
 
 Set  $V:={\cal L}_{v_1}$. 
We assign to a triple of decorated flags $\{{\cal F}_i\}$ in $V$ a pair $({\cal L}, \alpha)$, where the decoration $\alpha$ is given by the flag 
 ${\cal F}_i$ at the fibers  ${\cal L}_{v_i}$ for $i=1,2,3$, identified with $V$ as above. \end{proof}

A triple of decorated flags $({\cal A}, {\cal B}, {\cal C})$ gives rise to a twisted triple, setting
\be \la{TWTR}
({\cal A}, {\cal B}, {\cal C}) \lms ({\cal A}, {\cal B}, {\cal C}, {\cal A}', {\cal B}', {\cal C}' ), \ \ {\cal A}'=-{\cal A}, \ {\cal B}'=-{\cal B}, \ 
{\cal C}'=-{\cal C}.
\ee  

\paragraph{2. Describing  twisted  pairs of decorated flags in a two dimensional space.} Let ${\cal A}$, ${\cal B}$ 
  be two  decorated flags in a two dimensional $R$-vector space $V_2$. 
Let 
\be \la{12.10.11.1s}
\begin{split}
&A:=  {\cal A}^1, \ \ \ \ \ \ \ B:=  {\cal B}^1,\\
&\overline A:=  {\cal A}/{\cal A}^1, ~~~ \overline B:=  {\cal B}/{\cal B}^1\\
\end{split}
\ee
 be the pairs of associate graded lines of these flags. 
Each of the lines is equipped with a non-zero vector. We can describe decorated flags by these vectors, using the notation
\be \la{ABCvv}
{\cal A}= (a_1, a_2), ~~ {\cal B}= (b_1, b_2).
\ee 

Projecting the line $A$ 
to the quotient line $\overline B$ we get an isomorphism of lines 
\be \la{12.9.11.100}
\varphi_{a_1, b_2}: A  \lra \overline B.
\ee
Since the lines are equipped with non-zero vectors, it 
is described by an element $\Delta(a_1,  b_2) \in R^*$: 
\be \la{12.9.11.300}
\varphi_{a_1, b_2}: a_1\lms \Delta(a_1, b_2) b_2. 
\ee
The inverse map $\varphi_{b_2, a_1}:  \overline B \lra \rm A$,   $b_2 \lms \Delta(b_2, a_1) a_1$  is described by   the invariant
\be \la{OPNxx}
\Delta(b_2, a_1):= \Delta(a_1, b_2) ^{-1}. 
\ee 
Generic pairs of flags $({\cal A},  {\cal B})$ are described up to an isomorphism by the invariants $\{\Delta(a_1, b_2), \Delta(b_1, a_2) \}$.  

 \begin{figure}[ht]
\centerline{\epsfbox{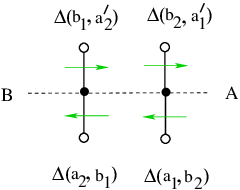}}
\caption{Invariants of a twisted pair of decorated flags   $({\cal A}, {\cal B}, -{\cal A}, -{\cal B})$ at the edges of a (disconnected)  bipartite ribbon graph. }  
\label{AB5}
\end{figure} 

Let us consider now      \underline{twisted}   pairs of decorated flags  $({\cal A}, {\cal B}, -{\cal A}, -{\cal B})$. Then the natural invariants are 
$$
\{\Delta(a_1, b_2), \Delta(b_1, a'_2)\}, \ \ \ \ a_2':=-a_2. 
$$

They are invariant under the twisted cyclic shift  $({\cal A}, {\cal B}, -{\cal A}, -{\cal B}) \to ({\cal B}, -{\cal A}, -{\cal B}, - {\cal A})$. 

 Note that on  the bipartite graph on Figure \ref{AB5},  consistently with (\ref{CONVa}),   there are two more invariants: 
\be \la{OPN}
\Delta(a_2, b_1) ^{-1} =-\Delta(b_1, a'_2), \ \ \ \ \Delta(b_2, a'_1):= - \Delta(a_1, b_2) ^{-1}. 
\ee 

\paragraph{3. Describing  twisted  triples of decorated flags in a two dimensional space.}
Let ${\cal A}$, ${\cal B}$, 
 ${\cal C}$ be three  decorated flags in a two dimensional $R$-vector space $V_2$. Similar to (\ref{12.10.11.1s}), we set 
\be \la{12.10.11.1}
\begin{split}
&A:=  {\cal A}^1, \ \ \ \ \ \ \ B:=  {\cal B}^1, \ \ \ \ \ \ \ C:=  {\cal C}^1,\\
&\overline A:=  {\cal A}/{\cal A}^1, ~~~ \overline B:=  {\cal B}/{\cal B}^1, ~~~ \overline C:=  {\cal C}/{\cal C}^1.\\
\end{split}
\ee
We   describe decorated flags by the  vectors in these lines, using the notation
\be \la{ABCv}
{\cal A}= (a_1, a_2), ~~ {\cal B}= (b_1, b_2), ~~{\cal C}= (c_1, c_2).
\ee
 The three flags provide six such line isomorphisms, which we 
organize into a diagram: 
\be \la{12.9.11.200}
\begin{gathered}
    \xymatrix{
    &   \overline B  &  & \overline A   &&    \overline C    & \\
A    \ar[ru]     &&      \ar[lu]       C   \ar[ru]   &&   \ar[lu]  B    \ar[ru] && \ar[lu]    A \\
        }
\end{gathered}
 \ee

 We compose the six maps  to a map $A  \to A$, where going against 
an arrow means inverting the isomorphism. 
By Lemma \ref{HEXTR}  the composition is equal to $-{\rm Id}$.  
This just means that
\be \la{12.9.11.400}
\Delta(a_1, b_2) \Delta(c_1, b_2)^{-1} \Delta(c_1, a_2) \Delta(b_1, a_2)^{-1} \Delta(b_1, c_2) \Delta(a_1, c_2)^{-1} =-1. 
\ee
 Using the   twisted triples of decorated flags  (\ref{TWTR}),   setting $a_1':=-a_1$, etc., and using (\ref{OPNxx}), we   write this as 
\be \la{12.9.11.400ax}
 \Delta(a_1, b_2) \Delta(b_2, c_1)  \Delta(c_1, a'_2) \Delta(a'_2, b'_1)  \Delta(b'_1, c'_2) \Delta(c'_2, a'_1)=1.   
\ee 
In fact this is just equivalent to 
\be \la{12.9.11.400axY}
 \Delta(a_1, b_2) \Delta(b_2, c_1)  \Delta(c_1, a_2) \Delta(a_2, b_1)  \Delta(b_1, c_2) \Delta(c_2, a'_1)=1. 
 \ee


 \begin{figure}[ht]
\centerline{\epsfbox{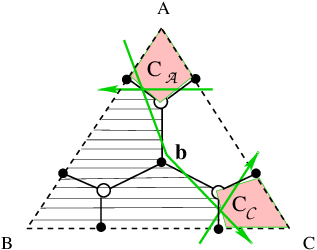}}
\caption{Sectors $C_{\cal A}$ and $C_{\cal C}$ (pink) and cosector $\widehat C_{\cal B}$ (shaded).}  
\label{bgc11}
\end{figure} 

 Let $\Gamma$ be the ${\rm GL_2}$-graph associated to the triangle $t$, see Figure \ref{bgc11}. 
The lines $A, B, C$ are assigned to the sectors $C_{\cal A}$, $C_{\cal B}$, $C_{\cal C}$. 
The lines $\overline A, \overline B, \overline C$ are assigned to the cosectors 
$\widehat C_{\cal A}$, $\widehat C_{\cal B}$, $\widehat C_{\cal C}$. The spectral surface of $\Gamma$ 
is a hexagon $h_t$ providing a $2:1$ cover of the triangle $t$, ramified at  
the central $\bullet$-vertex. \\

\begin{figure}[t]
\centerline{\epsfbox{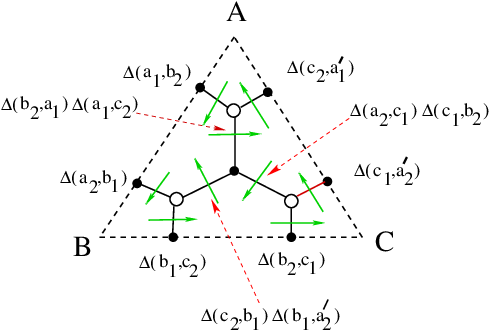}}
\caption{ The ${\cal A}$-coordinates on the bipartite ${\rm GL_2}-$graph assigned to a triangle. The  coordinates at the side $\rm C \rm A$ are assigned to the 
pair   $({\cal C}, {\cal A}')$ where ${\cal A}' := -{\cal A}$, and $(a'_1, a'_2) = (-a_1, -a_2)$.    Identity (\ref{12.9.11.400ax})  means that the counterclockwise product of the  coordinates at the edges at the central $\bullet-$vertex is   $1$.}  
\label{bgc11b0}
\end{figure}

\begin{figure}[t]
\centerline{\epsfbox{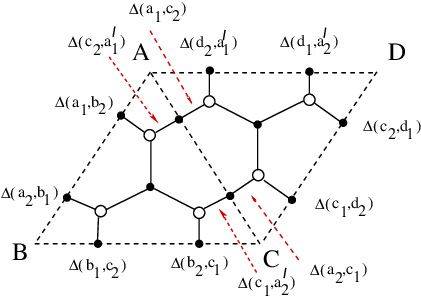}}
\caption{The ${\cal A}-$coordinates on the bipartite ${\rm GL_2}-$graph associated with a triangulated    rectangle. }  
\label{bgc11b}
\end{figure}

Using  this and notation (\ref{OPNxx}), we  assign the coordinates  to the edges of the  bipartite ${\rm GL_2}-$graph  related to  the counterclockwise oriented 
triangle  on Figure \ref{bgc11b0}. 
Then  we   assign   coordinates to the edges of a  bipartite ${\rm GL_2}-$graph associated with an ideal triangulation of the surface, as shown   on 
Figure  \ref{bgc11b}.

The cyclic product  of   elements at the edges  incident to a  $\circ-$vertex is    $-1$, and   for the   $\bullet-$vertices $(-1)^{{\rm va}l(v)-1}$.   
 For a $\bullet-$vertex this is the six-term identity (\ref{12.9.11.400ax}), or, equivalently, (\ref{12.9.11.400axY}).
 
\begin{figure}[t]
\centerline{\epsfbox{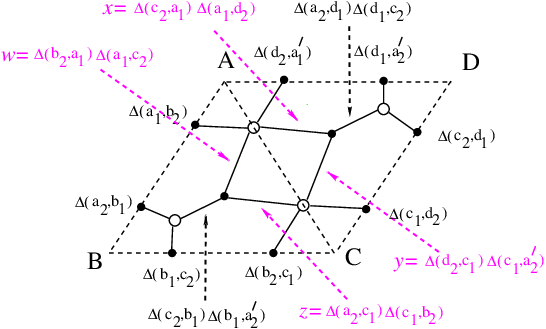}}
\caption{The ${\cal A}-$coordinates for the triangulation $AC$ of  counterclockwise oriented rectangle $ABCD$.}  
\label{bgc11c}
\end{figure}

The product of the two coordinates at the edges sharing a two-valent $\bullet-$vertex is $-1$. 
We  shrink the edge containing a two-valent $\bullet-$vertex into a   $\circ-$vertex. 
Then the cyclic product of the resulted coordinates at the edges of the new $\circ-$vertex is   $-1$, as needed. 
 The resulted coordinates for the two   triangulations of a rectangle are shown on 
Figures \ref{bgc11c}  - \ref{bgc11d}.

\paragraph{4. In coordinates.} Choose a basis $(e_1, e_2)$ in the two dimensional $R$-vector space 
$V_2$. Then the three pairs of vectors $(a_1, a_2; b_1, b_2; c_1, c_2)$ are described by a $2 \times 6$ matrix with the columns given by decompositions 
\be
\begin{split}
&a_1=a_{11}e_1+a_{12}e_2, \ \  a_2=a_{21}e_1+a_{22}e_2;\\
&b_1=b_{11}e_1+b_{12}e_2, \ \ \ b_2=b_{21}e_1+b_{22}e_2;\\
&c_1=c_{11}e_1+c_{12}e_2, \ \ \ c_2=c_{21}e_1+c_{22}e_2.\\
\end{split} 
\ee
The condition that $a_1 - \Delta(a_1, b_2) b_2 = \lambda b_1$ for some $\lambda \in R^*$ just means that 
 $$
(a_{11}e_1+a_{12}e_2) - \Delta(a_1, b_2) (b_{21}e_1+b_{22}e_2) = \lambda (b_{11}e_1+b_{12}e_2). 
$$
So we get a system of equations
\begin{align} \la{GRF1}
\left\{ \begin{array}{lll} 
a_{11} - \Delta(a_1, b_2) b_{21} = \lambda b_{11},\\ 
a_{12} - \Delta(a_1, b_2) b_{22} = \lambda b_{12}. 
\end{array}\right.
\end{align}
To solve the system, multiply the first equation   from the right by $b_{11}^{-1}$, the second by $b_{12}^{-1}$, and subtract.   
We get the unique   solution, which can be written in  two slightly different ways: 
\begin{align} \la{GRF}
\Delta(a_1, b_2) = (a_{11}-a_{12}b_{12}^{-1}b_{11}) (b_{21}-b_{22}b_{12}^{-1}b_{11})^{-1} =   (a_1, b_1)_{11}(b_1, b_2)_{21}^{-1};\\
\Delta(a_1, b_2) =  
(a_{11}b^{-1}_{11}b_{12} -a_{12}) (b_{21}b_{11}^{-1}b_{12}-b_{22})^{-1}= (a_1, b_1)_{12}(b_1, b_2)_{22}^{-1}. \nonumber
\end{align}

Here $(x_1, x_2)_{ij}$ is the quasideterminant \cite{GR} 
assigned to the $(ij)$ entry of a $2\times 2$ matrix $(x_{ij})$:
\begin{align} \la{GRREF}
(x_1, x_2)_{11}:= x_{11} - x_{12}x_{22}^{-1}x_{21}, ~~
(x_1, x_2)_{12}:= x_{12} - x_{11}x_{21}^{-1}x_{22}.\\
(x_1, x_2)_{21}:= x_{21} - x_{22}x_{12}^{-1}x_{11}, ~~
(x_1, x_2)_{22}:= x_{22} - x_{21}x_{11}^{-1}x_{12}. \nonumber
\end{align}

\paragraph{5. Remarks.}  1. Equations (\ref{GRF1}) reflect the parallel transform in the \underline{twisted} flat line bundle, see  Remark in the end of Section \ref{SEC2}. 

2. If $R$ is   commutative,    formulas (\ref{GRF}) look as follows. 
Pick a symplectic form $\omega\in {\rm det}(V_2^*)$. Then for any two vectors $v,w \in V_2$ we 
have an invariant $\omega(v,w)$, and formula (\ref{GRF}) can be written as: 
\be
\begin{split}
 \Delta(a_1, b_2) = & -\omega(a_1, b_1) \omega(b_1, b_2)^{-1}\\
=&\omega(b_1, a_1) \omega(b_1, b_2)^{-1}.\\
\end{split}
 \ee

\paragraph{6. Flips of triangulations and non-commutative Pl\"ucker relations.}   Given  a basis $(e_1, e_2)$ in the   $R$-vector space 
$V_2$,  the four pairs of vectors $(a_1, a_2; b_1, b_2; c_1, c_2; d_1, d_2)$ are described by a $2 \times 8$ matrix:
\be
\begin{split}
&a_1=a_{11}e_1+a_{12}e_2, \ \  a_2=a_{21}e_1+a_{22}e_2;\\
&b_1=b_{11}e_1+b_{12}e_2, \ \ \ b_2=b_{21}e_1+b_{22}e_2;\\
&c_1=c_{11}e_1+c_{12}e_2, \ \ \ c_2=c_{21}e_1+c_{22}e_2;\\
&d_1=d_{11}e_1+d_{12}e_2, \ \  d_2=d_{21}e_1+d_{22}e_2.\\
\end{split} 
\ee


Recall the  Pl\"ucker identity for   quasideterminants: 
  \be \la{PI}
 \begin{split}
  (a_1, c_1)_{11}  = &(a_1, b_1)_{11} (d_1, b_1)^{-1}_{11}  (d_1, c_1)_{11}  + 
  (a_1, d_1)_{11} (b_1, d_1)^{-1}_{11}  (b_1, c_1)_{11}. \\
 \end{split}
 \ee

 \bp There is a non-commutative Pl\"ucker relation for the elements $\Delta(*,*)$ from (\ref{GRF}):
\be \la{QPI}
 \begin{split}
 \Delta(a_1, c_2) = &\Delta(a_1, b_2) \Delta(d_1, b_2)^{-1}\Delta(d_1, c_2) + \Delta(a_1, d_2) \Delta(b_1, d_2)^{-1}\Delta(b_1, c_2).\\
 \end{split}
  \ee
 It is equivalent to the Pl\"ucker identity for  quasideterminants. 
 \ep

 \begin{proof} We use  formula  $\Delta(a_1, b_2) = (a_1, b_1)_{11}(b_1, b_2)_{21}^{-1}$ from (\ref{GRF}).    Then  we have to check that 
 \be
 \begin{split}
  (a_1, c_1)_{11} (c_1, c_2)_{21}^{-1}\stackrel{?}{=} &(a_1, b_1)_{11}    (b_1, b_2)^{-1}_{21}  \cdot  (b_1, b_2)_{21}   (d_1, b_1)^{-1}_{11}  \cdot (d_1, c_1)_{11}(c_1, c_2)_{21}^{-1} +\\
 &  (a_1, d_1)_{11}  (d_1, d_2)^{-1}_{21}  \cdot  (d_1, d_2)_{21}  (b_1, d_1)^{-1}_{11} \cdot (b_1, c_1)_{11}(c_1, c_2)_{21}^{-1}  \\
=&(a_1, b_1)_{11}      (d_1, b_1)^{-1}_{11}    (d_1, c_1)_{11}(c_1, c_2)_{21}^{-1} +\\
 &  (a_1, d_1)_{11}    (b_1, d_1)^{-1}_{11}  (b_1, c_1)_{11}(c_1, c_2)_{21}^{-1}. \\  
 \end{split}
 \ee 
 Multiplying this from the right by $(c_1, c_2)_{21}$ we get precisely   (\ref{PI}).  \end{proof} 

\paragraph{\it Preliminary comments.}  Recall  from Lemma \ref{CRTH} that ${\cal A}-$decorated ribbon graphs $\Gamma$ are in bijection with twisted one dimensional local systems on the corresponding spectral 
 surface $\Sigma = \bS_\Gamma$, trivialized at the boundary. By Theorem \ref{TH2.2},  proved in Section \ref{SECT1.1a},  the one dimensional 
 local systems on the spectral surface inject into twisted decorated local systems on the base surface $\bS= \bS_\Gamma$. In the case  considered in Theorem \ref{Th4.5}, the latter are just 
 twisted configurations of quadruples of  flags in $V_2$. In Theorem \ref{Th4.5} below we consider the case when 
 the ribbon graphs $\Gamma$ and $\Gamma'$ are dual to the two triangulations of the rectangle $ABCD$. If a certain quadruple of flags   is in the image of both injections, for the graphs $\Gamma$ and $\Gamma'$, then Theorem \ref{Th4.5} tells how  the   corresponding 
 ${\cal A}-$decorations of the graphs  $\Gamma$  and $\Gamma'$ are related. Note that a quadruple of flags is described by $12=8+4$ ${\cal A}-$coordinates. We stress that they are not coordinates, since there are monomial relations between them.\footnote{The coordinate at the internal 
 non-rhombus edge in  each  triangle  is recovered by the  monomial relation at its $\circ-$vertex.} The eight  ${\cal A}-$coordinates  
 at the external  edges, see Figure \ref{bgc11c}, remain intact. So we have  to explain how the ${\cal A}-$coordinates at the four rhombus edges  change under a flip. 
 The answer is given by formulas (\ref{FRM}).\footnote{Decorated flags in $V_2$ depend on $3$ parameters. So the space of isomorphism classes of quadruples of decorated flags in $V_2$ depends on  $8= 3 \cdot 4 -4$  parameters. 
 There are $4$ monomial relations between the  twelve "${\cal A}-$coordinates", one for each rhombus vertex.}
 
  \bt \la{Th4.5}
The  flip $AC \to BD$ of the triangulated  rectangle $ABCD$   amounts to the following transformation of the   ${\cal A}-$coordinates, as explained in the comment above,  depicted on Figures \ref{bgc11c}-\ref{bgc11d}: 
\be \la{FRM}
 \begin{split}
 &  \overline x  =   (1 + (zwxy)^{-1} ) z = 
 \ z (1 + (wxyz)^{-1} );  \\
  &   \overline y =  (1+   wxyz)^{-1}w \ \ = \  
 w(1+xyzw)^{-1};\\
 &  \overline z  =   (1 + (xyzw)^{-1} )x \ =  
~ x  (1 + (yzwx)^{-1} );\\
&  \overline w =  (1+yzwx)^{-1}y \ \ = \   
 ~ y(1+zwxy)^{-1}.\\
  \end{split}
 \ee
 Each of these  formulas  is equivalent to the non-commutative Pl\"ucker identity (\ref{QPI}).  
\et

   \begin{proof} The  clockwise products of the   rhombus edge coordinates on Figure \ref{bgc11c}, starting from a  $\bullet-$vertex, are given by: 
 \be \la{PIOP}
 \begin{split}
  & yzwx= - \Delta(d_2, c_1)\Delta(c_1, b_2)\Delta(b_2, a_1)\Delta(a_1, d_2),\\
  &wxyz=-\Delta(b_2, a_1)\Delta(a_1, d_2)\Delta(d_2, c_1)\Delta(c_1, b_2).\\
  \end{split}
  \ee
   \begin{figure}[ht]
 \centerline{\epsfbox{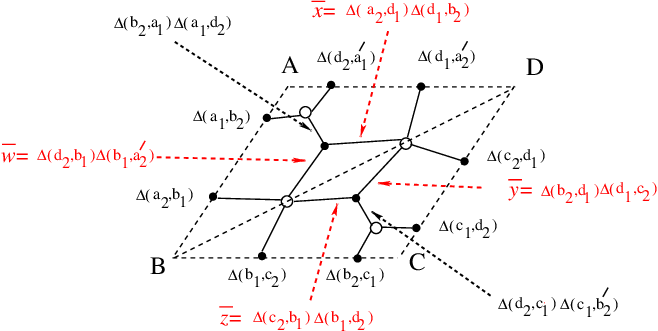}}
 \caption{  The ${\cal A}-$coordinates for the triangulation $BD$ of   counterclockwise oriented rectangle $ABCD$. 
We start with a configuration of $4$ decorated flags $({\cal A}, {\cal B}, {\cal C}, {\cal D})$. It  induces   triples  $({\cal A}, {\cal B}, {\cal C})$ and $({\cal B}, {\cal C}, {\cal D})$, which we promote   to  twisted configurations, and introduce the coordinates as on Figure \ref{bgc11b0}.}  
 \label{bgc11d}
 \end{figure} 
 
 \noindent
Clockwise products of   rhombus edge coordinates on Figure \ref{bgc11d}, starting from a  $\bullet-$vertex, are given by: 
 \be
 \begin{split}
 \overline x \ \overline y \ \overline z \ \overline w =
 -\Delta(a_2, d_1) \Delta(d_1, c_2)  \Delta(c_2, b_1) \Delta(b_1, a_2),\\  
   \overline z\ \overline w \ \overline x  \ \overline y    
 =- \Delta(c_2, b_1) \Delta(b_1, a_2)  \Delta(a_2, d_1) \Delta(d_1, c_2).\\
    \end{split}
  \ee
Now the proof  is deduced easily from the following Lemma. 

\bl
We have:
\be
 \begin{split}
 &y \overline w^{-1}=  1+ \Delta(d_2, c_1)\Delta(c_1, b_2)  \Delta(b_2, a_1)\Delta(a_1,d_2) =  1+yzwx;\\
 &w \overline y^{-1}=  1+ \Delta(b_2, a_1)\Delta(a_1, d_2)  \Delta(d_2, c_1)\Delta(c_1, b_2) =  1+wxyz;\\
 &z^{-1}\overline x=   1+ \Delta(b_2, c_1)  \Delta(c_1, d_2)  \Delta(d_2, a_1)\Delta(a_1, b_2) =  1+(wxyz)^{-1};\\
 &x^{-1}\overline z= 1+ \Delta(d_2, a_1)  \Delta(a_1, b_2)  \Delta(b_2, c_1)\Delta(c_1, d_2)= 1+( y z wx)^{-1}.\\
  \end{split}
 \ee 
  \el

\begin{proof} It is easy to see that we have the following:
 \be \la{44}
  \begin{split}
&w \overline y^{-1}=\Delta(b_2, a_1)\Delta(a_1, c_2)\Delta(c_2, d_1)\Delta(d_1, b_2),\\
&z^{-1}\overline x=\Delta(b_2, c_1)\Delta(c_1, a_2)\Delta(a_2, d_1)\Delta(d_1, b_2),\\
  \end{split}\ee
  and similar two identities obtained by the cyclic shift by two: $(a,b,c,d) \lms (c,d,a,b)$. 

Substituting   the non-commutative Pl\"ucker identity (\ref{QPI}) to   formula (\ref{44}) for $w\overline y^{-1}$   we get
\be
\begin{split}
 w \overline y^{-1}&=\Delta(b_2, a_1)\Bigl(\Delta(a_1, b_2) \Delta(b_2, d_1)\Delta(d_1, c_2) + \Delta(a_1, d_2) \Delta(d_2, b_1)\Delta(b_1, c_2)\Bigr)\Delta(c_2, d_1)\Delta(d_1, b_2)\\
&  = 1+ \Delta(b_2, a_1)\Delta(a_1, d_2) \Delta(d_2, b_1)\cdot \Delta(b_1, c_2) \Delta(c_2, d_1)\Delta(d_1, b_2)\\
&= 1+ \Delta(b_2, a_1)\Delta(a_1, d_2) \Delta(d_2, b_1)\cdot \Delta(b_1, d_2) \Delta(d_2, c_1)\Delta(c_1, b_2)\\
& = 1- \Delta(b_2, a_1)\Delta(a_1, d_2)  \Delta(d_2, c_1)\Delta(c_1, b_2) \\ 
&= 1 + wxyz.\\
\end{split}
\ee
Here to get the third  equality we use the following identity, equivalent to (\ref{12.9.11.400}), combined with (\ref{OPNxx}):
\be \la{88}
\Delta(b_1, c_2) \Delta(c_2, d_1)\Delta(d_1, b_2) =  -\Delta(b_1, d_2) \Delta(d_2, c_1)\Delta(c_1, b_2).
\ee
The identity for $y\overline w^{-1}$ is obtained from this by the cyclic shift by two. \\

 Similarly, to get  formula formula (\ref{44})  for $z^{-1}\overline x$ we need the non-commutative Pl\"ucker identity 
\be \la{QPII}
\Delta(c_1, a_2) = \Delta(c_1, d_2) \Delta(d_2, b_1)\Delta(b_1, a_2) + \Delta(c_1, b_2) \Delta(b_2, d_1)\Delta(d_1, a_2).
 \ee 
Substituting    identity (\ref{QPII}) to  formula (\ref{44}) for $z^{-1}\overline x$,   we get
\be
\begin{split}
z^{-1}\overline x&= \Delta(b_2, c_1)\Delta(c_1, a_2)\Delta(a_2, d_1)\Delta(d_1, b_2)  \\
& \stackrel{(\ref{QPII})}{=}\Delta(b_2, c_1)\Bigl( \Delta(c_1, d_2) \Delta(d_2, b_1)\Delta(b_1, a_2) + \Delta(c_1, b_2) \Delta(b_2, d_1)\Delta(d_1, a_2)\Bigr)\Delta(a_2, d_1)\Delta(d_1, b_2)  \\
& =\Delta(b_2, c_1)  \Delta(c_1, d_2) \Delta(d_2, b_1)\cdot \Delta(b_1, a_2) \Delta(a_2, d_1)\Delta(d_1, b_2)+ 1  \\
&= 1-\Delta(b_2, c_1)  \Delta(c_1, d_2) \Delta(d_2, b_1)\cdot \Delta(b_1, d_2) \Delta(d_2, a_1)\Delta(a_1, b_2) \\
&=1-\Delta(b_2, c_1)  \Delta(c_1, d_2)  \Delta(d_2, a_1)\Delta(a_1, b_2)\\
&\stackrel{(\ref{PIOP})}{=} 1+(wxyz)^{-1}.\\
  \end{split}
\ee
Here to get the fourth equality we use the identity similar to (\ref{88}). \end{proof} \end{proof}

\paragraph{7. Remark.}  Formulas (\ref{FRM})   imply that the   monodromies starting at the $\bullet-$vertex are related as follows: 
 \be \la{xyzw}
  \overline x \  \overline y \  \overline z  \  \overline w = (zwxy)^{-1}.
   \ee
 Indeed, the left hand side is equal to \be
  \begin{split}
   &z  (1 + (wxyz)^{-1} )(1+   wxyz)^{-1}wx  (1 + (yzwx)^{-1} )(1+yzwx)^{-1}y \\
  =&z  \cdot  (wxyz)^{-1} \cdot wx \cdot   (yzwx)^{-1}  \cdot y  =  (zwxy)^{-1}.\\ 
   \end{split}
    \ee
It is also clear from Figures \ref{bgc11c} - \ref{bgc11d} that each   side  of (\ref{xyzw}) is equal to   
 $
 \Delta(a_2, d_1) \Delta(d_1, c_2)  \Delta(c_2, b_1) \Delta(b_1, a_2).
 $
\\

We say that a twisted decorated two-dimensional $R-$local system on  $\bS$ is {\it generic}, 
if  for any   triangle $t$ of any ideal triangulation of $\bS$ the decorations  at the  vertices of $t$ are in generic position.

\bc 
The isomorphism classes of generic twisted decorated two-dimensional local systems on  $\bS$ can be described by a collection of $-3\chi(\bS)-1$ elements $a_i \in R^\times$ 
which satisfy a countable collection of inequalities $F_j(a_i) \not = 0$, where $F_j$ are non-commutative rational functions. 
\ec

For example,  generic quadruples of decorated flags in a two-dimensional $R-$vector space $V_2$ are parametrised by   $f_1, ..., f_8, x,y,z,w\in R^\times$ where  $xyzw \not = -1$ 
satisfying  four monomial relations. 
  
\paragraph{8. Proof of Theorem \ref{Th5.4}.} \begin{proof}

There are   two ways to assign the non-commutative cluster ${\cal A}-$coordinates  to  generic quadruples of decorated flags in a two-dimensional $R-$vector space $V_2$. 
They  correspond to the two 
 bipartite ${\rm GL}_2-$graphs assigned to the  two  triangulations of the rectangle.   
Formula (\ref{FRM}) proved in Theorem \ref{Th4.5} calculates the related coordinate transformation. Formula (\ref{FRM}) is the same  as  the one  (\ref{FRM1}).

 Denote by  ${\rm Conf}^{\rm df}_n(V_2)^\circ$ the isomorphism classes of $n$ generic decorated flags in   $V_2$. 
 Consider the cluster ${\cal A}-$coordinates on   ${\rm Conf}^{\rm df}_5(V_2)^\circ$, 
assigned to the bipartite ${\rm GL}_2-$graphs provided by the five different 
triangulations of the pentagon.  By Theorem \ref{Th4.5}, 
each flip of a triangulation of the pentagon is described by the   transformation (\ref{FRM1}) of  ${\cal A}-$coordinates for the corresponding bipartite ribbon graphs. 
By the  construction,  the composition of the five consecutive transformations  describes the identity map on  ${\rm Conf}^{\rm df}_5(V_2)^\circ$. Since 
the cluster ${\cal A}-$coordinates for a given triangulation of the pentagon provide an open embedding of ${\rm Conf}^{\rm df}_5(V_2)$ to the cluster torus, 
this implies the pentagon relation for the  ${\cal A}-$coordinates. 
\end{proof}

\section{Admissible dg-sheaves on surfaces and non-commutative cluster varieties} \la{SECCT7}

\subsection{Admissible dg-sheaves} \la{SEC7.1}

 \paragraph{1. DG-category of constructible dg-sheaves on a cell complex.} Let $X$ be a stratified topological space, such that the closure of each stratum is a ball. The strata  form a partially ordered set. We write $\sigma_1 < \sigma_2$ if $\sigma_1  \subset \overline \sigma_2$ and $\sigma_1 \not = \sigma_2$. It determines a quiver $Q$, whose  
  vertices   are the strata $\sigma$, and vertices $\sigma_1, \sigma_2$ are related by an arrow  $\sigma_1 \to \sigma_2$ if and only if $\sigma_1 <\sigma_2$.

 The quiver $Q$ determines  a nonunital $A_\infty$-category ${\cal Q}$ whose objects are the strata $\sigma$, and    

   \begin{equation} \label{Quiver}
{\rm Hom}_{\cal Q}(\sigma_1, \sigma_2) := 
\left\{ \begin{array}{lll} \Z
& \mbox{ if }    \sigma_1 <  \sigma_2, \\
0&  \mbox{ otherwise}.\end{array}\right.
\end{equation}
The composition of ${\rm Hom}$'s  is defined by the composition of arrows. 

\bd A constructible dg-sheaf on  $X$ is an  $A_\infty$-functor from the $A_\infty$-category ${\cal Q}$ to the dg-category ${\rm Vect}^\bullet$ of complexes of $R-$vector spaces. 
Constructible dg-sheaves  form a dg-category 
 $$
{\cal C}_X:= {\rm Fun}_{A_\infty}({\cal Q}, {\rm Vect}^\bullet).
$$
\ed

We consider it as a category linear over the underlying ground field - that is $\Q$ if $R$ is of characteristic zero, since  ${\rm Hom}_R$ and $\otimes_R$ do not have a natural $R-$linear structure. \vskip 2mm

Elaborating this definition, an $A_\infty$-functor ${\cal F} \in {\rm Ob}({\cal C}_X)$ is given by the following data:

\begin{itemize} \item For each nested collection of strata $\sigma_0 < \sigma_1 < \ldots < \sigma_n$, a map 
\be \la{MAPS}
f_{\sigma_0, \ldots , \sigma_n}:  {\cal F}(\sigma_0)  \otimes  {\rm Hom}_{\cal Q}(\sigma_0, \sigma_1) \otimes \ldots \otimes {\rm Hom}_{\cal Q}(\sigma_{n-1}, \sigma_n)   \lra {\cal F}(\sigma_n)[1-n].
\ee
These maps satisfy the well known quadratic equations.  
\end{itemize}

Since the spaces ${\rm Hom}_{\cal Q}(\sigma_0, \sigma_1)$ are canonically identified with $\Z$, one can write (\ref{MAPS}) just as 
\be \la{MAPSa}
f_{\sigma_0, \ldots , \sigma_n}:  {\cal F}(\sigma_0)     \lra {\cal F}(\sigma_n)[1-n].
\ee

\begin{figure}[ht]
\centerline{\epsfbox{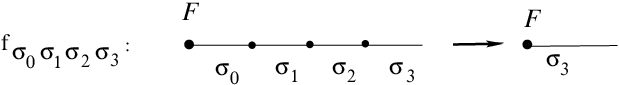}}
\caption{A   map $ {\cal F}(\sigma_0)  \otimes  {\rm Hom}(\sigma_0, \sigma_1)  \otimes  {\rm Hom}(\sigma_1, \sigma_2) \otimes {\rm Hom} (\sigma_{2}, \sigma_3)   \lra {\cal F}(\sigma_3)[-2].$ }
\label{ncls103}
\end{figure}

One  pictures a map (\ref{MAPS}) by a collection of points on the half-line, with the objects $\sigma_i$ at the intervals between the points, 
and the object ${\cal F}$ at the  end point, 
see Figure \ref{ncls103}. Then the terms of the quadratic equations match   all the ways to collapse a collection of consecutive points. 
 
 Elaborating  further,  an object of the category ${\cal C}_X$ is given by   the following data:
 
 A complex $C_\sigma$ assigned to each cell $\sigma$.

For each pair of strata $\sigma_0 < \sigma_1$,    a map of complexes 
$f_{\sigma_0, \sigma_1}: C_{\sigma_0} \lra C_{\sigma_1}$. 

For each triple of strata $\sigma_0 < \sigma_1 < \sigma_2$,  a homotopy $f_{\sigma_0, \sigma_1, \sigma_2}:  C_{\sigma_0} \to C_{\sigma_2}[-1]$:
$$
[d, f_{\sigma_0, \sigma_1, \sigma_2}] = f_{\sigma_1, \sigma_2}  \circ f_{\sigma_0, \sigma_1} - f_{\sigma_0, \sigma_2}.
 $$
 %
 And so on. 
 \vskip 3mm
 
 Given two $A_\infty$-functors ${\cal F}$ and ${\cal G}$, one defines   a graded vector space 
\be 
\begin{split}
&{\rm Hom}^\bullet({\cal F}, {\cal G}):= \\
&\prod_{n\geq 0} \prod_{\sigma_0< \ldots < \sigma_n } {\rm Hom}^\bullet\Bigl({\cal F}(\sigma_0) \otimes 
{\rm Hom}_{\cal Q}(\sigma_0, \sigma_1) \otimes \ldots \otimes {\rm Hom}_{\cal Q}(\sigma_{n-1}, \sigma_n),  {\cal G}(\sigma_n)\Bigl)[-n]=\\
&\prod_{n\geq 0} \prod_{\sigma_0< \ldots < \sigma_n } {\rm Hom}^\bullet\Bigl({\cal F}(\sigma_0),  {\cal G}(\sigma_n)\Bigl)[-n].\\
\end{split}
\ee
One equips it with  a natural differential,    getting a  complex ${\rm Hom}^\bullet({\cal F}, {\cal G})$.

 \paragraph{2. ${\cal H}$-{supported} complexes of constructible sheaves on a manifold.} Let $X$ be an oriented $n$-dimensional manifold. A {\it coorientation} of a hypersurface $H$ in $X$ is  a choice of a connected component 
 of the  conormal bundle to $H$ minus the zero section. 
  This component is called the  conormal bundle to the  {cooriented} hypersurface $H$.

  Recall that the {\it microlocal support}  
  of a complex of constructible sheaves ${\cal F}$ on $X$ \cite{KS}. It  is a closed conical Lagrangian subvariety of $T^*X$, defined as follows. 
Take a non-zero covector $\eta \in T_x^*X$, and 
  a little ball  passing through the point $x$ and tangent to the hyperplane $\eta=0$ in  $T_xX$. Move it slightly, getting a ball $B_\eta$ containing $x$; move it the other way, 
 getting a ball $B'_\eta$ which does not contain $x$. Let  $B''_\eta:= B'_\eta\cap B_\eta$.  Consider the cone of the restriction map from $B_\eta$ to $B_\eta''$: 
\be \la{MLS}
 {\rm Cone}\Bigl({\cal F}(B_\eta) \stackrel{\rm Res}{\lra} {\cal F}(B''_\eta)\Bigr).
\ee
Then $\eta$ does not enter 
 the microlocal support of ${\cal F}$ if and only if the complex (\ref{MLS}) is acyclic. Finally, if an open neighborhood $U$ of $x$ is contained in the support of ${\cal F}$, 
then  the zero section  of $T^*U$ belongs to the microlocal support of $ {\cal F}$.

\bd Let ${\cal H}$ be a collection   of    cooriented hypersurfaces in a manifold $X$. 
A complex of  constructible sheaves on $X$ is ${\cal H}$-supported    
 if its microlocal  support is contained in  the union of the zero section of $T^*X$ and the conormal bundles to the 
  cooriented hypersurfaces in  ${\cal H}$. 
  \ed
  
  The triangulated categories of complexes of constructible sheaves with microlocal support  on a given collection of hypersurfaces has been studied before, see \cite{STZ14}. 
  Our goal  is to describe a simple dg-model of its certain subcategory, that is a dg-category canonically equivalent to this subcategory, 
  which can be studied by elementary means, which also work in the non-commutative set up. 

\begin{figure}[ht]
\centerline{\epsfbox{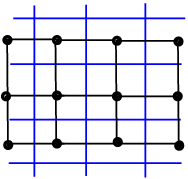}}
\caption{ A cubical complex dual to a generic stratification.}
\label{ncls105}
\end{figure}

 \paragraph{3. Generic cooriented stratification ${\cal H}$  of a manifold.} Take a collection   of cooriented hypersurfaces in generic position, that is  locally  isomorphic to 
  a collection of coordinate hyperplanes. It provides a stratification ${\cal H}$  of $X$. The 0-cells are the   intersection points of $n$ hypersurfaces.  
  The 1-cells are  components of the      intersections of $n-1$ hypersurfaces minus 0-cells, and so on. 
 Let $S$ be a   codimesion $k$ stratum 
  of the stratification ${\cal H}$. We can identify an open neighbourhood of $S$  with    $S \times \R^k$, where $\R^k$ is the   
  vector space with coordinates $\{x_i\}$, so that the  codimension $1$ strata containing $S$  are given by   equations  $x_i=0$, cooriented   by $dx_i>0$. 
   The conormal bundle $T^*_SX$ contains $2^k$  octants. 
Using an isomorphism  $T^*_SX = \R^k \times S$, they are given by inequalities  $\pm x_1 \geq 0, \ldots, \pm x_k \geq 0$.   
The octant   $x_1\geq0, \ldots , x_k\geq0$ is called  the {\it positive octant}, and denoted by $T^+_SX$. The dual  stratification  is described by an oriented cube  ${\cal K}_S$. 
   Its vertices  
 match the octants in $\R^k$.  The  oriented edges 
  match   the cooriented codinension 1 strata,  etc.
  The cubes  ${\cal K}_S$ for different strata $S$   are assembled into a cubical complex ${\cal K}$ dual to the stratification ${\cal H}$, see Figure \ref{ncls105}.   

 In particular, for each face $K$ of ${\cal K}$  there is  a vertex $v_-(K)$ assigned to the negative octant, and the 
  opposite vertex $v_+(K)$  assigned to the positive octant. They coincide if ${\rm dim}(K)=0$.

\bl The microlocal 
support of an ${\cal H}$-supported   constructible complex of sheaves ${\cal F} $ on $X$    is contained in the union of the {positive} octants   $ T_S^+X$, where $S$ are the strata of 
${\cal H}$. 
\el

\begin{proof}
The microlocal support of ${\cal F} $ is a union of certain octants in $T_S^*X$. Indeed, varying the convector $\eta$ inside of one octant we do not change the cohomology groups 
of  complex (\ref{MLS}).   
If the microlocal support intersects 
 one of the non-positive open octants, the closure of this octant intersects a negative component of the conormal bundle to one of the codimension $1$ strata, which  
contradicts to the condition that  ${\cal F} $ is ${\cal H}$-supported. \end{proof}

   \paragraph{4. ${\cal H}$-{supported} dg-sheaves on a manifold.}  A dg-enhancement of the derived category of  ${\cal H}$-{supported} complexes of constructible sheaves is given by the     
  ${\cal H}$-{supported} dg-sheaves. Let us define a more economic model for it for stratifications by generic cooriented hypersurfaces. 
    
  \bp  \la{OTHER} Let ${\cal H}$ be a stratification of a manifold $X$  by generic  cooriented hypersurfaces, such that the closures of all strata are balls. 
  Then ${\cal H}$-{supported dg-sheaves} can be  described as follows, using the dual cubical complex ${\cal K}$: 
  
\vskip 2mm
 {\bf Data:} 
  \begin{itemize}    
  \item A   graded vector space  $C^\bullet_v$ assigned to each vertex $v$ of  ${\cal K}$.
  
 \item  A map $\delta_K: C^\bullet_{v_-(K)} \lra C^\bullet_{v_+(K)}[1-k]$ for each $k$-dimensional face $K$ of  ${\cal K}$, where $k \geq 0$. 

  
  \end{itemize}

 {\bf Condition:}  
 For each face $K$ of ${\cal K}$, take the direct sum of the graded spaces $C_v^\bullet$ at the vertices $v$ of $K$, shifted   
 by the distance $d_v$ from $v$ to the  vertex  $v_-(K)$:\footnote{If one parametrizes the vertices of a cube $K$ by $(\varepsilon_1, \ldots, \varepsilon_k)$, where $\varepsilon_i \in \{0,1\}$, then the distance is  $\varepsilon_1 + \ldots + \varepsilon_k$.}
 $$
 C_K^\bullet:= \oplus_{v\in K} C^\bullet_{v}[-d_v].
 $$

 \begin{itemize}
 
\item Then the signed sum $d_K$ of the maps $\delta_F$ assigned to the faces $F$ of   $K$ in ${\cal K}$ satisfies $d_K^2=0$.  

\end{itemize}
 
  \ep
  
{\bf Example}.   Let us elaborate the cubical   description of an ${\cal H}$-supported dg-sheaf in the dimension 3:
  
    For each vertex $v$ we get a complex $(C^\bullet_v, \delta_v)$. 
 
 For each oriented edge $E = (0, 1)$ we get a morphism of complexes 
 $\delta_E: C^\bullet_{0} \lra C^\bullet_{1}$. 
 
 For each square  $K$ with the vertices $({00}, {01}, {10}, {11})$ we get a homotopy 
 $$
 h_K: C^\bullet_{{00}} \lra C^\bullet_{{11}}[-1], ~~~~  \delta_{11} h_K + h_K \delta_{00} = \delta_{01 \to 11}\circ \delta_{00 \to 01} - \delta_{10 \to 11}\circ \delta_{00 \to 10}.
 $$
 
  For each cube ${\rm C}$ with   vertices $\{{\varepsilon_1, \varepsilon_2, \varepsilon_3}\}$ we get a higher homotopy $H_{\rm C}: C^\bullet_{{000}} \lra C^\bullet_{{111}}[-2]$:
  \be
 \begin{split}
 &  \delta H_{\rm C} - H_{\rm C}\delta=  
 \delta_{110 \to 111}\circ h_{000 \to 110} + \delta_{101 \to 111}\circ h_{000 \to 101}  + \delta_{011 \to 111}\circ h_{000 \to 011}\\
 &-h_{100 \to 111}\circ \delta_{000 \to 100} - h_{010 \to 111}\circ \delta_{ 000\to 010}  - h_{001 \to 111}\circ \delta_{000 \to 001}.
 \end{split}
 \ee  \vskip 3mm

 \begin{proof}  Let us show how a data $(C^\bullet_v, \delta_K)$    determines a dg-sheaf on $X$. 
 We start with a simple  observation.  
  Let $H$ be a   component of a cooriented hypersurface strata. Denote by ${\cal D}^+$ and ${\cal D}^-$ the  connected domains of $X - {\cal H}$ 
sharing  $H$, so that   $H$ is cooriented towards  ${\cal D}^+$. 
Let ${\cal F}$ be an ${\cal H}$-supported  dg-sheaf  on $X$. Then 
 the natural map 
 \be \la{Cores1}
 {\cal F}_{|H} \lra {\cal F}_{|{\cal D}^-}
 \ee  
 is a quasi-isomorphism. Indeed, its cone   detects the microlocal support of ${\cal F}$ in the direction opposite 
 to the coorientation of $H$, and hence   is acyclic. 
 We can alter the dg-sheaf data on the component $H$ so that the map (\ref{Cores1}) is an isomorphism, getting a quasi-isomorphic   dg-sheaf.  
  
  Generalising this, let us construct an $A_\infty$-functor ${\cal F} $. Given  a stratum $\sigma$ of ${\cal H}$, 
let us define ${\cal F} (\sigma)$. Let us move a bit the point $p\in \sigma$ 
inside of the unique $n$-dimensional stratum $S_p$ in the negative direction from $p$, whose closure contains $p$.  
 The stratum $S_p$ is encoded by a vertex $v$ of ${\cal K}$. We set ${\cal F} (\sigma):= C^\bullet_v$. 
 The maps $f_{\sigma_1, ..., \sigma_n}$ are defined naturally.  \vskip 2mm

Let us describe  a dg-sheaf ${\cal C}^\bullet$ near   the intersection point $v$   of two cooriented lines in $\R^2$   on Figure \ref{ncls100a}. 
The cubical description of the dg-sheaf ${\cal C}^\bullet$ near $v$ is given by the following data, where $h:{C}^\bullet_{{00}} \lra {C}^\bullet_{{11}} [-1]$, see  also Figure \ref{ncls100a}:  
 \be \la{rs22ab} 
\begin{gathered}
    \xymatrix{ 
        {C}^\bullet_{{01}}\ar[r]^{{g_1}}    & {C}^\bullet_{{11}}  \\
         {C}^\bullet_{{00}}\ar[u]^{ {f_0}}  \ar[ru]^{h} \ar[r]_{ {g_0}}        & \ar[u]_{{f_1}}  {C}^\bullet_{{10}}}
\end{gathered}
 \ee 
The restriction of   ${\cal C}^\bullet$ to a ball  ${\cal B}$ around $v$ is given, by the   definition,    by the complex  $C_{00}^\bullet$.  
The restriction of   ${\cal C}^\bullet$ to the domain ${\cal D}$   on Figure \ref{ncls100a} is given by the following complex, where $C^0_{ij}$ is in the degree $i+j-1$: 
  \be \la{rs22abc} 
\begin{gathered}
    \xymatrix{ 
        {C}^\bullet_{{01}}\ar[r]^{{g_1}}    & {C}^\bullet_{{11}}  \\
                 & \ar[u]_{{f_1}}  {C}^\bullet_{{10}}}
\end{gathered}
\ee
Indeed, to calculate  ${\cal C}^\bullet_{|{\cal D}}$ we cover ${\cal D}$ by  open subsets ${\cal U}$ and ${\cal V}$, where  ${\cal U}$ is disjoint with the horizontal line, and ${\cal V}$ with the vertical. 
Then ${\cal C}^\bullet_{|{\cal U}}={C}^\bullet_{{01}}$ , ${\cal C}^\bullet_{|{\cal V}} = {C}^\bullet_{{10}}$. Their restrictions to ${\cal V}\cap {\cal U}$ are given by the maps $g_1, f_1$. 
\begin{figure}[t]
\centerline{\epsfbox{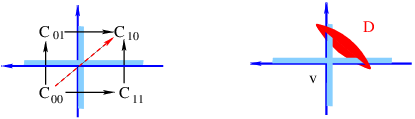}}
\caption{The   domain ${\cal D}$ intersects the two lines, but does not contain their intersection $v$.}
\label{ncls100a}
\end{figure}

Assume that ${\cal D}\subset {\cal B}$. Then, using (\ref{rs22abc}),  the restriction of  ${\cal C}^\bullet_{|{\cal B}}$     to  ${\cal D}$   is given by 
the   map of complexes: 
\be \la{MMK}
 (f_0,g_0,h): C_{00}^\bullet \lra    {\rm Cone}\Bigl(C_{01}^\bullet \oplus C_{10}^\bullet\stackrel{(g_1, f_1)}{\lra}  C_{11}^\bullet\Bigr)[-1].
\ee
The shifted cone of this map is  the total complex associated with (\ref{rs22ab}).
  \end{proof}

\bl
Given an ${\cal H}$-supported   complex of sheaves ${\cal F}$ on $X$, the  complex (\ref{MLS})  for any convector $\eta$ inside   the positive octant 
$T^+_SX$ is quasiisomorphic to the complex $(C^\bullet_K, d_K)$, where $K$ is the face of ${\cal K}$ dual to the strata $S$. 
 \el
 
 \begin{proof}   Let us elaborate the microlocal support at the vertex $v$ on Figure \ref{ncls100a}.  Then the map (\ref{MLS}) is given   by the map  (\ref{MMK}). 
 So the claim boils down  to the last sentence of  the proof of Proposition \ref{OTHER}. 
 \end{proof}
 
 \paragraph{5. ${\cal H}$-admissible dg-sheaves.} Among all  ${\cal H}$-supported  dg-sheaves  we distinguish a crucial for us subcategory of {\it ${\cal H}$-admissible dg-sheaves}.  
  
\bd Given a stratification ${\cal H}$ provided by a  collection of cooriented hypersufaces in a manifold $X$ whose Legendrians are disjoint, an ${\cal H}$-supported  dg-sheaf  as in Proposition \ref{OTHER}  
  is   ${\cal H}$-{\rm admissible} if

 \begin{itemize}
 
 \item   
 The complexes $(C_K^\bullet, d_K)$ satisfy the following conditions:
 \be \la{CONDZ}
 \begin{split}
 &i) ~~{\rm dim ~K} =0: ~~H^p(C_K^\bullet) =0 ~~\mbox{if $p\not \in \{0,1\}$},\\
 &ii) ~~{\rm dim ~K} =1: ~~H^p(C_K^\bullet) =0 ~~\mbox{if $p\not \in \{0\}$},\\
 &iii) ~~{\rm dim ~K} >1: ~~H^p(C_K^\bullet) =0~~\mbox{for all $p$}. 
 \end{split}
\ee
\end{itemize}
\ed

\bl The restriction of an ${\cal H}$-admissible dg-sheaf    to each hypersurface   of ${\cal H}$ is a local system. 
 \el
       
       \begin{proof} Follows the same idea as the proof of Lemma \ref{LE4.1}. \end{proof}

 {\bf Remark}.   Below we need subcategories of ${\cal H}-$admissible dg-sheaves with  microlocal support on the (conormal bundles to) hypersurfaces of an a priori given  collection of hypersurfaces ${\cal H}$. However we do not want to assume that  the induced stratification is topologically trivial. A typical example is a punctured surface $S$ with a collection of small loops near each puncture. So we consider in addition to  ${\cal H}$ a collection of dummy hypersurfaces, and require 
     that the microlocal support of  dg-sheaves on them is trivial, i.e. the dg-sheaves are non-singular along them. We denote by 
     $\widetilde {\cal H}$ 
     the bigger collection of hypersurfaces, and choose it so that the induced stratification is topologically trivial.  There is a  functor from 
     the category of  $\widetilde {\cal H}-$admissible dg-sheaves to  the triangulated category of  complexes of constructible sheaves with the microlocal support in ${\cal H}$, where the latter was defined in \cite{KS}, 
     which is is an equivalence. 
     Using these functors one sees that the categories of $\widetilde {\cal H}-$admissible dg-sheaves    are equivalent, and they all equivalent to the original 
     category of ${\cal H}-$admissible dg-sheaves. See also \cite[Section 3]{STZ14} where, as the referee pointed out, the related issues are discussed.

  \subsection{The bigon and triangle  moves} \la{Sec3}
  
 We consider moves of the hypersurfaces supporting the microlocal support of an admissible dg-sheaf such that 
 the corresponding Legendrians remain disjoint. We call them admissible moves. 
Focusing on the dimension 2, we prove that admissible moves preserve the category of ${\cal H}$-admissible dg-sheaves.   
 \begin{figure}[t]
\centerline{\epsfbox{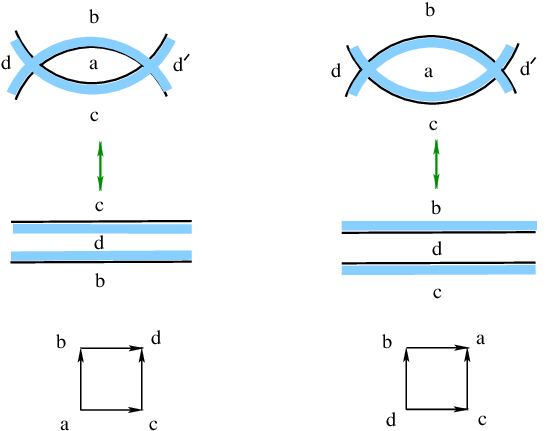}}
\caption{There are two different  admissible bigon moves.}
\label{cube3}
\end{figure}

\paragraph{The  bigon moves.} There are two   different admissible moves of a pair of cooriented lines in the plane, called  bigon moves and shown on Figure \ref{cube3}. 
 Complexes assigned to the faces and arrows between them corresponding to the cooriented 
 lines separating them are described by   the vertices and oriented edges of a square, shown at the bottom 
 of each move.

Recall that any object of the  triangulated category of representations of a Dynkin quiver is a direct sum of the shifts $R_i[n_i]$, $n_i\in \Z$ of  
indecomposable representations $R_i$. The isomorphism classes of the latter are  parametrized by  positive roots of the root system.

\begin{figure}[ht]
\centerline{\epsfbox{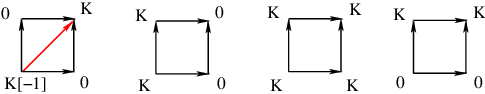}}
\caption{The squares for   bigon moves.  The non-trivial homotopy is shown by a red arrow. Each arrow $K \to K$ is an isomorphism. }
\label{cube4}
\end{figure}

The two arrows from the bottom left vertex of the square form an $A_3$-quiver. It has six indecomposable representations. Taking into account a $\Z/2\Z-$symmetry of the $A_3$-quiver, there are just four of them to consider, depicted on Figure \ref{cube4}. Each of them determines uniquely an acyclic homotopy square, as shown on Figure \ref{cube4}.\footnote{{Below we use the  fact that the classification of the indecomposable reprsentations of quivers of type $A_m$, $m\leq 3$ and $D_4$ over a skew field is exactly the same 
as over a field. Following the request  of the referee, we present the proof  in the $D_4$ case when  all arrows  to  a  vertex $v$. 
 In the paper we also need the case when they go out of $v$. Over a field one reduces to the other by the dualisation. Over a skew field dualisation does not work, but one 
 proceeds  by the analogy. So take a $D_4-$quiver  with  vertices $a, b, c, v$, and all arrows going  to $v$. 
The only fact we use is that a proper linear subspace of an $R-$vector space has a complement. So we  reduce to the case when all maps  are injective: otherwise we  split off a 
representation with $V_v=0$.   So from now on $V_a, V_b, V_c\subset V_v$. Similarly we can assume    $V_a + V_b +V_c = V_d$. If  
$W_{abc}:= V_a \cap V_b \cap V_c \not = 0$, we get a  non-decomposable  if  $V_a  =  V_b = V_c$ is of rank 1, which we can split off otherwise.     So we can 
assume  $W_{abc} = 0$. If 
 $V_c=0$, we get a representation of the quiver $\A_3$. If also $V_b=0$, we split off an indecomposable 
with the dimension vector $(1, 0, 0, 1)$. Otherwise we get an indecomposable only of $V_a=V_b$ is one dimensional. 
So we can assume  $V_a, V_b, V_c\not = 0$. Next, we can  assume that  any two of $V_a, V_b, V_c$ intersect by $0$: otherwise we split off a representation of an $\A_3$ quiver with the dimension vector $(1,1,1)$. If ${\rm dim}(V_a)>1$, we can split off a subspace from $V_a$. So we can assume $V_a, V_b, V_c$  are one-dimensional. We get the last indecomposable with the dimension vector $(1,1,1,2)$.}}

The quiver representation obtained by forgetting the  vertex $a$ of the square and the arrows from this vertex describe 
a dg-sheaf obtained by the   elementary move of the original dg-sheaf.

\bt \la{TRTB} The bigon moves transform admissible dg-sheaf 
to admissible ones. 
\et

\begin{proof} 
  We use the remark in the end of Section \ref{SEC7.1}  to extend to a bigger collection $\widetilde {\cal H}$ with topologically trivial strata of the induced stratification, to be able to restrict considerations to a small ball where we will perform the isotopy. Then the claim  is clear from   Figure \ref{cube4}.    
\end{proof}

\begin{figure}[ht]
\centerline{\epsfbox{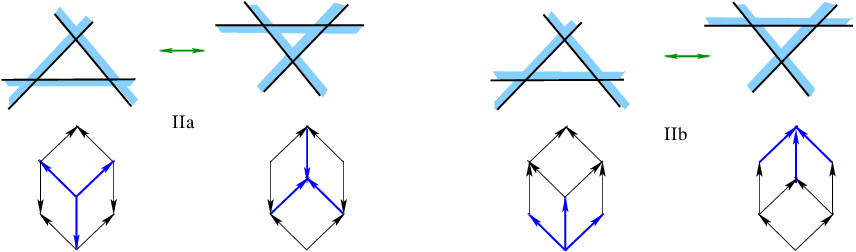}}
\caption{The triangle moves.  The $D_4$-quivers are shown by blue arrows.}
\label{cube2}
\end{figure}

\paragraph{The triangle moves.} 
There are two  moves of a triple of cooriented lines in the plane, called  IIa and IIb  moves, or just triangle moves,  shown on Figure \ref{cube2}.  
The complexes assigned to the faces and the arrows between them  are described by   the quivers  shown on the bottom.

\begin{figure}[ht]
\centerline{\epsfbox{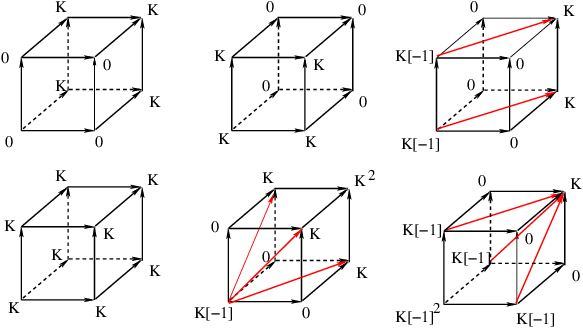}}
\caption{The cubes  for elementary moves IIa and IIb. Homotopies are shown by red arrows.}
\label{cube1}
\end{figure}

  The three arrows from the central   vertex   form an $D_4-$quiver. 
 A $D_4-$quiver has $12$ indecomposable representations. Pick an orientation of the edges of the $D_4-$quiver  such that the three arrows go out of the central vertex. 
The symmetry group of this  quiver is the group $S_3$. Modulo the $S_3-$symmety, there are $6$ indecomposable representations. They are 
 depicted on Figure \ref{cube1} for the $D_4-$quiver 
whose central vertex is the bottom left vertex of the cube. 

A  triangle move of three cooriented lines in a 2d space  is described in the 3d space-time by the three cooriented coordinate planes.  
The ${\cal H}$-supported dg-sheaves for the union of cooriented coordinate  planes are described by  data on a cube  shown on Figure \ref{cube1}, and their shifts.

Forgetting the top right vertex of the cube diagram we describe the original admissible dg-sheaf. Forgetting the bottom left vertex we describe the resulting one. 
This shows how  indecomposable dg-sheaves   transforms  under type IIa moves.  Type IIb moves are treated similarly, with the comment that this time we remove bottom right vertex instead of the  top right. 

\bt \la{TRT} Triangle moves transform admissible dg-sheafs 
to admissible ones. 
\et

\begin{proof} This is clear from   Figure \ref{cube1}.   Indeed, the  dg-sheaf in $\R^3$ assigned to the cube is   admissible,  
so forgetting any vertex we still get an admissible dg-sheaf. 
\end{proof}

 Triangle moves provide   a bijection between the isomorphism classes of indecomposable objects 
for two differently oriented $A_3$-quivers: the one oriented out of the central vertex to the one oriented towards the central vertex for the IIa move, and to the one 
with two edges oriented towards the central vertex for the   IIb move.

\subsection{Local systems on ideal bipartite  graphs   and  decorated surfaces} \la{SSec4}

\paragraph{1. Bipartite ribbon graphs and zig-zag paths.}

\begin{figure}[ht]
\centerline{\epsfbox{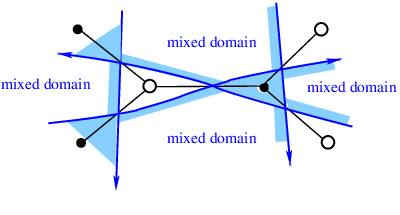}}
\caption{Zig-zag paths for a bipartite graph $\Gamma$ on a surface $\bS$ cut  the surface   into $\bullet$, $\circ$, and mixed domains. 
Each zig-zag path is coorientated   by a shaded area to the right of the path.}
\label{ncls104}
\end{figure}

A   bipartite graph $\Gamma$    gives rise to  a collection of  zig-zag paths on the decorated surface  $\bS$ associated with $\Gamma$.    
Denote  by    ${\cal Z}$   their union. The  complement $\bS - {\cal Z}$  
is a disjoint union of connected domains of three types:   $\bullet$, $\circ$, and mixed domains, see Figure \ref{ncls104}. The $\bullet-$domains are the ones containing a single $\bullet-$vertex of $\Gamma$. 
The $\circ-$domains contain a single $\circ-$vertex of $\Gamma$. The rest   are mixed domains.   The  $\bullet$ and $\circ-$domains are contractible. 

 Each zig-zag path is   orientated so that the $\bullet-$vertices are on the right.  
 This plus the surface orientation   provide  coorientation of zig-zag paths: moving along a zig-zag path   
the   coorientation points to the right.  All sides of a $\bullet-$domain are cooriented inside, all sides of a $\circ-$domain are cooriented outside, and  coorintations of the sides of a mixed domain alternate.

So the union of zig-zag paths give rise to  a cooriented stratification ${\cal Z}$ of   $\bS$.   
Its zero stratum is the set  ${\cal Z}^0$  of crossing points. So   
  there is a category of  ${\cal Z}$-admissible dg-sheaves on $\bS$.

\begin{figure}[ht]
\centerline{\epsfbox{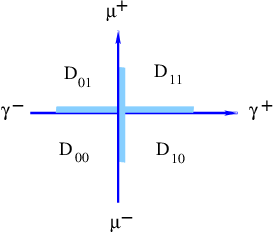}}
\caption{The stratification near a crossing point of zig-zag paths $\gamma$ and $\mu$ has one $0$-dimensional  stratum, four  $1$-dimensional strata 
$\gamma^\pm, \mu^\pm$ and four $2$-dimensional strata ${\cal D}_{\ast\ast}$.}
\label{ncls102}
\end{figure}

\paragraph{2. ${\cal Z}$-admissible dg-sheaves.} 
A ${\cal Z}$-admissible  dg-sheaf  on $\bS$ is given by the following data:

\begin{itemize}
  
  \item A   complex of local systems    ${\cal F}^\bullet_{\cal D}$ on each domain  ${\cal D}$ of  $\bS-{\cal Z}$, concentrated in   degrees $[0, 1]$:

 \item  A map  of complexes $
{\varphi}_{\gamma}:  {\cal F}^\bullet_{{\cal D}^-} \to {\cal F}^\bullet_{{\cal D}^+}
$  for each  component $\gamma$  of ${\cal Z}-{\cal Z}^0$, such that 
 \be \la{CONEGx}
H^i{\rm Cone}(\varphi_{\gamma}) =0~~\mbox{unless  $i=0$}.
\ee

 \item For each crossing point of   zig-zags $\gamma$ and $\mu$, a homotopy - see Figure \ref{ncls102}: 
 $$
 h: {\cal F}^\bullet_{{\cal D}_{00}}  \lra {\cal F}^\bullet_{{\cal D}_{11}}[-1], ~~~~dh  + h d =  \varphi_{\mu^+} \circ \varphi_{\gamma^-} -\varphi_{\gamma^+} \circ \varphi_{\mu^-}.
  $$
 \be \la{rs22a} 
\begin{gathered}
    \xymatrix{ 
        {\cal F}^\bullet_{{\cal D}_{01}}\ar[r]^{\varphi_{\mu^+}}    & {\cal F}^\bullet_{{\cal D}_{11}}  \\
         {\cal F}^\bullet_{{\cal D}_{00}}\ar[u]^{\varphi_{\gamma^-}}  \ar[ru]^{h} \ar[r]_{\varphi_{\mu^-}}        & \ar[u]_{\varphi_{\gamma^+}}  {\cal F}^\bullet_{{\cal D}_{10}}}
\end{gathered}
 \ee

\item  Let    ${\rm D}$   be   the signed sum of the differentials on ${\cal F}^\bullet_{{\cal D}_{ij}}$,  the maps  $\varphi_{\mu^\pm}$,  
 $\varphi_{\gamma^\pm}$, and $h$.  The data above satisfies the condition that the following complex with the differential ${\rm D}$ is acyclic:
 \be \la{CONDg}
  \begin{split}
 & {\cal F}^\bullet_{{\cal D}_{00}} \oplus {\cal F}^\bullet_{{\cal D}_{01}}[-1] \oplus {\cal F}^\bullet_{{\cal D}_{10}}[-1]\oplus {\cal F}^\bullet_{{\cal D}_{11}}[-2]. \\
 \end{split}
 \ee
  \end{itemize}
 
The last condition is equivalent to   either of the following    conditions:   
 
 1. The   natural map of complexes $
 {\rm Cone}( \varphi_{\gamma^-}) \lra {\rm Cone}( \varphi_{\gamma^+})$  is a quasi-isomorphism 

 2. The    natural map of complexes $
 {\rm Cone}( \varphi_{\mu^-}) \lra {\rm Cone}( \varphi_{\mu^+}) 
$  is a quasi-isomorphism.

 \paragraph{a) ${\cal Z}$-admissible  dg-sheaves on $\bS$ $\lra$  local systems on zig-zag paths.}  Condition   (\ref{CONEGx}) just means that 
  on each   component $\gamma$  of  ${\cal Z} - {\cal Z}^0$ there is  a local system 
$
 {\rm Cone}(\varphi_{\gamma}).
$
\bl \la{LE4.1}
Given a  zig-zag path $\gamma$,   local systems  ${\rm Cone}(\varphi_{\gamma})$ on  $\gamma - \{\mbox{\rm crossing points}\}$   
extend   to a local system on  $\gamma$,  denoted  by ${\rm Cone}(\varphi_{\gamma})$. 
\el

\begin{proof} Condition      (\ref{CONDg}) implies that for each crossing point on  $\gamma$ separating  two  germs of zig-zags $\gamma^-$ and $\gamma^+$,
 there is an isomorphism  
 $
 {\rm Cone}(\varphi_{\gamma^-}) \stackrel{\sim}{\lra} {\rm Cone}(\varphi_{\gamma^+}).
 $
 It  allows  to extend   local systems on  components of $\gamma -\{ \mbox{crossing points}\}$ to a local system on $\gamma$.
 \end{proof}

\paragraph{b) Local systems on $\Gamma$ $\lra$ ${\cal Z}$-admissible  dg-sheaves on $\bS$, acyclic on mixed strata.}   A local system  of vector 
spaces  on a graph $\Gamma$ is given by a collection of vector spaces 
$V_v$ at the vertices $v$ of $\Gamma$ and linear maps $\varphi_{(v, w)}: V_v\to V_w$ for each edge $ (v,w)$ of $\Gamma$, so that   $\varphi_{(w,v)}= \varphi_{(v,w)}^{-1}$.

Given    a local system ${\cal V}$ on a bipartite graph $\Gamma$ on $\bS$, let us define a ${\cal Z}$-admissible  dg-sheaf ${\cal F}_{\cal V}$ 
on $\bS$.   Since  the $\bullet$ and $\circ$   domains of $\bS - {\cal Z}$ are contractible,  
  we can talk about the fibers  of a complex of  ${\cal Z}$-constructible sheaves     over them.  Then we define
  
\begin{figure}[t]
\centerline{\epsfbox{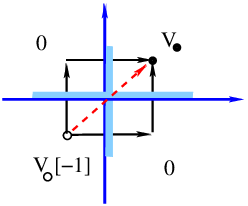}}
\caption{The admissible dg-sheaf ${\cal F}_{\cal V}$ on $\bS$ assigned to a local system ${\cal V}$ on $\Gamma$.}
\label{ncls100}
\end{figure}

  \begin{itemize}

\item The fiber of ${\cal F}_{\cal V}$ at the $\circ-$domain   with a  vertex $\circ$  
of $\Gamma$ is   the shifted fiber ${\cal V}_\circ[-1]$ of  ${\cal V}$ at   $\circ$. 

\item  The fiber of ${\cal F}_{\cal V}$ at the  $\bullet-$domain   containing a   vertex $\bullet$ 
of $\Gamma$ is   the fiber ${\cal V}_\bullet$ of  ${\cal V}$.

\item The fibers of ${\cal F}_{\cal V}$ over mixed domains are zero. 
So  the maps $\varphi_{\gamma}$ for   zig-zags $\gamma$ are zero. 

\item The   homotopy map for each edge $\circ - \bullet$ of   $\Gamma$, as illustrated on 
Figure \ref{ncls100}:
\be \la{HOM} 
h: {\cal V}_\circ[-1] \lra {\cal V}_\bullet[-1]. 
\ee
This map is  the shift of the   parallel transport map  along the edge $\circ - \bullet$. 
\end{itemize}

So we get the following diagram near each      crossing of two zig-zag paths:

\be \la{rs22} 
\begin{gathered}
    \xymatrix{
        0\ar[r]    & {\cal V}_\bullet   \\
         {\cal V}_\circ[-1] \ar[u]  \ar[ru]^{h} \ar[r]      & \ar[u]0}
\end{gathered}
 \ee
 
  \paragraph{c) ${\cal Z}$-admissible  dg-sheaves on $\bS$, acyclic on mixed strata $\stackrel{}{\lra}$ local systems on $\Gamma$. }  
 Construction b) can be inverted. Namely, let ${\cal F}$ be a  ${\cal Z}$-admissible  dg-sheaf ${\cal F}$ on $\bS$, acyclic on mixed strata. 
  Then ${\cal Z}$-admissibility implies that
   \begin{equation} \label{LOCALS}
  H^i ({\cal F}_{{\cal D}_\circ}) =  
 \begin{array}{lll} 0
& \mbox{ if }    i\not = 1 ,
\end{array}  ~~~~~~
  H^i ({\cal F}_{{\cal D}_\bullet}) =  
  \begin{array}{lll}0
& \mbox{ if }    i\not = 0.
\end{array}\end{equation}
Furthermore, the homotopy $h$ provides an isomorphism
 $$
 h:  H^1 ({\cal F}_{{\cal D}_\circ})  \stackrel{=}{\lra}   H^0 ({\cal F}_{{\cal D}_\bullet})[-1]. 
  $$
 So we get a local system ${\cal L}_{\cal F}$ on the bipartite graph $\Gamma$ by assigning the spaces $ H^1 ({\cal F}_{{\cal D}_\circ})[1]$ to the $\circ-$vertices $\circ$ of $\Gamma$, 
      the ones       $ H^0 ({\cal F}_{{\cal D}_\bullet})$  to the $\bullet-$vertices $\bullet$ of $\Gamma$, and using the isomorphism $h$ 
      as the parallel transport along the edge  $\circ\to \bullet$.  
      
      \paragraph{d) Framed local systems on $\bS$ as  ${\cal Z}_m$-admissible dg-sheaves.} Given a decorated surface $\bS$, 
  consider a collection of $m$ nested simple clockwise loops    around each puncture, and 
  a collection of $m$  simple clockwise arcs    around
   each special point, which do not intersect and do not self-intersect.  They are  cooriented    "inside the surface", see Figure \ref{ncls101}. 
  Denote by  ${\cal Z}_m$ their union.  
  
  \bl \la{FRAML} A ${\cal Z}_m$-admissible dg-sheaf on a decorated surface $\bS$, acyclic   
 near the punctures and special points, is equivalent to an $m$-dimensional local system on $\bS$ with a framing. 
  \el 
  
  \begin{proof} Take a punctured disc ${\cal D}^*$ with $m$ concentric circles going clockwise around the puncture, see Figure \ref{ncls101}.
  Take a radius of the disc and restrict  to it a ${\cal Z}_m$-admissible complex of sheaves   on ${\cal D}^*$. The restriction  is described by    complexes 
  $V_0^\bullet, \ldots V_m^\bullet$ and maps $\varphi_{i, i+1}:  V_i^\bullet \lra V_{i+1}^\bullet$ between them. To proceed further we need the following simple observation.
    
  \begin{figure}[ht]
\centerline{\epsfbox{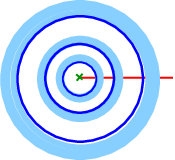}}
\caption{An admssible dg-sheaf  on  a punctured disc  stratified  by $m$ concentric outward cooriented circles, which is zero near the puncture, is a filtered vector space.}
\label{ncls101}
\end{figure}

Let $f: W^\bullet_1 \lra W^\bullet_2$ be a map of complexes of vector spaces with the following properties: 
  
  1. ${\rm Cone}(f)$    is  quasi isomorphic to a one dimensional vector space in the degree $0$.  
  
  2. $H^i(W^\bullet_1)=0$ if $i\not =0$. 
  
  Then $H^i(W^\bullet_2)=0$ if $i\not =0$,  and the induced map $H^0(W^\bullet_1) \lra H^0(W^\bullet_2)$ is injective. 
  
    We apply  this to the maps $\varphi_{i, i+1}$. Since  $H^\bullet(V_{0}^\bullet) = 0$, we get a flag of vector spaces:   
  $$
  V_0 \hra V_1 \hra V_2 \hra \ldots \hra V_m, ~~~~{\rm dim}V_i=i, ~~~~V_i := H^0(V_i^\bullet).
  $$
The monodromy around the puncture provides   an operator in the vector space $V_m$ preserving the flag. This just means that we get a framed $m$-dimensional local system of vector spaces over the punctured disc ${\cal D}^*$. Applying this argument at all punctures and special points we get an $m$-dimensional local system on the surface $\bS$ equipped with an invariant flag near each puncture and each special point, that is a framed local system on the decorated surface $\bS$.  
  \end{proof}

 \paragraph{Framed local systems on $\bS$ $\longleftrightarrow$  flat line bundles on an ideal bipartite graph $\Gamma$.} Let ${\cal L}$ be a  one dimensional local system  
  on a    ${\rm GL}_m-$bipartite graph $\Gamma$ on $\bS$. It determines a  ${\cal Z}$-admissible dg-sheaf  ${\cal F}_{\cal L}$ on $\bS$ via  Construction b). 
By the very definition, the zig-zag paths of $\Gamma$ can be deformed  using the elementary moves   to  a collection ${\cal Z}_m$,  
formed by $m$  nested   loops/arcs near each marked points. 
By Lemma \ref{FRAML}, the resulting   ${\cal Z}_m$-admissible dg-sheaf is just a framed $m$-dimensional local system on $\bS$.  
The elementary   moves  provide birational equivalences between the corresponding categories of ${\cal Z}$-admissible dg-sheaves.    

\vskip 2mm

Conversely, let ${\cal L}$ be  a framed $m$-dimensional local system on $\bS$.   Construction d) assigns  to   it   a ${\cal Z}_m$-admissible dg-sheaf on $\bS$. 
   Deforming back using the elementary  moves  the collection ${\cal Z}_m$ of loops/arcs around the marked points
 to the collection ${\cal Z}_\Gamma$ of zig-zag paths of the ideal bipartite graph $\Gamma$, we get a ${\cal Z}_\Gamma$-admissible 
dg-sheaf ${\cal E}_{\cal L}$ on $\bS$, which is zero near the marked points.  
\begin{figure}[ht]
\centerline{\epsfbox{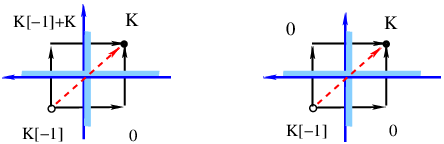}}
\caption{Two admissible dg-sheaves near a crossing point. Each provides a rank one local system on the microlocal support.}
\label{ncls107}
\end{figure}
  As   Figure \ref{ncls107} shows, 
this admissible dg-sheaf may not be acyclic on all mixed domain. However for a  generic  ${\cal Z}_m$-admissible dg-sheaf on $\bS$ we get 
an admissible dg-sheaf   acyclic on the mixed domain. 
The microlocal   support of the obtained ${\cal Z}_\Gamma$-admissible dg-sheaf is described by  flat line bundles on  zig-zag paths $\gamma$ of $\Gamma$.  
They are the same line bundles  we had on $\gamma$ before the deformation.    Construction c) gives us a flat line bundle on the graph $\Gamma$. 

\vskip 2mm
So     
  we get a  regular injective   map from the moduli space  ${\rm Loc}_1(\Gamma)$ of flat line bundles on the graph $\Gamma$ to the moduli space  ${\cal X}_{{\rm GL_m}, \bS}$ 
of framed $m$-dimensional local systems on $\bS$,   called the reconstruction map: 
\be \la{CLM}
R_\Gamma: {\rm Loc}_1(\Gamma) \hra {\cal X}_{{\rm GL}_m, \bS}.
\ee
It provides an open domain  in the moduli space space  ${\cal X}_{{\rm GL}_m, \bS}$, which is a cluster Poisson torus.
This is achieved by the following chain of    constructions: 
\be \la{1717}
\begin{split}
&\mbox{Flat line bundles  on a  rank $m$ ideal  bipartite graph $\Gamma$ on $\bS$} \stackrel{b)+ c)}{\longleftrightarrow}\\
&\mbox{${\cal Z}_\Gamma$-admissible 
dg-sheaves   on $\bS$, acyclic on mixed domains} \stackrel{\mbox{\rm Sec. \ref{Sec3}}}{\subset} \\
&\mbox{${\cal Z}_m$-admissible 
dg-sheaves, vanishing near the marked points} \stackrel{d)}{\longleftrightarrow}\\
&  \mbox{Framed $m$-dimensional local systems  on $\bS$.}\\
\end{split}
\ee
The middle map in (\ref{1717}) is an open embedding. The others are isomorphisms of moduli spaces.

 \begin{figure}[t]
\centerline{\epsfbox{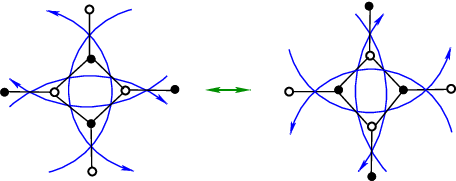}}
\caption{A two by two move of bipartite graphs.}
\label{ncls125}
\end{figure}

\paragraph{A Lagrangian surface in ${\cal L}_\Gamma \subset T^*\bS$.}  Given a bipartite graph $\Gamma$, note the following: 

1. Each $\bullet$ or $\circ-$domain ${\cal D}\subset \bS$ gives rise to a domain ${\cal D}\subset T^* \bS$ -  the zero section on ${\cal D}$.

2. Each cooriented zig-zag path $\gamma$ gives rise to the conormal bundle $T^*_\gamma\bS \subset T^* \bS$.

\bd  The surface ${\cal L}_\Gamma \subset T^*\bS$ is   the union of all these domains in $T^*\bS$:
\be
{\cal L}_\Gamma := \cup_\bullet {\cal D}_\bullet \cup_{\circ} {\cal D}_\circ \cup_\gamma T^*_\gamma
\ee
\ed

\begin{figure}[ht]
\centerline{\epsfbox{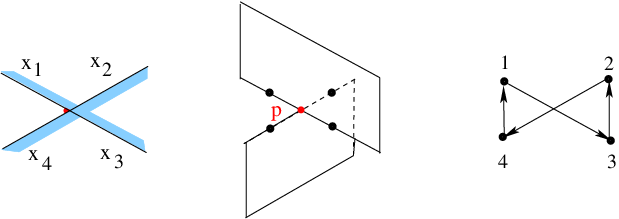}}
\caption{The link of the surface ${\cal L}_\Gamma$ near a crossing point $p$ is a circle $1\to 3\to 2\to 4\to 1$.}
\label{ncls130}
\end{figure}

\bl The surface ${\cal L}_\Gamma$ is homeomorphic to a smooth oriented surface,   identified with the oriented spectral surface $\Sigma_\Gamma$. 
\el

\begin{proof} The surface ${\cal L}_\Gamma$ is evidently homeomorphic to a smooth surface everywhere except the crossing point of the conormal bundles 
$p:=T^*_{\gamma_1}\bS \cap T^*_{\gamma_2}\bS$ given by the zero cotangent vector at the crossing point $\gamma_1 \cap \gamma_2$ of two zig-zag paths.  
To prove that it is homeomorphic to a smooth surface at $p$, we calculate the intersection of a small sphere at $p$ with 
$T^*_{\gamma_1}\bS \cup T^*_{\gamma_2}\bS$. We claim that it is homeomorphic to a circle. This is proved on Figure \ref{ncls130}, where we showed the link of a crossing point, observing at the same time that the projection ${\cal L}_\Gamma \lra \bS$ reversing the orientation on one of the triangles. (Indeed, both arrows $3\to 2$ and $4 \to 1$ look up on  Figure \ref{ncls130}). \end{proof}

\begin{figure}[ht]
\centerline{\epsfbox{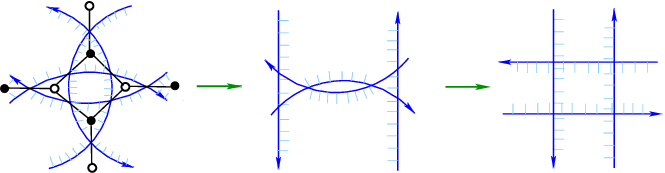}}
\caption{From the left to the right: two type II triangle  moves, followed by a bigon move.}
\label{ncls106}
\end{figure}

\begin{figure}[ht]
\centerline{\epsfbox{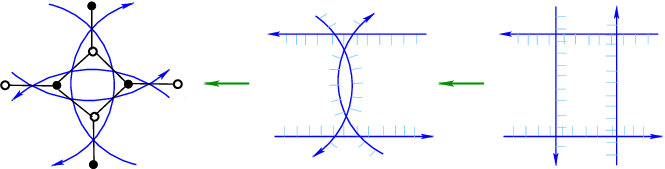}}
\caption{From the right to the left:  a bigon move, followed by two triangle  moves.}
\label{ncls124}
\end{figure}

\paragraph{Two by two moves of bipartite graphs.} Recall the two by two move $\Gamma_1 \to \Gamma_2$ of bipartite graphs on Figure \ref{ncls125}. It is decomposed into a sequence of  bigon and triangle moves 
on Figures \ref{ncls106} - \ref{ncls124}. Note that the collections of cooriented lines which we get in the process  are no longer associated with   bipartite graphs.


Starting from a   local system ${\cal L}_1$ on the graph $\Gamma_1$, we interpret it as an admissible dg-sheaf, and transform 
accordantly along the sequence of admissible moves above. 
In the process we get admissible dg-sheaf not associated with any bipartite graphs. However, assuming that we started from a generic local system ${\cal L}_1$ on the graph $\Gamma_1$, 
the resulting admissible dg-sheaf in the end does correspond to a local system ${\cal L}_1$ 
on the   graph $\Gamma_2$. So we get a birational isomorphism of stacks of local systems 
on the graphs $\Gamma_1$ and  $\Gamma_2$. In particular, we get a birational isomorphism between the stacks   of flat line bundles:
$$
{\rm Loc}_1(\Gamma_1) \lra {\rm Loc}_1(\Gamma_2).
$$
Our next goal is to get an explicit expression for this map for   any skew field $R$. 

Let   $B_1, B_2, B_3, B_4$ be the fibers of the local system over the $\bullet-$domains and   $W_1, W_3$ 
 the fibers over the $\circ-$domains on the left on Figure \ref{ncls106}. Then in the domain   on the Figure  the local system and hence the corresponding dg-sheaf are described by  
 the following diagram of isomorphisms:

\be \la{rs22} 
\begin{gathered}
    \xymatrix{
    &&B_2 &&\\
  B_1 &\ar[l]_{\sim} W_1 \ar[ru]^{\sim}  \ar[rd]_{\sim} &&\ar[lu]_{\sim} \ar[ld]^{\sim} W_3 \ar[r]^{\sim}    & B_3   \\
        &&B_4     & & }
\end{gathered}
 \ee

Let us set
$$
{\cal W}_1:=\frac{B_1\oplus B_2 \oplus B_4}{W_1}, ~~~~ {\cal W}_3:= \frac{B_2 \oplus B_3 \oplus B_4}{W_3}.
 $$
 Here we   identify $W_1$   with its image under the map $W_1 \hra B_1\oplus B_2 \oplus B_4$, and similarly for $W_3$. 
Then  after the triangle  moves at the two $\circ-$vertices   we get a dg-sheaf with the microlocal support shown in the middle picture on Figure \ref{ncls106}.  
Since an admissible dg-sheaf related to a bipartite graph has zero fibers over the mixed domains, the only relevant picture describing a single type IIa  triangle move  is the 
  fifth picture   on Figure \ref{cube1}. It follows that the resulting dg-sheaf is described by the following   diagram:
\be \la{rs22y} 
\begin{gathered}
    \xymatrix{
    &&\ar[ld] B_2 \ar[rd]&&\\
  B_1 \ar[r]&{\cal W}_1  & &  {\cal W}_3 &\ar[l]     B_3   \\
        &&\ar[lu] B_4   \ar[ru]  & & }
\end{gathered}
 \ee  
These maps are no longer isomorphisms   since  ${\rm rk}({\cal W}_1) = {\rm rk}({\cal W}_3)=2 {\rm rk}(B_i)$, and 
the diagram is no longer  assigned to a bipartite graph. 
There are  canonical isomorphisms
 $$
 {\cal W}_1    \stackrel{\sim}{\longleftarrow} B_2 \oplus B_4 \stackrel{\sim}{\lra}   {\cal W}_3.  
 $$
 The left one is  the embedding $B_2 \oplus B_4 \hra B_1\oplus B_2 \oplus B_4$, followed by the projection to  ${\cal W}_1$. The right one is similar. 
 Combining them we  get an isomorphism $f: {\cal W}_1 \stackrel{\sim}{\lra} {\cal W}_3$. Setting $X:= {\cal W}_1\stackrel{f}{=} {\cal W}_3$, we get a diagram
\be \la{rs23}  
\begin{gathered}
    \xymatrix{
    & B_2 \ar[d]&\\
  B_1 \ar[r]&X & \ar[l]     B_3   \\
        & B_4  \ar[u]   &  }
\end{gathered}
 \ee
 It describes the dg-sheaf on the right of Figure \ref{ncls106}. The transition   (\ref{rs22}) $\to$ (\ref{rs23}) describes the composition of the  three moves 
 on Figure \ref{ncls106}. 
Then we apply the  three moves  on Figure \ref{ncls124}. 
The corresponding transformations of  quiver representations are recorded on diagram (\ref{rs22c}). 

\be \la{rs22c} 
\begin{gathered}
    \xymatrix{
    &&B_2&&\\
    && \ar[ld]    W_2 \ar[u]\ar[rd] &&\\
   & B_1 && B_3       &   \color{red}{ \longleftarrow}  \\
        &&\ar[lu]    W_4 \ar[d] \ar[ru]  & &\\
        &&B_4  &&}
\end{gathered}
\begin{gathered} 
    \xymatrix{
    &&B_2\ar[d] &&\\
    &&    {\cal W}_2 &&\\
   & B_1 \ar[ru]\ar[rd] && \ar[lu]  B_3     \ar[ld]    &  \color{red}{ \longleftarrow}   ~~~~(\ref{rs23}).\\
        &&   {\cal W}_4   & & \\
        &&B_4  \ar[u]&&}
\end{gathered}
 \ee
 
 As the result   we altered the dg-sheaf in the domain ${\cal D}$  inside of the four $\bullet-$vertices, but did not alter it outside. 
 Let us calculate how the corresponding local system inside of the domain ${\cal D}$ changes. 
 Consider the following diagram on the left, where all maps   are isomorphisms:
 \be \la{rs22ca} 
\begin{gathered} 
    \xymatrix{
    &&   {{\rm L}}_2 &  \\
 & {\rm L}_1  & \ar[l]_{i_1} {\rm L}  \ar[u]^{i_2}  \ar[r]^{i_3} &  {\rm L}_3     \\
          }
\end{gathered}
\begin{gathered} 
    \xymatrix{
    &&    \ar[d] {{\rm L}}_2  &\\
 & {\rm L}_1 \ar[r] &  {\cal L}  &   \ar[l]  {{\rm L}}_3   \\
   }
\end{gathered}
 \ee
 Setting  ${\cal L}:= ({\rm L}_1  \oplus {\rm L}_2  \oplus {\rm L}_3)/{\rm L}$, it determines the diagram on the right. So there are three subspaces in ${\cal L}$. 
 Then the third subspace ${\rm L}_3$ provides a linear  
  isomorphism 
 $$
\varphi_{{\rm L}_1}: {\rm L}_2  \lra {\rm L}_3.
 $$ 
 We claim that  
 $ 
 \varphi_{{\rm L}_1} =   - i_3 \circ i_2^{-1}.
   $ 
 Indeed, for any   $l \in {\rm L}$, let $l_k:= i_k(l)$. Then there are three vectors in ${\cal L}$ which sum to   zero:
  $$
   (l_1, 0, 0)+ (0, l_2, 0)+ (0, 0, l_3) = 0 ~~\mbox{ in} ~~{\cal L}.
   $$  
 Now return to   diagram (\ref{rs23}), providing four subspaces in $X$. Then  there are two  isomorphisms:
 $$
\varphi_{{\rm B}_3}: {\rm B}_1  \lra {\rm B}_2, ~~~~\varphi_{{\rm B}_4}: {\rm B}_1  \lra {\rm B}_2.
$$  
Let  $M_{{\rm B}_2}: {\rm B}_2 \lra  {\rm B}_2$ be the \underline{minus} counterclockwise monodromy  in   diagram (\ref{rs22}), given  by  
$$
{\rm B}_2 \lra  {\rm W}_1 \lra  {\rm B}_4 \lra   {\rm W}_3 \lra  {\rm B}_2.
$$
\bl One has 
 $$
 \varphi_{{\rm B}_4} =  (1+ M_{{\rm B}_2}) \varphi_{{\rm B}_3}.
 $$
  \el
 
\begin{proof}   If our local system is a 
 rank one flat  $R-$line bundle,  it is  proved by the calculation done in proof of Proposition \ref{10.12.11.5X}. \end{proof}

 \vskip 2mm Changing $\{$$\bullet-$vertices$\}$ $\longleftrightarrow$ $\{$$\circ-$vertices$\}$  and reversing the arrows amounts to dualisation of the diagrams. Precisely,  let   $W_1, W_2, W_3, W_4$ be the fibers of the local system over the black domains, and   $B_1, B_3$ 
 the fibers over the white domains. Then  the local system is described by  
 the following diagram where all maps are isomorphisms:

\be \la{rs22x} 
\begin{gathered}
    \xymatrix{
    &&\ar[ld]_\sim W_2 \ar[rd]^\sim &&\\
  W_1 \ar[r]^\sim  &B_1  &&B_3   &  \ar[l]_\sim  W_3   \\
        && \ar[lu]^\sim W_4   \ar[ru]_\sim  & & }
\end{gathered}
 \ee
Let us set
$$
{\cal B}_1:={\rm Ker}\Bigl(W_1\oplus W_2 \oplus W_4 \lra B_1\Bigr), ~~~~ {\cal B}_3:={\rm Ker}\Bigl(W_2\oplus W_3 \oplus W_4 \lra B_3\Bigr).
 $$
After the  triangle  moves at the  two  $\bullet-$vertices   we get a   dg-sheaf  described by the   diagram 
\be \la{rs22xx} 
\begin{gathered}
    \xymatrix{
    &&  W_2  &&\\
  W_1 &\ar[l]\ar[ur] \ar[rd]{\cal B}_1  & & \ar[lu] {\cal B}_3 \ar[ld]\ar[r] &    W_3   \\
        && W_4     & & }
\end{gathered}
 \ee   
The composition of  canonical isomorphisms
 $
 {\cal B}_1   \stackrel{\sim}{\lra} W_2 \oplus W_4 \stackrel{\sim}{\longleftarrow}    {\cal B}_3 
 $
 provides an isomorphism $g: {\cal B}_1 \stackrel{\sim}{\lra} {\cal B}_3 $. Setting $Y:= {\cal B}_1\stackrel{}{=} {\cal B}_3$, we get a diagram
\be \la{rs23x} 
\begin{gathered}
    \xymatrix{
    & W_2  &\\
  W_1& \ar[l]\ar[u] Y  \ar[r] \ar[d]  &  W_3   \\
        & W_4   &  }
\end{gathered}
 \ee
 
  \paragraph{Conclusion.} Section \ref{SSec4} gives an alternative approach to  main definitions $\&$ results  of Section \ref{SEC5.1}. 
  
  \paragraph{ Twisted and untwisted dg-sheaves.} Given a sheaf of categories ${\cal C}$ on a space $X$, and a cocycle $c_2$ representing a class in 
  $H^2(X, {\rm Cent}({\cal C})^*)$, where   ${\rm Cent}({\cal C})^*$ is given by the invertible elements of   ${\cal C}$, 
  one can make a new sheaf of categories $\widetilde {\cal C}$,  by twisting ${\cal C}$ by the 2-cocyle   $c_2$. Precisely, take an open cover  $\{U_i\}$ of $X$. 
  Then an object of $\widetilde {\cal C}$ over $X$ is given by  objects $A_i$ of   categories ${\cal C}(U_i)$,  related by isomorphisms $f_{ij}: A_i \lra A_j$ on $U_i \cap U_j$ which 
  satisfy the following  relations on the triple intersections: 
  $$
   f_{jk}  \circ f_{ij} = c_2(U_i, U_j, U_k) f_{ik}
  $$  
  
Below we  use  twisting by the second Stiefel-Whitney class ${\rm sw}_2(\bS)\in H^2(\bS, \Z/2\Z)$.  
  We can twist both the Fukaya category and the category of constructible or dg-sheafs. The Fukaya category of objects with the microlocal support on $T^*_Z\bS$ is naturally equivalent to the twisted constructible category on $\bS$, and vice versa, see \cite[Remark 5.2]{GPS}.   
  Note that ${\rm sw}_2(\bS)=0$  unless 
  $\bS$ is a compact non-oriented surface. So in our case the twisted category is equivalent to the original one. 
Yet the twisted version is somewhat better. It is the one   providing positivity phenomenon in all formulas.  

\section{Non commutative  stacks of framed Stokes data} \la{sec8.4}

At Section \ref{sec8.4} we adopt a  faster pace of the exposition than in the rest of the paper. 

\subsection{Stokes data: local and global pictures   } \la{sec9.1}
  
 \paragraph{1. Admissible sheaves and local  Stokes data.} We start with   two  topological  definitions. 
  
  \bd \la{DEFSDI} A collection ${\cal L}$ of smooth oriented loops $L_1, \ldots , L_n$ in the   punctured disc ${\rm D}^*\subset {\Bbb C}^*$ is called {\em ideal} if 
   it enjoys the following properties: 
   
   \begin{enumerate}
    \item 
 The  (self)intersection $L_i \cap L_j$ of any two   loops is   transversal.

  \item Each  loop  $L_i$ rotates to the right, that is the argument function $\varphi$ on ${\rm D}^*$ 
  strictly decreases as we go along the loop following its orientation. We coorient the loops out of the puncture.

  \end{enumerate}
  \ed 
  
  Cooriented loops $L_i$ are the connected components of the Legendrian  link in ${\rm ST}({\rm D}^*)$ provided by   ${\cal L}$.

  We  refer to ${\cal L}$  as   an {\it ideal  Legendrian link}, although 
   it is only a projection of a  Legendrian. 
   
  \bl 
  An ideal Legendrian link on the punctured disc ${\rm D}^*$ is the same thing as an element of the cyclic envelope $[{\rm Br}_n^+]$ of the braid semigroup, $n \geq 0$. 
  \el 
  
\begin{proof} Cut the collection of curves ${\cal L}$ by a generic ray $r$ from the origin. 
  Scanning the ${\cal L}$ by rotating the ray to the right we get a sequence of permutations, and hence an element of the braid semigroup ${\rm Br}_n^+$, given by their product. The $n=0$ case 
is the empty link.    This construction can evidently be inverted. Scanning stating from a different ray produces the same element of the cyclic braid semigroup. 
  \end{proof}

     Consider  the subcategory ${\cal C}_0({{\rm D}^*}, {\cal L})$ of the category  ${\cal C}({{\rm D}^*}, {\cal L})$ of admissible dg-sheaves   on  ${\rm D}^*$ for   an ideal Legendrian  ${\cal L}$ consisting of admissible dg-sheaves which  are    
   acyclic near the puncture. By Proposition \ref{SSD}, they are automatically concentrated in the degree zero, so we refer to them just as admissible sheaves. 
   We denote by   $[{\cal L}]$  the homotopy class of an ideal Legendrian ${\cal L}$ on ${\rm D}^*$ under Legendrian isotopies keeping Legendrians  ideal and disjunct,   
   called admissible isotopies. 
       Categories assigned to admissibly isotopic  Legendrians ${\cal L}$ are equivalent by Theorems \ref{TRTB}-\ref{TRT}. So the category 
   ${\cal C}_0({{\rm D}^*}, {\cal L})$ is determined  by the admissible isotopy class  $[{\cal L}]$.

    \bd \la{DEFSD} An ideal Legendrian ${\cal L}$ in the   punctured disc ${\rm D}^*\subset {\Bbb C}^*$ is of {\em  Stokes type} if 
   it  
   
   \begin{itemize}
  
  \item 
 The  intersection $L_i \cap L_j$ of any two   different components   of ${\cal L}$ is non-empty.   
 
  \end{itemize}
  \ed 
  
    \bd \la{DEFLSD}  A local {\em Stokes data} is a local system ${\cal V}$ of vector spaces on the circle   of   rays in the  tangent space   at zero, 
   with an increasing  filtration ${\cal F}_\bullet$ 
  everywhere except   the fibers over a 
  finite number of 
  rays, called {\it Stokes rays}, 
  related near   them as follows. 
  
  Denote by ${\cal F}^\pm_\bullet({\cal V})$    filtrations   to the right and   to the left of a Stokes ray 
  $r$. Then  ${\cal F}^+_i({\cal V}) = {\cal F}^-_i({\cal V})$  if   $i\not =a$, and   the    subquotient 
 $
   {\rm gr}_{[a+1, a]}{\cal F}({\cal V})
 $ is the direct sum of  its 
  subspaces ${\rm gr}_{a}{\cal F}^\pm({\cal V})$:
  \be \la{GKOND}
   {\rm gr}_{[a+1, a]}{\cal F}({\cal V}):= {\cal F}_{a+1}({\cal V})/  {\cal F}_{a-1}({\cal V}) =  {\rm gr}_{a}{\cal F}^+({\cal V}) \oplus {\rm gr}_{a}{\cal F}^-({\cal V}).
   \ee
  \ed

 \bp \la{SSD} The following categories are canonically equivalent:
 
 \begin{enumerate} 
 
 \item  The category of   local Stokes data from Definition \ref{DEFLSD}.

 \item The  direct sum  of categories  $\bigoplus_{[{\cal L}]}{\cal C}_0({{\rm D}^*}, {\cal L})$     over all  homotopy types $[{\cal L}]$ of   Stokes Legendrians.  \end{enumerate}
  
 \ep
  
\begin{proof} Before we proceed to the proof, let us state the   following simple Lemma. 

\bl \la{ACP} An admissible sheaf near a  simple  crossing point 
 \begin{figure}[ht]
\centerline{\epsfbox{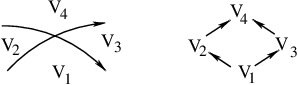}}
\caption{An admissible sheaf near a crossing point.}
\label{quiver2}
\end{figure}  
is determined by a commutative diagram of four vector spaces  $V_i$  on the right of Figure \ref{quiver2}, 
 where all arrows are injective, and the   map 
\be
V_2/V_1 \oplus V_3/V_1 \stackrel{\sim}{\lra} V_4  
\ee 
 is an isomorphism. 
 \el

Now let us prove the Proposition. \vskip 2mm

$2) \lra 1)$. The circle $S^1$ of rays minus the Stokes rays is a union of a finite number of intervals:
$$
S^1 - \{\mbox{Stokes rays}\} = {\rm I}_1 \cup  \ldots \cup  {\rm I}_m.
$$
 Since the sheaf is zero near the puncture,  the admissibility just means that on each of the  
intervals ${\rm I}_k$ we get   a local system of filtered vector spaces. 
It follows from   Lemma \ref{ACP} that the admissibility implies that this  local system extends across the Stokes rays, while the filtration is altered:  
the filtration in  ${\rm gr}_{[a+1, a]}$ flips, 
and the rest of the filtration data is preserved.

 $1) \lra 2)$.  We start from a local Stokes data on $S^1$. For each   ray $r\in {\rm I}_k$, we assign  the  associate graded quotients ${\rm gr}_{i}{\cal F}^+({\cal V})$  to some points of the ray so that the index $i$ grows when we go along the ray, 
 and  when the ray moves in the interval ${\rm I}_k$, we get a collection of curves together with local systems of 
   vector spaces on them for each sector between the Stokes rays. We erase the curves with the zero local systems on them.  
   
Let us glue these   curves with   local systems, which we have so far  in the sectors of the punctured disc determined by the Stokes rays. 
 For each  $i \not =a, a+1$, at a Stokes ray $r$, we glue the curves which carry  the local systems   ${\rm gr}_{i}$,  and then 
 glue the local systems on them using the local system structure on 
   ${\cal V}$ near the ray. 
  Finally, we glue the curve which carries ${\rm gr}_{a}$ on the one side of $r$ with the one with ${\rm gr}_{a+1}$ on the other side, and   glue   local systems on them via the isomorphisms     provided by 
  condition (\ref{GKOND}):    \be
  {\rm gr}_{a}{\cal F}^-({\cal V}) \stackrel{\sim}{\lra} {\rm gr}_{a+1}{\cal F}^+({\cal V}), \ \ \ \ \ \    {\rm gr}_{a+1}{\cal F}^-({\cal V}) \stackrel{\sim}{\lra}   {\rm gr}_{a}{\cal F}^+({\cal V}).
    \ee
We get a collection ${\cal L}$ of loops   on the   disc, which never turn back, with   simple crossings at the Stokes rays. 
  Each of the loops carries   a local system of vector spaces.

The claim that we get an admissible sheaf follows  from Lemma \ref{ACP}. Proposition \ref{SSD} is proved. 
  \end{proof}

\paragraph{Remark.} One can relax condition 1) in Definition \ref{DEFSD} allowing   $m>2$ curves   intersect at a point,   with 
different tangents.  Let us describe   the corresponding analog of condition (\ref{GKOND}) in Definition \ref{DEFLSD}. 
Recall that a pair of length $m$ flags ${\cal F}_\bullet$ and ${\cal G}_\bullet$ in a vector space $V$  is complimentary if 
$V = {\cal F}_a \oplus {\cal G}_b$ if $a+b=m$. 

\begin{figure}[ht]
\centerline{\epsfbox{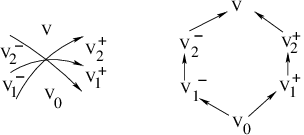}}
\caption{An admissible sheaf near a triple crossing point: the following two flags  in $V/V_0$  are complimentary to each other: $V_1^-/V_0 \subset V_2^-/V_0$ and $V_1^+/V_0 \subset V_2^+/V_0$. }
\label{stokes3}
\end{figure}

\bl \la{ACP1} An admissible sheaf near an $m-$fold  crossing point with different tangents,  see  Figure \ref{stokes3},  
is determined by a pair of vector spaces $V_0\subset V$ and two  length $m$ complimentary flags   in   $V/V_0$.   
 \el
 
 \begin{proof} Given such an admissible sheaf, going from the bottom to the top on the left of the crossing points gives one flag, and going from the right gives the other flag. 
 \end{proof}

Then there is an analog of Proposition \ref{SSD} if instead of (\ref{GKOND}) we have the following condition: 
$
\mbox{\it The left and right flags in the subquotient associated with the crossing point are opposite to each other. } 
$ \\
 
Our next goal is to explain how  Stokes Legendrian links  ${\cal L}$   and   the categories ${\cal C}_0({\rm D}^*, {\cal L})$ 
arise in complex geometry, as topological  invariants   of   holomorphic vector bundles  with connections on a punctured   disc 
with any, possibly irregular, singularity at zero.  For this we need to recall 
 classification of   connections on a formal punctured   disc   ${\rm Spec}~\C((z))$,  due to Hukuhara, Turrittin,  Levelt \cite{M2}. 

\paragraph{2. Formal classification of irregular singularities.} A Puiseux series is a Laurent series in $z^{1/d}$, $d \in \Z_{\geq 1}$. 
The field of all Puiseux series   
  is the algebraic closure $ \overline {\C((z))}$ of the field $\C((z))$ of Laurent series in $z$, with the Galois group $\widehat \Z$. 
The integral closure of the ring of Taylor series $\C[[z]]$ is the local ring ${\cal O}_{\overline {\C((z))}} $ of this field. It is given by Taylor series in $z^{1/d}$, $d \geq 1$. 

Consider  the quotient of the field of  Puiseux series by  its local ring:  
\be \la{Q}
  \overline {\C((z))}\  / \ {\cal O}_{\overline {\C((z))}} = \lim_{{d \geq 1}}\ \C((z^{1/d})) \ /  \ \C[[z^{1/d}]].
\ee

\bd \la{DEFIND} The space ${\cal G}$ is the space of orbits  of the Galois action   on the elements $f(z)$ in  (\ref{Q}). 
\ed

So an element $f(z) \in   {\cal G}$  is presented as a Puiseux polynomial in $z^{-1/d}$ for some positive integer $d$:  
\be
f(z)= \sum_{k \geq 1} \lambda_{k} z^{-k/d}, \ \ \ \ \lambda_{k} \in \C, 
\ee
 where we  identify   $\{f(z)\}\sim \{f(e^{2\pi i/d}z)\}$, and  assume $d$ is the smallest such integer. 

Consider the   $d:1$ cover   of the formal punctured  disc:  
$$
\pi_d: {\rm Spec}\ \C((w))\lra {\rm Spec}\ \C((z)), \ \ \ \ z=w^d, \ \ d\geq 1.
$$
Take a holonomic ${\cal D}-$module on this cover, given by the tensor product of the ${\cal D}-$module generated by  $e^{f(z)}$, and 
a local system  ${\cal V}_w$  on     
${\rm Spec}~\C((w))$:
$$
e^{f(z)} \otimes {\cal V}_w.
$$
The ${\cal V}_w$ is determined by a local system on the circle  of the rays in the tangent space to $0$ in the $w-$plane.  
We consider the push-forward  of this ${\cal D}-$module under the map $\pi_d$:
$$
\pi_{d*}(e^{f(z)}\otimes {\cal V}_w).
$$

\bt \la{DM} The category of holonomic  ${\cal D}-$modules on the formal punctured disc is decomposed into blocks, parametrised by the elements of ${\cal G}$. So 
  any holonomic ${\cal D}-$modules on ${\rm Spec} ~ \C((z))$ is isomorphic to a direct sum of the standard ${\cal D}-$modules 
\be \la{DIRS}
\bigoplus_{j=1}^m\pi_{d_j*}(e^{f_j(z)} \otimes {\cal V}_{w_j}), \ \ \ \ f_j(z) \in {\cal G}.
\ee  
\et

\paragraph{3. Legendrian links ${\cal L}$ of Stokes type, and local Stokes data.} Take a flat connection  on a disc with a single singular point. By Theorem \ref{DM}, its restriction 
to  the formal punctured disc  is   a direct sum   (\ref{DIRS}). So it determines 
a collection of different Puiseux polynomials with negative exponents:\footnote{An irreducible analytic flat connection on a punctured disc, restricted to the formal punctured disc, can be a sum of several ${\cal D}-$modules from different blocks.} 
\be
f_j(z)= \sum_{k \in \Z_{\geq 1}} \lambda_{j,k} z^{-k/d_j}, \ \ \ \ \lambda_{j,k} \in \C, \ \ d_j \in \Z_{\geq 1}, \ \ \ \ j=1, ..., m.
\ee
 They define  multivalued functions $f_j(z)$. Restricting them  to a certain concrete but sufficiently small radius $\varepsilon$ circle $\{\varepsilon e^{i \theta}\}$ centered at $0$, 
  we get a collection of curves $\{L_j\}$ given parametrically by  
\be
Z_j(\theta):=    R_j(\theta) \cdot e^{i\theta}, \  \ \ \ \ R_j(\theta) := e^{{\rm Re} f_j(\varepsilon e^{i \theta})} \ \ \ \ \ j=1, ..., m. 
\ee
Taking the 1-jet of the radial part $R_j(\theta)$ of the function $Z_j(\theta)$,  we get   a Legendrian knot $L_j$ in the space ${\cal J}^1(S^1)$ of 1-jets of real functions on the circle. 
Explicitly,  
  $$
\theta \lms (\theta,  R_j(\theta), d  R_j(\theta)).
$$
 Their union is a Legendrian link ${\cal L}$ in ${\cal J}^1(S^1)$. We call it a Legendrian link  of {\it Stokes stype}.  For a generic small $\varepsilon$ 
 it fits Definition \ref{DEFSD}.

 The push-forward of the local systems ${\cal V}_{w_j}$ provides a local system ${\cal V}_j$ of vector spaces on the loop $L_j$. 
   
  The  pairs $(L_j, {\cal V}_j)$  describe the microlocal support of   an admissible sheaf.

  The filtration by the order of growth provides the filtration at   all but Stokes rays.  
  
  {Therefore we get a local Stokes data from 
  Definition \ref{DEFLSD}.} \\
  
     Legendrian links ${\cal L}$ assigned to different sufficiently small   $\varepsilon$ are isotopic, although considered as collections of curves they may look differently. Therefore the corresponding categories are canonically equivalent by Theorem \ref{TRT}.    \\
     
     One can  consider, as Deligne and Malgrange do 
  \cite{M1}, \cite{M2},   the limiting case when $\varepsilon \to 0$. Then the corresponding Legendrian link may have 
  non-transversal crossings, with different tangents. Note that Deligne-Malgrange Stokes rays are determined by the singularity, while our  Stokes rays depend also on the choice of a small radius $\varepsilon$.

\begin{figure}[t]
\centerline{\epsfbox{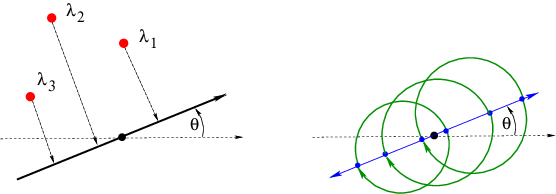}}
\caption{On the left: projecting red complex points $\lambda_1, \lambda_2, \lambda_3$ onto a line. On the right: rotating their exponentials - blue points on the  axis. 
The curve ${\cal L}$ is the union of three green loops obtained this way.}
\label{stokes}
\end{figure}

\paragraph{4. Example.}  Given an $n \times n$ matrix $A$ with distinct   eigenvalues $\{ \lambda_j\}$, consider a differential equation  
\be \la{SSS}
df = \frac{A}{z^2} f(z) dz.
\ee
Its solutions  in the eigenbasis of $A$  are linear combinations of   functions $f_j(z) = e^{-\lambda_j/z}$. To draw the curve ${\cal L}$, we project orthogonally   points $\{\lambda_j\}$  
onto the oriented line obtained by rotating the real axis by an  angle $\theta$,  exponentiate  their coordinates  on the line, and put them onto the positive ray in the   direction  $\theta$.
Rotating $\theta$ from $0$ to $2\pi$, these points move along  the curve ${\cal L}$. 
Rotating 
$\theta \to \theta+\pi$, we change the order of the points on the ray   to the opposite one, see Figure \ref{stokes}.  So any two of these loops intersect. 
Therefore we get a collection of loops as in Definition \ref{DEFSD}. 

The Stokes rays are the rays 
${\rm Re}(\lambda_i-\lambda_j)=0$ in $\C$. 
The local Stokes data  from  Definition \ref{DEFLSD} is given by filtration of the space of solutions in a sector around a non-Stokes ray $r$ by the growth of a solution. It satisfies the tranversality assumption (\ref{GKOND}) if no three $\lambda_i$'s are on the same line. 

In a more general situation, solutions of type $f_j(z) = z^{\alpha} e^{-\lambda_j/z}$ describe one dimensional local systems on these circles with the monodromies   $e^{2\pi i \alpha}$. 

 \paragraph{5. Legendrian links of analytic Stokes type.} Only some  Stokes Legendrians   describe the singular type  of a connection  on a disc. We call them 
 Legendrian links  of  {\it analytic Stokes} type.  
 
 One can describe  Legendrian links  of  {analytic Stokes} type  as follows. 
Take a    radius $\varepsilon$ sphere $S^3\subset \C^2$,  and delete from it the circle  in the $x=0$ plane.  The obtained contact threefold $S^3_{(0)}$ is identified with 
 ${\cal J}^1(S^1)$. 
 Let us consider a germ of an algebraic curve  $C$ in $\C^2$, reduced, but possibly reducible, 
 with a singularity at $(0,0)$. So it is the set of zeros of a product of  different irreducible polynomials in two variables. 
 Intersecting  the curve $C$ with   $S_{(0)}^3 \subset \C^2$ we get a Legendrian link. All Legendrian links of analytic Stokes type are obtained this way.  \vskip 2mm
 
 We are now return back to arbitrary  Stokes Legendrinans,  and describe the corresponding category of admissible sheaves in terms of representations of a certain quiver.

\begin{figure}[t]
\centerline{\epsfbox{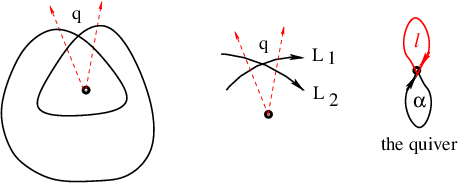}}
\caption{At the crossing point $q$, the local arc $L_1$ is above the one $L_2$,  producing a loop $\alpha$ of the quiver. The loop $l$ on the diagram, shown red, 
corresponds to the loop itself}
\label{quiver1}
\end{figure}

\paragraph{6. Quiver description.} Consider  an 
   oriented  quiver $Q_{\cal L}$ with the vertices $v_i$ assigned to the oriented loops $L_i$ of   ${\cal L}$,  and the following  arrows  
   assigned to  the loops and crossing points of ${\cal L}$, see Figure \ref{quiver1}:
  
    \begin{enumerate}
   \item 
 For each vertex $v$, there is a simple loop $l_v: v \lra v$     at this vertex.   
  
  \item For each crossing  $q\in L_i \cap L_j$,   we have an arrow $\alpha_q: v_i \to v_j$ if  $L_i$ is above  $L_j$   to the right  of  $q$.

  \end{enumerate}
  
    \bt \la{SD1} \begin{enumerate}

  \item  The category ${\cal C}({{\rm D}^*},  {\cal L})$   is abelian.        It  is equivalent to the category ${\cal R}(Q^*_{\cal L})$of representations   of the  quiver $Q_{\cal L}$, with invertible arrows  $l_v$.

\item  If a Legendrian link ${\cal L}$ is of analytic Stokes type, that is arises from a   singular point   of a holomorphic vector bundle with connection on ${\rm D}^*$,  possibly irregular, 
       then the stack ${\cal M}_{{\rm D}^*}({\cal L})$  is equivalent to    the  Betti stack of given irregular type, also known as  the stack of   Stokes data.
       
       \end{enumerate}      
         
         \et

{\it Remark.} Since the category of representations of the quiver $Q_{\cal L}$ is  abelian,   the category ${\cal C}_{{\rm D}^*}( {\cal L})$   is abelian. 
The category ${\cal C}({{\rm D}^*},  {\cal L})$ assigned to a general ideal Legendrian is not   abelian.  
Indeed, if   ${\cal L}$ is isotopic to a collection of concentric circles, we get the category of filtered vector spaces.

  \begin{proof} 
  
  i) Pick a reference point $p_i$ at each loop $L_i$ of ${\cal L}$.  
  
  For each intersection/selfintersection point $q\in L_i \cap L_j$,   let $\kappa_1, \kappa_2$ be little intervals on the loops $L_i$, $L_j$ containing  $q$. 
  They are   ordered so that the interval $\kappa_1$ goes above   $\kappa_2$ to the right of $q$, see Figure \ref{quiver1}.   Pick  paths    in the positive direction  connecting the point 
  $q$ with   reference points   $p_i$ and $p_j$. \\

   Let us assign to a representation ${\cal R}$ of the quiver $Q_{\cal L}$ a sheaf ${\cal S}_{\cal R}$ on the punctured disc ${\rm D}^*$. 
   
   First, we define a local system ${\cal V}_i$ on each   loop  $L_i$. Recall the vertex $v$ of the quiver assigned to the loop $L_i$. 
   Put the vector space  
  ${\cal R}(v)$  to the  reference point $p_i$ on $L_i$. Then the operator ${\cal R}(l_v): {\cal R}(v)\lra {\cal R}(v)$ for the loop $l_v$ desrcibes   
   the monodromy of the local system   on the loop $L_i$ in the positive direction.

   Next, we define the fiber of the  sheaf  ${\cal S}_{\cal R}$ at a point  $d\in {\rm D}^* - \{\mbox{Stokes rays}\}$.  
  Denote by $\gamma_1, \ldots , \gamma_k$ the little intervals on the curve ${\cal L}$ containing the intersection points of the  segment  connecting $0$ and $d$ with ${\cal L}$. 
  Then 
  fiber of the  sheaf ${\cal S}_{\cal R}$ at $d$ is    the direct sum of the fibers  of the local systems ${\cal V}_i$   at the   intervals:
  $$
  {\cal S}_{\cal R}(d) :=   \oplus_{i=1}^k {\cal V}_i(\gamma_i).
  $$  
 We order the  intervals   away from $0$, and filter  the fiber   by      subspaces $\oplus_{i=1}^p {\cal V}_i(\gamma_i)$. \\

  We extend the  sheaf to ${\rm D}^* - \{\mbox{crossing points}\}$ in the obvious way.  Note that the filtration above is defined only 
  away from the Stokes rays. 
  Let us use the maps $  {\cal R}(\alpha_q)$ assigned to the arrows $\alpha_q$ of the quiver to extend the   sheaf  to   the crossing points. 
  
   Denote by ${\cal V}_{\kappa_i}$   the fiber  of the local system on ${\cal L}$ over the interval    $\kappa_i$, $i=1,2$. 
   The stalk ${\cal S}_{\cal R}(q)$ of the sheaf at a crossing point $q$ is mapped by the restriction map to the stalks ${\cal S}_{\cal R}(q_\pm)$   on the right/left of $q$. 
  These maps are isomorphisms, since the sheaf 
  does not have microlocal support in those directions. So an extension of the sheaf to the point $q$ provides  isomorphisms
  $$
{\cal V}_{<q} \oplus {\cal V}_{\kappa_1} \oplus {\cal V}_{\kappa_2}=: {\cal S}_{\cal R}(q_-) \stackrel{\sim}{\longleftarrow} {\cal S}_{\cal R}(q)  \stackrel{\sim}{\lra} {\cal S}_{\cal R}(q_+)  :=   {\cal V}_{<q} \oplus {\cal V}_{\kappa_1} \oplus {\cal V}_{\kappa_2}.
$$
 So to   extend  the sheaf to the point $q$ one just needs to define a linear automorphism of ${\cal V}_{<q} \oplus {\cal V}_{\kappa_1} \oplus {\cal V}_{\kappa_2}$, which is 
 the identitity on  ${\cal V}_{<q}$, and  preserves   ${\cal V}_{\kappa_1} \oplus {\cal V}_{\kappa_2}$. We define the latter    by setting 
    \be
  (v_1, v_2) \lra (v_1, {\cal R}(\alpha_q)(v_1) + v_2).\ \ \ \ \ \ 
   {\cal R}(\alpha_q): V_{\kappa_1}\lra V_{\kappa_2}.
\ee
Then  it also preserving the filtration $V_{\kappa_1} \subset V_{\kappa_1} \oplus V_{\kappa_2}$, as required in our situation. \\

      This functor is an equivalence. This is a non-trivial result, and we do not present here a proof.    \\

ii) This is a reformulation of the known description \cite{DMR}, \cite{M2} of the analytic Stokes data in terms of the filtration of the space of solutions given by the order of growth, and the way it changes when we cross   Stokes rays. \end{proof}
  
       \paragraph{7. Global Stokes data.} Below $S$ is an   oriented topological surface, possibly  with punctures. 

\bd  \la{D9/11}A Legendrian ${\cal L}$ on   $S$ is {\em ideal} /   {\em  Stokes type},    if it  
 admissibly isotopic to a disjoint union  of  Legendrians  ${\cal L}_p$   from Definition \ref{DEFSDI} / \ref{DEFSD}  near  punctures $p$ on $S$:
 $$
 {\cal L} = \cup_p  {\cal L}_p, \ \ \ \ {\cal L}_p\cap {\cal L}_q=\emptyset.
  $$
  \ed

\bd A {\em Stokes data} on   $S$   is   a local system ${\cal V}$ of vector spaces on $S$, plus a local Stokes data at the  specialization of ${\cal V}$ 
  to each puncture.  \ed
  
  The following Proposition is an immediate consequence of Proposition \ref{SSD}. 
 \bp Given an ideal Legendrian ${\cal L}$ on $S$, 
  the category  of Stokes data on   $S$ of the topological type ${\cal L}$ is equivalent to the category ${\cal C}_0(S; {\cal L})$  
  of admissible sheaves on $S$ with the microlocal support on  ${\cal L}$ and the zero section of $T^*S$,  which vanish near the marked points. 
  \ep

 \paragraph{8. Non-commutative stacks of Stokes data.}    Denote by ${\cal S}({\cal L}, d)$  the  stack of Stokes data 
  of a given topological type $({\cal L}, d)$ on $S$. The topological types are described by the pairs $({\cal L}, d)$, where ${\cal L}$ is an ideal Legendrian links    on ${\rm ST}(S)$ from  
  Definition \ref{D9/11}, 
  modulo    admissible    Legendrian isotopies, equipped with local systems of vector spaces on the connected components ${\cal L}_\alpha$ of the  Legendrian   ${\cal L}$.  Dimensions of   local systems are encoded by the dimension vector $d = \{d_\alpha\}$.  When ${\cal L}$ is empty, we recover the  stack ${\rm Loc}_m(S)$.   
    
\bd \la{NCSSD} The non-commutative stack of Stokes data ${\cal S}^{\rm nc}({\cal L}, d)$ of the topological type $({\cal L}, d)$ is   the stack  ${\cal M}_0({\cal L}, d)$ of ${\cal L}-$admissible dg-sheaves of $R-$vector spaces on $S$   vanishing near the punctures of $S$,  whose microlocal support on ${\cal L}_\alpha$ is  a   local system on the loop  ${\cal L}_\alpha$  of dimension $d_\alpha$:
$$
{\cal S}^{\rm nc}({\cal L}, d):= {\cal M}_0({\cal L}, d).
$$
  \ed

         \subsection{The Riemann-Hilbert correspondence and the wild mapping class group}  \la{SECT9.2} 
         
         In Section \ref{SECT9.2} we work in the commutative set-up. 
         
     \paragraph{1. The classical Riemann-Hilbert correspondence.} Given a Riemann surface $\Sigma$, recall     the de Rham stack ${\cal M}^{\rm reg}_{\rm DR}(m, \Sigma)$  of rank $m$ 
 vector bundles with  meromorphic connections  with regular singularities on   $\Sigma$. 
Its Betti variant     is   the stack ${\rm Loc}_m(S)$      of    local systems of    $m-$dimensional vector spaces on         the  topological surface $S$  underlying $\Sigma$.

 The classical Riemann-Hilbert correspondence 
         \cite{D} identifies  the associated complex analytic stacks: 
       $$
       {\cal R}{\cal H}^{\rm reg}_\Sigma: {\cal M}^{\rm reg}_{\rm DR}(m, \Sigma)(\C)    \stackrel{\sim}{\lra}     {\rm Loc}_m(S)(\C).
       $$
 However the algebraic structures of these stacks are  different.        
       The algebraic structure of ${\cal M}^{\rm reg}_{\rm DR}(m, \Sigma)$  depends on               
     the complex structure on   $\Sigma$.  When $\Sigma$ varies, we get the algebraic de Rham  stack ${\cal M}^{\rm reg}_{\rm DR}(m)$   over the moduli space ${\cal M}_{g,n}$. 
  Contrary to this, the algebraic structure of the Betti stack  depends only on               
     the topology of  $\Sigma$. 
    When $\Sigma$ varies, we get a local system of  the Betti stacks   over   ${\cal M}_{g,n}$. 
       Equivalently,  the mapping class group $\Gamma_S = \pi_1({\cal M}_{g,n})$ of $S$ acts by automorphisms of the stack $ {\rm Loc}_m(S)$. 
       
       Summarising, we get the canonical universal Riemann-Hilbert equivalence: 
   \be \la{BeDRa} 
 \begin{gathered}
     \xymatrix{
   {\cal M}^{\rm reg}_{\rm DR}(m) (\C) \ar[dr]   \ar[rr]^{ {\cal R}{\cal H} } & & \ar[dl]{\cal M}^{\rm reg}_{\rm B}(m)(\C)\\
         &    {\cal M}_{g,n}    &   }
 \end{gathered}
  \ee    
To generalize it  to irregular singularities, we need   the Betti side, given by    moduli spaces of  Stokes data.

\paragraph{2.         The   Riemann-Hilbert correspondence.} 
 
 We start with a pair $(S, {\cal L})$, where  ${\cal L} = \cup{\cal L}_p$ is  a  Legendrian link of  analytic
  Stokes type on  a surface $S$, considered  up to   admissible isotopies. 
 \bd
 The  moduli space 
         ${\cal P}_{S, {\cal L}}$  parametrises data $(\Sigma, \{f_j\})$, where $\Sigma$ is a Riemann surface   of the topological type $S$, and 
         $\{f_j\in {\cal G}_p\}$ is a collection of distinct Puiseux polynomials\footnote{The space ${\cal G}_p$ is defined via the local ring  ${\cal O}_p$ of functions  at a point $p\in \Sigma$. It  does not depend on  a local parameter $z$.}  from Definition \ref{DEFIND} near the  punctures $p$ on   $\Sigma$, 
           with a given admissible isotopy class    of the corresponding Legendrian  ${\cal L}$. 
           \ed
           
      Forgetting   Puiseux polynomilas, we get   a   projection, where $S$ is a genus $g$ surface with $n$ punctures: 
      \be \la{MLP}
      {\cal P}_{S, {\cal L}} \lra {\cal M}_{g,n}.      \ee      

    Given a Legendrian link ${\cal L}$ of Stokes type on $S$, there are the following categories and stacks:\vskip 2mm
    
1. The category of admissible sheaves  ${\cal C}(S, {\cal L})$,  and the stack  
  of its objects ${\cal M}(S, {\cal L})$.     

2. The  DeRham stack  ${\cal M}_{\rm DR}(S, {\cal L})$ of holomorphic vector bundles with connection over a Riemann surface $\Sigma$ of topological type $S$, with the given 
homotopy class   of the Legendrian link ${\cal L}$  assigned to its formal equivalence class. It is fibered over the stack  ${\cal P}_{S, {\cal L}}$:
    \be
 {\cal M}_{\rm DR}(S, {\cal L}) \lra      {\cal P}_{S, {\cal L}}.      
      \ee                
           
    3.      The Betti stack    ${\cal M}_{\rm B}(S, {\cal L})$, which is a local system  of  stacks  ${\cal M}_{\rm B}(S, {\cal L})$ over  
                  ${\cal P}_{S, {\cal L}}$. 
                  
                  Its fiber over a point $(\Sigma, \{f_j\}) \in {\cal P}_{S, {\cal L}}$ is the stack ${\cal M}(S, {\cal L})$ for 
                   the Legendrian link ${\cal L}$ of $\{f_j\}$.

                   The fact that it is indeed a well defined local system of stacks, e.g. does not depend on the choice of radius $\varepsilon$ 
                   in the assignment $\{f_j\} \lms {\cal L}$, and locally constant over the base,  is guaranteed by  Theorem \ref{TRT}.  \\

          The wild Riemann-Hilbert correspondence  \cite{M2} is     a complex analytic equivalence of stacks over   ${\cal P}_{S, {\cal L}}$:
 \be \la{rs22xxys} 
 \begin{gathered}
     \xymatrix{
   {\cal M}_{\rm DR}(S, {\cal L}) \ar[dr]   \ar[rr]^{ {\cal R}{\cal H} } & & \ar[dl]{\cal M}_{\rm B}(S, {\cal L})  \\
         &    {\cal P}_{S, {\cal L}}     &   }
 \end{gathered}
  \ee                    
  The algebraic structure of  the  fibers of the De Rham stack over ${\cal P}_{S, {\cal L}}$ varies, while for  the Betti stack  
    it     does not.  
          The  Betti stack form a local system of stacks over ${\cal P}_{S, {\cal L}}$.   In particular,  we can define  the analog of  isomonodromic deformation for  irregular holomorphic  connections.     
          
    \paragraph{3. The wild mapping class group.} 
          Comparing  the classical Riemann-Hilbert correspondence (\ref{BeDRa}) with the one (\ref{rs22xxys}), we see that   the right analog    of       ${\cal M}_{g,n}$ is the moduli space ${\cal P}_{S, {\cal L}}$, and arrive at   
       
\bd 
 The {\em  wild mapping class group}   is the fundamental group  $\Gamma_{S, {\cal L}}:= \pi_1({\cal P}_{S, {\cal L}})$.
 \ed          
          
  The wild mapping class group  $\Gamma_{S, {\cal L}}$ acts by automorphisms of the  stack ${\cal M}(S, {\cal L})$.      
          
The map (\ref{MLP}) provides a  surjection $\Gamma_{S, {\cal L}}\lra \Gamma_S$.     The   group $\Gamma_{S, {\cal L}}$ is an extension of the mapping class group   $\Gamma_S$ by the product over the punctures $p$ of $S$ 
       of the local wild mapping class groups 
        ${\Gamma}_{{\rm D}^*, {\cal L}_p}$:
        \be
        1 \lra \prod_{p} {\Gamma}_{{\rm D}^*, {\cal L}_p} \lra      \Gamma_{S, {\cal L}} \lra \Gamma_S \lra 1.    \ee

\paragraph{Example.}  In  example (\ref{SSS}) 
in Section \ref{sec9.1},  the space    ${\cal P}_{{\rm D}^*, {\cal L}}$ for the punctured disc is given by  collections of distinct complex numbers $\{\lambda_j\}$. 
                  The   wild mapping class  group $\Gamma_{S, {\cal L}}$ is the braid group ${\rm Br}_n$.      \vskip 2mm

 Below we will discuss topological ramifications of this definition.

            \paragraph{4. Three flavors of   Legendrian mapping class groups.} There are  three flavors of the space of Legendrian links ${\cal L}\subset 
              {\rm ST}^*(S)$ and  Hamiltonian isotopies betwen them. 
            Therefore there are three flavors of  Legendrian mapping class group,  defined as the   fundamental group of the connected component of ${\cal L}$. \vskip 2mm

             1. The space of \underline{all} Legendrians  ${\cal L}\subset {\rm ST}^*({S})$ and Hamiltonian isotopies. 
              The {\it Legendrian mapping class group} ${\Gamma}_{{\rm leg}, S, {\cal L}}$ is the fundamental group of the connected component  containing   ${\cal L}$. 
              By \cite{GKS}, it acts by   autoequivelences of the triangulated/dg category of constructible sheaves on $S$.              \vskip 2mm 
                           
             2.  Legendrians given by  \underline{smooth} loops in $S$, and \underline{admissible} Hamiltonian isotopies between them.      
               The {\it admissible Legendrian mapping class group} $ {\Gamma}_{{\rm adm}, S, {\cal L}}$ is the fundamental group of the   component  of   ${\cal L}$.  
               By Theorems \ref{TRTB} - \ref{TRT},  it acts by autoequivelences of the   category ${\cal C}({S}, {\cal L})$ of admissible dg-sheaves.          \vskip 2mm            
                        
             
             3.   The \underline{ideal} Legendrians ${\cal L}$  which are, by the definition,  admissibly isotopic to a union of local  {ideal} Legendrians ${\cal L}_p$ near the punctures $p$ of $S$, modulo   admissible Hamiltonian isotopies.                     
             The {\it wild mapping class group} $ {\Gamma}_{S, {\cal L}}$ is the fundamental group of the connected component  of  ${\cal L}$.       
            It acts by preserving   cluster structures of the  moduli spaces of Stokes data, see Section \ref{Sect9.5}. \vskip 2mm

             These spaces are nested:  $(3) \subset (2) \subset (1)$. Therefore  there are     maps of the corresponding groups: 
             \be
 {\Gamma}_{S,  {\cal L}} \lra  {\Gamma}_{ {\rm adm}, S, {\cal L}}   \lra   {\Gamma}_{{\rm leg}, S, {\cal L}}.
    \ee

  \paragraph{5. A combinatorial description of the wild mapping class
   group.} Given a Weyl group $W$, consider the following  stratified space ${\cal C}_W$. Its strata ${\cal S}_\beta$ are given by    
           collections $\beta=\{w_1, ..., w_n\}$ of  the cyclically ordered elements  of $W$  assigned to 
           distinct points on a circle, so that the  cyclic order of the points on the circle   coincides  with the one of  $\beta$. So   ${\cal S}_\beta$ is homeomorphic to 
           $S^1\times \R^{n-1}$.            
           
           If $l(w_iw_{i+1})= l(w_i) + l(w_{i+1})$,  the stratum   
           ${\cal S}_{w_1, ..., w_iw_{i+1}, ... ,w_n}$        lies on the  boundary of         ${\cal S}_{w_1, ..., w_i, w_{i+1}, ... w_n}$.         
           Denote by ${\cal C}_\beta$ the connected component of the cell complex      ${\cal C}_W$  containing the strata ${\cal S}_{ w}$.    
           The space   ${\cal C}_\beta$ is determined by the  image  $[\beta]$ of the element $\beta$ in the cyclic envelope of the braid semigroup  ${\rm Br}^+_W$  assigned to   $W$, discussed in Section \ref{Sect9.6}.      
           
           Let $W=S_{N}$.   Then $[\beta]$     determines a cyclic braid. Denote the corresponding Legendrian by ${\cal L}_{[\beta]}$.
                 \bp
            The local wild  mapping class
   group  $\Gamma_{{\cal L}_{[\beta]}}$ is the fundamental group of    ${\cal C}_{\beta}$:
             $$
       \Gamma_{ {\cal L}_{[\beta]}}:=      \pi_1({\cal C}_{\beta}).
             $$      
             \ep

   \subsection{Non-commutative cluster structure of  spaces of Stokes data} \la{Sect9.5}

            By a non-commutative space ${\cal X}$ we mean a functor from a category of skew fields to the category of sets 
            which assigns to a skew field $R$ a set ${\cal X}(R)$. We say that ${\cal X}$ is defined over ${\rm Spec}(\Q)$ if we consider all skew fields $R$ of 
            characteristic zero. \vskip 2mm

            By a non-commutative stack ${\cal X}$ we just mean 
            the quotient of a non-commutative space ${\cal Y}$, equipped with an action of the group ${\rm GL}_N(R)$,  by the action of this group. So for any skew field $R$ we have a set ${\cal X}(R):= {\cal Y}(R)/{\rm GL}_N(R)$. \vskip 2mm

           We say, e.g. in Theorem \ref{TSD},  that a non-commutative space ${\cal X}$ admits a structure of a {\it non-commutative cluster Poisson variety} if it is birationally isomorphic to 
            a cluster Poisson variety associated to a bipartite ribbon graph $\Gamma$, as discussed in paragraph 3 of Section \ref{SECT1.1} and  in Section \ref{SEC5S}.
            
            This means that there exists a non-empty collection of bipartite ribbon graphs $\{\Gamma\}$ such that     any two of the graphs   are related by a sequence of two by two moves,  for each $\Gamma$ from the collection 
              there is a 
            birational isomorphism
            $$
           \alpha:  {\rm Loc}_1(\Gamma)\lra {\cal X}, 
            $$
            and for any pair of the graphs related by a two by two move $\mu: \Gamma \lra \Gamma'$ there is a commutative diagram
                      \be \la{NCXCV**} 
 \begin{gathered}
     \xymatrix{
   {\rm Loc}_1(\Gamma) \ar[d]^{\mu_*} \ar[r]^{\quad \alpha}    & {\cal X}    \ar[d]^{=}  \\
        {\rm Loc}_1(\Gamma')   \ar[r]^{ \quad \alpha'}  &  {\cal X}      }
 \end{gathered}
  \ee                
where the birational isomorphism $\mu_*$ was described in Section \ref{SEC5.1}. 
        
            Next, let ${\rm Mod}_{\cal X}$ be  a discrete group acting by birational 
            automorphisms 
            of ${\cal X}$. We say that ${\cal X}$ admits a structure of a {\it non-commutative  ${\rm Mod}_{\cal X}-$equivariant cluster Poisson variety}     if     
            for any    $g\in {\rm Mod}_{\cal X}$ there is a bipartite ribbon graph $\Gamma'$ together with a birational isomorphism 
              $$
           \alpha':  {\rm Loc}_1(\Gamma') \lra {\cal X}, 
            $$
            and a composition of two by two moves $\mu_g: \Gamma'\lra \Gamma'$,  which induces 
            a birational 
            non-commutative cluster Poisson transformation $\mu_{g*} : {\rm Loc}_1(\Gamma) \lra {\rm Loc}_1(\Gamma')$ providing a commutative diagram 
            \be \la{NCXCV} 
 \begin{gathered}
     \xymatrix{
   {\rm Loc}_1(\Gamma)(R)  \ar[d]^{\mu_{g*}} \ar[r]^{\qquad \alpha}    & {\cal X}    \ar[d]^{g}  \\
        {\rm Loc}_1(\Gamma')(R)   \ar[r]^{ \qquad \alpha'}  &  {\cal X}      }
 \end{gathered}
  \ee

             \paragraph{1. Non-commutative stacks of framed Stokes data.}    Just like in the regular case, to quantize the stack  of  Stokes data 
  ${\cal S}({\cal L}, d)$,  we   introduce  a larger stack  ${\cal X}({\cal L}, d)$ of {\it framed} Stokes data. Namely, 
  we add a framing, defined as a complete filtration of the local system ${\cal L}_\alpha$ on each loop $\alpha$ of the Legendrian ${\cal L}$, invariant under the monodromy around the loop. \vskip 2mm

\begin{figure}[ht]
\centerline{\epsfbox{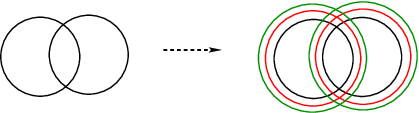}}
\caption{Inflating a Legendrian link.}
\label{stokes7}
\end{figure}

We interpret a  filtration on a $d_j-$dimensional local system on the loop $\alpha_j$ of the Legendrian   ${\cal L}$ as follows. 
We replace the loop $\alpha_j$ by a collection of $d_j$   loops $\alpha_j^{(1)}, \ldots , \alpha_j^{(d_j)}$ obtained by expanding the   loop $\alpha_j$   out of the puncture 
by $\varepsilon, 2 \varepsilon, ..., d_j\varepsilon$, see Figure \ref{stokes7}. 
We equip each loop $\alpha_j^{(p)}$ with the local system   given by ${\rm gr}^p$ of the   filtration on the 
 loop $\alpha_j$. So get a new link ${\cal L}_d$   supporting one dimensional local systems. 
  
 Set ${\bf 1}:= (1, ..., 1)$. 
  Consider the stack ${\cal M}_0({\cal L}_d, {\bf 1})$ of ${\cal L}_d-$admissible  dg-sheaves on   $S$,   vanishing near the punctures, whose microlocal support is described by  the union of the zero section and  arbitrary 1-dimensional local systems on 
  the  Legendrian ${\cal L}_d$.   The next   observation   generalizes Lemma \ref{LE1}.

\bl \la{NCSDI}
The non-commutative stack   ${\cal M}_0({\cal L}_d, {\bf 1})$  is  canonically equivalent to the one of framed Stokes data of     type $({\cal L}, d)$:
\be
{\cal M}_0({\cal L}_d, {\bf 1}) \stackrel{}{=} {\cal X}({\cal L}, d).
\ee
\el

 There is a  projection  given   by  assigning  the microlocal support  local system  on the   Legendrian ${\cal L}_d$:
   \be \la{111}
 \mu:  {\cal M}_0({\cal L}_d, {\bf 1}) \lra {\rm Loc}_1(S^1)^{|d|}, \ \ \ \ \ \ |d|:= \sum_\alpha d_\alpha.  
  \ee

The product of  the symmetric groups $\prod_{\alpha}S_{d_{\alpha}} $ acts by birational automorphisms of   ${\cal X}({\cal L}, d)$.  
To define it, note that  given a generic  $p-$dimensional $R-$local system on a circle, there are $p!$ complete filtrations of the local system. Indeed, a generic linear automorphism $L$ of an $p-$dimensional vector space $V_p$ is decomposed uniquely into a direct sum of  
 $L-$invariant one-dimensional subspaces. The filtrations are parametrised by orderings of these subspaces. Therefore the symmetric group $S_{p}$ acts on the fltrations. 
 
 This just means that the canonical projection, obtained by forgetting the filtrations
 \be
 p: {\cal X}({\cal L}, d) \lra {\rm Loc}({\cal L}, d)
  \ee
 is a Galois cover at the generic point with the Galois group $\prod_{\alpha}S_{d_{\alpha}} $.

   \bt \la{TSD} For any decorated surface $\bS$,   and   any topological   data $({\cal L}, d)$ where ${\cal L}$ is an  ideal Legendrian on $\bS$, the 
stack $ {\cal X}({\cal L}, d)$ of framed non-commutative Stokes data   has a non-commutative cluster Poisson structure equivariant under the 
           wild cluster mapping class group $ \Gamma_{{\cal L}, d}$.  

In the commutative case, the subalgebra $\mu^*{\cal O}({\Bbb G}_m^{|d|})$, see (\ref{111}), is   the Poisson center of    ${\cal O}( {\cal X}({\cal L}, d))$. 
  \et
  
  We prove Theorem \ref{TSD} assuming that in the case when $\bS= S^2-\{\infty\}$ the Legendrian ${\cal L}$ is not very degenerate, see Theorem \ref{MTHSSCL1} for the precise statement. 
  
  \paragraph{2. Idea of the  proof of Theorem \ref{TSD}.}   
    Lemma \ref{NCSDI} tells that any stack   of framed Stokes data has a canonical realization as a stack of admissible sheaves 
with  the following two features:
\vskip 2mm

i) The admissible sheaves  vanish  near the punctures. 

ii) The microlocal  support  local systems     on the    Legendrian   ${\cal L}_d$ are \underline{one-dimensional}. 
\vskip 2mm


  We show that the Legendrian link ${\cal L}_d$ has an admissible deformation to  (a slight modification  of) the web ${\cal Z}_\Gamma$ 
  of zig-zag strands for a certain bipartite graph $\Gamma$ on $S$. 
 By the invariance of the stack of admissible dg-shaves under admissible deformations,   
 ${\cal M}({\cal L}_d, {\bf 1}) = {\cal M}({\cal Z}_\Gamma, {\bf 1})$. 
   By Lemma \ref{LE2}, the stack ${\rm Loc}_1(\Gamma)$ of 1-dimensional local systems on  $\Gamma$  is realised as the open substack 
       defined by the   condition of vanishing on mixed domains: 
 $
      {\rm Loc}_1(\Gamma) \subset  {\cal M}({\cal Z}_\Gamma, {\bf 1}).       
 $
   Since vanishing on mixed domains implies vanishing near the punctures, we   get an open embedding
      $$
      {\rm Loc}_1(\Gamma) \subset     {\cal M}_0({\cal L}_d, {\bf 1}) \stackrel{ }{= } {\cal X}({\cal L}, d).
        $$
 The last equivalence is given by Lemma \ref{LE2}. This way we get an open cluster Poisson 
 torus assigned to the graph $\Gamma$. The tori assigned to different graphs $\Gamma$ provide the cluster atlas. 
 The elements of the cluster wild mapping class group transform such a   graph $\Gamma$  to another  one $\Gamma'$,  related to $\Gamma$ by   
  two by two moves. 
 Therefore we get a cluster Poisson atlas equivariant under the wild mapping class group.  \vskip 2mm

          Before we start the proof, let us recall the following. 
          
  \paragraph{3. Bipartite graphs and triple crossing diagrams.}    
Let us recall    triple crossing diagrams \cite{Th}.

    \begin{figure}[ht]
\centerline{\epsfbox{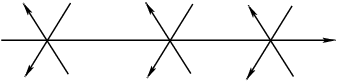}}
\caption{A triple crossing diagram.}
\label{stokes9}
\end{figure}

 \bd  A {\it   triple crossing diagram} is a collection  ${\cal I}$ of smooth oriented strands on an oriented decorated surface $\bS$ such that:

1.    Each intersection point of ${\cal I}$ is a triple   point.
       
2.  Strand orientations   are consistent with an orientation of   each  component  of $\bS-{\cal I}$.

3. Each strand is either a loop, or ends on the boundary of $\bS$. Endpoinds of strands are disjunct.
   
              \ed

\begin{figure}[ht]
\centerline{\epsfbox{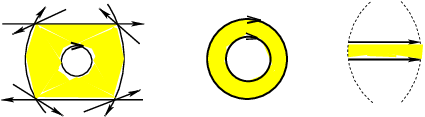}}
\caption{Strand orientations   are not consistent with any orientations of   internal (yellow)   components.}
\label{stokes16}
\end{figure}

See counterexamples to condition 2)   on  Figure \ref{stokes16}). Note that condition  2) implies   the following:

$\bullet$ For each intersection point,   passing strands are oriented as   on Figure \ref{stokes9}. 

$\bullet$   Going along a strand of ${\cal I}$, the orientations of the intersecting it strands   alternate, as   on Figure \ref{stokes9}.

Any bipartite graph can be transformed to a bipartite graph with 3-valent $\bullet-$vertices obtained by splitting the $\bullet-$vertices of valency $>3$ by   adding  edges, and 
 introducing 2-valent    $\circ-$vertices in the middle of these edges. 
A bipartite ribbon graph $\Gamma$ with 3-valent $\bullet-$vertices gives rise to a triple crossing diagram on the associated surface $S_\Gamma$, provided by its zig-zag strands, 
 \cite[Section 2.6]{GK}. Namely, we first push zig-zag strands out of the graph, and then shrink the triangles near the  $\bullet-$vertices to points, see Figure \ref{ncls104}. 
 Then all original double crossings collide into triple crossings corresponding to the $\bullet-$vertices.

This construction can be inverted. Any  triple crossing diagram  ${\cal I}$ on $\bS$ has a canonical up to isotopy resolution defined as follows \cite[Section 2.3]{GK}. 
   For each crossing, move slightly one of the strands, inflating the crossing into a little triangle, whose strands are oriented  consistently with  the orientation 
  of $\bS$ -  counterclockwise -  getting a web ${\cal I}'$.  There are three types of components of $\bS-{\cal I}'$, see Figure \ref{ncls104}: 
  
  \vskip 1mm
  i)  {\it Black} components, whose boundary strands are oriented clockwise. 
  
  ii)  {\it White} components, whose boundary strands are oriented counterclockwise. 
  
  iii)     {\it Mixed} components, whose boundary strands   orientations alternate.   
  
  \vskip 1mm
  
  We shrink all black (respectively white) components to black (respectively white) vertices, and connect two vertices by an edge if they are separated by a crossing of the web ${\cal I}'$. 
  We get a bipartite graph $\Gamma_{\cal I}$.

 \paragraph{4. Ideal Legendrians $\to$  triple crossing diagrams $\to$ non-commutative cluster Poisson varieties.} 
 Take an ideal Legendrian ${\cal L}$ on a punctured surface $S$,  whose local components ${\cal L}_p$   are  ideal Legendrians 
   near  punctures $p$. So   each    ${\cal L}_p$ 
  satisfies conditions of Definition \ref{DEFSDI},   different ${\cal L}_p$  are disjunct, and ${\cal L} = \cup_p{\cal L}_p$. The number of crossings  of a generic ray $r$ from  $p$ with ${\cal L}_p$ is the {\it local rank} of 
    ${\cal L}_p$ at $p$. An ideal Legendrian ${\cal L}$ is of {\it rank} $m$ if for each local     Legendrian ${\cal L}_p$ its local {rank} at $p$ is equal to $m$.   
    
           \bt \la{MTHSSCL1} Any rank $m$ ideal Legendrian ${\cal L}$ on a punctured surface $S$ different from $S^2-\{\infty\}$, 
           or  a rank $m$ ideal Legendrian   on   $S^2-\{\infty\}$  corresponding to a cyclic closure of the braid,  $w_0xw_0y$, for any $x,y \in {\rm Br}_m$, 
                    can be admissibly isotoped  to a triple crossing diagram with the following  property: \vskip 2mm
           
           The black and white components are either contractible, or 
           homeomorphic to an annulus. 
           
           The union of all annuli is a disjoint union of punctured discs around  (some of the) punctures. 
           
           \et
           
              Before we proceed to the proof of Theorem \ref{MTHSSCL1}, let us discuss its main implications. \vskip 2mm

               The triple crossing diagrams ${\cal L}$ which appear in our story have   black and white components which are either contractible, or 
            annuli. Furthermore, let us assume that the union of all annuli is given by a union of discs around the punctures. Assuming these conditions, 
            consider the category ${\cal C}(S; {\cal L}, {\bf 1})$ of ${\cal L}-$admissible dg-sheaves 
           such that the local systems of the microlocal support  supported on  the components of ${\cal L}$ are 1-dimensional. Denote by ${\cal M}_0(S; {\cal L}, {\bf 1})$ the substack   
           of the   stack of objects ${\cal M}(S; {\cal L}, {\bf 1})$  of the category ${\cal C}(S; {\cal L}, {\bf 1})$           given by   
            ${\cal L}-$admissible dg-sheaves which are zero on the mixed components.   
         
           \bl The stack ${\cal M}_0(S; {\cal L}, {\bf 1})$ is a torus.       
           \el

       \begin{proof}     According to our assumption, the annuli surrounding a given puncture 
       \end{proof}

           \begin{corollary} Any rank $m$ ideal Legendrian ${\cal L}$ on $S$ gives rise to a Zariski open torus            ${\cal M}_0(S; {\cal L}, {\bf 1})$     in the non-commutative 
           moduli stack    ${\cal M}(S; {\cal L}, {\bf 1})$ 
           of   ${\cal L}-$admissible dg-sheaves with one-dimensional local systems on the components of ${\cal L}$. \ec

We will   proceed to the proof of Theorem \ref{MTHSSCL1}, after discussing   two   examples.

   \paragraph{5. The simplest example.} The stack of   $m-$dimensional local systems on $S^2-\{0, \infty\}$ with complete filtrations near the punctures  
  is just the stack of admissible sheaves    vanishing near the punctures      for the oriented  web ${\cal L}_m$, illustrated for $m=3$ on the left of Figure \ref{stokes17}. 
   An admissible deformation of ${\cal L}_m$   shown on the right describes the following constructible sheaves. Denote by ${\cal A}_1, {\cal A}_2, ..., {\cal A}_{2m+1}$ the open  
   components of $S^2 - {\cal L}_m$, counted from    $0$. Then    the sheaves are zero on  ${\cal A}_{2k+1}$, 
   and are 1-dimensional local systems on  
     $\overline {\cal A}_{2k}$. Their moduli space   is   $(\Bbb G_m)^m$,  parametrising  the monodromies. 
        
                 \begin{figure}[ht]
\centerline{\epsfbox{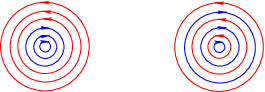}}
\caption{Left: the link ${\cal L}_3$ describing the  stack of rank 3  local systems on $\C^*$ with complete filtrations near the punctures. Right: admissible deformed link providing an open locus 
identified with ${\Bbb G}_m^3$.}
\label{stokes17}
\end{figure}    
  
\paragraph{6. Another simple example.}   Let $S^2-\{0, \infty\}$ be a sphere  a regular puncture at $\infty$, and an irregular one at $0$. 
     The    Legendrian ${\cal L}_0$   goes  clockwise around   $0$. 
           It bounds a  domain ${\cal D}_{{\cal L}_0}$.  Precisely, take a point $x \in {\rm D}^*$.        
 Let  $r_x$ be  the ray  from   $p$ to $x$. 
           Then  $x \in {\cal D}_{{\cal L}_0}$ if and only if it is in the convex hull of       $\{r_x \cap   {\cal L}_0\}$.

        \begin{figure}[ht]
\centerline{\epsfbox{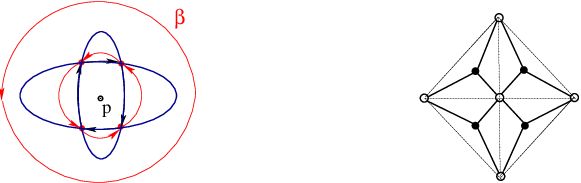}}
\caption{Left: a triple crossing diagram from two intersecting clockwise oriented loops around the puncture. Right: the corresponding bipartite graph.}
\label{stokes8}
\end{figure} 

Now we follow \cite{FG3}. Call   the selfintersection points of $ {\cal L}_0$ {\it vertices}. 
       Take a   component  ${\cal V}$ of   ${\cal D}_{{\cal L}_0} -  {\cal L}_0$.       
       Connect the  rightmost vertex  on the boundary of ${\cal V}$   with the   leftmost   one  by a counterclockwise path $\gamma_{\cal V}\subset {\cal V}$,  transversal to ${\cal L}_0$.       
       Take the union of  these  paths, a counterclockwise arc $\beta$ near $\infty$, and   $  {\cal L}_0$, see   Figure \ref{stokes8}: 
       \be
       {\cal I}:= \cup_{\cal V} \gamma_{\cal V} \cup \beta \cup {\cal L}_0.
       \ee   
       It is a triple point diagram. Its  number of   counterclockwise arcs is the same as the local rank of  ${\cal L}_0$ at $0$.

  \paragraph{7. Proof of Theorem \ref{MTHSSCL1}.}  
Take  a  triangulation ${\cal T}$ of $S$   with   vertices at the   punctures. Inflate each side  of   ${\cal T}$ to a pair of thin {\it bigons}, see Figure \ref{stokes13}. 
 Then $S$ -   $\{$the bigons$\}$  =  a union of  {\it internal triangles}. 
 Performing   admissible isotopies,  we can assume that  ${\cal L}$ crosses internal triangles $t$ and bigons $b$ as follows.

1. For each  vertex $v$ of   $t$:  the ${\cal L}$ induces    $m$    non-intersecting,    clockwise  (red) arcs  around  $v$. 
 
 2. For each bigon $b$: near one of its  angles, called the {\it red angle},   ${\cal L}$ induces $m$   non-intersecting   clockwise arcs near its vertex.
  We put no conditions  on the restriction of ${\cal L}$ to the other one,   called the {\it blue angle}.
 
 3. For each vertex $v$ of each   edge    of ${\cal T}$:   just one of the  bigons at this edge has a red angle at $v$.

      \begin{figure}[ht]
\centerline{\epsfbox{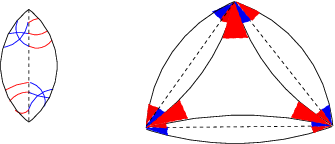}}
\caption{Inflating an edge to a pair of   bigons.  Inflating edges of a triangle to pairs of colored bigons. }
\label{stokes13}
\end{figure} 


Let $b$ be a  bigon. Denote by  ${\cal L}_b$  the Legendrian in the bigon. 
 
\bl \la{9.17} One can    admissibly isotope the   Legendrian ${\cal L}_b$  to a   Legendrian 
${\cal L}_b'$ so that 

   \begin{enumerate}
   
   \item   On each   side  of the bigon $b$, the  endpoints of the red and blue strands of ${\cal L}_b'$ alternate. 
       
  \item  The strands of ${\cal L}_b'$ form a triple crossing diagram.

  \end{enumerate} 
  \el
  
          \begin{figure}[ht]
\centerline{\epsfbox{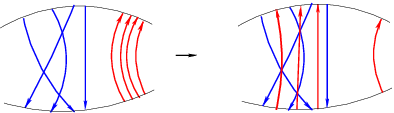}}
\caption{From concentric red and blue arcs near the vertices of a bigon, to a triple crossing diagram.}
\label{stokes14}
\end{figure} 

  \begin{proof} Move the red arcs towards the blue angle of $b$ so that condition 1) holds. 
 Then there is a unique up to isotopy way to deform them to a triple crossing diagram, with the minimal number of crossings. 
 
 Namely, forget the red arcs. Then for each component of the complement to the blue arcs pick the leftmost vertex/boundary point, and the similar rightmost one, and 
 connect them by the unique up to isotopy red arc inside the component, see Figure \ref{stokes14}. We get a triple crossing diagram. 
  \end{proof}

Lemma \ref{9.17}  provides   triple crossing diagrams in   the bigons. 
For each   internal triangle $t$, see   Figure \ref{stokes13},
 there is a standard triple crossing diagram provided by the $m-$trianulation of the triangle, which induces 
on the boundary   the    given alternating sequence of   intersection points of ${\cal L}'_b$ provided by Lemma \ref{9.17}. Indeed, take the $m-$trianulation of a triangle, see 
Figure \ref{stokes18}, and move each of its  strands  towards the parallel to it side of the triangle, so that the centers of looking down small triangles become the triple crossing points. 
Then the endpoints   alternate, and we can match them with the ones  on the bigons ${\cal L}_b'$. We denote this triple crossing diagram by ${\cal L}_t'$. 
Take the union of these triple point diagrams:
\be
{\cal L}_{\cal T}:= \cup_t {\cal L}_t' \bigcup \cup_b{\cal L}'_{b}.
\ee
Theorem \ref{MTHSSCL1} is proved. 
           \begin{figure}[ht]
\centerline{\epsfbox{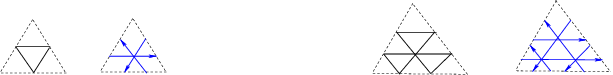}}
\caption{A triple crossing diagram assigned to an $m$-triangulation of the triangle $t$, for $m=1,2$.}
\label{stokes18}
\end{figure}

   \subsection{Cluster structure of stacks of $\G-$Stokes data on surfaces} \la{Sect9.6}

\paragraph{1. The  cyclic envelope ${\rm Br}_{\mathfrak{g}}^+$ of the braid semigroup.} Let $\G$ be a split semi-simple algebraic group  over $\Q$.
Denote by  ${\rm I}$    the set of  positive simple roots,  and by $C_{ij}$   the Cartan matrix of the root system of the Lie algebra $\mathfrak{g}$ of $\G$. 
 Let ${\rm Br}_{\mathfrak{g}} $ (respectively ${\rm Br}_{\mathfrak{g}}^+$) 
be the braid group (respectively semigroup) of  the root system of  $\mathfrak{g}$.  It is generated by the elements $s_i, i\in {\rm I}$, subject to the relations 
\begin{equation} \label{1}
\begin{array}{lcl}
s_is_j = s_js_i &\mbox{if}& C_{i j}=C_{ji}=0,\\
s_is_js_i = s_js_is_j &\mbox{if}& C_{i j}=C_{ji}=-1,\\
s_is_js_is_j = s_js_is_js_i &\mbox{if}& C_{ij}=2C_{ji}=-2,\\
s_is_js_is_js_is_j = s_js_is_js_is_js_i &\mbox{if}& C_{ij}=3C_{ji}=-3.
\end{array}
\end{equation}
 The 
Weyl group $W$ is the quotient of the group ${\rm Br}_{\mathfrak{g}}$ by the relations $s_i^2=1$. Denote by  $l(*)$   the length function on $W$. There is a   set theoretic section $\mu: W \to {\rm Br}_{\mathfrak{g}}^+$, 
such that $\mu(s_i)=s_i$ and 
$ 
\mu(ss') = \mu(s)\mu(s') \quad \mbox{if $l(ss') = l(s)+l(s')$}.
$ 
Abusing notation, we denote by $s_i$   the elements $\mu(s_i) \in {\rm Br}_{\mathfrak{g}}^+$.  

Denote by $[{\rm Br}_{\mathfrak{g}}^+]$ the set of coinvariants of the cyclic shift $s_{j_1} \ldots s_{j_{n-1}}  s_{j_{n}}  \lms  s_{j_{n}} s_{j_1} \ldots s_{ j_{n-1}}$  on  ${\rm Br}_{\mathfrak{g}}^+$,  called   the {\it cyclic envelope} of  the braid semigroup.  

 We denote by   $[s_{j_1} \ldots s_{j_n}]\in [{\rm Br}_{\mathfrak{g}}^+]$  the cyclic element  corresponding to   $s_{j_1} \ldots s_{j_n}\in {\rm Br}_{\mathfrak{g}}^+$. \\

Given a pair of flags $\B, \B'$, we define  $w(\B, \B') \in W$  via the isomorphism  $\G\backslash ({\cal B} \times {\cal B}) = W$. 

Take a  defining relation  $s_{l_1}s_{l_2} ... s_{l_m} = s_{r_1}s_{r_2} ... s_{r_m}$ in (\ref{1}).   
For any    flags $(\B_0, ..., \B_m)$ with   $w(\B_{a-1}, \B_{a}) =s_{j_a}$, where  $a=1, ..., m$, 
 there is a unique collection of flags $(\B'_0, ..., \B'_m)$   with   $w(\B'_{a-1}, \B'_{a}) = s_{r_a}$.   
 
 Given a reduced decomposition $s_{j_1}s_{j_2} ... s_{j_k}$ of a cyclic element $\beta \in [{\rm Br}_{\mathfrak{g}}^+]$ there is the moduli space of 
 $\G-$orbits of flags $(\B_1, ..., \B_k)$ such that $w(\B_{i-1}, \B_i)= s_{j_i}$, where $i \in \Z/k\Z$. The moduli spaces assigned to different reduced decompositions of   $\beta$ are canonically isomorphic.

      \paragraph{ The  moduli space ${\cal X}({\G}, \beta; S)$.} 
     Let  ${\cal B}$ be the flag variety of $\G$. A $\G-$local system ${\cal L}$ on $S$ gives rise to  
      the associated  local system of flags ${\cal L}_{\cal B}:= {\cal L}\times_\G {\cal B}$.   
    Denote by  $S^1_p$ the circle of rays at $T_pS$.  Given a reduced decomposition   $s_{j_1}  \ldots   s_{j_k}$ of the cyclic word ${\beta}_p \in [{\rm Br}^+_{\mathfrak{g}}]$,  take  distinct points 
   $x_1 , ... , x_k \in S^1_p$ going cyclically around  the circle, and 
   put the elements  $s_{j_i}$ at the points $x_i $. 
   
             \begin{figure}[ht]
\centerline{\epsfbox{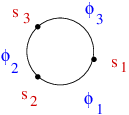}}
\caption{Cyclic framing of type $[s_1s_2s_3]$.}
\label{stokes19}
\end{figure} 
     
  \bd \la{11}
      Let     $\beta:= \{{\beta}_p\}$ be a collection of  elements of   $[{\rm Br}^+_{\mathfrak{g}}]$ assigned to the punctures $\{p\}$ of a  surface $S$.       
      A {\it   ${\beta}_p-$framing at a puncture $p$} on a $\G$-local system ${\cal L}$ on $S$ is  given by flat sections 
  $\varphi = (\varphi_1, \ldots , \varphi_k)$     of the   local system of flags ${\cal L}_{\cal B}$
   over  intervals $x_1x_2$, ..., $x_kx_1$ on $S^1_p$ so that  $w(\varphi_i, \varphi_{i+1})=s_{\alpha_i}$:  
      $$
   \varphi_1 \stackrel{s_{j_1}}{\lra}  \varphi_2  \stackrel{s_{j_2}}{\lra}   \ldots  \varphi_k  \stackrel{s_{j_k}}{\lra}   \varphi_1.
   $$     
   A $\beta-$framing on  ${\cal L}$   is  a collection of   local framings $\{\beta_p\}$ at the punctures of $S$.   
   
   The moduli space ${\cal X}({\G}, \beta; S)$ parametrises $\beta-$framed $\G$-local systems ${\cal L}$ on $S$. \ed
   
The stack ${\cal X}({\G}, \beta; S)$ depends   on the elements $\beta_p\in [{\rm Br}^+_{\mathfrak{g}}]$, rather than their reduced decompositions.

 \paragraph{Examples of stacks ${\cal X}({\G}, \beta; S)$.} 
 
 \begin{enumerate}
 
 \item If $\G = {\rm GL}_m$, a local \underline{framed} 
 Stokes data 
 is just a $\beta-$framing at the puncture on an $m-$dimensional  local system of vector spaces on the punctured disc 
 for an element $\beta \in [{\rm Br}_{\mathfrak{gl_m}}^+]$. 
 
\item  If $\beta_p=e$ is the unit, then a local $\beta_p-$framing at the puncture $p$ is just a monodromy invariant flag near the puncture. 
 
 \item   If all  ${\beta}_p=e$ are the unit elements, then  ${\cal X}({\G}, \beta; S)$  is  the stack   ${\cal X}({\G}, S)$   \cite{FG1} 
  parametrising       $\G$-local system ${\cal L}$ on $S$ 
   with  a flat section of the local system of flags ${\cal L}_{\cal B}$ near each puncture of $S$.  
   
  \item    Let $\bS$ be a decorated surface  with 
   $k$ boundary components, with $b_1, ..., b_k$ special points  on   them. 
   
    If ${\beta}_p = [\mu(w_0) \cdot \ldots \cdot \mu(w_0)] \in [{\rm Br}^+_{\mathfrak{g}}]$ is the cyclic product  
    of $b_p$ elements $\mu(w_0)$ assigned to the longest element  $w_0\in W$, then   ${\cal X}({\G}, \beta; S)$    is the stack ${\cal X}_{\G, \bS}$ for  from \cite{FG1}, as we show in Section \ref{Sec10.1}.
   
 The cluster Poisson structure on moduli spaces 1) and 2) was defined  for $\G = {\rm PGL}_m$ in \cite{FG1}, 
for   classical groups and ${\rm G}_2$ by case-by case study in \cite{Le1}-\cite{Le2}, and in the   general case in \cite{GS19}.  

\item      If $S =  S^2 - \{0, \infty \}$ is a punctured disc, and ${\beta}_0=e$,  ${\beta}_\infty = \beta$,  then the stack ${\cal X}({\G}, \{e, \beta\}; {\rm D}^*)$   carries a
 cluster Poisson variety structure 
 assigned to the  cyclic braid semigroup element  
   $\beta$  in \cite{FG4}. 

 \item     If $S = S^2 - \{0, \infty \}$, then the stack  ${\cal X}({\G}, \beta_0, \beta_p; {\rm D}^*)$   contains as a Zariski open part  the double Bott-Samelson variety \cite{SW}  assigned to a pair 
   cyclic braid semigroup  elements 
   $\beta_0, \beta_\infty \in [{\rm Br}^+_{\mathfrak{g}}]$. 
   \end{enumerate}
   
         \bt \la{MTHSSCL222}  
             Any stack of framed $\G$-Stokes data on a  surface $S$ is the  stack   ${\cal X}({\G}, \beta; {S})$    for some $\beta$.
              \et
            \begin{proof}         This can be deduced from the formal classification in the $\G$-case just as  in Section \ref{sec9.1}.    \end{proof}

           \bt \la{MTHSSCL111} Let $S$ be a punctured surface, and $S \not = S^2-\{\infty\}$. Then            
           the stack ${\cal X}({\G}, \beta; {S})$     has a   cluster Poisson structure, equivariant under the wild mapping clsass group. 
             Therefore the stack  ${\cal X}({\G}, \beta; {S})$  admits  a cluster quantization, equivariant under the wild mapping class group.       
             
             For each puncture $p$, the braid group ${\rm Br}_{\mathfrak{g}} $ acts by automorphisms  of the stack  ${\cal X}({\G}, \beta; {S})$, preserving the cluster Poisson structure.     \et
           
           \begin{proof} Assume first that $S\not = S^2-\{0, \infty\}$. Then there is an ideal triangulation ${\cal T}$ of $S$. 
           
            Take an ideal triangulation ${\cal T}$ of $S$, and inflate each of its edges $E$ to a pair of thin bigons $b_E$, as we did in Section \ref{Sect9.5}.  
           For each puncture $p$, put all the points $x_1, ..., x_{b_p}$ inside of the arc on the circle $S^1_p$ located inside of the left bigon at one of the edges sharing $p$. 
           Then the moduli space ${\cal X}({\G}, \beta; {S})$ is birationally equivalent to the result of amalgamation of the moduli spaces 
           ${\cal P}_{\G, t}$  from \cite{GS19} assigned to the triangles $t$ of the triangulation, and the moduli spaces     ${\cal P}_{\G, r_E}$     assigned to the rectangles corresponding to the bigons, that is the edges of ${\cal T}$. The independence of the choice of  ${\cal T}$ follows from  \cite{FG4} $\&$ \cite{GS19}. 
   \vskip 2mm
           
         Now let $S  = S^2-\{0, \infty\}$ be a cylinder. Recall the   decorated flag variety ${\cal A}:=\G/\rm U$. Recall the elementary moduli space ${\cal P}_{\G, s}$ paramterising triples   $(\B_1, \A_2, \A_3)\in \G\backslash ({\cal B} \times {\cal A} \times {\cal A})$,     such that the pairs $(\B_1, \A_2)$ and 
         $(\B_1, \A_2)$ are generic, and   
         $w(\B_2, \B_3)=s$  for the pair of flags $(\B_2, \B_3)$ underlying the  decorated flags $(\A_2, \A_3)$. We picture such a triple by a wedge with the flag $\B_1$ at the  vertex, and two  arrows pointing out to the decorated flags at the   vertices on the short edge. 
         The space ${\cal P}_{\G, s}$ carries a   cluster Poisson structure \cite[Section 6]{GS19}. The same space with     the opposite  Poisson structure is denoted by ${\cal P}_{\G, \overline s}$.                 
         
         Take   reduced decompositions for the braid semigroup elements $\beta_0 = s_{i_1} ... s_{i_a}$ and $\beta_\infty 
         = \overline s_{j_1} ... \overline s_{j_b}$. Amalgamate the corresponding moduli spaces   ${\cal P}_{\G, s_i}$ and   ${\cal P}_{\G, \overline s_j}$   so that the order 
         of the spaces    ${\cal P}_{\G, s_i}$ (respectively ${\cal P}_{\G, \overline s_j}$) follows a reduced decomposition $\beta_0$ (respectively $\beta_\infty$),   see Figure \ref{stokes21}.      Then glue 
       the left and the right pairs $(\A, \B)$ for the resulting moduli spaces if they oriented the same way. Otherwise alter one of them by   $({\B_1}, {\A_2}) \lra \overline w_0({\B_1}, {\A_2})$ where $\overline w_0$ is the canonical lift of $w_0$ to $\G$, and then glue.

 This way we get  a Zariski open part of ${\cal X}({\G}, \beta; {S})$. 
          Therefore the result follows from \cite{GS19}.         \vskip 2mm
          
              \begin{figure}[t]
\centerline{\epsfbox{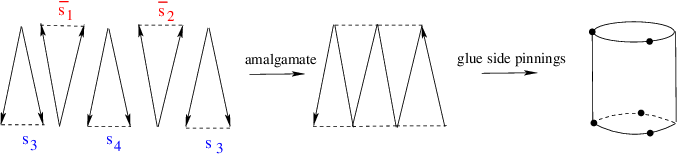}}
\caption{Getting a cluster chart with $\beta_0=s_1s_2$, $\beta_\infty=s_3s_4s_3$ on a cylinder by the amalgamation.}
\label{stokes21}
\end{figure} 
        The second claim  follows from this by applying the cluster Poisson quantization machine \cite{FG4}.       \\  
        
        Let us define the action $a_p$ of the braid group ${\rm Br}_{\mathfrak{g}} $ on the stack  ${\cal X}({\G}, \beta; {S})$.    
        Pick a generator $s_i$ of the  the braid group ${\rm Br}_{\mathfrak{g}}$. Take a cyclic word ${\bf b}= [s_{j_1} \cdot \ldots \cdot s_{j_n}]$  assigned to the cyclic collection of points $x_1, ..., x_n$ on the circle $S^1_p= T_pS-\{0\}/\R_{>0}$. We write this schematically as $[s_{j_1}|x_1,  \cdot \ldots \cdot s_{j_n}|x_n]$. Then we find its representative of $[{\bf b}]$ starting from $s_1$, and shift  the cyclic word by one step to the right:
        $$
        [s_{1}|x_1,  s_{j_2}|x_2, \cdot \ldots \cdot s_{j_n}|x_n]      \lms    [s_{1}|x_2,  s_{j_2}|x_3, \cdot \ldots \cdot s_{j_n}|x_1]. 
          $$
        which starts from $s_i$. 
        \end{proof}

        \paragraph{7. Cluster quantization of  ${\cal M}_{\rm DR}(\G, \beta; S)$.}  To quantize the  De Rham stack ${\cal M}_{\rm DR}(\G, \beta; S)(\C)$, 
  we apply  the $\G-$analog of the universal Riemann-Hilbert equivalence (\ref{rs22xxys}), and then 
  develope  the $ {\Gamma}_{\G, \beta; S}-$equivariant quantization of the   Betti  
    stack $  {\cal M}_{\rm B}(\G, \beta; S)$.      
   
    To do this,  we   quantize first the bigger stack ${\cal X}({\G}, \beta; S)$ of framed Stokes data. 
  There is a finite group  $W_{\G, S, {\beta}}$  which 
     acts birationally on ${\cal X}({\G}, \beta; S)$ by altering   filtrations.  Forgetting the filtrations 
   we get  a  map, which is   a finite Galois cover at the generic point:
$$
  p: {\cal X}({\G},  \beta; S)  \lra   {\cal M}_{\rm B}([S, {\cal L}]).
  $$

  The space ${\cal X}({\G}, \beta; S)$ has a ${\Gamma}_{\G, \beta; S}\times W_{\G, \beta; S}-$equivariant Poisson structure. 
   Therefore   \cite{FG4}   gives  its   
  ${\Gamma}_{\G, \beta; S}\times W_{\G, \beta; S}-$equivariant cluster Poisson quantization,  given by:\vskip 2mm
  
  1. A $q-$deformed algebra of functions ${\cal O}_q({\cal X}({\G},  \beta; S))$.
  
  2. A unitary projective representation of the group ${\Gamma}_{{\G, \beta; S}}\ \times W_{\G, \beta; S}$ in a Hilbert space ${\cal H}$.
      
  3. A representation   of the $\ast-$algebra given by the modular double  
 \be\la{HHH}
   {\cal O}_q({\cal X}({\G},  \beta; S)) \otimes {\cal O}_{q^\vee}({\cal X}({\G},  \beta; S)), \ \ \ \ q=e^{i\pi \hbar}, \ \ q^\vee = e^{i \pi / \hbar}. 
 \ee
  by unbounded operators in the Hilbert space ${\cal H}$, 
  intertwining the   action of the group    $\Gamma_{\G, \beta; S} \times W_{\G, \beta; S}$ on the $\ast-$algebra (\ref{HHH}) with its  projective unitary action  in ${\cal H}$ from 2). 
  Here $\hbar \in \R_{>0}$ or $|\hbar|=1$. \vskip 2mm

The   group $W_{\G, \beta; S}$ acts by cluster Poisson transformations (the regular case is done in \cite{GS19}). 
 Set 
  $$
  {\cal O}_q( {\cal M}_{\rm B}({\G},  \beta; S)):= {\cal O}_q({\cal X}({\G},  \beta; S))^{W_{\G, \beta; S}}. 
     $$
 We consider the   representation of the modular double of this subalgebra in the Hilbert space ${\cal H}$. It  commutes with the unitary action of the group $W_{\G, \beta; S}$ in ${\cal H}$.

           \subsection{More examples}     \la{Sec10.1}
                                    
                         \paragraph{1. Non-commutative stacks ${\cal X}_m(\bS)$ and ${\cal A}_m(\bS)$ as stacks of   Stokes data}  
                        
                          Below $\bS$ is a decorated surface with boundary. We  assume that the local systems on the strands of the webs, 
                           provided by the  microlocal supports of admissible sheaves, are   one dimensional.  \vskip 2mm

                  We assign to each boundary component $\pi$ of $\bS$, with $b_\pi$ special   points, an ideal Legendrian link ${\cal L}_{\pi,m}$  
                   in a little cylinder $C_\pi\subset \bS$ containing 
                   $\pi$,  obtained by   concatenating $b_\pi$ copies 
                   of the   Legendrian   with $m$ strands  on a rectangle  realising the longest permutation $w_0: (1, 2, ..., m) \lms (m, ..., 2, 1)$, see Figure \ref{stokes4}. 
                   
                   Let us assign to each puncture $p$ on $\bS$ a collection ${\cal L}_{p,m}$ of $m$ small disjoint loops   around the puncture.

                   Denote by ${\cal L}_{\bS, m}$  the union of these Legendrians:
                   \be
                   {\cal L}_{\bS, m}:=\bigcup_{\mbox{ $\pi \in \pi_0(\partial \bS)$}} {\cal L}_{\pi, m}\ \ \cup \ \ \bigcup_{\mbox{punctures $p$ }}{\cal L}_{p,m}.
                                       \ee
                   \bd The non-commutative stack ${\cal M}_{\cal X}^{(1)}(\bS, {\cal L}_{\bS, m})$ parametrises admissible sheaves on $\bS$ from ${\cal C}(\bS, {\cal L}_{\bS, m})$,
                     vanishing near marked points, whose microlocal support restricted to    the Legendrian ${\cal L}_{\bS, m}$ 
                   is given by   one dimensional      non-commutative local systems.        
                    
                       The stack ${\cal M}_{\cal A}^{(1)}(\bS, {\cal L}_{\bS, m})$     parametrises     similar data  
                       with   trivialized   local systems on the Legendrian. \ed
                   
 \bt The stack ${\cal M}_{\cal X}^{(1)}(\bS, {\cal L}_{\bS, m})$ is birationally isomorphic to the stack ${\cal X}_m(\bS)$.                    
It  has a  non-commutative $\Gamma_{\bS, {\cal L}_{\bS, m}}-$equivariant cluster   Poisson structure. 
The stack ${\cal M}_{\cal A}^{(1)}(\bS, {\cal L}_{\bS, m})$ is birationally isomorphic to the stack ${\cal A}_m(\bS)$.                    
It has a   $\Gamma_{\bS, {\cal L}_{\bS, m}}-$equivariant non-commutative cluster ${\cal A}-$structure.                 
 \et

\begin{proof} 
The complement $\bS-{\cal L}_{\bS, m}$ contains the unique internal component. It is homotopic to $\bS$. Its boundary is a union of circles $\{S^1_p\}$ and $\{S^1_\pi\}$ 
parametrised by   punctures $p$ and boundary components $\pi$. 
 
Restricting  an admissible sheaf ${\cal S} \in {\cal M}^{(1)}(\bS, {\cal L}_{\bS, m})$ to 
the internal component of $\bS-{\cal L}_{\bS, m}$  we get an  $m-$dimensional local system ${\cal V}$.  
  For each puncture $p$, its restriction to the circle $S^1_p$  has a complete filtration, provided by the  
 Legendrian ${\cal L}_p$. 
 
 \begin{figure}[ht]
\centerline{\epsfbox{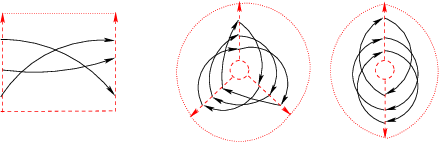}}
\caption{A Legendrian realising permutation $(1,2,3)\to (3,2,1)$, and examples of concatenation.}
\label{stokes4}
\end{figure} 

Given a boundary component $\pi$,   consider $b_\pi$ generic normal rays   cutting the cylinder $C_\pi$ into $b_\pi$ rectangles, see Figure \ref{stokes4}, 
so that the link ${\cal L}_{\pi, m}$ is 
 a  concatenation of  Legendrians in these  rectangles  realizing the longest permutation $w_0$. 
   For each of the  rays $r$,  the   admissible sheaf ${\cal S}$ provides   a   flag   ${\cal F}_{\pi, r}$ 
  in the fiber   of   ${\cal V}$ over  $S^1_\pi \cap r$.  
  The local system ${\cal V}$ with the   flags on the circles $\{S^1_\pi\}$  and filtrations on the circles  $\{S^1_p\}$ 
is exactly the data which determines a framed local system on $\bS$ of dimension $m$. 
 
Evidently,  this construction can be inverted   starting with a framed local systems with the flags  ${\cal F}_{\pi, r}$  in generic position. So   we get a birational isomorphism.  
The ${\cal A}-$case is completely similar. 

The rest follows from the already established results about the spaces ${\cal X}_m(\bS)$ and ${\cal A}_m(\bS)$.   

Alternatively,  we can deform a Legendrian ${\cal L}_{\bS, m}$ to a web of zig-zags of a bipartite graph of type $\bS$, say by using an ideal 
 triangulation of $\bS$ for a start. So we get  a cluster   torus of each of the two flavors. 
                   \end{proof}

 \paragraph{2. Non-commutative Grassmannians as stacks of Stokes data} 

 Consider a web ${\cal W}_{k,n}$ given by a collection of oriented arcs on a disc, cooriented away from the boundary, 
 as shown on Figure \ref{stokes2}. We can describe it   as follows. Take $n$ special points on the boundary, 
 and consider arcs which encircle $k-1$ subsequent points,  so that for each  point there is just one arc starting a bit to the right of a special points, 
 and ending a bit to the left of another special point. 

 Consider the stack ${\rm C}_{\cal X}^\circ({\cal W}_{k,n})$ of admissible dg-sheaves with the following two properties: 
 
 i) They   vanish in the middle of each interval between the special points. 
 
 ii) The local system on each arc describing the microlocal support is one dimensional.
 
  In fact  these dg-sheaves are  sheaves: they are  concentrated in the degree zero. 
  
  The ${\cal A}-$flavor of this space is the stack ${\rm C}_{\cal A}^\circ({\cal W}_{k,n})$ of 
 admissible sheaves as above  such that  the one dimensional local systems on each curve are trivialized. 
 
  \bl i) The stack ${\rm C}_{\cal X}^\circ({\cal W}_{k,n})$ describes  cyclically ordered   collections of $n$ one dimensional subspaces in a  $k-$dimensional $R-$vector space such that any 
 subsequent $k$ of them 
 generate the space. It has a non-commutative cluster Poisson variety structure.  
 
 ii) The stack ${\rm C}_{\cal A}^\circ({\cal W}_{k,n})$ describes  cyclic  collections of $n$ vectors in a  $k-$dimensional $R-$vector space such that any 
 subsequent $k$ of them 
 generate the space. It has a cluster ${\cal A}-$variety structure. 
 \el

   \begin{figure}[ht]
\centerline{\epsfbox{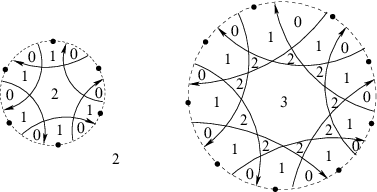}}
\caption{The webs ${\cal W}_{2,5}$ and ${\cal W}_{3,9}$ on a disc  describe configurations of $5$ and $9$ lines/vectors in  non-commutative vector spaces of dimension $2$ and $3$ repectively, as stacks of admissible  sheaves. The numbers tell  dimensions of the fibers of  admissible sheaves in each component of the complement to the web.}
\label{stokes2}
\end{figure}

 \begin{proof} The ambient vector space  is reconstructed as the fiber of an admissible sheaf in the center of the disc. 
 The one-dimensional subspaces are the fibers of the admissible sheaf near   special boundary points.  
 The non-zero vectors generating these one dimensional subspaces are given by trivializations of the local systems sitting on the arcs of the web. 
 The admissibility   at a crossing point  is given by   Lemma \ref{ACP}. This   implies that the  admissibility at the crossing points 
 is equivalent to the condition that each $k$ subsequent   subspaces generate the space.

 To get a cluster description, it is sufficient to find a bipartite graph whose zig-zag web is admissibly isotopic to a web ${\cal W}_{k,n}$. 
Then our general result implies that the torus ${\rm Loc}_1(\Gamma)$ sits inside of the stack we consider.  
This kind of zig-zag webs are the ones considered by Postnikov \cite{P}.  \end{proof}

\section{Theorem \ref{THEOREM4.12} and the  Steinberg relation for non-commutative 2-forms} \la{sec9}

In Section \ref{sec9} we  prove Theorem \ref{THEOREM4.12}, and then observe that one of the key steps of the calculation can be 
viewed as an identity between non-commutative 2-forms generalizing the Steinberg relation. It would be very interesting to find similar identities for the non-commutative $n-$forms generalizing the Steinberg identity $d\log (z_1) \wedge \ldots \wedge d\log (1+x) \wedge \ldots \wedge d\log (x) \wedge \ldots \wedge d\log(z_{n-2})=0$.

\subsection{The proof of Theorem \ref{THEOREM4.12}}

Recall the transformation formulas, where $A_1=a_1a_2a_3a_4$ etc.:
\be \la{FRM+}
 \begin{split}
 &  b_1  =   (1 + A_3^{-1} ) a_3 \ = \ \ \ 
 a_3  (1 + A_4^{-1} ), \\
 &   b_2 =  (1+   A_4)^{-1}a_4 = \ \ \ 
 a_4(1+A_1)^{-1},\\
&  b_3  =   (1 + A_1^{-1} )a_1\ = \ \ \ 
a_1  (1 + A_2^{-1} ),\\
  & b_4 =  (1+A_2)^{-1}a_2 = \ \ \ 
 a_2(1+A_3)^{-1}.\\  \end{split}
 \ee
 Recall the identities 
\be  \la{F1QQ}
\begin{split}
&\{a,b\} =  -\{b^{-1}, a^{-1}\},\\
& \{b, c\} - \{ab, c\} + \{a, bc\} - \{a,b\} = 0,\\
&\{xy, y^{-1}\}= - \{x, y\}, \ \ \ \ \ \ \{x^{-1}, xy\}=-\{x,  y\}.\\
\end{split}
\ee

 \bl \la{5.15}
\be \la{L5.14}
\begin{split}
  &\{b_2, b_3\} = 
      \{a_2a_3a_4, (1 + A_1)^{-1}   \}  + \{a_2a_3, a_4\}   + \{1 + A_1, a_4^{-1}\}.\\
 \end{split}
\ee
\el

\begin{proof} By the definition (\ref{FRM}), we have
\be
\begin{split}
 \{b_2, b_3\} =   &\{a_4(1 + A_1)^{-1}, (1 + A_1^{-1}) a_1\}.   \\
\end{split}
\ee
By the cocycle identity in (\ref{F1QQ}) for the triple $ (a_4, (1 + A_1)^{-1}, (1 + A_1^{-1}) a_1)$, 
\be
\begin{split}
 \{b_2, b_3\}      
= & \{(1 + A_1)^{-1}, (1 + A_1^{-1}) a_1\} + \{a_4, A_1^{-1}a_1\} - \{a_4, (1 + A_1)^{-1}\}. \\ \end{split}
\ee
By the cocycle identity  for the triple $  ((1 + A_1)^{-1}, 1 + A_1^{-1}, a_1)$ we can rewrite the first term, and get 
\be
\begin{split}
  \{b_2, b_3\} =   &
  \{A_1^{-1}, a_1\} 
 - \{ 1 + A_1^{-1}, a_1\} +  \{(1 + A_1)^{-1},  1 + A_1^{-1}\} + 
    \{a_4, A_1^{-1}a_1\}  - \{a_4, (1 + A_1)^{-1}\}.\\
 \end{split}
\ee
Writing the middle term   as $-  \{1 + A_1,  A_1^{-1}\}$ by using the last identiy in (\ref{F1QQ}), 
applying the cocycle identity for the triple   $  (1 + A_1,  A_1^{-1}, a_1)$, and observing that 
$
 \{a_4, A_1^{-1}a_1\}  =     \{a_2a_3, a_4 \}     
$  
 and     $ \{a_3, A_4^{-1}\}  =  -\{A_4, a_3^{-1}\} =  \{a_4a_1a_2, a_3\} $ using the last line in (\ref{F1QQ}), we get  the claim.  
 \end{proof}

 \bl \la{5.15}
\be \la{L5.15}
\begin{split}
 -\{b_1, b_2\} =     &\{   a_3  A_4^{-1}, (1 + A_4) \}      -
 \{ a_4, a_1a_2  \}     
  -\{a^{-1}_4,   1 + A_4  \}.   \\
 \end{split}
\ee
\el

\begin{proof}  By the definition (\ref{FRM}), we have
\be
\begin{split}
 \{b_2^{-1}, b_1^{-1}\} =   &\{a_4^{-1} (1 + A_4), (1 + A_4^{-1})^{-1}a^{-1}_3\}.   \\
\end{split}
\ee
Using the cocycle identity   for the triple $  (a_4^{-1} (1 + A_4), (1 + A^{-1}_4)^{-1},  a_3^{-1})$ we  write this as 
\be
\begin{split}
  -\{b_1, b_2\} =  &-\{ (1 + A^{-1}_4)^{-1},  a^{-1}_3\}  + \{a^{-1}_4  A_4,  a_3^{-1}\}  + \{a_4^{-1} (1 + A_4), (1 + A^{-1}_4)^{-1} \}\\
 = &\{ a_3, 1 + A^{-1}_4   \}  - \{a_1a_2,  a_3\}- \{1 + A^{-1}_4,   (1 + A_4)^{-1}a_4\}. \\ 
\end{split}
\ee
The cocycle identity in (\ref{F1QQ}) allows to write the last term  as
\be
\begin{split}
& 
\{ 1 + A_4, (1 + A^{-1}_4)^{-1}  \}      
-\{a^{-1}_4,  A_4  \}  -  \{a^{-1}_4,    1 + A_4   \}.\\
 \end{split}
\ee
Since the first term  here is equal to 
$\{ A_4^{-1},  1 + A_4   \}$ and $\{a_3,     A^{-1}_4a_4\}     = \{a_1a_2, a_3\}$, we get 
\be \la{L5.15*}
\begin{split}
 -\{b_1, b_2\} =     &\{   a_3,  A_4^{-1}(1 + A_4) \}  - \{a_1a_2,a_3 \}    +
 \{ A^{-1}_4, (1 + A_4)  \}     
-\{a_4,  a_1a_2a_3  \}  -\{a^{-1}_4,   1 + A_4  \}.   \\
 \end{split}
\ee
Using the cocycle relation for the triple $( a_3,  A_4^{-1}, 1 + A_4)$, and then the cocycle relation for $(a_4, a_1a_2, a_3)$ we get the claim.  \end{proof}

Denote by ${\rm Cycl}_2$ the operator on functions in $a_1, ..., a_4$ given by  the identity + the cyclic shift by two:
$$
{\rm Cycl}_2\rm F(a_1, a_2, a_3, a_4)= \rm F(a_1, a_2, a_3, a_4) + \rm F(a_3, a_4, a_1, a_2).
$$
  We split the sum to  calculate  as A) + B), where:\\
       
  A)  It is the total contribution of the middle terms in the Lemmas is: 
      \be
\begin{split}
   &{\rm Cycle}_2\Bigl( -\{a_4,a_1a_2 \}  +  \{a_2a_3, a_4\} \Bigr)   \\
    =&-\{a_4,a_1a_2 \}  +  \{a_2a_3, a_4\}   -\{a_2,a_3a_4 \}  +  \{a_4a_1, a_2\}    \\
    =&\{a_1, a_2 \}  -\{a_2,a_3 \}   +\{a_3, a_4 \}  -  \{a_4, a_1 \}.\\\end{split}
\ee 
The last equality uses   cocycle identities for  $(a_4,a_1, a_2)$ and $(a_2,a_3, a_4 )$. We get   exactly what we wanted. \\

B)   It is the rest:
   \be \la{249}
   \begin{split}
  & {\rm Cycle}_2\Bigl( \{1+A_1, a^{-1}_4\}  +  \{a_2a_3a_4, (1 + A_1)^{-1}   \}    
     -\{a^{-1}_4, 1+A_4\}  
        -\{   (1 + A_4)^{-1}, a_4a_1a_2 \}.   \\
       \end{split}
       \ee
       This  can be rewritten as 
     \be
\begin{split}   
 {\rm Cycle}_2\Bigl( &\Bigl(    - dA_1 a_4^{-1}da_4     +  
      dA_1 (a_2a_3a_4)^{-1} d(a_2a_3a_4)     \Bigr)   (1 + A_1)^{-1} + \\ 
  &   \Bigl(   da_4  a^{-1}_4 dA_4    - d(a_4a_1a_2)  (a_4a_1a_2)^{-1} dA_4    \Bigr) (1 + A_4)^{-1} \Bigr).\\
  \end{split}
\ee 
Expanding this,   we get 
       \be
\begin{split}
&         
     d(a_1 a_2a_3a_4) (a_2a_3a_4)^{-1} d(a_2a_3)a_4       (1 + A_1)^{-1}   \\ 
&       +  d(a_3 a_4a_1a_2) (a_4a_1a_2)^{-1} d(a_4a_1)a_2       (1 + A_3)^{-1}    \\ 
       & -     
a_4d(a_1a_2)  (a_4a_1a_2)^{-1}  dA_4   (1 + A_4)^{-1} \\
       & -       
a_2d(a_3a_4)  (a_2a_3a_4)^{-1}  dA_2   (1 + A_2)^{-1}.\\
\end{split}
\ee 
 Interchanging lines 2 and 3, writing $1+A_4 = a_4(1+A_1)a^{-1}_4$, and similarly for $1+A_3$ and $1+A_2$,  we get:
        \be \la{234}
\begin{split}
&         
     d(a_1 a_2a_3a_4) (a_2a_3a_4)^{-1} d(a_2a_3)a_4       (1 + A_1)^{-1}    \\ 
       & -     
d(a_1a_2)  (a_4a_1a_2)^{-1}  d(a_4a_1a_2a_3)  a_4 (1 + A_1)^{-1} \\
&   +   (a_3a_4)^{-1}   d(a_3 a_4a_1a_2) (a_4a_1a_2)^{-1} d(a_4a_1)a_2a_3a_4       (1 + A_1)^{-1}   \\ 
       & -       
(a_3a_4)^{-1} d(a_3a_4)  (a_2a_3a_4)^{-1}  d(a_2a_3a_4a_1)  a_2a_3a_4 (1 + A_1)^{-1}.\\
\end{split}
\ee 

\bl The expresssion (\ref{234}) is  equal to zero. 
\el

 \begin{proof} Let us write the first two lines as 
          \be \la{2345}
\begin{split}
&   \Bigl( d(a_1 a_2) a_2^{-1} d(a_2a_3)a_4            
 +a_1 a_2 d(a_3a_4) (a_2a_3a_4)^{-1} d(a_2a_3)a_4          \\ 
            & -     
d(a_1a_2)  (a_4a_1a_2)^{-1}  d(a_4a_1)a_2a_3  a_4    -     
d(a_1a_2)  a_2^{-1}  d(a_2a_3)  a_4\Bigr)(1+A_1). \\
\end{split}
\ee
The  terms 1 and 4 cancel, so we get 
\be
\begin{split}
(\ref{2345})&             
 = \Bigl(a_1 a_2 d(a_3a_4) (a_2a_3a_4)^{-1} d(a_2a_3)a_4          
             -     
d(a_1a_2)  (a_4a_1a_2)^{-1}  d(a_4a_1)a_2a_3  a_4\Bigr)(1+A_1).\\
\end{split}
\ee 

\vskip 2mm
A similar calculation with the lines 3 and 4 in (\ref{234}) gives
        \be \la{23456}
\begin{split}
&   \Bigl(  (a_3a_4)^{-1}   d(a_3 a_4)a_4^{-1} d(a_4a_1)       
   +   d(a_1a_2) (a_4a_1a_2)^{-1} d(a_4a_1)       \\ 
       & -       
(a_3a_4)^{-1} d(a_3a_4)  (a_2a_3a_4)^{-1}  d(a_2a_3)a_4a_1    -       
(a_3a_4)^{-1} d(a_3a_4)  a_4^{-1}  d(a_4a_1) \Bigr) a_2a_3a_4(1+A_1).\\
\end{split}
\ee
The terms 1 and 4  cancel, so we get 
\be
\begin{split}
(\ref{23456}) &     =\Bigl(  d(a_1a_2) (a_4a_1a_2)^{-1} d(a_4a_1)        -       
(a_3a_4)^{-1} d(a_3a_4)  (a_2a_3a_4)^{-1}  d(a_2a_3)a_4a_1\Bigr) a_2a_3a_4(1+A_1).\\
\end{split}
\ee

Therefore  (\ref{2345}) + (\ref{23456}) is equal to   
              \be
\begin{split}
&         
     \Bigl( a_1 a_2d(a_3a_4) (a_2a_3a_4)^{-1} d(a_2a_3)a_4           \\ 
       & -     
d(a_1a_2)  (a_4a_1a_2)^{-1}  d(a_4a_1)a_2a_3   a_4  \\
&   +   d(a_1a_2) (a_4a_1a_2)^{-1} d(a_4a_1)a_2a_3a_4           \\ 
       & -       
(a_3a_4)^{-1} d(a_3a_4)  (a_2a_3a_4)^{-1}  d(a_2a_3)a_4A_1\Bigr)(1+A_1).\\
\end{split}
\ee  
Terms 2 and 3 evidently cancel. Terms 1 and 4 cancel   if we permute $A_1$ and $1+A_1$  in the  term  4. 
 \end{proof} 
 
 Therefore the sum B) is equal to zero. Theorem \ref{THEOREM4.12} is proved.

\subsection{An analog of the Steinberg relation   for non-commutative 2-forms.}  \la{SEC11.1}
       
       Note  that 
$$
\{x, y^{-1}\}= -  x^{-1}dx y^{-1}dy, \quad \{x^{-1}, y\}= - dx x^{-1}dy y^{-1}.
$$

\bt \la{T11.4} Let us introduce the notation 
\be
\begin{split}
&\langle a | b\rangle:=  
\{1+ab, b^{-1}\}-  (b^{-1}, 1+ba\}.\\
\end{split}
\ee
Then one has 
   \be \la{432}
\begin{split}
&  
{\rm Cycle}_2\Bigl(\langle a_1a_2a_3 | a_4\rangle -\langle a_1| a_2a_3a_4\rangle    \Bigr)  =0. \\
\end{split}
   \ee\et

\begin{proof} We can write (\ref{249}) as 
\be \la{249s*}
   \begin{split}
        &{\rm Cycle}_2\Bigl( \{1+A_1, a^{-1}_4\}  -  \{1 + A_1, (a_2a_3a_4)^{-1}   \}    
     -\{a^{-1}_4, 1+A_4\}  
        +\{ (a_4a_1a_2)^{-1}, 1 + A_4  \}  
       \Bigr).  \\
       \end{split}
       \ee
  This coincides with (\ref{432}).  But we already proved that  the sum  (\ref{249}) is equal to zero. 
\end{proof}
 
\paragraph{Remark.} 
 If  variables $a_i$ commute, then $\langle a | b\rangle = -2 \cdot d\log(1+ab) \wedge d\log b$, and identity (\ref{432}) reads as
$$
d\log(1+a_1a_2a_3a_4) \wedge d\log (a_1a_2a_3 a_4)=0.
$$
However for a general matrix $X$ we have  $d\log (1+X) \wedge d\log X \not = 0$. The identity (\ref{432}) is the analog 
of the Steinberg relation $d\log (1+z) \wedge d\log (z) =0$ for the non-commutative 2-forms. 
It is interesting that it involves four non-commutative variables $a_1,a_2,a_3,a_4$ rather than one. 

It would be interesting to find a similar non-commutative generalization of the de Rham realization of the defining relation  in  the Milnor K-group:
$d\log (1+z) \wedge d\log (z) \wedge d\log (z_1) \wedge \ldots \wedge d\log (z_n) = 0$.

\end{document}